%% file: annales-toulouse.tex
\documentclass[10pt,twoside,leqno]{amsart}
\usepackage{amssymb,amsbsy,amsmath,amsfonts,amssymb,amscd,times,graphics,color,xypic}
\usepackage[french]{babel}
\usepackage[T1]{fontenc} 
\newcommand{\C}                 {\mathbb{C}}

\newcommand{\K}                 {\mathbb{K}}

\newcommand{\R}                 {\mathbb{R}}
\newcommand{\N}                 {\mathbb{N}}

\newcommand{\V}                 {\mathbb{V}}
\def\dl{{[\![}}
\def\dr{{]\!]}}
\def\1{{\text{\bf 1}}}

\newtheorem{theorem}{Th\'eor\`eme}
\newtheorem{lemma}{Lemme}
\newtheorem{proposition}{Proposition}
\newtheorem{corollary}{Corollaire}
\newtheorem{definition}{D\'efinition}

\newtheorem{assertion}{Assertion}

\begin{document}


\title[R\'egularit\'e analytique
de l'application de r\'eflexion CR formelle
]{
\'Etude de la r\'egularit\'e analytique
de l'application de r\'eflexion CR formelle}

\author{Jo\"el Merker}

\address{
CNRS, Universit\'e de Provence, LATP, UMR 6632, CMI, 
39 rue Joliot-Curie, F-13453 Marseille Cedex 13, France}

\email{merker@cmi.univ-mrs.fr} 


\date{\number\year-\number\month-\number\day}

\begin{abstract}
La recherche de formes normales\footnotemark[1] pour les
sous-vari\'et\'es analytiques r\'eelles de $\C^n$ soul\`eve la
question de la convergence des normalisations formelles. En 1983,
J.K.~Moser et M.W.~Webster ont donn\'e des exemples de surfaces
analytiques r\'eelles dans $\C^2$ \`a tangente complexe isol\'ee et
hyperbolique au sens de E.~Bishop, qui sont formellement mais non
holomorphiquement normalisables (\`a cause d'un ph\'enom\`ene de
petits diviseurs), m\^eme lorsque la forme normale est elle-m\^eme
analytique ou alg\'ebrique. En revanche, il appara\^{\i}t qu'un tel
ph\'enom\`ene ne se produit pas pour les sous-vari\'et\'es dont la
dimension CR est localement constante, d'apr\`es des r\'esultats
r\'ecents dus \`a S.M.~Baouendi, P.~Ebenfelt et L.-P.~Rothschild, et
qui sont \'enonc\'es avec des hypoth\`eses de
non-d\'eg\'en\'eres\-cence relativement simples, mais satisfaites en
un point Zariski-g\'en\'erique\footnotemark[2]. Ces auteurs
\'etablissent notamment que toute application CR formelle inversible
entre deux sous-vari\'et\'es de $\C^n$ g\'en\'eriques, analytiques
r\'eelles, finiment non-d\'eg\'en\'er\'ees et minimales (au sens de
J.-M.~Tr\'epreau et A.E.~Tumanov) est convergente. Nous d\'emontrons
ici un th\'eor\`eme de convergence plus g\'en\'eral, valable sans
aucune hypoth\`ese de non-d\'eg\'en\'erescence, et qui confirme la
rigidit\'e de la cat\'egorie CR ({\it voir} Th\'eor\`eme~1.23). Ce
r\'esultat s'interpr\`ete alors comme un principe de sym\'etrie de
Schwarz formel pour les applications CR. Nous en d\'eduisons que toute
\'equivalence CR formelle entre deux sous-vari\'et\'es de $\C^n$
g\'en\'eriques, analytiques r\'eelles et minimales est convergente si
et seulement si les deux sous-vari\'et\'es sont holomorphiquement
non-d\'eg\'en\'er\'ees (au sens de N.~Stanton). Enfin, nous
\'etablissons que deux sous-vari\'et\'es de $\C^n$ g\'en\'eriques,
analytiques r\'eelles et minimales sont formellement CR \'equivalentes
si et seulement si elles sont biholomorphiquement \'equivalentes.

\smallskip
\noindent
{\sc Abstract} {\sf in English}~: {\sl after the bibliography}.

\end{abstract}

\maketitle

\begin{center}
\begin{minipage}[t]{11cm}
\baselineskip =0.35cm
{\scriptsize

\centerline{\bf Table des mati\`eres}

\smallskip

{\bf Partie~I \dotfill 1.}

\ \
{\bf 1.~Introduction \dotfill 2.}

\ \
{\bf 2.~Pr\'eliminaire~: 
s\'eries formelles, analytiques et alg\'ebriques
\dotfill 20.}

\ \
{\bf 3.~G\'eom\'etrie locale des paires de feuilletages minimales
\dotfill 25.}

\ \
{\bf 4.~Jets de sous-vari\'et\'es de Segre et invariance de
l'application de r\'eflexion \dotfill 34.}

{\bf Partie~II \dotfill 41.}

\ \
{\bf 5.~Convergence d'applications CR formelles finiment 
non-d\'eg\'en\'er\'ees \dotfill 41.}

\ \
{\bf 6.~Convergence d'applications CR formelles 
Segre non-d\'eg\'en\'er\'ees \dotfill 47.}

\ \
{\bf 7.~Convergence de l'application de r\'eflexion CR formelle
\dotfill 55.}

{\bf R\'ef\'erences \dotfill 82.}

\smallskip

{\footnotesize\tt [Avec 2 figures]}

}\end{minipage}
\end{center}

\bigskip

\section*{\S1.~Introduction}

\bigskip
\bigskip

\footnotetext[1]{
{\it Voir}~les travaux fondateurs de S.-S.~Chern, J.K.~Moser~\cite{
cm1974} et de J.K.~Moser, S.M.~Webster~\cite{ mw1983} ainsi que les
articles plus r\'ecents de S.M.~Webster~\cite{ we1992}, de X.~Huang, 
S.G.~Krantz~\cite{ hk1995}, de X.~Gong~\cite{ go1996} et de
P.~Ebenfelt~\cite{ eb1998}.}

\footnotetext[2]{
{\it Voir}~\cite{ ber1997}, \cite{ za1997}, 
\cite{ber1999a}, \cite{ ber1999b}, \cite{ brz2001}. 
}

Cette {\sl Introduction} \'etendue est n\'ecessaire \`a la
pr\'esentation des r\'esultats principaux de ce m\'emoire.
Dans les Sections~3 et~4, qu'il est conseill\'e de lire en
parall\`ele, le lecteur trouvera une exposition autonome,
\'el\'ementaire et d\'etaill\'ee des concepts de base qui y sont
utilis\'es.

Notre objectif principal est d'\'elaborer un point de vue conceptuel
et synth\'etique sur le principe de r\'eflexion dit <<analytique>> en
plusieurs variables complexes. Notamment, nous insisterons sur quatre
concepts majeurs~: la minimalit\'e, les jets de sous-vari\'et\'es de
Segre, l'application de r\'eflexion CR et les conditions
CR-horizontales de non-d\'eg\'en\'erescence.

Nous esp\'erons qu'un tel point de vue permettra de mieux
appr\'ehender l'abondance des th\'eor\`emes d\'emontr\'es r\'ecemment
dans cette direction de recherche. Au lieu de c\'eder, comme certains
auteurs, \`a la tentation d'exprimer tous les r\'esultats et toutes
les cons\'equences possibles, ce qui nuirait \`a l'unit\'e du sujet,
nous s\'electionnerons rigoureusement trois \'enonc\'es~:
Th\'eor\`emes~1.11, 1.23 et~1.39, ainsi que deux applications~:
Corollaires~1.45 et~1.46. Puisqu'une telle s\'election ne saurait se
passer d'explications motiv\'ees, nous d\'etaillerons nos angles
d'attaque \`a partir d'une analyse de la litt\'erature r\'ecente, en
formulant des principes combinatoires organisateurs.

\subsection*{1.1.~\'Equivalences CR formelles entre sous-vari\'et\'es
de $\C^n$ analytiques r\'eelles et g\'en\'eriques} Soit $M$ une
sous-vari\'et\'e locale de $\C^n$, analytique r\'eelle, passant par
l'origine et envisag\'ee au voisinage de ce point de r\'ef\'erence. Si
elle est de codimension $1$, on l'appelle une {\sl hypersurface}. En
g\'en\'eral, nous supposerons que $M$ est de codimension $d \geq 1$
strictement positive et qu'elle est {\sl g\'en\'erique},
c'est-\`a-dire que $T_0 M + J T_0 M = T_0 \C^n$, o\`u $J$ est la
structure complexe standard de $T\C^n$. Nous supposons aussi que sa
dimension CR, \'egale \`a $m := n-d$, est strictement positive.

Dans les coordonn\'ees holomorphes $t: = (t_1,\,\dots,\, t_n)$
canoniques de $\C^n$, la sous-vari\'et\'e $M$, entendue comme
sous-ensemble de $\C^n$, peut \^etre repr\'esent\'ee comme l'ensemble
des $t\in \C^n$ o\`u s'annulent exactement $d$ s\'eries enti\`eres
analytiques $\rho_j (t,\, \bar t)\in\C \{ t,\, \bar t\}$, pour
$j=1,\,\dots,\, d$~; de telles s\'eries enti\`eres doivent satisfaire
aux conditions de r\'ealit\'e $\rho_j (t,\, \bar t) \equiv \overline{
\rho}_j (\bar t,\, t)$ pour $j=1,\, \dots,\, d$, s'annuler pour $t=0$
et poss\'eder des diff\'erentielles r\'eelles $d\rho_1,\, \dots,\, d
\rho_d$ qui sont ind\'ependantes \`a l'origine. La g\'en\'ericit\'e de
$M$ \'equivaut alors au fait que les diff\'erentielles complexes
$\partial \rho_1,\,\dots,\, \partial \rho_d$ sont elles aussi, de
surcro\^{\i}t, ind\'ependantes \`a l'origine.

De m\^eme, soient $\rho_1' (t' ,\, \bar t') = 0,\, \dots,\, \rho_{
d'}' (t',\, \bar t') = 0$ des \'equations cart\'esiennes d\'efinissant
une autre sous-vari\'et\'e locale de $\C^{ n'}$, analytique r\'eelle,
g\'en\'erique, passant par l'origine, qui est de codimension $d' \geq
1$ et de dimension CR \'egale \`a $m':= n'- d'\geq 1$, toutes deux
strictement positives.

Soit $h(t) := (h_1 (t),\, \ldots,\, h_n (t))$ une collection de
s\'eries formelles $h_i (t) \in \C \dl t \dr$, $i=1,\, \dots,\, n$,
dont les termes constants s'annulent. Par d\'efinition, $h$ induit une
{\sl application CR formelle entre $M$ et $M'$} s'il existe une
matrice de taille $d' \times d$ de s\'eries formelles $b(t , \,\bar
t)$ telle que l'on a l'identit\'e formelle vectorielle $\rho'( h(t)
,\, \overline{ h} ( \bar t )) \equiv b (t,\, \bar t) \, \rho(t ,\,\bar
t)$, interpr\'et\'ee dans le produit $\C \dl t,\, \bar t\dr^d := \C\dl
t,\, \bar t \dr \times \cdots \times \C\dl t,\, \bar t \dr$ contenant
exactement $d$ facteurs. On notera $h: (M,\, 0)
\longrightarrow_{\mathcal{ F}} (M',\, 0)$ une telle application, avec
en indice la lettre $\mathcal{ F}$, initiale du mot <<formel>>, car
une telle application n'a rien d'une application ponctuelle. 

Dans le langage abstrait de la g\'eom\'etrie analytique locale,
les applications formelles s'identifient \`a des morphismes entre
anneaux locaux de s\'eries enti\`eres convergentes, le tout \'etant
ins\'er\'e dans une th\'eorie architectur\'ee, mais nous pensons qu'il
est encore trop t\^ot pour entreprendre le travail, certes
d\'esirable, d'abstraction de tous les concepts de g\'eom\'etrie CR
qui interviendront dans ce m\'emoire. 

De notre point de vue, c'est \`a partir de l'identit\'e formelle de
d\'epart $\rho' \left( h(t) ,\,
\overline{ h} ( \bar t ) \right) \equiv b \left( t,\, \bar t \right)
\, \rho \left( t ,\,\bar t \right)$ que tous nos calculs se
d\'eploieront de mani\`ere absolument concr\`ete et structur\'ee.

Nous dirons que $h$ est une {\it \'equivalence} CR formelle si $n' =
n$ et si le d\'eterminant $\hbox{ det} \, \left( { \partial h_{i_1}
\over \partial t_{i_2} } (0) \right)_{ 1 \leq i_1 ,\,i_2 \leq n}$ ne
s'annule pas~; on d\'emontre que ceci implique $d'=d$ et $m'= m$.

Ces d\'efinitions g\'en\'eralisent au cas formel la notion
d'application CR continue, lisse, analytique ou alg\'ebrique entre
sous-vari\'et\'es de $\C^n$, analytiques r\'eelles et g\'en\'eriques.
Avant de poursuivre, dressons un bref historique r\'ecent du principe
de r\'eflexion afin justifier le fait que nous ne travaillerons pas
avec des sous-vari\'et\'es alg\'ebriques.

\subsection*{1.2.~Principe de r\'eflexion alg\'ebrique et variations} 
L'\'etude de la r\'egularit\'e des applications CR formelles est une
variation du principe dit <<de r\'eflexion>>, historiquement initi\'e
dans le cas d'applications de classe $\mathcal{ C}^1$ entre
hypersurfaces strictement pseudoconvexes par S.~Pinchuk
dans~\cite{ pi1975}
et plus tard, mais ind\'ependamment, par H.~Lewy dans~\cite{
le1977}. Le paragraphe de~\cite{ cm1974} consacr\'e \`a la
convergence des formes normales formelles de J.K.~Moser contient une
d\'emonstration implicite du fait que toute \'equivalence CR formelle
entre hypersurfaces de $\C^n$ analytiques r\'eelles et Levi
non-d\'eg\'en\'er\'ees est convergente. Dans le premier article~\cite{
ber1997} consacr\'e \`a la param\'etrisation des applications CR par
un jet d'ordre fini en un point ainsi qu'\`a la r\'egularit\'e des
applications CR formelles, S.M.~Baouendi, P.~Ebenfelt et
L.-P.~Rothschild d\'emontrent que toute application CR formelle entre
deux hypersurfaces de $\C^n$ analytiques r\'eelles et finiment
non-d\'eg\'en\'er\'ees est convergente. Ce r\'esultat r\'epondait \`a
une question soulev\'ee par F.~Treves~; la non-d\'eg\'en\'erescence
finie g\'en\'eralise la Levi non-d\'eg\'en\'erescence par passage aux
d\'eriv\'ees d'ordre sup\'erieur. Peu de temps apr\`es, ces
\'enonc\'es furent transf\'er\'es \`a la codimension quelconque par
les m\^emes auteurs dans~\cite{ ber1998}, en supposant la
sous-vari\'et\'e $M$ minimale \`a l'origine (au sens de
J.-M.~Tr\'epreau et A.E.~Tumanov). La preuve est bas\'ee sur les
it\'erations de sous-vari\'et\'es de Segre, appel\'ees <<ensembles de
Segre>> dans l'article~\cite{ ber1996}, suite de~\cite{
br1995}. Gr\^ace \`a cet outil complexe, ces trois auteurs
\'etablissent une condition suffisante optimale pour l'alg\'ebricit\'e
des biholomorphismes entre sous-vari\'et\'es de $\C^n$ alg\'ebriques
r\'eelles et minimales de codimension quelconque
({\it cf.}~la recension~\cite{ tr2000}).

\'Evidemment, les r\'esultats pr\'ecit\'es ne s'appuient
pas sur une g\'en\'eralisation au cas Levi d\'eg\'en\'er\'e des formes
normales de J.K.~Moser, car nous sommes bien loin actuellement de
poss\'eder des indications valables pour \'edifier une th\'eorie
satisfaisante, malgr\'e quelques tentatives incompl\`etes effectu\'ees
dans cette direction pour la dimension complexe trois ({\it cf.} par
exemple~\cite{ eb1998}).

Par ailleurs, gr\^ace \`a la notion de {\sl degr\'e de transcendance}
d'une application holomorphe, telle qu'elle appara\^{\i}t pour la
premi\`ere fois dans le travail~\cite{ pu1990} de Y.~Pushnikov, bas\'e
sur des indications de S.~Pinchuk, les trois auteurs B.~Coupet,
F.~Meylan et A.~Sukhov sont parvenus \`a \'eliminer dans~\cite{
cms1999} toute hypoth\`ese de rang sur une application holomorphe $h$
entre sous-vari\'et\'es alg\'ebriques r\'eelles. Cependant, en toute
rigueur, on pourrait se dispenser de la terminologie <<degr\'e de
transcendance>>, car le point-cl\'e de~\cite{ pu1990} et de \cite{
cms1999} consiste \`a raisonner avec des polyn\^omes minimaux pour les
relations alg\'ebriques entre les composantes de $h$ et
celles de $\overline{ h}$. Concr\`etement,
cette approche offre une simplification \'el\'egante des calculs
majeurs d\'evelopp\'es pour les sous-vari\'et\'es essentiellement
finies dans~\cite{ bjt1985}, qui avaient \'et\'e illustr\'es
auparavant sur un exemple par M.~Derridj dans~\cite{ de1985}.

Aussi l'approche <<degr\'e de transcendance>> a-t-elle conduit \`a une
s\'erie de travaux r\'ecents~: \cite{ cps2000}, \cite{ me2001a},
\cite{ da2001}, \cite{ cdms2002}, \cite{ mmz2002}, \cite{ mmz2003},
tous bas\'es sur le m\^eme proc\'ed\'e alg\'ebrique
d'\'elimination. En effet, il a \'et\'e rapidement remarqu\'e que le
sch\'ema de d\'emonstration s'\'etendait, sans obstacle majeur et modulo
quelques variations mineures, \`a la situation (l\'eg\`erement plus
g\'en\'erale) o\`u la sous-vari\'et\'e image $M'$ est suppos\'ee
alg\'ebrique r\'eelle, tandis que la sous-vari\'et\'e source $M$ est
suppos\'ee analytique r\'eelle.

Conjecturalement, on s'attend \`a ce que presque tous les principes de
r\'eflexion connus soient valides dans le cas o\`u $M$ et $M'$ sont
{\it toutes deux analytiques r\'eelles}. Malheureusement, la finitude
intrins\`eque au concept d'alg\'ebricit\'e, qui est fortement
utilis\'ee dans ces travaux, fait d\'efaut dans le cas g\'en\'eral o\`u
$M'$ est analytique r\'eelle. On peut alors se demander si une
sous-vari\'et\'e analytique r\'eelle peut \^etre rendue alg\'ebrique
dans un syst\`eme de coordonn\'ees locales ad\'equat.

R\'ecemment, en collaboration avec H.~Gaussier, l'auteur a \'etabli
dans~\cite{ gm2004} que la plupart (au sens de R.~Baire), des
sous-vari\'et\'es de $\C^n$ minimales et finiment
non-d\'eg\'en\'er\'ees dont le groupe d'automorphismes holomorphes
locaux est commutatif et de dimension $n$, ne peuvent \^etre rendues
alg\'ebriques dans aucun syst\`eme de coordonn\'ees holomorphes
locales, quelle que soit la transcendance relative du changement de
coordonn\'ees. On dira que de telles sous-vari\'et\'es ne sont pas
localement {\sl alg\'ebrisables}. Par le biais d'arguments
heuristiques qui extrapoleraient les th\'eor\`emes sp\'ecifiques
de~\cite{ gm2004}, on pourrait se convaincre que dans toute classe de
sous-vari\'et\'es analytiques r\'eelles g\'en\'eriques locales non
homog\`enes dont le groupe d'automorphismes holomorphes poss\`ede une
structure fix\'ee, la plupart d'entre elles ne sont pas localement
alg\'ebrisables. Quant aux sous-vari\'et\'es homog\`enes, le fait
qu'elles soient localement alg\'ebrisales doit se discuter au cas par
cas~; cette caract\'eristique ne d\'epend en effet que de la nature de
leur groupe transitif de transformation, lequel permet bien entendu de
les reconstruire sans ambigu\"{\i}t\'e comme l'orbite d'un point
donn\'e, du reste quelconque.

Un tel ph\'enom\`ene g\'en\'eral de <<non-alg\'ebrisabilit\'e>> locale
g\'en\'erique (au sens de R.~Bai\-re) laisse entrevoir que les
r\'esultats pr\'ecit\'es, qui utilisent fortement l'alg\'ebricit\'e,
sont d'une port\'ee restreinte. C'est pourquoi nous pr\'ef\`ererons
raisonner avec des outils purement analytiques, par exemple les
espaces de jets (\S1.7) ou ce que nous appellerons l'application de
r\'eflexion associ\'ee \`a $h$ (\S1.12). Commen\c cons par pr\'esenter
en r\'esum\'e la g\'eom\'etrie de la complexification extrins\`eque de
$M$.

\subsection*{1.3.~Minimalit\'e} 
C'est un principe g\'en\'eral, d\'ej\`a \`a l'{\oe}uvre dans le
th\'eor\`eme de C.F.~Gauss sur l'existence de coordonn\'ees
isothermes, omnipr\'esent dans les math\'ematiques de la fin du
dix-neuvi\`eme si\`ecle, dans les travaux d'\'E.~Cartan et dans le
d\'eveloppement contemporain de la g\'eom\'etrie CR, qu'il est
judicieux d'introduire d\`es que possible des variables
suppl\'ementaires. En suivant S.-S.~Chern, J.K.~Moser et S.M.~Webster
(\cite{ cm1974}, \cite{ we1977}, \cite{ we1978}, \cite{ mw1983}),
rempla\c cons donc $\bar t$ par une variable $\tau \in \C^n$
ind\'ependante~; on obtient la relation vectorielle $\rho' \left(
h(t),\, \overline{ h} (\tau) \right) \equiv b(t,\, \tau) \, \rho (t,\,
\tau)$, qui exprime que l'application formelle {\sl complexifi\'ee}
$h^c(t,\, \tau):= \left( h (t),\, \overline{ h} (\tau) \right)$ de
l'application CR formelle $h$ induit une application formelle entre la
{\sl complexification de $M$}, qui est la sous-vari\'et\'e analytique
complexe de $\C^{2n}$ d\'efinie par $\mathcal{ M}:= \{ (t,\,\tau)\in
\C^n \times \C^n : \, \rho_j(t, \, \tau)= 0,\, j =1,\, \dots,\, d\}$,
et la complexification $\mathcal{ M}' := \{ (t',\,\tau') \in \C^{n'}
\times \C^{n'} : \, \rho_{j'}'(t', \, \tau') = 0,\ j'=1,\, \dots,\,
d'\}$. Dor\'enavant, nous raisonnerons le plus souvent possible
directement avec les objets g\'eom\'etriques complexifi\'es.

On v\'erifie ({\it voir}\, le~\S3.2 pour des explications) que dans
tout syst\`eme de coordonn\'ees holomorphes locales $(z ,\, w) \in
\C^{m} \times \C^{ d}$ telles que $T_0 M + \left( \{ 0\} \times \C^{ 
d} \right) = T_0 \C^{ n }$, la sous-vari\'et\'e complexifi\'ee
$\mathcal{ M }$ peut \^etre repr\'esent\'ee par $d$ \'equations
analytiques complexes graph\'ees 
de la forme $\xi_{ j} = \Theta_{ j}( \zeta,\,
t)$, pour $j=1,\, \dots,\, d$, o\`u $\tau = ( \zeta,\, \xi ) \in \C^{
m } \times \C^{ d}$. Dans ce syst\`eme de coordonn\'ees, la collection
de s\'eries enti\`eres $\Theta_{ j} \in \C \{ \zeta,\, t\}$ est alors
unique. 

Commen\c cons par pr\'esenter en r\'esum\'e le concept de {\sl
minimalit\'e}. On renvoie \`a la Section~3 pour une exposition de
sa signification g\'eom\'etrique. Soit $p$ un point $\mathcal{ M}$ de
coordonn\'ees $(z_p,\, w_p,\, \zeta_p,\, \xi_p) \in \C^m
\times \C^d \times \C^m \times \C^d$.
Soit $z_1 \in \C^m$ et soit
$\zeta_1 \in \C^m$. On d\'efinit les deux applications
\def\theequation{1.4}\begin{equation}
\left\{
\aligned
\mathcal{ L}_{z_1}(z_p,\, w_p,\, \zeta_p,\, \xi_p) := 
& \
\left(
z_p+z_1,\, \overline{\Theta}(z_p+z_1,\, \zeta_p,\, \xi_p),\, 
\zeta_p,\, \xi_p \right) 
\ \ \ \ \ 
{\rm et} \\
\underline{\mathcal{L}}_{\zeta_1}(z_p,\,w_p,\,\zeta_p,\,\xi_p):=
& \
\left( z_p,\,w_p,\,\zeta_p+\zeta_1,\,
\Theta(\zeta_p+\zeta_1,\,z_p,\,w_p) \right). \\
\endaligned\right.
\end{equation} 
Notons que $\mathcal{ L}_{ z_1} (p)$ et $\underline{ \mathcal{ L}}_{
\zeta_1} (p)$ appartiennent \`a $\mathcal{ M}$. It\'erons ces
applications en les composant alternativement l'une avec l'autre~: en
partant de l'origine $0 \in \mathcal{ M}$, d\'efinissons $\underline{
\Gamma}_1 (z_1) := \underline{ \mathcal{ L} }_{ z_1} (0)$, puis
\def\theequation{1.5}\begin{equation}
\underline{\Gamma}_2(z_1,\,z_2):=
\mathcal{L}_{z_2}(\underline{\mathcal{L}}_{z_1}(0)), 
\end{equation}
puis $\underline{ \Gamma}_3 (z_1,\, z_2,\, z_3) :=
\underline{ \mathcal{ L}}_{ z_3}( \mathcal{ L}_{ z_2} (\underline{
\mathcal{ L} }_{ z_1} (0 )))$, et encore
\def\theequation{1.6}\begin{equation}
\underline{ \Gamma}_4 (z_1,\,
z_2,\, z_3,\, z_4) := \mathcal{ L}_{ 
z_4}( \underline{ \mathcal{ L}}_{
z_3}( \mathcal{ L}_{ z_2} (
\underline{ \mathcal{ L}}_{ z_1} (0 )))),
\end{equation}
et ainsi de suite. Pour tout $k\in \N\backslash \{ 0\}$, on obtient
des applications holomorphes locales $\underline{ \Gamma}_k$ de
$\C^{mk}$ \`a valeurs dans $\mathcal{ M}$ et satisfaisant $\underline{
\Gamma}_k (0) = 0$.

Nous dirons que $M$ est {\sl minimale} (au sens de J.-M.~Tr\'epreau et
A.E.~Tumanov) s'il existe un entier $\mu_0$ tel que l'image par
$\underline{ \Gamma}_{\mu_0}$ d'un voisinage de l'origine
arbitrairement petit dans
$\C^{m\mu_0}$ contient un voisinage de $0$ dans $\mathcal{ M}$. On
d\'emontre que cette condition est invariante par changement de
coordonn\'ees holomorphes locales. C'est tout ce que nous aurons
besoin de savoir au sujet de la g\'eom\'etrie de la sous-vari\'et\'e
source $\mathcal{ M}$.

\subsection*{ 1.7.~Jets de sous-vari\'et\'es de Segre}
Pour le principe de r\'eflexion CR formel, ce sont les \'equations de
la sous-vari\'et\'e image $\mathcal{ M}'$ qui joueront le r\^ole le
plus important. Cependant, le d\'efaut majeur des \'equations
analytiques de la forme $\rho_{j'}'(t',\, \tau') = 0$ r\'eside dans la
non-unicit\'e des s\'eries enti\`eres $\rho_{ j'}'$. En effet, la
sous-vari\'et\'e $\mathcal{ M}'$ est aussi bien repr\'esent\'ee
comme le lieu d'annulation de tout autre jeu de $d'$ s\'eries
enti\`eres analytiques d\'efinissantes qui sont de la forme
$\widetilde{ \rho }_{ j'}' (t',\, \tau'):= \sum_{l'=1}^{d'} c_{j',\,
l' }'(t',\, \tau') \, \rho_{l'}' (t',\, \tau')$, o\`u la matrice de
taille $d' \times d'$ des s\'eries enti\`eres $c_{ j',\, l'}' \in \C\{
t',\, \tau'\}$ est inversible. Pour \'eliminer ce d\'efaut qui
pourrait s'av\'erer g\^enant ({\it cf.}~les
commentaires qui suivent l'\'enonc\'e
du th\'eor\`eme principal~1.23), il convient de repr\'esenter comme un
graphe aussi bien $\mathcal{ M}'$ que $\mathcal{ M}$. Soit donc un
syst\`eme de coordonn\'ees holomorphes locales $(z' ,\, w') \in
\C^{m'} \times \C^{ d'}$ avec $T_0 M' + \left( \{ 0\} \times \C^{ 
d'} \right) = T_0 \C^{ n' }$, dans lequel la sous-vari\'et\'e
complexifi\'ee $\mathcal{ M }'$ est repr\'esent\'ee par $d'$
\'equations analytiques complexes de la forme $\xi_{ j'}' = \Theta_{
j'}'( \zeta',\, t')$, pour $j'=1,\, \dots,\, d'$, o\`u $\tau' = (
\zeta',\, \xi' ) \in \C^{ m' } \times \C^{ d'}$.

Fixons $t'\in \C^{n'}$ et $k \in \N$. La {\sl sous-vari\'et\'e de
Segre complexifi\'ee conjugu\'ee} est la sous-vari\'et\'e analytique
complexe de $\C^{ n'}$ d\'efinie par $\underline{ \mathcal{ S }}_{t'}'
:= \{ (\zeta',\, \xi ') \in \C^{ n'} : \, \xi' = \Theta ' (\zeta',\,
t')\}$. D\'efinissons alors l'application de ses jets d'ordre $k$
explicitement par
\def\theequation{1.8}\begin{equation}
\varphi_k'(\zeta',\, t')
:=
\left( \zeta',\
\left({1\over\beta'!}\, \partial_{\zeta'}^{\beta'}
\Theta_{j'}'(\zeta',\, t')\right)_{
1\leq j'\leq d',\, \vert \beta' \vert \leq
k}\right).
\end{equation}
Elle est \`a valeurs dans $\C^{m'+ N_{d',\, m',\, k}}$, pour un
certain entier $N_{ d',\, m',\, k}$. Comme l'ont compris K.~Diederich
et S.M.~Webster dans~\cite{ dw1980}, ce sont les propri\'et\'es de
cette application analytique complexe locale, d\'efinie au voisinage
de l'origine dans $\C^{ m'} \times \C^{ n'}$, qui gouvernent les
divers principes de r\'eflexion possibles, y compris pour les
applications CR formelles. Nous dirons que $\mathcal{ M}'$ est

\begin{itemize}
\item[{\bf (nd1)}]
{\sl Levi non-d\'eg\'en\'er\'ee \`a l'origine} si $\varphi_1 '$ est de
rang $m'+n'$ en $(\zeta',\, t')= (0,\, 0)$~;
\item[{\bf (nd2)}]
{\sl finiment non-d\'eg\'en\'er\'ee \`a l'origine} s'il existe un
entier $k_0$ tel que $\varphi_k'$ est de rang $n'+ m'$ en $(\zeta',\,
t')= (0,\, 0)$, pour tout $k\geq k_0$~;
\item[{\bf (nd3)}]
{\sl Essentiellement finie \`a l'origine} si 
$\varphi_k'$ est est une application holomorphe
finie en $(\zeta',\, t')= (0,\, 0)$, pour tout 
$k\geq k_0$~; 
\item[{\bf (nd4)}]
{\sl Segre non-d\'eg\'en\'er\'ee \`a l'origine} s'il existe un entier
$k_0$ tel que la restriction de $\varphi_k'$ \`a la sous-vari\'et\'e
de Segre complexifi\'ee conjugu\'ee $\underline{ \mathcal{ S }}_0'$
(qui est de dimension complexe \'egale \`a $m'$) est de rang
g\'en\'erique \'egal \`a $m'$, pour tout $k\geq k_0$~;
\item[{\bf (nd5)}]
{\sl holomorphiquement non-d\'eg\'en\'er\'ee} s'il existe un entier
$k_0$ tel que l'application $\varphi_k'$ est de rang g\'en\'erique
maximal possible, \'egal \`a $m' +n'$, pour tout $k\geq k_0$.
\end{itemize}

On d\'emontre que ces cinq conditions ne d\'ependent pas du syst\`eme
de coordonn\'ees $(z',\, w')$ dans lequel on a repr\'esent\'e
$\mathcal{ M}'$ sous la forme $\xi ' = \Theta'( \zeta',\, t')$ ({\it
voir}~\cite{ me2003} pour les d\'etails).

Bien entendu, par souci de puret\'e, on pourrait en toute rigueur
\'eviter d'introduire une terminologie sp\'ecifique pour de telles
conditions de non-d\'eg\'en\'erescence. N\'eanmoins, puisque l'usage a
d\'ej\`a consacr\'e les quatre conditions que sont {\bf (nd1)}, {\bf
(nd2)}, {\bf (nd3)} et {\bf (nd5)}, nous adopterons ces
d\'enominations \'etablies, et nous en introduirons de nouvelles
ult\'erieurement.

La n\'ecessit\'e d'introduire la condition {\bf (nd4)} est due \`a un
exemple int\'eressant que l'on trouve page~721 de~\cite{ ber2000}.

On v\'erifie ({\it voir} par exemple~\cite{ me2003}) les quatre
implications
\def\theequation{1.9}\begin{equation}
\text{\bf (nd1)} \Rightarrow 
\text{\bf (nd2)} \Rightarrow 
\text{\bf (nd3)} \Rightarrow 
\text{\bf (nd4)} \Rightarrow 
\text{\bf (nd5)},
\end{equation}
dont seule la troisi\`eme est non triviale.

Les conditions {\bf (nd2)} et {\bf (nd3)} apparaissent explicitement
dans~\cite{ dw1980}. C.K.~Han obtient dans~\cite{ ha1983} un principe
de r\'eflexion avec une formulation diff\'erente de {\bf (nd2)}. Il
est important de noter que {\bf (nd1)}, {\bf (nd2)} sont des
hypoth\`eses de rang constant sur l'application de jets~\thetag{
1.4}~; c'est pourquoi elles sont des plus ais\'ees \`a manipuler ({\it
cf.}~\cite{ ber1996}, \cite{ ber1997}, \cite{ ber1998}, \cite{
ber1999a}, \cite{ ber1999b}, \cite{ brz2001}). Quant \`a la condition
{\bf (nd3)}, elle n'est pas tr\`es \'eloign\'ee de {\bf (nd2)} du
point de vue de la th\'eorie des singularit\'es~; elle a \'et\'e
r\'ep\'etitivement pos\'ee comme hypoth\`ese du principe de
r\'eflexion analytique, suite au travail~\cite{ bjt1985} qui
d\'eveloppait une technologie ad\'equate ({\it cf.}~\cite{ br1988},
\cite{ br1990}, \cite{ br1995}, \cite{ ber1999a}, \cite{ cps1999},
\cite{ ber2000} et \cite{ da2001} pour une synth\`ese r\'ecente~; voir
aussi~\cite{ df1988} pour une approche alternative, plus
g\'eom\'etrique).

Il est important de remarquer que les quatre premi\`eres conditions
{\bf (nd1)}, {\bf (nd2)}, {\bf (nd3)} et {\bf (nd4)} sont ponctuelles,
tandis que la derni\`ere ne l'est pas, puisqu'il s'agit d'un rang
g\'en\'erique. Il existe des exemples de
sous-vari\'et\'es g\'en\'eriques qui sont holomorphiquement
non-d\'eg\'en\'er\'ees, mais qui ne sont ni finiment
non-d\'eg\'en\'er\'ees, ni essentiellement finies, ni m\^eme Segre
non-d\'eg\'en\'er\'ees en certains points appartenant \`a un
sous-ensemble analytique r\'eel non vide de $M'$. Ainsi, la notion de
non-d\'eg\'en\'erescence holomorphe est la plus fine des cinq.

La condition interm\'ediaire {\bf (nd4)}, qui appara\^{\i}t
dans~\cite{ me2000}, est d\'ej\`a plus d\'elicate. En v\'erit\'e, pour
les sous-vari\'et\'es $M'$ dites {\sl rigides} dont les fonctions
d\'efinissantes $\Theta_{ j'}' \equiv \Theta_{ j' }'( \zeta',\, z')$
ne d\'ependent pas de $w'$, on v\'erifie que les conditions {\bf
(nd4)} et {\bf (nd5)} sont \'equivalentes. Comme la condition {\bf
(nd5)}, la condition {\bf (nd4)} autorise que les fibres $( \varphi_k'
)^{ -1} ( \varphi_k' ( \zeta',\, t' ))$ de l'application de jets
soient de dimension non localement constante au voisinage de
l'origine. Il est bien connu alors que les concepts standard de
g\'eom\'etrie diff\'erentielle sont insuffisants~: on entre dans le
domaine de la th\'eorie des singularit\'es analytiques complexes.

Dans~\cite{ hi1973}, H.~Hironaka met au point un proc\'ed\'e
d'\'eclatements locaux successifs qui permet, par transformations
strictes successives, de remplacer tout morphisme analytique local par
un morphisme <<redress\'e>> qui satisfait la condition alg\'ebrique
dite de <<platitude>> introduite par J.-P.~Serre. Ce th\'eor\`eme
dit d'<<aplatissement local>> implique la constance locale de la
dimension des fibres, lorsque les espaces d'arriv\'ee et de but du
morphisme <<aplati>> sont lisses. L'existence de ce proc\'ed\'e
sugg\`ere de l'appliquer \`a l'\'etude des applications CR formelles
({\it cf.} \cite{ te1996}). Nous avons constat\'e que cette approche
aboutit dans le cas des applications CR d\'efinies par des s\'eries
enti\`eres, mais dans cet article, nous n'utiliserons pas la th\'eorie
de H.~Hironaka. En effet, gr\^ace \`a un th\'eor\`eme dit
d'approximation d\^u \`a M.~Artin ({\it voir}\, Th\'eor\`eme~2.5
ci-dessous), nous pourrons r\'esumer en partie la complexit\'e
caus\'ee par les singularit\'es de l'application de jets d'ordre
infini~\thetag{ 1.4}, pour $k=\infty$~; l'utilisation de ce
th\'eor\`eme dans le sujet remonte \`a M.~Derridj dans~\cite{ de1986},
d'apr\`es une suggestion de A.~Douady. La complexit\'e de la preuve du
r\'esultat principal (Th\'eor\`eme~1.23) ci-dessous demeurera
substantielle, car elle implique un grand nombre de collections
infinies d'identit\'es formelles.

Dans~\cite{ me2002}, en utilisant la technique g\'eom\'etrique dite
des <<disques analytiques>>, nous avons \'etabli un principe de
r\'eflexion pour les diff\'eomorphismes CR de classe $\mathcal{ C
}^\infty$ entre hypersurfaces analytiques r\'eelles holomorphiquement
non-d\'eg\'en\'er\'ees. Afin d'obtenir une version de ce r\'esultat
en codimension quelconque, nous pensons que le th\'eor\`eme
d'aplatissement local de H.~Hironaka, M.~Lejeune et B.~Tessier devrait
\^etre appliqu\'e au morphisme de jets~\thetag{ 1.8}.

\subsection*{1.10.~Non-d\'eg\'en\'erescence holomorphe} 
Dans~\cite{ st1995} et \cite{ st1996}, N.~Stanton a introduit la
non-d\'eg\'en\'erescence holomorphe. Elle dit que $M'$ est {\sl
holomorphiquement non-d\'eg\'en\'er\'ee} s'il n'existe pas de champ de
vecteurs $X' = \sum_{ i' = 1}^{ n'} \, a_{i'} ' (t') \frac{ \partial
}{ \partial t_{ i' }'}$ {\it non nul}\, et \`a coefficients
holomorphes qui est tangent \`a $M'$. On d\'emontre que cette
d\'efinition \'equivaut \`a celle que nous avons formul\'ee ({\it
voir}~le \S3.4 dans~\cite{ me2003}).

La condition de non-d\'eg\'en\'erescence holomorphe est connue pour
\^etre la condition n\'ecessaire la plus naturelle pour qu'un principe
de r\'eflexion soit valide ({\it cf.}~le premier article \cite{
br1995} contenant une telle observation dans un cadre alg\'ebrique,
{\it cf.}~\cite{ me2001a} pour un r\'esultat sans condition de rang
dans le cadre alg\'ebrique, et {\it cf.}~\cite{ me2002} pour un
r\'esultat r\'ecent dans le cadre $\mathcal{ C
}^\infty$). Rappelons-en le principe. Soit $X' := \sum_{ i'=1}^{ n'}
\, a_{i'} ' (t') \, \frac{ \partial }{ \partial t_{ i'}'}$ un champ de
vecteurs non nul et \`a coefficients holomorphes qui est tangent \`a
$M'$. Soit $(s',\, t')\longmapsto \exp (s' X')(t')$ le flot local de
$X'$, o\`u $s'\in \C$ et $t'\in \C^{n'}$. Gr\^ace \`a la condition de
tangence, ce flot induit une famille \`a un param\`etre
d'automorphismes holomorphes locaux de $M'$. Dans le \S2.16
ci-dessous, nous v\'erifierons qu'il existe de nombreuses s\'eries
enti\`eres {\it non convergentes}\, $\varpi' ( t') \in \C\dl t' \dr$
dont le terme constant est nul, telles que la composition du flot avec
la substition du temps complexe $s'$ par $\varpi'(t')$, c'est-\`a-dire
$t' \mapsto_{ \mathcal{ F}} \exp (\varpi'(t') X') (t')$, est une
auto-application CR formelle inversible de $M'$ non convergente. Une
obstruction similaire au principe de r\'eflexion se produit lorsque
$M'$ est holomorphiquement non-d\'eg\'en\'er\'ee, dans la
cat\'egoris alg\'ebrique, $\mathcal{ C}^\infty$ ou $\mathcal{ C}^0$.

La simplicit\'e de cette obstruction laisse \'evidemment deviner que
la non-d\'eg\'en\'erescen\-ce holomorphe pourrait \^etre une condition
n\'ecessaire et suffisante \`a la convergence de $h$. En supposant $M$
minimale \`a l'origine, nous \'etablissons la r\'eciproque
attendue. C'est notre premier r\'esultat principal, annonc\'e
dans~\cite{ me2001c}.

\def\thetheorem{1.11}\begin{theorem}
Soit $h:\, (M,\, 0) \to_{ \mathcal{ F}} (M',\, 0)$ une \'equivalence
CR formelle entre sous-vari\'et\'es de $\C^n$ analytiques r\'eelles,
g\'en\'eriques, de m\^eme codimension $d \geq 1$ et de m\^eme
dimension CR \'egale \`a $m := n-d \geq 1$. Si $M$ est minimale \`a
l'origine et si $M'$ est holomorphiquement non-d\'eg\'en\'er\'ee,
l'application CR formelle $h$ est convergente.
\end{theorem}

Reste \`a s'interroger sur la n\'ecessit\'e de supposer $M$ minimale
\`a l'origine. On pourrait formuler et d\'emontrer une version formelle
du Th\'eor\`eme~2.7 de~\cite{ me2001a} \'enonc\'e dans un cadre
alg\'ebrique, ce qui donnerai~: si $M$ n'est nulle part minimale et s'il
existe un groupe \`a un param\`etre r\'eel d'auto-applications
holomorphes locales stabilisant $M$, il existe de nombreuses
auto-applications CR formelles de $M$ qui ne sont pas
convergentes. Par ailleurs, il a \'et\'e d\'emontr\'e dans~\cite{
ber1996} que le principe de r\'eflexion alg\'ebrique pour les
biholomorphismes est valide en supposant seulement que la
sous-vari\'et\'e source $M$ est minimale en tout point hors d'un
sous-ensemble analytique r\'eel strict. On parle alors de
minimalit\'e {\sl en un point Zariski-g\'en\'erique}. Il existe donc
une conjecture <<folklorique>> d'apr\`es laquelle le Th\'eor\`eme~1.11
devrait \^etre vrai en supposant seulement que $M$ est minimale en un
point Zariski-g\'en\'erique, ce qui
signifie qu'elle n'est pas forc\'ement minimale
\`a l'origine, mais qu'elle est minimale
en des points arbitrairement proches de l'origine.

Il s'agit l\`a d'une hypoth\`ese fine qui pourrait \`a nouveau
impliquer la th\'eorie des singularit\'es. Pour l'instant, bien que
quelques r\'esultats soient connus pour le principe de r\'eflexion
$\mathcal{ C }^\infty$ entre hypersurfaces non minimales ({\it
cf.}~\cite{ eb2002}), aucun analogue formel n'est connu. Les id\'ees
font d\'efaut, mais nous consid\'erons que cette question m\'erite
d'\^etre \'etudi\'ee.

Une version du Th\'eor\`eme~1.11, valable avec l'hypoth\`ese {\bf
(nd4)} exprimant que $M'$ est Segre non-d\'eg\'en\'er\'ee \`a
l'origine, a \'et\'e d\'emontr\'ee dans le travail non publi\'e~\cite{
me1999}. Pour les applications entre hypersurfaces, le
Th\'eor\`eme~1.11 a \'et\'e d\'emontr\'e d'abord dans la premi\`ere
version de \cite{ me2001b}, et ind\'ependamment ensuite par N.~Mir
dans \cite{ mi2000} (publi\'e auparavant), dont la premi\`ere version
contenait aussi une d\'emonstration d'un \'enonc\'e plus g\'en\'eral~:
la convergence de l'{\sl application de r\'eflexion} pour les
\'equivalences CR formelles entre hypersurfaces.

\subsection*{1.12.~Application de r\'eflexion}
Dans le cas d'une variable complexe, le principe de r\'eflexion
classique d\^u \`a K.H.A.~Schwarz est valable sans hypoth\`ese de
non-d\'eg\'en\'erescence sur la fronti\`ere. Cela n'est pas
\'etonnant, puisque tout arc analytique r\'eel est localement
biholomorphe \`a un segment ouvert de l'axe r\'eel. D\`es qu'il y a
plus de deux variables complexes, il est bien connu que le principe de
r\'eflexion pour les applications CR entre hypersurfaces serait faux
sans une hypoth\`ese forte de non-d\'eg\'en\'erescence, telle que la
Levi non-d\'eg\'en\'erescence ou telle que la finitude essentielle.
Intuitivement parlant, de telles hypoth\`eses expriment que toutes les
variables horizontales $z' \in \C^{ m'}$ ainsi que les conjugu\'ees
complexifi\'ees $\zeta' \in \C^{ m'}$ sont (tr\`es~!)
pr\'esentes dans les s\'eries d\'efinissantes $\Theta_{ j' }'( \zeta'
,\, z',\, w')$ de $\mathcal{ M}'$. On peut toutefois se demander
s'il n'existe pas une g\'en\'eralisation du principe de
r\'eflexion de Schwarz, qui soit valide sans faire aucune hypoth\`ese
sur la fronti\`ere. Dans un travail ant\'erieur ({\it voir}~\cite{
me1997b} et aussi la remarque page~1098 de l'article~\cite{ mm1999}),
l'auteur a effectivement trouv\'e un invariant dont les propri\'et\'es
de r\'egularit\'e sont plus g\'en\'erales que celles de l'application
CR.

Pour une application CR formelle $h :\, (M,\, 0) \longrightarrow_{
\mathcal{ F}} (M',\, 0)$, cet invariant se pr\'esente comme suit. Soient
$\xi_{j'}' - \Theta_{ j'}'( \zeta',\, t')=0$, $j'= 1,\, \dots,\, d'$,
des \'equations complexes arbitraires pour $\mathcal{ M}'$ dans un
voisinage de l'origine. Alors l'{\sl application de r\'eflexion}
associ\'ee \`a $h$ et \`a ce syst\`eme de coordonn\'ees s'exprime par
une s\'erie formelle vectorielle \`a deux variables 
$\tau' \in \C^{n'}$ et $t \in \C^n$~:
\def\theequation{1.13}\begin{equation}
\mathcal{ R}_h' (\tau',\, t):= 
\xi' - \Theta '( \zeta',\, h(t)) \in 
\C \dl \tau',\, t\dr^{d'}. 
\end{equation}
En toute rigueur, cette application $\mathcal{ R}_h'$ d\'epend du
syst\`eme de coordonn\'ees dans lequel on a repr\'esent\'e $\mathcal{
M}'$, mais gr\^ace \`a l'invariance biholomorphe des sous-vari\'et\'es
de Segre, on d\'emontre que la convergence de $\mathcal{ R}_h'$,
c'est-\`a-dire la propri\'et\'e $\mathcal{ R }_h'( \tau' ,\, t) \in \C
\{ \tau',\, t\}^{d'}$, est une propri\'et\'e invariante ({\it voir}~la
Section~4 ci-dessous). Dans ce m\'emoire, nous \'enoncerons et
d\'emontrerons un r\'esultat de convergence de $\mathcal{ R}_h'$,
valable sans aucune hypoth\`ese de non-d\'eg\'en\'erescence sur $M'$~:
le Th\'eor\`eme principal~1.23 ci-dessous.

\subsection*{1.14.~R\'esultats r\'ecents}
Au cours de la d\'emonstration des r\'esultats principaux des
articles~\cite{ bjt1985} et~\cite{ br1988} (consacr\'es au principe de
r\'eflexion pour les applications CR de classe $\mathcal{ C
}^\infty$), une expression voisine de~\thetag{ 1.13} appara\^{\i}t,
{\it mais on n'a pas pas vu qu'elle devrait jouir d'une propri\'et\'e
de r\'egularit\'e analytique, m\^eme lorsque la sous-vari\'et\'e image
$M'$ est tr\`es d\'eg\'en\'er\'ee}. Pour les applications CR entre
hypersurfaces essentiellement finies, le prolongement holomorphe de
$\mathcal{ R}_h'$ s'identifie au prolongement de $h$ en tant que
correspondance, tel qu'il est trait\'e dans~\cite{ df1988}, \cite{
dfy1994}, \cite{ dp1995}, \cite{ dp1998}, \cite{ sh2000}, \cite{
pv2001}, \cite{ dp2003}, \cite{ sh2003}, parfois avec l'hypoth\`ese
plus forte que $M'$ ne contienne pas de courbe holomorphe.

En 1998, inspir\'e par nos conjectures sur l'application de
r\'eflexion dans~\cite{ me1997b} et par la remarque de la page~1098 du
preprint de~\cite{ mm1999}, N.~Mir a d\'emontr\'e de mani\`ere
ind\'ependante dans~\cite{ mi1998} que l'application de r\'eflexion
associ\'ee \`a un biholomorphisme entre hypersurfaces alg\'ebriques
est alg\'ebrique. Sa d\'efinition de l'application de r\'eflexion
\'elimine la variable $\xi'$ dans~\thetag{ 1.13}~; celle-ci a pourtant
un sens g\'eom\'etrique, puisque la d\'efinition de $\mathcal{ R}_h'$
est intrins\`equement reli\'ee aux sous-vari\'et\'es de Segre
complexifi\'ees. Ni dans ce travail, ni dans d'autres travaux
ult\'erieurs~\cite{ mi2000}, \cite{ mi2002}, \cite{ bmr2002}
consacr\'es \`a l'application de r\'eflexion, N.~Mir ne traite
l'invariance biholomorphe de cette application. Or, pour peu que l'on
\'etablisse que l'alg\'ebricit\'e de $\mathcal{ R}_h'$ est une
propri\'et\'e invariante par changement de coordonn\'ees alg\'ebriques
({\it cf.}~\cite{ me2001a} et \cite{ me2002}), le r\'esultat principal
de~\cite{ mi1998} devient un corollaire \'el\'ementaire de~\cite{
br1995}. En effet ({\it cf.}~le \S11 de~\cite{ me2001a}), en d\'epla\c
cant l\'eg\`erement le point de r\'ef\'erence en un point o\`u
l'application de jets~\thetag{ 1.8} est de rang localement
constant\,--\,ce qui est autoris\'e puisque l'application
consid\'er\'ee est d\'ej\`a holomorphe dans un ouvert\,--\,et en
\'eliminant des variables muettes, on se ram\`ene \`a un
biholomorphisme local entre deux sous-vari\'et\'es alg\'ebriques
r\'eelles, g\'en\'eriques et finiment non-d\'eg\'en\'er\'ees contenues
dans des espaces euclidiens complexes de dimensions inf\'erieures~;
alors les r\'esultats de~\cite{ br1995} (codimension $1$) ou de~\cite{
ber1996} (codimension quelconque) s'appliquent directement. En
r\'esum\'e, pour le principe de r\'eflexion alg\'ebrique,
l'alg\'ebricit\'e de la fonction de r\'eflexion \'equivaut \`a
l'alg\'ebricit\'e d'une application holomorphe entre sous-vari\'et\'es
g\'en\'eriques de dimension inf\'erieure ({\it cf.} le
Th\'eor\`eme~11.4 de~\cite{ me2001a} qui \'etablit cette \'equivalence
sans aucune hypoth\`ese de rang).

En revanche, pour les applications CR formelles (ou $\mathcal{
C}^\infty$), il est vraiment impossible de d\'eplacer la situation
locale en un point o\`u les singularit\'es de l'application de
jets~\thetag{ 1.8} disparaissent. 

Avant d'\'enoncer notre r\'esultat principal, pr\'esentons une
deuxi\`eme liste hi\'erarchis\'ee, comportant cinq conditions
CR-horizontales de non-d\'eg\'en\'erescence.

\subsection*{1.15.~Conditions CR-horizontales 
de non-d\'eg\'en\'erescence} Pour les sous-vari\'et\'es g\'en\'eriques
$M'$ quelconques, sans aucune condition de non-d\'eg\'en\'erescence,
le principe de r\'eflexion implique d\'ej\`a des questions
d\'elicates, m\^eme en supposant, pour simpli\-fier, que l'application
$h$ est inversible. Mais puisqu'un grand nombre de raffinements ont eu
cours durant la derni\`ere d\'ecennie et qu'il est presque toujours
possible d'arguer de la <<nouveaut\'e>> d'un r\'esultat qui suppose
l'application $h$ toujours un peu plus d\'eg\'en\'er\'ee\,--\,pourvu
que les techniques connues s'appliquent encore\,--, nous pensons qu'il
est n\'ecessaire d'exposer un principe organisateur pour
pr\'esenter les conditions de non-d\'eg\'en\'erescence sur
l'application $h$.

Pour exprimer ces conditions, travaillons d'embl\'ee avec les
\'equations complexes graph\'ees du \S1.3 pour $\mathcal{ M}'$, ou
plut\^ot avec les \'equations conjugu\'ees $w_{j'} = \overline{
\Theta}_{j'}' (z',\, \tau')$, $j'=1,\, \dots,\, d'$, qui sont
\'equivalentes, d'apr\`es le \S3.2. Repr\'esentons aussi la
complexification $\mathcal{ M }$ de $M$ par des \'equations complexes
de la forme $w_j = \overline{ \Theta}_j (z,\, \tau )$, $j = 1,\,
\dots,\, d$ dans des coordonn\'ees adapt\'ees $t = (z,\, w) \in \C^m
\times \C^d$. Posons $\overline{ r}_j (\tau,\, t) := w_j - \overline{
\Theta }_j (z,\, \tau)$ et $\overline{ r}_{j'}' (\tau',\, t') :=
w_{j'}'- \overline{ \Theta}_{j'}' (z',\, \tau')$. Par hypoth\`ese, il
existe une matrice de taille $d' \times d$ de s\'eries formelles
$\overline{ b} (\tau,\, t)$ telle que $\overline{ r}' \left(
\overline{ h} (\tau),\, h (t) \right) \equiv \overline{ b}(\tau ,\, t)
\, \overline{ r} (\tau,\, t)$ dans $\C \dl t,\, \tau
\dr^d$. D\'ecomposons les composantes de l'application d'une mani\`ere
compatible avec le scindage $(z',\, w')\in \C^{ m'} \times \C^{ d'}$
des coordonn\'ees, ce qui donne $h(t) =: ( f(t),\, g(t))\in \C \dl t
\dr^{ m'} \times \C \dl t \dr^{ d'}$. En rempla\c cant $w$ par
$\overline{ \Theta} (z,\, \tau)$ dans l'identit\'e fondamentale
$\overline{ r}' \left( \overline{ h} (\tau ),\, h (t) \right) \equiv
\overline{ b}( \tau,\, t) \, \overline{ r} (\tau,\, t)$, le second
membre s'annule identiquement et nous obtenons les identit\'es
formelles suivantes, valables dans $\C \dl z,\, \tau \dr$~:
\def\theequation{1.16}\begin{equation}
g_{j'} \left(z,\, \overline{ \Theta} (z,\, \tau) \right) \equiv
\overline{ \Theta}_{ j'}' \left( f(z,\, \overline{
\Theta} (z,\, \tau),\, \overline{ h}( \tau)) \right),
\end{equation}
pour $j'=1,\, \dots,\, d'$. En posant $\tau = 0$ dans ces
identit\'es, on obtient les identit\'es
\def\theequation{1.17}\begin{equation}
g_{j'} \left(z,\, \overline{ \Theta} (z,\, 0) \right) \equiv
\overline{ \Theta}_{ j'}' \left( f(z,\, \overline{
\Theta} (z,\, \tau),\, 0 ) \right). 
\end{equation}
Classiquement, on les interpr\`ete en exprimant que la restriction de
l'application formelle $h$ \`a la sous-vari\'et\'e de Segre $S_0$
passant par l'origine, d\'efinie par $\{ (z,\, w) \in \C^n : \, w =
\overline{ \Theta} (z,\, 0)\}$, induit une application formelle \`a
valeurs dans la sous-vari\'et\'e de Segre $S_0'$ de l'espace image
d\'efinie par $\{ (z',\, w') \in \C^{ n'} : \, w' = \overline{ \Theta}
' (z',\, 0)\}$. Alors la restriction de $h$ \`a $S_0$ co\"{\i}ncide
avec l'application formelle
\def\theequation{1.18}\begin{equation}
\C^m \ni z \longmapsto_{\mathcal{ F}} \
\left(
f\left(
z,\, \overline{ \Theta} (z,\, 0)
\right), \ \
\overline{ \Theta}'\left(
f\left(
z,\, \overline{ \Theta} (z,\, 0),\, 0
\right)
\right)
\right)\in \C^{m'} \times \C^{ d'}.
\end{equation}
Ce sont les propri\'et\'es de non-d\'eg\'en\'erescence de cette
application formelle induite qui gouvernent les divers raffinements
possibles du principe de r\'eflexion analytique. Par projection sur
le sous-espace $\C^{ m'} \times \{ 0\}$, on peut \'evidemment
identifier cette application avec sa {\em partie CR-horizontale}
d\'efinie par
\def\theequation{1.19}\begin{equation}
\C^m \ni z \longmapsto_{\mathcal{ F}} \
f\left(
z,\, \overline{ \Theta} (z,\, 0)
\right)\in \C^{m'}.
\end{equation}
Avec ces notations, nous pouvons formuler tr\`es concr\`etement cinq
conditions sur l'application~\thetag{ 1.19}, que nous ordonnons par
ordre croissant de g\'en\'eralit\'e. L'application CR formelle $h$
sera dite

\begin{itemize}
\item[{\bf (cr1)}]
{\sl CR-inversible} \`a l'origine si $m'=m$ et si sa partie
CR-horizontale est une \'equivalence formelle en $z=0$~;
\item[{\bf (cr2)}]
{\sl CR-submersive} \`a l'origine si $m' \leq m$ et si sa partie
CR-horizontale est une submersion formelle en $z=0$~;
\item[{\bf (cr3)}]
{\sl CR-finie} \`a l'origine si $m' = m$ et si sa partie
CR-horizontale est une application formelle finie en $z=0$~;
\item[{\bf (cr4)}]
{\sl CR-dominante} \`a l'origine 
si $m' \leq m$ et si sa partie CR-horizontale est
dominante en $z=0$~;
\item[{\bf (cr5)}]
{\sl CR-transversale} \`a l'origine si sa partie CR-horizontale est
transversale en $z=0$~;
\end{itemize}

Nous renvoyons le lecteur au \S2.3 pour des d\'efinitions compl\`etes
de ces cinq conditions de non-d\'eg\'en\'erescence, valides dans la
cat\'egorie des applications formelles quelconques~; la condition {\bf
(cr5)}, qui n'implique {\it aucune in\'egalit\'e entre $m'$ et $m$},
est exprim\'ee {\it in extenso}\, au d\'ebut du \S1.21. Bien entendu,
en utilisant l'invariance biholomorphe des sous-vari\'et\'es de Segre,
on d\'emontre que ces d\'efinitions ne d\'ependent pas des syst\`emes
de coordonn\'ees dans lesquels on repr\'esente $M$ et $M'$.

Les cinq conditions du \S2.3 sont tout \`a fait classiques en
g\'eom\'etrie analytique locale. On v\'erifie les quatre implications
\def\theequation{1.20}\begin{equation}
\text{\bf (cr1)} \Rightarrow
\text{\bf (cr2)} \Rightarrow
\text{\bf (cr3)} \Rightarrow
\text{\bf (cr4)} \Rightarrow
\text{\bf (cr5)},
\end{equation}
pourvu que $m'=m$ dans la deuxi\`eme et dans la troisi\`eme. Dans le
contexte CR, les conditions {\bf (cr1)} et {\bf (cr2)} apparaissent
dans~\cite{ za1997}~; la condition {\bf (cr3)}, maintenant classique,
est une condition naturelle pour les applications CR entre
hypersurfaces essentiellement finies~; elle appara\^{\i}t dans~\cite{
df1988}, \cite{ br1988} et dans d'autres r\'ef\'erences. La condition
{\bf (cr4)} appara\^{\i}t dans~\cite{ br1990}. Enfin, la condition
{\bf (cr5)} appara\^{\i}t dans~\cite{ ber2000}, avec une appellation
diff\'erente. Mais dans cette r\'ef\'erence, les auteurs supposent
la sous-vari\'et\'e
$M'$ essentiellement finie~: ils travaillent avec la condition
{\bf (nd3)}, bien comprise depuis le travail fondateur~\cite{
bjt1985}.

Le pr\'efixe commun aux cinq conditions <<CR->> se justifie de la
mani\`ere suivante~: puisque l'espace tangent au point $0 \in
\mathcal{ S}_0$ co\"{\i}ncide avec l'espace tangent complexe \`a $M$
en $0$, lequel absorbe la {\sl structure CR infinit\'esimale} de $M$
en $0$, on peut penser que l'application formelle induite $h
\vert_{\mathcal{ S}_0} : \, (\mathcal{ S}_0,\, 0) \longrightarrow_{
\mathcal{ F}} (\mathcal{ S}_0',\, 0)$ est un <<prolongement>> de
l'application {\sl CR tangente} $dh : \, T_0^c M \to T_0^c M'$.

S\'electionnons maintenant la condition {\bf (cr5)}, puisque c'est la
plus g\'en\'erale.

\subsection*{1.21.~R\'esultat principal} 
Par d\'efinition ({\it cf.}~le \S2.3 ci-dessous), une application CR
formelle $h$ comme dans le \S1.15 est {\sl CR-transversale \`a
l'origine} s'il n'existe pas de s\'erie formelle $F'(z_1',\, \dots,\,
z_{ m'}') \in \C \dl z_1',\, \dots,\, z_{ m'}' \dr$ non nulle telle
que l'on a l'identit\'e 
\def\theequation{1.22}\begin{equation}
F'\left( f_1 \left(z,\, \overline{
\Theta} (z,\, 0)\right),\, \dots,\,
f_{m'} \left(z,\, \overline{ \Theta} 
(z,\, 0)\right) \right) \equiv 0,
\end{equation}
dans $\C \dl z \dr$. Le r\'esultat principal de cet article, dont le
Th\'eor\`eme~1.11 d\'ecoule en v\'erit\'e comme corollaire, est le
suivant.
\def\thetheorem{1.23}\begin{theorem}
Soit $h:\, (M,\, 0) \to_{ \mathcal{ F }} (M',\, 0)$ une application CR
formelle entre deux sous-vari\'et\'es de $\C^n$, $\C^{n'}$ analytiques
r\'eelles, g\'en\'eriques, de codimensions $d \geq 1$, $d' \geq 1$ et
de dimensions CR \'egales \`a $m:= n-d \geq 1$, $m' := n' - d' \geq
1$. Si $M$ est minimale \`a l'origine et si $h$ est CR-transversale,
pour tout syst\`eme de coordonn\'ees $(z',\, w') \in \C^{ m'} \times
\C^{ d'}$ dans lequel la complexification $\mathcal{ M}'$ est
repr\'esent\'ee par $\xi ' = \Theta '( \zeta',\, t')$, l'application
de r\'eflexion CR formelle associ\'ee $\mathcal{ R }_h' (\tau',\, t) := 
\xi' - \Theta '( \zeta',\, h(t))$ est convergente.
\end{theorem}

Dans le \S4.30, nous d\'emontrerons que si cette propri\'et\'e de
convergence est satisfaite dans un tel syst\`eme de coordonn\'ees
$(z',\, w')$, alors pour tout autre syst\`eme de coordonn\'ees
$(z'',\, w'')$ centr\'ees \`a l'origine dans lesquelles la
complexification de la sous-vari\'et\'e transform\'ee est
repr\'esent\'ee par des \'equations similaires $\xi_{j'}''- \Theta_{
j'}''( \zeta'',\, t'') = 0$, $j'=1,\, \dots,\, d'$, l'application de
r\'eflexion associ\'ee est elle aussi convergente.

La force principale de ce th\'eor\`eme r\'eside dans le fait qu'il ne
requiert aucune condition de non-d\'eg\'en\'erescence sur $M'$. Comme
pour le Th\'eor\`eme~1.11, nous pensons bien entendu qu'il demeure
valide en supposant seulement que $M$ est minimale en un point
Zariski-g\'en\'erique.

Attention, il y a un pi\`ege~! Par souci de g\'en\'eralit\'e, on
pourrait \^etre tent\'e comme dans~\cite{ bmr2002} de raisonner avec
des \'equations d\'efinissantes analytiques r\'eelles arbitraires
$\rho_{ j'} '( t',\, \bar t') = 0$ pour $M'$, telles qu'introduites
dans le \S1.1. L'application de r\'eflexion associ\'ee serait alors
d\'efinie par $\widehat{ \mathcal{ R}}_h' (\tau',\, t):= \rho_{ j'}'
(h(t),\, \tau') \in \C \dl t,\, \tau' \dr^{ d'}$, et le
Th\'eor\`eme~1.23 exprimerait, sous les m\^emes hypoth\`eses, qu'elle
est convergente. Mais en 1997, J.-M.~Tr\'epreau nous a fait remarquer
qu'un tel \'enonc\'e serait trivialement faux.

En effet, choisissons une s\'erie enti\`ere non convergente $\varpi(
z_2) \in \C \dl z_2\dr$ telle que $\varpi (z_2) = z_2 + {\rm
O}(z_2^2)$ et consid\'erons l'application formelle d\'efinie par
$h(z_1,\, z_2,\, w ):= (z_1,\, \varpi(z_2),\, w)$. C'est une
\'equivalence CR formelle entre l'hypersurface alg\'ebrique $M$ de
$\C^3$ d\'efinie par $w= \bar w + i z_1\bar z_1$ et (la m\^eme~!)
l'hypersurface de $\C^3$ d\'efinie par $r ' =0$, o\`u $r' := \bar w' -
w' + i z_1' \bar z_1'$. Notons que $M'$ est holomorphiquement
d\'eg\'en\'er\'ee, puisque le champ holomorphe $\frac{ \partial }{
\partial z_2'}$ lui est tangent. Il est vrai que l'application de
r\'eflexion $\mathcal{ R}_h' (\tau',\, t)$ \'egale
\`a $\xi' - w + i z_1 \,
\zeta_1 '$ est convergente. Par contre, il est vraiment faux que
l'application de r\'eflexion associ\'ee \`a une \'equation
d\'efinissante arbitraire pour $M'$ est convergente. En effet, prenons
par exemple la fonction $\rho'(t',\, \bar t') : = [1+ z_1' \bar z_1' +
z_2' \bar z_2'] \, r'(t',\, \bar t')$~; son lieu d'annulation
co\"{\i}ncide avec $M'$. Si l'application de r\'eflexion
\def\theequation{1.24}\begin{equation}
\widehat{ \mathcal{ R}}_h' ( \tau',\, t) := 
\left[
1+ z_1 \zeta_1 ' + \varpi(z_2)\, \zeta_2']\cdot
[\xi' - w + iz_1\, \zeta_1'
\right]
\end{equation}
\'etait convergente par rapport aux six variables $(z_1,\, z_2,\, w,\,
\zeta_1',\, \zeta_2',\, \xi')$, on en d\'eduirait en consid\'erant
$\left. \frac{ \partial^2 }{ \partial \zeta_2 ' \partial \xi '}
\widehat{ \mathcal{ R}}_h ' (\tau',\, t) \right\vert_{ \tau' =0}$ que
la s\'erie formelle $\varpi (z_2)$ est convergente, ce qui contredirait
notre choix initial. Par cons\'equent, il n'est pas anodin de choisir
d'embl\'ee des \'equations complexes graph\'ees $\xi' = \Theta '(
\zeta',\, t')$ pour la repr\'esentation de la complexifi\'ee
$\mathcal{ M }'$ ainsi que pour la d\'efinition de l'application de
r\'eflexion CR formelle.

Comme nous l'avons mentionn\'e \`a la fin du \S1.11, N.~Mir a obtenu
dans~\cite{ mi2000} une d\'emonstration du Th\'eor\`eme~1.23 pour les
\'equivalences formelles dans le cas $d= d' =1$, mais la m\'ethode,
astucieuse, achoppe d\`es que la codimension $d$ de $M$ est
sup\'erieure ou \'egale \`a $2$. La premi\`ere version de~\cite{
me2001b}, qui a circul\'e avant~\cite{ mi2000}, contenait seulement le
Th\'eor\`eme~1.11 dans le cas $d = d' =1$. Dans cette r\'ef\'erence,
l'existence de paires d'identit\'es de r\'eflexion conjugu\'ees
apparaissait clairement, bien qu'exploit\'ee de mani\`ere
insuffisante. En fait, dans la version publi\'ee~\cite{ me2001b}, il a
suffi d'inclure le court \S9 pour obtenir le Th\'eor\`eme~1.23 pour
les \'equivalences formelles dans le cas $d= d' = 1$. Cette {\sl
paire d'identit\'es de r\'eflexion conjugu\'ees} \'etant absolument
cruciale pour la d\'emonstration du Th\'eor\`eme~1.23, nous allons
l'exposer dans le \S1.25 ci-dessous.

En Mai 2000, une d\'emonstration compl\`ete du Th\'eor\`eme~1.23 pour
les \'equivalences CR formelles a \'et\'e annonc\'ee dans~\cite{
me2000}. Cette annonce \'electronique a donn\'e lieu \`a la
publication r\'esum\'ee~\cite{ me2001c}. Sept mois plus tard, en
d\'ecembre 2000, S.M.~Baouendi, N.~Mir et L.-P.~Rothschild ont
annonc\'e \'electroniquement le m\^eme type de r\'esultats, avec les
raffinements attendus sur le rang de l'application $h$, lesquels ne
s'\'el\`event pourtant que jusqu'au niveau {\bf (cr4)}. Un examen de
la publication~\cite{ bmr2002} \`a laquelle a donn\'e lieu ce travail
(qui ne contient plus les r\'ef\'erences \`a nos travaux pr\'esentes
dans la version \'electronique) montre que ces auteurs utilisent les
paires d'identit\'es de r\'eflexion conjugu\'ees, ce que seul
l'ultra-sp\'ecialiste peut d\'eceler dans le c{\oe}ur technique de la
d\'emonstration principale ({\it voir}~les \'equations~\thetag{ 5.2}
et~\thetag{ 5.3}, la Proposition~6.1 et le Lemme~7.1 de~\cite{
bmr2002}). Par ailleurs, ces auteurs, qui n'emploient pas la
terminologie <<application de r\'eflexion>> (utilis\'ee pourtant
dans~\cite{ mi2000}, \cite{ mi2002}), introduisent une notion
alternative d'<<id\'eal de Segre>>, laquelle est d\'efinie au moyen
d'\'equations analytiques r\'eelles arbitraires $\rho_{ j'} '( t',\,
\bar t') =0$ pour $M'$. Ce choix pour \'enoncer leurs th\'eor\`emes de
convergence, les contraint \`a quelques circonlocutions, puisque la
convergence de $\widehat{ \mathcal{ R }}_h' (\tau',\, t) := \rho_{
j'}' (h(t),\, \tau') \in \C \dl t,\, \tau' \dr^{ d'}$ n'est satisfaite
que pour les repr\'esentants graph\'es de l'id\'eal engendr\'e par les
s\'eries enti\`eres complexifi\'ees $\rho_{ j'} ' (t',\, \tau')$,
comme nous venons de le voir.

En conclusion de ce paragraphe, la prolixit\'e et le raffinement des
r\'esultats pr\'esent\'es dans~\cite{ bmr2002} confinent \`a un
certain herm\'etisme auquel nous n'adh\`ererons jamais, puisque les
quatre concepts analytico-g\'eom\'etriques qui sont impliqu\'es dans
le sujet sont relativement simples~:

\begin{itemize}
\item[{\bf (1)}]
Minimalit\'e locale comme propri\'et\'e des orbites de champs de
vecteurs CR complexifi\'es, dont les sous-vari\'et\'es int\'egrales
co\"{\i}ncident avec les sous-vari\'et\'es de Segre complexifi\'ees
(Section~3)~;
\item[{\bf (2)}]
jets d'ordre $k$ des sous-vari\'et\'es de Segre complexifi\'ees
et diverses conditions de non-d\'eg\'en\'erescence (Section~4)~;
\item[{\bf (3)}]
application de r\'eflexion CR comme invariant fondamental 
qui jouit de propri\'et\'es de r\'egularit\'e
(Sections~1, 5 et 6)~;
\item[{\bf (4)}]
conditions de non-d\'eg\'en\'erescence CR-horizontales
(Section~1)~;
\end{itemize}

Ce sont ces hypoth\`eses multiples, combin\'ees souvent \`a
l'alternative entre cat\'egorie alg\'ebrique et cat\'egorie
analytique, qui sont responsables de la combinatoire de th\'eor\`emes
possibles publi\'es r\'ecemment sur le principe de r\'eflexion
analytique. Toutefois, cette diversit\'e s'exerce au d\'etriment de
r\'esultats plus rares o\`u une difficult\'e substantielle a \'et\'e
surmont\'ee et elle occulte leur rep\'erage.

Exposons maintenant le point-cl\'e qui est \`a la base de la
d\'emonstration du Th\'eor\`eme~1.23.

\subsection*{1.25.~Paire d'identit\'es de r\'eflexion conjugu\'ees} 
Soient $w_j = \overline{ \Theta}_j (z,\, \tau)$, $j= 1,\, \dots,\, d$
un syst\`eme de $d$ \'equations complexes graph\'ees pour la
complexification $\mathcal{ M}$ de $M$ et soient $\xi_j = \Theta_j
(\zeta,\, t)$, $j=1,\,\dots,\, d$, les \'equations conjugu\'ees. On
consid\`ere la paire de syst\`emes de $m$ champs de vecteurs
holomorphes tangents \`a $\mathcal{ M}$ d\'efinis comme suit~:
\def\theequation{1.26}\begin{equation}
\left\{
\aligned
\mathcal{L}_k:= & \ 
{\partial \over\partial z_k}+
\sum_{j=1}^d\,
{\partial\overline{\Theta}_j\over \partial z_k}
(z,\, \tau)
\, {\partial\over\partial w_j}, \ \ \ \ 
k=1,\dots,\,m,\\
\underline{\mathcal{L}}_k:= & \
{\partial\over\partial \zeta_k}+\sum_{j=1}^d\,
{\partial \Theta_j\over \partial \zeta_k}(\zeta,\, t)\,
{\partial\over\partial \xi_j}, \ \ \ \
k=1,\dots,\,m.
\endaligned\right.
\end{equation} 
Par hypoth\`ese, l'application CR formelle $h(t)= (f(t),\, g(t)) \in
\C \dl t \dr^{ m'} \times \C \dl t \dr^{ d'}$ satisfait les $d'$
identit\'es formelles~\thetag{ 1.16} et leurs conjugu\'ees
complexifi\'ees, que nous \'ecrirons ensemble comme suit~:
\def\theequation{1.27}\begin{equation}
\left\{
\aligned
\overline{ g}_{ j'} \left(
\zeta, \, \Theta (\zeta,\, t) \right) 
& \
\equiv
\Theta_{ j'} '\left(
\overline{ f} (\zeta,\, \Theta (\zeta,\, t)),\, 
h(t)
\right), \\
g_{j'} \left(z,\, \overline{ \Theta} (z,\, \tau) \right) 
& \
\equiv
\overline{ \Theta}_{ j'}' \left( f\left(
z,\, \overline{
\Theta} (z,\, \tau) \right),\, \overline{ h}( \tau) \right).
\endaligned\right.
\end{equation}
On les abr\`egera en les \'ecrivant $g(t) = \overline{ \Theta} '
\left(
f(t),\, \overline{ h} ( \tau)\right)$ 
et $\overline{ g} ( \tau) = \Theta ' \left(
\overline{ f} (\tau),\, h(t) \right)$, 
\'etant entendu que $(t ,\, \tau) \in
\mathcal{ M}$. D\'eveloppons les fonctions $\Theta_{ j'}'$ par rapport
aux puissances de $\zeta'$~: on obtient des expressions de la forme
$\Theta_{ j'} '(\zeta',\, t') = \sum_{ \gamma ' \in \N^{ m'}} \,
(\zeta')^{ \gamma'} \, \Theta_{j',\, \gamma'} ' (t')$, avec des
s\'eries enti\`eres convergentes $\Theta_{j',\, \gamma'}' (t')\in \C
\{ t' \}$ qui satisfont bien s\^ur une estim\'ee de Cauchy, puisque
les fonctions $\Theta_{ j'} '( \zeta',\, t')$ sont holomorphes par
rapport aux deux variables $\zeta'$ et $t'$. En utilisant ce
d\'eveloppement, nous pouvons tout d'abord r\'e\'ecrire l'application
de r\'eflexion~\thetag{ 1.13} sous la forme plus explicite
\def\theequation{1.28}\begin{equation}
\mathcal{ R}_h' (\tau',\, t) =
\xi' - 
\sum_{ \gamma ' \in \N^{m'}} \, 
(\zeta')^{\gamma '} \, 
\Theta_{ \gamma '} '( h(t)). 
\end{equation}
Dans la Section~5 ci-dessous, nous \'etablirons que la convergence de
l'application de r\'eflexion est \'equivalente \`a la convergence de
la collection infinie de s\'eries formelles $\Theta_{ j',\, \gamma'} '
(h(t))$, pour tous $j'= 1,\, \dots,\, d'$ et tous $\gamma ' \in
\N^{m'}$. Nous appellerons {\sl composantes de l'application de
r\'eflexion} ces s\'eries formelles $\Theta_{ j',\, \gamma'} '( h(t))$,
parfois not\'ees sous la forme vectorielle abr\'eg\'ee $\Theta_{
\gamma '} ' (h(t))$.

En utilisant le m\^eme d\'eveloppement partiel en s\'eries enti\`eres
des $\Theta_{ j',\, \gamma'}'(\zeta',\, t')$, on peut aussi
r\'e\'ecrire les relations fondamentales~\thetag{ 1.27} sous une forme
plus explicite, dont le m\'erite principal est de faire clairement
appara\^{\i}tre toutes les composantes de l'application de
r\'eflexion~:
\def\theequation{1.29}\begin{equation}
\left\{
\aligned
\overline{ g} (\tau) 
& \
=
\sum_{ \gamma' \in \N^{m'}}\, 
\overline{ f} (\tau)^{\gamma '} \,
\Theta_{\gamma' } ' (h(t)), \\
g (t) 
& \
=
\sum_{\gamma' \in\N^{ m'}} \, 
f(t)^{\gamma '} \, 
\overline{ \Theta}_{ 
\gamma'} ' \left(\overline{ h} (\tau) \right),
\endaligned\right.
\end{equation}
o\`u $(t,\, \tau) \in\mathcal{ M}$.
Venons-en maintenant aux identit\'es de r\'eflexion. Pour un
multiindice arbitraire $\beta = (\beta_1,\, \beta_2,\, \dots,\,
\beta_m) \in \N^m$, on note $\mathcal{ L }^\beta$ et $\underline{
\mathcal{ L} }^\beta$ les d\'erivations holomorphes et antiholomorphes
d'ordre $\vert \beta \vert$ d\'efinies par
\def\theequation{1.30}\begin{equation}
\left\{
\aligned
\mathcal{ L}^\beta :=
& \ 
( \mathcal{ L }_1 )^{ \beta_1}(
\mathcal{ L }_1 )^{ \beta_2} \cdots ( \mathcal{L }_m)^{ \beta_m}
\ \ \ \ \
{\rm et} \\
\underline{ \mathcal{ L}}^\beta :=
& \ 
(\underline{ \mathcal{ L }}_1 )^{ \beta_1}
(\underline{ \mathcal{ L
}}_1 )^{ \beta_2} \cdots ( 
\underline{ \mathcal{L }}_m)^{ \beta_m}.
\endaligned\right.
\end{equation}

Classiquement, on applique les d\'erivations antiholomorphes
$(\underline{ \mathcal{ L } }_1 )^{ \beta_1 }( \underline{ \mathcal{ L
}}_1 )^{ \beta_2} \cdots ( \underline{ \mathcal{L } }_m)^{ \beta_m }$
au premier jeu d'\'equations~\thetag{ 1.29}. De mani\`ere
\'equivalente, \`a une conjugaison pr\`es, on pourrait appliquer les
d\'erivations conjugu\'ees $( \mathcal{ L }_1 )^{ \beta_1}( \mathcal{
L }_1 )^{ \beta_2} \cdots ( \mathcal{L }_m)^{ \beta_m}$ au second
jeu d'\'equations~\thetag{ 1.29}. Au total, les deux proc\'ed\'es
reviennent \`a choisir une fois pour toutes les variables $t$ ou les
variables $\bar t$ pour \'ecrire les identit\'es de r\'eflexion.
C'est le point de vue qui est adopt\'e dans tous les travaux
consacr\'es au principe de r\'eflexion analytique que sont \cite{
pi1975}, \cite{ le1977}, \cite{ we1978}, \cite{ we1982}, \cite{
dw1980}, \cite{ ha1983}, \cite{ de1985}, \cite{ bjt1985}, \cite{
br1988}, \cite{ pu1990}, \cite{ br1990}, \cite{ br1995}, \cite{
ss1996}, \cite{ ber1996}, \cite{ ber1997}, 
{\small \cite{ mi1998}}, \cite{
ber1999a}, \cite{ cms1999}, \cite{ cps1999}, \cite{ ber1999b}, \cite{
cps2000}, \cite{ ber2000}, \cite{ me2001a}, \cite{ brz2001}, \cite{
me2002}, \cite{ eb2002}, \cite{ cdms2002}, \cite{ mmz2002} et \cite{
mmz2003}. Nous l'analyserons dans le \S1.33 ci-dessous.

Respectons une exigence de compl\'etude~: nous disposons de deux
jeux de $d'$ \'equations formelles fondamentales~\thetag{ 1.29} et de
deux jeux infinis de d\'erivations fondamentales~\thetag{ 1.30}. Au
total, ce ne sont donc pas deux mais quatre identit\'es de r\'eflexion
que nous devrions obtenir.

Pour les \'ecrire, on observe que $\underline{ \mathcal{ L}}_k (h)
\equiv 0$ et que $\mathcal{ L }_k \left( \overline{ h } \right) \equiv
0$ pour $k = 1,\, \dots,\, m$, ce qui est \'evident d'apr\`es les
formules~\thetag{ 1.26}. Il en d\'ecoule que $\underline{ \mathcal{
L}}^\beta ( h ) \equiv 0$ et que $\mathcal{ L }^\beta \left(
\overline{ h } \right) \equiv 0$ puis aussi $\underline{ \mathcal{ L}
}^\beta \Theta_{\gamma'} '( h) \equiv 0$ et $\mathcal{ L }^\beta
\left( \overline{ \Theta} ' \left( \overline{ h } \right) \right)
\equiv 0$, pourvu bien s\^ur que $\beta \neq 0$.

Ainsi, en appliquant les d\'erivations~\thetag{1.30} pour $\beta \neq
0$ aux identit\'es~\thetag{ 1.29}, on obtient quatre familles infinies
d'identit\'es de r\'eflexion. Disposons-les en deux paires
conjugu\'ees comme suit~: premi\`ere paire~:
\def\theequation{1.31}\begin{equation}
\left\{
\aligned
\underline{ \mathcal{ L}}^\beta \, 
\overline{ g} (\tau) 
& \
=
\sum_{ \gamma' \in \N^{m'}}\,
\underline{ \mathcal{ L}}^\beta \left[ 
\overline{ f} (\tau)^{\gamma '}\right] \,
\Theta_{\gamma' } ' (h(t)), \\
0
& \
=
\sum_{\gamma '\in \N^{ m'}} \, 
f(t)^{\gamma '} \, 
\underline{ \mathcal{ L}}^\beta \left[
\overline{ \Theta}_{ 
\gamma'} ' \left(\overline{ h} (\tau) \right)\right];
\endaligned\right.
\end{equation}
seconde paire, conjugu\'ee (modulo transposition) de la premi\`ere~:
\def\theequation{1.32}\begin{equation}
\left\{
\aligned
\mathcal{ L}^\beta 
g (t)
& \
=
\sum_{\gamma' \in \N^{ m'}} \,
\mathcal{ L}^\beta \left[ 
f(t)^{\gamma '}\right] \, 
\overline{ \Theta}_{ 
\gamma'} ' \left(\overline{ h} (\tau) \right), 
\\
0
& \
=
\sum_{ \gamma' \in \N^{m'}}\, 
\overline{ f} (\tau)^{\gamma '} \,
\mathcal{ L}^\beta \left[
\Theta_{\gamma' } ' (h(t))\right]. \\
\endaligned\right.
\end{equation}
Bien entendu, $\beta \neq 0$ et $(t,\, \tau) \in \
\mathcal{ M}$. Les deux paires~\thetag{ 1.31}
et~\thetag{ 1.32} sont donc conjugu\'ees terme \`a terme (modulo une
transposition)~; elles ne sont donc pas essentiellement
distinctes. Mais dans chacune des deux paires, une diff\'erence
importante est \`a noter~: tandis que ce sont 
{\it a priori}\, toutes les composantes
$\overline{ f}_{ k'}$ et $\overline{ g}_{ j'}$ de l'application
$\overline{ h}$ que l'on diff\'erentie dans la premi\`ere
identit\'e~\thetag{ 1.31}, ce sont les composantes conjugu\'ees de
l'application de r\'eflexion $\overline{ \Theta}_{ \gamma'} ' \left(
\overline{ h} (\tau) \right)$ que l'on diff\'erentie \`a la seconde
ligne.

La diff\'erence a son importance pour la raison suivante. D'apr\`es la
propri\'et\'e mentionn\'ee apr\`es~\thetag{ 1.28}, le
Th\'eor\`eme~1.23 \'enonce essentiellement que les composantes
$\Theta_{ \gamma'}'( h(t))$ de l'application de r\'eflexion sont des
s\'eries convergentes. Il n'\'enonce nullement que toutes les
composantes de l'application $h(t)$ sont convergentes. L'exemple
\'el\'ementaire discut\'e apr\`es le Th\'eor\`eme~1.23 (ou d'autres
analogues) montre qu'en g\'en\'eral, aucune contrainte de
convergence n'est exerc\'ee sur les composantes de $h$ qui
n'apparaissent pas dans les composantes de l'application de
r\'eflexion. C'est pourquoi la premi\`ere identit\'e de
r\'eflexion~\thetag{ 1.31} a le d\'efaut majeur de faire intervenir
in\'evitablement les d\'eriv\'ees d'\'eventuelles <<mauvaises>>
composantes de $h$, tout du moins celles qui ne sont pas
intrins\`equement li\'ees \`a l'application invariante $\mathcal{
R}_h'$. Au contraire, dans la seconde identit\'e~\thetag{ 1.31} (tout
aussi bien que dans la premi\`ere identit\'e~\thetag{ 1.32}), on
diff\'erentie les vrais objets invariants que sont les composantes
$\overline{ \Theta }_{ \gamma'}' \left( \overline{ h } \right)$ (ou
leurs conjugu\'ees $\Theta_{ \gamma '} '( h )$).

\`A notre connaissance, les seuls travaux dans lesquels on consid\`ere
aussi les seconde identit\'es de r\'eflexion~\thetag{ 1.31} (ainsi que
sa conjugu\'ee, la premi\`ere identit\'e de~\thetag{ 1.32})
sont~: \cite{ me2001b}, \cite{ me2001c} et \cite{ bmr2002}. Avant de
poursuivre le commentaire, passons en revue ce qu'il est possible
d'\'enoncer au moyen de la premi\`ere identit\'e de
r\'eflexion~\thetag{ 1.31}.

\subsection*{1.33.~Conditions de non-d\'eg\'en\'erescence sur $h$}
Au lieu de formuler conjointement des conditions de
non-d\'eg\'en\'erescence sur $M'$ avec des conditions de
non-d\'eg\'en\'erescence CR-horizontales sur $h$ ({\it cf.}~\cite{
br1988}, \cite{ ber1997}, \cite{ ber1998}, \cite{ za1997}, \cite{
ber1999a}, \cite{ ber1999b}, \cite{ ber2000}, \cite{ brz2001} et
\cite{ bmr2002}), formulons des conditions de non-d\'eg\'en\'erescence
sur les premi\`eres identit\'es de r\'eflexion~\thetag{ 1.31} ({\it
cf.}~\cite{ ha1990}, \cite{ ss1996}, \cite{ cps1999}, \cite{ me1999},
\cite{ da2001}). Pour cela, introduisons une collection infinie de
s\'eries enti\`eres formelles d\'ependant des trois variables $t\in
\C^n$, $\tau \in \C^n$ (avec $(t,\, \tau) \in \mathcal{ M}$) et $t'
\in \C^{n'}$, laquelle est d\'efinie en rempla\c cant $h(t)$ par $t'$
dans la premi\`ere ligne de~\thetag{ 1.31}, ce qui donne~:
\def\theequation{1.34}\begin{equation}
\Psi_{ j',\, \beta} '(t,\, \tau,\, 
t') := \underline{\mathcal{L}}^\beta \,
\overline{ g}_{j'}- \sum_{\gamma' \in\N^{ m '}}\,
\underline{\mathcal{L}}^\beta \left[
\overline{ f}^{\gamma'}\right]\,
\Theta_{j',\, \beta}'(t'),
\end{equation}
pour $j'= 1,\, \dots,\, d'$ et $\beta \in \N^m$. Dans~\thetag{ 1.36}
ci-dessous, on notera aussi ces s\'eries $\Psi_{ j', \, \beta} '( z,\,
w,\, \zeta,\, \xi,\, t')$.

Soit $k\in \N$. En ne consid\'erant que les multiindices $\beta \in
\N^m$ de longueur $\vert \beta \vert \leq k$, et en posant $(t,\,
\tau):= 0$, d\'efinissons l'application formelle suivante, qui est \`a
valeurs dans $\C^{ N_{ d',\, n',\, k}}$ pour un certain entier $N_{
d',\, n',\, k}$~:
\def\theequation{1.35}\begin{equation}
\psi_k': \ 
t' \longmapsto 
\left(
\Psi_{ j',\, \beta} ' (0,\, 0,\, t')
\right)_{1 \leq j' \leq d',\, 
\vert \beta \vert \leq k}.
\end{equation}
Puisque les termes $\underline{ \mathcal{ L}}^\beta \, \overline{ g}_{
j'}$ et $\underline{ \mathcal{ L}}^\beta \left[ \overline{ f}^{
\gamma'}\right]$ sont constants lorsque l'on pose $(t,\, \tau) = (0,\,
0)$, l'application $\psi_k'$ est holomorphe au voisinage de l'origine.

L'application CR formelle $h$ sera dite

\begin{itemize}
\item[{\bf (h1)}]
{\sl Levi non-d\'eg\'en\'er\'ee} \`a l'origine 
si $\psi_1'$ est de rang $n'$ en $t' = 0$~;
\item[{\bf (h2)}]
{\sl finiment non-d\'eg\'en\'er\'ee} \`a l'origine s'il existe un
entier $\ell_0$ tel que $\psi_k '$ est de rang $n'$ en $t' = 0$ pour tout
$k\geq \ell_0$~;
\item[{\bf (h3)}]
{\sl essentiellement finie} \`a l'origine s'il existe un 
entier $\ell_0$ tel que $\psi_k'$ est une application 
holomorphe finie pour tout $k \geq \ell_0$~;
\item[{\bf (h4)}]
{\sl Segre non-d\'eg\'en\'er\'ee} \`a l'origine s'il existe des
entiers $j'(1),\,\dots,\, j'(n')$ satisfaisant $1\leq j'(i_1') \leq
d'$ pour $i_1'=1,\, \dots,\, n'$, et des multiindices distincts
$\beta(1),\, \dots,\, \beta (n')\in \N^m$ tels que le 
d\'eterminant suivant~:
\def\theequation{1.36}\begin{equation}
{\rm det} \ \left( \frac{ \partial 
\Psi_{ j'( i_1'), \, \beta(
i_1')}'}{ \partial t_{ i_2'}'} 
\left( z,\, \overline{ \Theta} (z,\,
0), \, 0,\, 0,\, h \left( z,\, 
\overline{ \Theta} (z,0) \right) \right) \right)_{ 1\leq
i_1',\, i_2' \leq n'}
\end{equation}
ne s'annule pas identiquement dans $\C \dl z \dr$.
\end{itemize}

On v\'erifie que ces quatre conditions ne d\'ependent pas du syt\`eme
de coordonn\'ees holomorphes locales dans lequel la sous-vari\'et\'e
$M'$ est repr\'esent\'ee par $\xi' = \Theta ' ( \zeta',\, t')$ ({\it
cf.}~\cite{ me2003}). La condition {\bf (h1)} appara\^{\i}t
dans~\cite{ ss1996} ainsi que dans d'autres r\'ef\'erences
pr\'ec\'edentes~; la condition {\bf (h2)} appara\^{\i}t dans~\cite{
la2000}~; la condition {\bf (h3)} appara\^{\i}t dans~\cite{ cps1999},
\cite{ me1999}, \cite{ da2001} et \cite{ dm2002}~; ces r\'ef\'erences
expriment la condition {\bf (h3)} en disant, de mani\`ere
\'equivalente, que la {\sl vari\'et\'e caract\'eristique} $\V_0' := \{
t ': \, \Psi_{ j',\, \beta} '( 0,\, 0,\, t') = 0, \ j'= 1,\, \dots,\,
d', \ \beta \in \N^m\}$, qui est un sous-ensemble analytique complexe
de $\C^{ n'}$ passant par l'origine, est de dimension z\'ero en
$t'=0$~; enfin, la condition nouvelle {\bf (h4)} appara\^{\i}t
dans~\cite{ me1999}.

Avant de pr\'esenter les quatre principes de r\'eflexion auxquels
donnent naissance les quatre conditions {\bf (h1)}, {\bf (h2)}, {\bf (h3)} et
{\bf (h4)}, nous voudrions mentionner qu'elles sont satisfaites si
l'on effectue une hypoth\`ese de non-d\'eg\'en\'erescence sur $M'$
combin\'ee \`a une hypoth\`ese de non-d\'eg\'en\'erescence
CR-horizontale sur $h$.

\def\theproposition{1.37}\begin{proposition}

Supposons $h$ CR-transversale \`a l'origine. 
\begin{itemize}
\item[{\bf (1)}]
Si $M'$ est Levi non-d\'eg\'en\'er\'ee \`a l'origine, $h$ est
finiment non-d\'eg\'en\'er\'ee \`a l'origine.
\item[{\bf (2)}]
Si $M'$ est finiment non-d\'eg\'en\'er\'ee \`a l'origine, $h$
est finiment non-d\'eg\'en\'er\'ee \`a l'origine.
\item[{\bf (3)}]
Si $M'$ est essentiellement finie \`a l'origine, 
$h$ est essentiellement finie \`a l'origine.
\item[{\bf (4)}]
Si $M'$ est Segre non-d\'eg\'en\'er\'ee \`a l'origine, 
$h$ est Segre non-d\'eg\'en\'er\'ee \`a l'origine.
\end{itemize}
\end{proposition}

Puisque la condition {\bf (cr5)} est la plus g\'en\'erale, cette
proposition se d\'emultiplie en quatre autres propositions
\'enonc\'ees avec {\bf (cr1)}, avec {\bf (cr2)}, avec {\bf (cr3)} ou
avec {\bf (cr4)} \`a la place de {\bf (cr5)}. Une partie de ces vingt
assertions, mais pas la totalit\'e, se trouve implicitement
d\'emontr\'ee dans les travaux de S.M.~Baouendi, L.P.~Rothschild et
divers co-auteurs. Pour la d\'emonstration relativement technique de
ce r\'esultat que nous n'utiliserons pas, nous renvoyons le lecteur 
au Th\'eor\`eme~4.3.1 de~\cite{ me2003}.

Mentionnons toutefois les quatre observations~:
\begin{itemize}
\item[{\bf (1)}]
$h$ Levi non-d\'eg\'en\'er\'ee $\Rightarrow$ 
$M'$ Levi non-d\'eg\'en\'er\'ee~;
\item[{\bf (2)}]
$h$ finiment non-d\'eg\'en\'er\'ee 
$\Rightarrow$ $M'$ finiment non-d\'eg\'en\'er\'ee~;
\item[{\bf (3)}]
$h$ essentiellement finie $\Rightarrow$
$M'$ essentiellement finie~;
\item[{\bf (4)}]
$h$ Segre non-d\'eg\'en\'er\'ee $\Rightarrow$ 
$M'$ Segre non-d\'eg\'en\'er\'ee.
\end{itemize}
Chacune de ces implications est stricte, comme le montreraient des
exemples \'el\'ementaires analogues \`a ceux qui sont d\'evelopp\'es
dans~\cite{ me2001d}. La condition de CR-transversalit\'e sur $h$ est
la condition la plus fine qu'il faut ajouter pour garantir les
r\'eciproques de ces implications.

\subsection*{ 1.38.~Quatre principes de r\'eflexion CR formels pour $h$}
Ainsi, travaillons directement avec les quatre conditions de
non-d\'eg\'en\'erescence sur $h$.

\def\thetheorem{1.39}\begin{theorem}
{\rm (\cite{ me1999}, Theorem~1.2.1)} Soit $h:\, (M,\, 0) \to_{
\mathcal{ F }} (M',\, 0)$ une application CR formelle entre deux
sous-vari\'et\'es de $\C^n$, $\C^{n'}$ analytiques r\'eelles,
g\'en\'eriques, de codimensions $d \geq 1$, $d' \geq 1$ et de
dimensions CR \'egales \`a $m:= n-d \geq 1$, $m' := n' - d' \geq
1$. Si $M$ est minimale \`a l'origine et si
\begin{itemize}
\item[{\bf (i)}]
$h$ est Levi ou finiment non-d\'eg\'en\'er\'ee \`a l'origine~; ou si
\item[{\bf (ii)}]
$h$ est essentiellement finie \`a l'origine~; ou si
\item[{\bf (iii)}]
$h$ est Segre
non-d\'eg\'en\'er\'ee \`a l'origine,
\end{itemize}
$h$ est convergente.
\end{theorem}

La partie {\bf (i)} est essentiellement d\'emontr\'ee dans~\cite{
ber1997} et~\cite{ ber1998} ({\it voir} aussi \cite{ la2000}). La
partie {\bf (ii)} est essentiellement d\'emontr\'ee dans~\cite{
ber2000}, gr\^ace \`a une variation sur les arguments classiques
de~\cite{ bjt1985}, mais avec des conditions de
non-d\'eg\'en\'erescence s\'epar\'ees sur $M'$ et sur $h$. On trouve
une copie conforme de cette d\'emonstration relativement complexe
dans~\cite{ mi2002}, o\`u la condition {\bf (h3)} est emprunt\'ee
\`a~\cite{ cps1999} et \`a~\cite{ me1999}. Le Th\'eor\`eme~1.2
de~\cite{ dm2002} co\"{\i}ncide avec ce m\^eme r\'esultat, mais
gr\^ace \`a l'argument \'el\'ementaire de tranchage tir\'e de~\cite{
me1999}, on peut affirmer que la complexit\'e de la preuve donn\'ee
dans~\cite{ ber2000} ou dans~\cite{ mi2002} n'a pas lieu
d'\^etre. R\'ep\'etons que la condition {\bf (h3)} est paradigmatique
depuis le travail classique~\cite{ bjt1985}.

Lorsque {\bf (h1)} ou {\bf (h2)} est satisfaite, en appliquant le
th\'eor\`eme des fonctions implicites, on d\'emontre qu'il existe un
entier $\ell_0 \geq 1$ et une application $\Phi$ holomorphe locale
telle que les identit\'es de r\'eflexion apparaissant \`a la
premi\`ere ligne de~\thetag{ 1.31} se r\'esolvent 
par rapport \`a $h(t)$, gr\^ace au
th\'eor\`eme des fonctions implicites, sous la
forme~:
\def\theequation{1.40}\begin{equation}
h(t) \equiv 
\Phi \left( t,\, \tau,\, J_\tau^{\ell_0} \, 
\overline{ h} (\tau)\right),
\end{equation} 
pour $(t,\, \tau ) \in \mathcal{ M}$ (le \S5.2 ci-dessous fournit les
d\'etails). Ici, la notation $J_\tau^{\ell_0} \overline{ h} (\tau)$
d\'esigne le jet d'ordre $\ell_0$ de $\overline{ h}(\tau)$. Une telle
relation remonte aux travaux fondateurs de S.~Pinchuk~\cite{ pi1975}
et H.~Lewy~\cite{ le1977}. Elle est appel\'ee {\sl identit\'e de
r\'eflexion basique} dans~\cite{ ber1999a}. Afin de rendre plus
accessible la d\'emonstration du Th\'eor\`eme~1.23, nous
reconstituerons la d\'emonstration tr\`es simple de la convergence de
$h$ sous cette hypoth\`ese dans la Section~5, en utilisant notre
propre formalisme.

Lorsque {\bf (h3)} est satisfaite, on d\'emontre gr\^ace \`a un
proc\'ed\'e d'\'elimination alg\'ebrique standard ({\it cf.}~les
travaux \cite{ bjt1985}, \cite{ ber1999a} reprenant les d\'etails de
ce proc\'ed\'e bien connu depuis le trait\'e de Van der Waerden~\cite{
vdw1970}), qu'il existe un entier $\ell_0 \geq 1$ et $n'$ polyn\^omes
tels que les identit\'es de r\'eflexion qui apparaissent \`a la
premi\`ere ligne de~\thetag{ 1.31} fournissent une
<<quasi-r\'esolution polynomiale>> de $h(t)$ par rapport au jet
d'ordre $\ell_0$ de $\overline{ h}(\tau)$, c'est-\`a-dire que l'on a
$n'$ relations polynomiales de la forme
\def\theequation{1.41}\begin{equation}
h_{i'}(t)^{N(i')}+ 
\sum_{ 1\leq k' \leq N(i')} \, 
A_{ i',\, k'}' \left( t,\, \tau,\, 
J_\tau^{\ell_0} \overline{ h} (\tau) \right)
\, 
h_{i'}(t)^{N(i')-k'}
\equiv 0,
\end{equation}
pour $i'=1,\, \dots,\, n'$ et pour $(t,\, \tau) \in \mathcal{ M}$, 
o\`u les $A_{ i',\, k'}'$ sont des applications holomorphes locales.
Nous renvoyons au Th\'eor\`eme~1.2 de~\cite{ dm2002} pour une
d\'emonstration \'epur\'ee de la convergence de $h$ sous cette
hypoth\`ese.

Enfin, dans la Section~6, nous exposerons la d\'emonstration
non publi\'ee du Th\'eor\`eme~1.39 {\bf (iii)}
contenue dans~\cite{ me1999}. 

\subsection*{ 1.42.~Absence de finitude relative dans le cas
analytique r\'eel} Il y a une explication au fait que la seconde
collection d'identit\'es de r\'eflexion n'est g\'en\'eralement pas
consid\'er\'ee dans les r\'ef\'erences pr\'ecit\'ees. En effet, les
hypoth\`eses de non-d\'eg\'en\'erescence du type {\bf (h1)}, {\bf
(h2)}, {\bf (h3)} ou {\bf (h4)} permettent toutes de ramener
l'infinit\'e des identit\'es de r\'eflexion \`a un nombre fini de
relations de d\'ependance entre les composantes $h_{ i'} (t)$ et un
jet d'ordre fini $J_\tau^{\ell_0} \overline{ h} (\tau)$ de
l'application formelle conjugu\'ee. L'essence m\^eme du principe de
r\'eflexion de Schwarz \`a une ou plusieurs variables complexes se
joue dans la possibilit\'e d'exprimer $h$ en fonction de $\overline{
h}$\,--\,cela est bien connu. Mais lorsque $M'$ n'est ni finiment
non-d\'eg\'en\'er\'ee ni essentiellement finie, une telle r\'esolution
est impossible~: aucune version du th\'eor\`eme des fonctions
implicites n'est valide lorsque le morphisme des $k$-jets de
sous-vari\'et\'es de Segre~\thetag{ 1.8} est d'une complexit\'e
arbitraire. En g\'en\'eral, il est donc n\'ecessaire de consid\'erer
l'infinit\'e des s\'eries $\Theta_{ j',\, \gamma'}'( h(t) )$.

Il y a pourtant un cas o\`u une r\'esolution finie est possible~:
c'est lorsque $M'$ est alg\'ebrique. En fait, dans~\cite{ mi2002},
N.~Mir obtient le Th\'eor\`eme~1.23 pour les
applications CR-dominantes en supposant $M'$ alg\'ebrique, car dans ce
cas (subrepticement simplifi\'e), toutes les s\'eries $\{
\Theta_{j',\, \gamma'}' (t' )\}_{1\leq j' \leq d',\, \gamma' \in \N^{
m'} }$ sont alors {\it alg\'ebriquement d\'ependantes} par rapport \`a
un nombre {\it fini} d'entre elles. Gr\^ace \`a cette propri\'et\'e
cruciale de finitude, on peut obtenir des identit\'es de r\'eflexion
analogues \`a~\thetag{ 1.41}~: pour toute composante $\Theta_{ j',\,
\gamma'}' ( h(t))$ de l'application de r\'eflexion, il existe une
relation polynomiale, satisfaite pour $(t,\, \tau) \in \mathcal{ M}$~:
\def\theequation{1.43}\begin{equation}
\sum_{ 0\leq k' \leq N_0'} \, 
A_{ k'} ' \left(
t,\, \tau,\, 
J_\tau^{ \ell_0} \overline{ h} (\tau)
\right) 
\, 
\left[
\Theta_{ j',\, \gamma'}' (h(t))
\right]
= 0, 
\end{equation}
o\`u chaque $\ell_0$, chaque $N_0'$ et chaque s\'erie analytique $A_{
k'}'$ d\'epend de $j'$ et de $\gamma'$.

Au contraire, dans la cat\'egorie analytique r\'eelle, il est
absolument faux qu'\'etant donn\'e un nombre infini de s\'eries
enti\`eres convergentes $\varphi_k ( {\sf x} ) \in \C \{ { \sf x} \}$,
$k=1,\, 2,\, 3,\, \dots,\, \infty$, ${\sf x} \in \C^n$, $n\geq 1$, il
en existe un entier $N\geq 1$ tel que pour tout $k \geq N+1$, il
existe une application holomorphe locale $G_k ({\sf x},\, {\sf X}_1,\,
\dots,\, {\sf X}_N,\, Y_k)$ telle que $G_k \left( {\sf x},\, h_1({\sf
x}),\, \dots,\, h_N ({\sf x} ),\, h_k( {\sf x})\right)$ $\equiv 0$. Ce
ph\'enom\`ene est reli\'e \`a un exemple classique d\^u \`a Osgood,
aux de travaux de A.M.~Gabrielov, de E.~Bierstone, P.D.~Milman, de
B.~Malgrange~; il exhibe une diff\'erence majeure avec la
g\'eom\'etrie alg\'ebrique locale. Aussi, l'utilisation des conditions
de non-d\'eg\'en\'erescence dans les travaux ant\'erieurs \'etait-elle
cruciale, puisqu'\`a chaque \'etape de la d\'emon\-stration, on peut
substituer \`a la consid\'eration de l'{\it infinit\'e} de s\'eries
formelles $\{ \Theta_{ j',\, \gamma'}'( h(t) ) \}_{1\leq j'\leq d',\,
\gamma' \in \N^{ m'}}$ celle des $n$ composantes de $h$ seulement.
Dans ce cas, toutes les $\{ \Theta_{ j',\, \gamma'} '( h(t) )\}_{1\leq
j'\leq d',\, \gamma' \in \N^{m'}}$ sont clairement holomorphes
respectivement \`a $(h_1 (t),\, \ldots,\, h_{n'}(t ))$. Mais ici, dans
le cas g\'en\'eral, on travaillera directement avec cette collection
infinie $\{ \Theta_{ j',\, \gamma'}' (h(t)) \}_{1\leq j'\leq d',\,
\gamma' \in \N^m}$, en utilisant de mani\`ere cruciale les deux paires
d'identit\'es de r\'eflexion conjugu\'ees~\thetag{ 1.31} et~\thetag{
1.32} ({\it voir} le~\S7.61 ci-dessous).

\subsection*{ 1.44.~Applications} 
Terminons cette introduction par l'\'enonc\'e de deux applications
principales. D'apr\`es le Th\'eor\`eme~1.11, les composantes $\Theta_{
j',\, \gamma '} '( h(t))$ de l'application de r\'eflexion associ\'ee
\`a une \'equivalence formelle sont convergentes~: elles s'identifient
\`a des s\'eries convergentes $\theta_{ j',\, \gamma'} '( t) \in \C \{
t\}$. En appliquant le th\'eor\`eme d'approximation de M.~Artin aux
\'equations analytiques $\Theta_{j',\, \gamma '}' (h(t)) - \theta_{
\gamma '}' (t) \equiv 0$ satisfaites par l'application formelle $h(t)
\in \C \dl t \dr^{n'}$, on d\'eduit l'existence d'une application
convergente ${\sf H}(t) \in \C\{ t \}^{ n'}$ telle que $\Theta_{ j',\,
\gamma '}' ({\sf H} (t)) - \theta_{ j',\, \gamma '}' (t) \equiv
0$. Dans le \S7.146, nous v\'erifierons qu'une telle application ${\sf
H} (t)$ \'etablit un biholomorphisme local entre $M$ et $M'$. Ainsi~:

\def\thecorollary{1.45}\begin{corollary}
Deux sous-vari\'et\'es de $\C^n$ analytiques r\'eelles g\'en\'eriques
et minimales sont formellement \'equivalentes si et seulement si elles
sont biholomorphes.
\end{corollary}

Ce r\'esultat a \'et\'e obtenu dans~\cite{ ber1997} par S.M.~Baouendi,
P.~Ebenfelt et L.P.~Rothschild avec l'hypoth\`ese simple de
non-d\'eg\'en\'erescence finie. Ici, nous l'obtenons sans hypoth\`ese
de non-d\'eg\'en\'erescence sur $M'$, mais en utilisant fortement la
minimalit\'e. Nous pensons qu'il devrait \^etre vrai sans aucune
hypoth\`ese sur les sous-vari\'et\'es $M$ et $M'$ formellement
\'equivalentes, except\'e le fait qu'elles sont g\'en\'eriques.

Gr\^ace au m\^eme argument d'approximation, on peut aussi d\'eduire du
Th\'eor\`eme~1.23 que pour tout entier $N \geq 1$, il existe une une
application ${\sf H}^N (t) \in \C\{ t \}^N$ dont la s\'erie de Taylor
co\"{\i}ncide avec celle de $h$ jusqu'\`a l'ordre $(N-1)$ compris,
telle que ${\sf H}^N (t)$ \'etablit une application holomorphe locale
de $M$ \`a valeurs dans $M'$ ({\it voir} le Corollaire~7.147).

Enfin, dans le \S2.20 nous \'etablirons le corollaire suivant qui donne
un crit\`ere g\'en\'eral pour la convergence de $h$.
\def\thecorollary{1.46}\begin{corollary}
Sous les hypoth\`eses du Th\'eor\`eme~1.23, supposons de plus qu'il
existe des entiers $j'(1),\, \dots,\, j'(n')$ tels que $1\leq j'(
i_1') \leq n'$ pour $i_1'= 1,\, \dots,\, n'$ et des multiindices
distincts $\gamma' ( 1),\, \dots,\, \gamma '( n') \in \N^{ m'}$ tels
que le d\'eterminant
suivant~: 
\def\theequation{1.47}\begin{equation}
{\rm det} \
\left(
\frac{ \partial \Theta_{j'(i_1'), \ \gamma' (i_1')}'}{
\partial t_{ i_2'}'} (h(t))
\right)_{1\leq i_1',\, i_2'\leq n'}
\end{equation}
ne s'annule pas identiquement dans $\C \dl t \dr$. 
Alors $M'$ est holomorphiquement non-d\'eg\'en\'e\-r\'ee
et $h$ est convergente. 
\end{corollary} 

Notre premier Th\'eor\`eme~1.11 d\'ecoule en v\'erit\'e comme
cons\'equence directe du Th\'eor\`eme~1.23 et de ce corollaire ({\it
voir}~\S2.20).

\subsection*{ 1.48.~Remarque finale}
Le lecteur aura remarqu\'e que cette introduction contient de
nombreuses r\'ef\'erences aux travaux de S.M.~Baouendi, de
P.~Ebenfelt, de F.~Meylan, de N.~Mir, de L.P.~Rothschild et de
D.~Zaitsev. Durant la p\'eriode 1998--2004, ces auteurs ont
r\'eguli\`erement suivi l'\'evolution de nos travaux, publi\'es dans
des revues sp\'ecialis\'ees ou pr\'epubli\'es \'electroniquement. Tous
nos travaux sur les applications CR font r\'ef\'erence \`a leurs
travaux. En revanche, pour la p\'eriode 1998--2004, force
est de constater qu'il n'existe qu'une seule publication de ces
auteurs, group\'ee ou individuelle, dont la {\it bibliographie}\,
contienne une r\'ef\'erence \`a l'un de nos travaux~: il s'agit du
livre~\cite{ ber1999a}, qui cite notre travail de th\`ese paru en
1997, portant sur les singularit\'es \'eliminables pour les
fonctions CR (un autre sujet de recherche) ainsi que l'article~\cite{
mm1999}, \'ecrit en 1997. Insistons sur le fait que cette absence de
citation est constatable non seulement dans les travaux que nous
citons ici, mais aussi dans {\it tous les autres travaux {\rm
(}pr\'e{\rm )}publi\'es par ces auteurs durant cette p\'eriode}. Les
travaux de N.~Mir, tr\`es proches des n\^otres, constituent le cas le
plus frappant d'absence de citation bibliographique. Pour cette
raison, nous nous devions de d\'etailler dans cette introduction la
chronologie pr\'ecise de l'apparition des r\'esultats r\'ecents sur
les applications CR formelles.

\subsection*{ 1.49.~Remerciement}
Je remercie vivement F.~Panigeon pour ses relectures minutieuses
sur \'ecran.

\section*{\S2.~Pr\'eliminaire~: 
s\'eries formelles, analytiques et alg\'ebriques}

\subsection*{2.1.~S\'eries formelles, analytiques, alg\'ebriques}
Dans ce paragraphe liminaire et \'el\'ementaire, destin\'e seulement
\`a fixer fermement nos notations et \`a pr\'esenter le th\'eor\`eme
d'approximation de M.~Artin, la lettre $\K$ d\'esigne ou bien le corps
$\R$ des nombres r\'eels, ou bien le corps $\C$ des nombres
complexes. Soit $n\in \N$ et ${\sf x} :=({\sf x}_1, \dots,\, {\sf
x}_n) \in \K^n$ des ind\'etermin\'ees. Soit $\K \dl {\sf x} \dr$
l'anneau local des s\'eries formelles en les variables $({\sf x}_1,
\dots,\, {\sf x}_n)$. Par d\'efinition, un \'el\'ement $\varphi( {\sf
x}) \in \K\dl {\sf x} \dr$ s'\'ecrit sous la forme $\varphi( {\sf x})
= \sum_{\alpha \in \N^n}\, \varphi_\alpha \, {\sf x }^\alpha$, o\`u
${\sf x}^\alpha$ est le mon\^ome ${\sf x }_1^{ \alpha_1 } {\sf x }_2^{
\alpha_2} \cdots {\sf x }_n^{ \alpha_n}$ et o\`u les coefficients
$\varphi_\alpha$, du reste arbitraires, appartiennent \`a $\K$ pour
tout multiindice $\alpha := (\alpha_1, \dots,\, \alpha_n)\in
\N^n$. Une telle s\'erie est {\sl identiquement nulle} si tous ses
coefficients $\varphi_\alpha$ sont nuls. Nous \'ecrirons cette
propri\'et\'e $\varphi( {\sf x}) \equiv 0$ (dans $\K\dl {\sf x}\dr$).
Cette relation, bien que triviale, sera fr\'equemment utilis\'ee dans
ce m\'emoire. Lorsque $\K = \C$, la variable ${\sf x}$ et les
coefficients $\varphi_\alpha$ sont complexes et on d\'efinit
$\overline{ \varphi} ( {\sf x} ) := \sum_{ \alpha \in \N^n} \,
\overline{ \varphi }_\alpha \, {\sf x }^\alpha$ en ne conjuguant que
les coefficients de la s\'erie, de telle sorte que l'on a $\overline{
\varphi ( {\sf x} ) } \equiv \overline{ \varphi } (\overline{ {\sf x}
})$, la barre de conjugaison se distribuant sur la s\'erie
et sur la variable.

La {\sl longueur} du multiindice $\alpha$ est l'entier
$\vert \alpha \vert:= \alpha_1+ \cdots + \alpha_n$. La d\'eriv\'ee
partielle correspondante sera not\'ee $\partial_{\sf x}^\alpha :=
\partial_{ { \sf x }_1 }^{ \alpha_1} \partial_{ { \sf x }_2 }^{
\alpha_2 } \cdots \partial_{ { \sf x }_n }^{ \alpha_n}$ et parfois
$\partial^{ \vert \alpha \vert}\varphi ( {\sf x}) / \partial
{\sf x}_1^{\alpha_1} \partial {\sf x}_2^{ \alpha_2 } \cdots \partial
{\sf x }_n^{ \alpha_n}$. \'Evidemment, on a $\varphi_\alpha = [1 /
\alpha!] \, \partial_{\sf x }^\alpha \varphi( {\sf x}) \vert_{ { \sf
x} =0}$, o\`u le symbole $\alpha! := \alpha_1 ! \alpha_2! \cdots
\alpha_n!$ est le produit des factorielles des $\alpha_j$.

Sur $\K^n$, la norme $\vert {\sf x} \vert := \max \, \{ \vert {\sf
x}_1\vert ,\, \vert{\sf x}_2 \vert,\dots,\, \vert {\sf x}_n \vert \}$
sera la plus commode. Si les coefficients satisfont une estim\'ee de
Cauchy de la forme $\vert \varphi_\alpha \vert \leq C \rho^{ -\vert
\alpha \vert}$, o\`u $C >0$ et $\rho>0$, on dira que la s\'erie
formelle $\varphi( {\sf x})$ {\sl converge normalement} dans le {\sl
cube ouvert} $\square_n (\rho) := \{{\sf x} \in \K^n: \, \vert x \vert
< \rho\}$. Bien entendu, sous cette condition, une valeur num\'erique
$\varphi({\sf x})\in \K$ peut \^etre assign\'ee univoquement \`a
$\varphi$ en tout point ${\sf x}\in \square_n(\rho)$. On dira que
$\varphi$ est {\sl $\K$-analytique} et on \'ecrira $\varphi \in \K \{
{\sf x} \}$, les constantes $\rho$ et $C$ \'etant de peu d'importance
pour les probl\`emes que nous \'etudierons. Si de plus il existe un
polyn\^ome {\it non nul}\, $P({\sf X}_1,\dots,{\, \sf X}_n,\, \Phi)
\in \K[{\sf X}_1, \dots,\, {\sf X}_n,\, \Phi]\backslash \{0\}$ tel que
$P({\sf x}_1,\dots,\, {\sf x}_n ,\, \varphi({\sf x}_1, \dots,\, {\sf
x}_n))\equiv 0$, on dira que $\varphi$ est {\sl $\K$-alg\'ebrique} (au
sens de J.~Nash) et on \'ecrira $\varphi( {\sf x}) \in \mathcal{ A}_\K
\{{\sf x}\}$. \'Evidemment, on a les deux inclusions strictes
\def\theequation{2.2}\begin{equation}
\K\dl {\sf x} \dr \supset
\K \{ {\sf x} \} \supset
\mathcal{ A}_\K \{ {\sf x} \}.
\end{equation}
Ces trois ensembles $\K \dl {\sf x} \dr$, $\K \{ {\sf x} \}$ et
$\mathcal{ A}_\K \{ {\sf x} \}$ sont des anneaux dits de {\sl s\'eries
enti\`eres} qui sont locaux, noth\'eriens, factoriels et qui satisfont
les th\'eor\`emes de pr\'eparation et de division de
K.~Weierstrass. Ils sont stables par composition et par
diff\'erentiation~; le th\'eor\`eme des fonctions implicites y est
valide.

\subsection*{2.3.~Application formelles} 
Soient $n$ et $n'$ deux entiers strictement positifs. Une {\sl
application formelle de $\K^n$ dans $\K^{n'}$} consiste en la donn\'ee
d'un $n'$-uplet de s\'eries enti\`eres formelles $\varphi ({\sf x}) :=
(\varphi_1 ({\sf x}),\,\dots,\, \varphi_{ n'} ( {\sf x}))$ appartenant
\`a $\K \dl {\sf x} \dr$ et sans terme constant, {\it i.e.}
satisfaisant $\varphi_{i'} (0) = 0$. Une telle application est dite

\begin{itemize}
\item[{\bf (1)}]
{\sl inversible} si $n' = n$ et si ${\rm det}\, 
([\partial \varphi_{i_1}/ \partial {\sf x}_{ i_2}](0))_{1 \leq i_1,\,
i_2\leq n}\neq 0$~;
\item[{\bf (2)}]
{\sl submersive} si $n' \leq n$ et s'il existe des entiers $1\leq i(1)
< \cdots < i( n')\leq n$ tels que ${\rm det }\, ([ \partial 
\varphi_{ i_1'}
/ \partial {\sf x}_{i(i_2') }] (0))_{ 1 \leq i_1',\, i_2'\leq n'}\neq
0$~;
\item[{\bf (3)}]
{\sl finie} si l'id\'eal engendr\'e par les composantes $\varphi_1(
{\sf x}),\, \dots,\, \varphi_{n'}( {\sf x})$ est de codimension finie
dans $\K \dl {\sf x} \dr$, ce qui implique que $n' \geq n$~;
\item[{\bf (4)}]
{\sl dominante} si $n' \leq n$ et s'il existe des entiers $1\leq i(1)
< \cdots < i( {n'}) \leq n$ tels que le d\'eterminant ${\rm
det}([\partial \varphi_{ i_1'}/ \partial {\sf x }_{ i(i_2' )}] (
{\sf x}))_{1
\leq i_1',\, i_2' \leq n'} \not \equiv 0$ ne s'annule pas
identiquement dans $\K \dl {\sf x} \dr$~;
\item[{\bf (5)}]
{\sl transversale} s'il n'existe pas de s\'erie formelle $F' ({\sf
x}_1',\, \dots,\, { \sf x}_{n'}') \in \C \dl {\sf x}_1',\, \dots,\,
{\sf x}_{n'}'\dr$ non nulle telle que $F '(\varphi_1( {\sf x}),\,
\dots,\, \varphi_{n'} ( {\sf x}))\equiv 0$ dans $\K \dl {\sf x} \dr$.
\end{itemize}

On d\'emontre de mani\`ere \'el\'ementaire que ces conditions sont
ordonn\'ees par ordre croissant de g\'en\'eralit\'e ({\it voir} \cite{
me2003}). La derni\`ere condition, de loin la plus g\'en\'erale,
n'implique aucune in\'egalit\'e entre les dimensions $n$ et $n'$.

\subsection*{2.4.~Approximation} 
Le principal outil (non trivial) de g\'eom\'etrie analytique que nous
utiliserons {\it ad nauseam}\, \'enonce que la s\'erie de Taylor de
toute solution purement formelle d'\'equations $\K$-analytiques peut
\^etre corrig\'ee \`a l'infini de mani\`ere \`a la rendre convergente,
de telle sorte que la s\'erie modifi\'ee demeure solution des
\'equations analytiques donn\'ees.

\def\thetheorem{2.5}\begin{theorem}
{\rm ({\sc M.~Artin}~: \cite{ ar1968}, \cite{ ar1969})} 
Soit $n \in \N$ avec $n \geq
1$, soit ${\sf x} = ({\sf x}_1, \dots,\, {\sf x}_n) \in \K^n$, soit $m
\in \N$, avec $m \geq 1$, soit ${\sf y} =({\sf y}_1, \dots ,\, {\sf y
}_m ) \in \K^n$, soit $d \in \N$ avec $d \geq 1$ et soit $R_1 ( {\sf
x}, {\sf y }), \dots, \, R_d({\sf x},\, {\sf y})$ une collection
arbitraire de s\'eries enti\`eres appartenant \`a $\K \{{\sf x},\,
{\sf y} \}$ ou \`a $\mathcal{A }_\K \{ {\sf x} ,\, {\sf y} \}$ qui
s'annulent \`a l'origine, c'est-\`a-dire $R_j(0,\, 0) =0$ pour $j=1,
\dots,\, d$. Supposons qu'il existe une application formelle $h( {\sf
x}) = (h_1( {\sf x }), \dots,\, h_m( {\sf x})) \in \K\dl {\sf x}
\dr^m$ avec $h(0) =0$ telle que
\def\theequation{2.6}\begin{equation}
R_j\left({\sf x} ,\, h({\sf x})\right)\equiv 0 \ \ 
{\rm dans} \ \ \K\dl {\sf x} \dr, \ \ \ \
{\rm pour} \ \
j=1,\dots,\, d.
\end{equation}
Soit $\mathfrak{m}( {\sf x}) := {\sf x}_1 \K\dl {\sf x} \dr + \cdots+
{\sf x}_n \K \dl {\sf x} \dr$ l'id\'eal maximal de $\K \dl {\sf x}
\dr$. Pour tout entier $N \geq 1$, il existe une s\'erie enti\`ere
convergente $h^N( {\sf x})$ qui appartient \`a 
$\K\{ {\sf x} \}^m$ ou \`a $\mathcal{A}_\K \{{\sf
x}\}^m$ telle que
\def\theequation{2.7}\begin{equation}
R_j\left({\sf x},\, h^N({\sf x})\right)\equiv 0 \ \ 
{\rm dans} \ \ \K \dl {\sf x} \dr, \ \ \ \
{\rm pour} \ \ 
j=1,\dots,\, d,
\end{equation}
et qui approxime $h({\sf x})$ \`a l'ordre $N-1$, 
c'est-\`a-dire qui satisfait
\def\theequation{2.8}\begin{equation}
h^N({\sf x}) \equiv h({\sf x}) \ \
{\rm mod} \, \left(\mathfrak{m}(x)^N \right).
\end{equation}
\end{theorem}

Dans la d\'emonstration du th\'eor\`eme principal~1.23, nous verrons
tr\`es fr\'equemment appara\^{\i}tre des solutions formelles d'un
nombre infini d'\'equations analytiques qu'il faudra transformer en
solutions convergentes. Heureusement, la consid\'eration d'un nombre
fini d'\'equations analytiques $R_j({\sf x},\, {\sf y})=0$, $j=1,
\dots,\,d$ dans le th\'eor\`eme de M.~Artin n'est en rien
restrictive. En effet, s'il l'on se donne au contraire un nombre
infini de telles \'equations $R_j({\sf x},\, {\sf y})=0$, pour $j=1,
\,2, \,3, \dots,\, \infty$, gr\^ace \`a la noeth\'erianit\'e de $\K\{
{\sf x},\, {\sf y}\}$ ou de $\mathcal{ A}_\K \{ {\sf x},\, {\sf y}
\}$, il existe un entier $d\geq 1$ tel que l'id\'eal engendr\'e par
tous les $R_j$ co\"{\i}ncide avec l'id\'eal engendr\'e par $R_1,\,
R_2, \dots,\, R_d$. En d'autres termes, pour tout entier $l=1,\,2,\,
\dots,\, \infty$ et tout $j =1,\, \dots,\, d$, il existe des
coefficients $\lambda_{ l,\, j} ( {\sf x},\, {\sf y})$ appartenant \`a
$\K\{ {\sf x},{\sf y}\}$ ou \`a $\mathcal{ A}_\K( {\sf x},\, {\sf y})$
tels que $R_l ({\sf x},\, {\sf y}) \equiv \sum_{ j=1 }^d\, \lambda_{
l,\, j}( {\sf x},\,{\sf y})\, R_j( {\sf x},\,{\sf y})$. On d\'eduit
imm\'ediatement que les \'equations finies $R_1\left( {\sf x},\, h({
\sf x}) \right) \equiv \cdots \equiv R_d \left( {\sf x},\, h({ \sf x})
\right) \equiv 0$ sont satisfaites si et seulement si les \'equations
infinies $R_l( {\sf x},\, h({\sf x})) \equiv 0$ pour $l = 1,\, 2,\dots
,\, \infty$ le sont~; on d\'eduit de m\^eme qu'une solution
convergente $h^N ( {\sf x}) \in \K\{ {\sf x} \}$ ou $h^N ({\sf x}) \in
\mathcal{ A}_\K \{ {\sf x}\}$ des \'equations finies $R_1\left( {\sf x
},\, h^N({\sf x}) \right) \equiv \cdots \equiv R_d\left( {\sf x} ,\,
h^N ( {\sf x}) \right) \equiv 0$ est automatiquement une solution des
\'equations infinies $R_l \left( {\sf x},\, h^N({\sf x}) \right)
\equiv 0$, $l=1, \,2, \dots,\,\infty$. En d\'efinitive, le
Th\'eor\`eme~2.5 est tout aussi valide pour un nombre infini
d'\'equations analytiques.

Notons qu'il n'est pas possible d'appliquer directement le
Th\'eor\`eme~2.5 pour d\'emontrer le Corollaire~1.45 ou le corollaire
plus g\'en\'eral du Th\'eor\`eme~1.23 cit\'e avant l'\'enonc\'ee du
Corollaire~1.46. En effet, sous les hypoth\`eses dudit th\'eor\`eme,
l'application CR formelle complexifi\'ee $h^c:\, ( \mathcal{ M},\, 0)
\longrightarrow_{ \mathcal{ F}} (\mathcal{ M}',\, 0)$ satisfait les
$d'$ \'equations analytiques complexes de la premi\`ere ligne
de~\thetag{ 1.27}, dont les variables sont ${\sf x} := (\zeta,\, t)
\in \C^{ m+n}$. En appliquant le Th\'eor\`eme~2.5, on trouve pour tout
$N$ une application $(H^N(\zeta,\, t),\, \Phi^N(\zeta,\, t))\in \C\{
\zeta,\, t\}^{ 2n}$ approximant $h^c$ jusqu'\`a l'ordre $(N-1)$ et
satisfaisant les m\^emes \'equations analytiques. Malheureusement,
rien ne permet de s'assurer que la variable $\zeta$ n'appara\^{\i}t
pas dans $H^N$ et plus encore, que $\Phi^N = \overline{ H}^N$. Toutes
les diff\'erentiations que nous effectuerons dans les Sections~5, 6 et
7 auront pour v\'eritable but de s\'eparer les variables $t$ de leurs
conjugu\'ees complexifi\'ees $\tau$.

\'Enon\c cons maintenant un corollaire direct (et connu) du
Th\'eor\`eme~2.5 qui fournit un crit\`ere puissant pour la convergence
d'applications formelles et rappelons-en une d\'emonstration,
\'el\'ementaire parmi d'autres.

\def\thecorollary{2.9}\begin{corollary}
Sous les hypoth\`eses du Th\'eor\`eme~2.5, si $d=m$ et si
le d\'eterminant
\def\theequation{2.10}\begin{equation}
{\rm det}\, 
\left(
{\partial R_j \over \partial {\sf y}_k}
\left({\sf x} ,\, h({\sf x})\right)
\right)_{1\leq j,\, k\leq m}
\end{equation}
ne s'annule pas identiquement dans
$\K \dl {\sf x} \dr$, la s\'erie formelle $h( {\sf x})$ est
en fait convergente analytique ou alg\'ebrique, {\it i.e.} appartient
\`a $\K\{ {\sf x}\}^m$ ou \`a $\mathcal{ A}_\K \{ {\sf x} \}^m$.
\end{corollary}

\proof
Gr\^ace \`a la formule de Taylor sous forme int\'egrale, on trouve
ais\'ement des s\'eries enti\`eres $S_{j, \, k} ({\sf x},\, {\sf y},\,
{\sf y}')$ analytiques ou alg\'ebriques satisfaisant les formules de
division
\def\theequation{2.11}\begin{equation}
R_j\left( {\sf x},\, {\sf y}\right) -R_j\left(
{\sf x},\, {\sf y}' \right)\equiv 
\sum_{k=1}^m\, 
S_{j,k}\left( {\sf x},\,
{\sf y},\, {\sf y}'\right)
\, \left[ {\sf y}_k-{\sf y}_k'\right],
\end{equation}
pour $j=1,\dots, \, d$. Notons qu'en faisant tendre ${\sf y}_k'$ vers
${\sf y}_k$, il est clair qu'on obtient
\def\theequation{2.12}\begin{equation}
S_{j,\, k}({\sf x},\, 
{\sf y},\,{\sf y})
\equiv [\partial R_j/\partial
y_k]({\sf x},\, {\sf y}). 
\end{equation}
Gr\^ace au Th\'eor\`eme d'approximation~2.5, on trouve pour tout
entier $N \geq 1$ des s\'eries enti\`eres convergentes analytiques ou
alg\'ebriques $h^N ( {\sf x})$ qui sont solutions des \'equations $R_j
\left( {\sf x} ,\, h^N ( {\sf x } ) \right) \equiv 0$, pour $j =1,
\dots,\, m$, et qui satisfont $h^N( {\sf x }) \equiv h ({\sf x}) \ \
{\rm mod} \, \left( \mathfrak{ m}( {\sf x })^N \right)$.

Rempla\c cons maintenant ${\sf y}$ par $h({\sf x})$ et ${\sf y}'$ par
$h^N \left( {\sf x} \right)$ dans les \'equations~\thetag{ 2.11}, ce
qui donne
\def\theequation{2.13}\begin{equation}
\left\{
\aligned
0
& \
\equiv R_j\left( {\sf x},\, h({\sf x})\right)- R_j\left({\sf x} ,\,
h^N({\sf x})\right)
\\
& \
\equiv
\sum_{k=1}^m\, 
S_{j,k}\left(
{\sf x},\,h({\sf x}),\,h^N({\sf x})\right)\, 
[h_k({\sf x})-h_k^N({\sf x})],
\endaligned\right.
\end{equation}
pour $j=1,\dots, \, m$. En tenant compte de l'identit\'e~\thetag{
2.12}, l'hypoth\`ese que le d\'eterminant~\thetag{ 2.10} ne s'annule
pas identiquement s'\'ecrit de mani\`ere \'equivalente
\def\theequation{2.14}\begin{equation}
{\rm det} 
\left(
S_{j,k}\left({\sf x},\,h({\sf x}),\,h({\sf x})\right) 
\right)_{1\leq j,\, k\leq m}\not\equiv 0 \ \
{\rm dans} \ \K \dl {\sf x} \dr.
\end{equation}
Il en d\'ecoule ais\'ement que si $N$ est suffisamment grand,
le d\'eterminant
\def\theequation{2.15}\begin{equation}
{\rm det}\, 
\left(
S_{j,\, k}\left({\sf x},\,h({\sf x}),\,h^N({\sf x})\right) 
\right)_{1\leq j,\, k\leq m}
\end{equation}
ne s'annule (lui non plus) pas identiquement dans $\K \dl {\sf
x}\dr$. Finalement, en interpr\'etant~\thetag{ 2.13} comme un
syst\`eme lin\'eaire homog\`ene dont les inconnues sont les $h_k (
{\sf x}) - h_k^N ({\sf x})$, pour $k=1, \dots ,\, m$, la
non-annulation du d\'eterminant~\thetag{ 2.15} implique
imm\'ediatement que $h_k( {\sf x}) \equiv h_k^N({\sf x})$ pour
$k=1,\dots ,\, m$. En conclusion, l'application formelle $h( {\sf x
})\equiv h^N( {\sf x})$, qui s'identifie donc \`a l'une des
applications {\it convergentes}\, qui sont fournies par
l'approximation de M.~Artin pour $N$ assez grand, est effectivement
convergente, analytique ou alg\'ebrique.
\endproof

\subsection*{ 2.16.~\'Equivalences CR formelles non convergentes
entre sous-vari\'et\'es holomorphiquement d\'eg\'en\'er\'ees} Dans ce
paragraphe, nous v\'erifions l'assertion faite au \S1.10 de
l'Introduction. Soit donc $(s',\, t') \longmapsto \exp( s' X')(t')=:
\varphi'( s',\, t')$ le flot local d'un champ de vecteurs holomorphe
$X' = \sum_{ i'= 1}^{ n'} \, a_{i'}'(t')\, \frac{ \partial }{\partial
t_{i'}'}$ non nul tangent \`a une sous-vari\'et\'e $M'$ de $\C^{ n'}$
analytique r\'eelle, g\'en\'erique et passant par l'origine. Par
d\'efinition, cette application holomorphe d\'efinie au voisinage de
$(0,\, 0)$ dans $\C\times \C^{n'}$ est uniquement d\'etermin\'ee par
la condition initiale $\varphi '( 0,\, t') \equiv t'$ et par le
syst\`eme d'\'equations diff\'erentielles ordinaires $\partial_{ s'}
\, \varphi_{ i'}'( s', \, t') \equiv a_{ i'}' (\varphi'( s',\, t'))$,
pour $i'= 1,\, \dots,\, n'$.

Supposons $M'$ repr\'esent\'ee par les \'equations analytiques
r\'eelles $\rho_{ j'}' \left( t',\, \bar t' \right) = 0$, pour $j' =
1,\, \dots,\, d'$. La condition de tangence de $X'$ \`a $M'$ implique
\'evidemment que le flot de $X'$ stabilise $M'$, c'est-\`a-dire qu'il
existe une matrice inversible de taille $d' \times d'$ de s\'eries
enti\`eres convergentes $c'(s',\, t',\, \bar s',\, \bar t')$ telle que
l'on a l'identit\'e vectorielle
\def\theequation{2.17}\begin{equation}
\rho ' \left( \varphi' ( s',\, t'), \, \overline{ \varphi}' ( \bar
s',\, \bar t')\right) \equiv c' \left( s',\, t',\, 
\bar s',\, \bar t'\right) \,
\rho' \left( t',\, \bar t' \right).
\end{equation}
De plus, comme $X'$ ne s'annule pas identiquement, apr\`es une
renum\'erotation \'eventuelle, on peut supposer que $a_1' (t') \not
\equiv 0$. Il en d\'ecoule que $\partial_{ s'} \, \varphi_1' (0,\,
t') \equiv
a_1'\left(\varphi'(0,\, t')\right) 
\not \equiv 0$ dans $\C\{ t'\}$. Autrement dit,
dans le d\'eveloppement en s\'erie enti\`ere de
$\partial_{ s'} \, \varphi_1' (s',\, t')$ par rapport aux puissances
de $s'$, qui s'\'ecrit 
$\sum_{ k = 0}^\infty \, (s')^k \, \varphi_{1,\, k}'
(t')$, on a
$\varphi_{ 1,\, 0}'( t') \not \equiv 0$. On v\'erifie alors
facilement que pour toute s\'erie formelle $\varpi' (t') \in \C \dl t'
\dr$ dont l'ordre d'annulation est suffisamment \'elev\'e, {\it i.e.}
$\varpi'( t') \in [\mathfrak{ m} (t')]^N$ pour $N$ assez grand, la
composition $\partial_{ s'} \varphi_1'( \varpi'( t'),\, t')$ ne
s'annule pas identiquement dans $\C \dl t' \dr$. Le lemme suivant
ach\`eve d'\'etablir l'\'enonc\'e d\'esir\'e.

\def\thelemma{2.18}\begin{lemma}
Pour toute s\'erie enti\`ere $\varpi ' ( t' ) \in \C \dl t' \dr$ {\rm
non convergente} qui satisfait $\partial_{ s'} \varphi_1' \left(
\varpi'( t'),\, t' \right) \not \equiv 0$, l'application $t' \mapsto_{
\mathcal{ F}} \, \varphi ' \left( \varpi '( t'),\, t' \right)$ induit
une auto-application CR formelle de $M'$ qui n'est pas convergente.
\end{lemma}

\proof
En effet, si l'on remplace $s'$ par $\varpi' ( t')$ dans~\thetag{
2.17}, on voit imm\'ediatement que $t' \longmapsto_{\mathcal{ F}} \,
\varphi '( \varpi' (t'),\, t')$ induit une auto-application CR
formelle de $M'$.

Supposons par l'absurde qu'elle est convergente. En particulier, sa
premi\`ere composante $\varphi_1'( \varpi' (t' ),\, t')$ s'identifie
\`a une s\'erie convergente. Notons-la $\alpha_1' (t') \in\C \{ t'\}$
et consid\'erons l'identit\'e formelle
\def\theequation{2.19}\begin{equation} 
\varphi' \left( \varpi'( t'),\, t'\right) -
\alpha_1 '(t') \equiv 0,
\end{equation}
qui exprime que $\varpi' (t')$ est solution d'\'equations analytiques.
Puisque l'on suppose que $\partial_{ s'} \varphi_1 ' \left( \varpi'
(t'),\, t'\right) \not \equiv 0$, l'hypoth\`ese principale du
Corollaire~2.9 est exactement satisfaite, ce qui implique que $\varpi'
( t')$ est convergente. Cette contradiction conclut le raisonnement
par l'absurde. En conclusion, l'auto-application CR formelle $t'
\mapsto_{ \mathcal{ F}} \, \varphi '\left( \varpi '( t'),\, t'
\right)$ n'est pas convergente.
\endproof

\subsection*{ 2.20.~D\'emonstrations du Corollaire~1.46 et
du Th\'eor\`eme~1.11} D'apr\`es le Th\'eor\`eme 1.23, toutes les
composantes $\Theta_{j',\, \gamma'} '( h(t) )$ s'identifient \`a des
s\'eries enti\`eres convergentes $\theta_{ j',\, \gamma'} ' (t) \in
\C\{ t\}$, pour $j'= 1,\, \dots,\, d'$ et $\gamma ' \in
\N^{m'}$. Consid\'erons les \'equations analytiques suivantes, en
nombre infini~:
\def\theequation{2.21}\begin{equation}
\Theta_{ j',\, \gamma'}' (h(t)) - 
\theta_{ j',\, \gamma '}' (t) \equiv 0~;
\end{equation}
elles sont satisfaites par l'application formelle $h(t)$. Les
hypoth\`eses du Corollaire~1.46 sont exactement celles qui assurent
que le Corollaire~2.9 s'applique. Donc $h(t) \in \C\{ t\}^{ n'}$
converge. Le Corollaire~1.46 est d\'emontr\'e.

Pour \'etablir le Th\'eor\`eme~1.11, rappelons que $M'$ est
holomorphiquement non-d\'eg\'en\'e\-r\'ee si et seulement si il existe
des entiers $j'(1),\, \dots,\, j'(n')$ tels que $1\leq j'(i_1')\leq
d'$ pour $i_1 ' = 1,\, \dots,\, n'$ et des multiindices distincts
$\gamma' (1),\, \dots,\, \gamma' ( n') \in \N^{ m'}$ tels le
d\'eterminant suivant~:
\def\theequation{2.22}\begin{equation}
{\rm det} \
\left(
\frac{ \partial \Theta_{j'(i_1'), \ \gamma' (i_1')}'}{
\partial t_{ i_2'}'} ( t' )
\right)_{1\leq i_1',\, i_2'\leq n'}
\end{equation}
ne s'annule pas identiquement dans $\C \dl t' \dr$ ({\it cf.}~\cite{
st1996}, \cite{ ber1999a} ou le Lemme~3.2.49 {\bf (5)} dans \cite{
me2003}). Si $h(t)$ est inversible \`a l'origine elle est \`a la fois
CR-transversale et transversale \`a l'origine. Le Th\'eor\`eme~1.23
s'applique~: les composantes $\Theta_{j',\, \gamma'} '( h(t) )$
s'identifient \`a des s\'eries enti\`eres convergentes $\theta_{ j',\,
\gamma'} ' (t) \in \C\{ t\}$, pour $j'= 1,\, \dots,\, d'$ et $\gamma '
\in \N^{m'}$.

Puisque $h$ est transversale \`a l'origine, le d\'eterminant~\thetag{
1.47} pour les m\^emes entiers $j'( i_1')$ et les m\^emes multiindices
$\gamma ' ( i_1')$ ne peut pas s'annuler identiquement. L'hypoth\`ese
principale du Corollaire~1.46 est satisfaite, donc $h(t) \in \C\{
t\}^{ n'}$ converge. Nous en d\'eduisons m\^eme que le
Th\'eor\`eme~1.11 est vrai avec l'hypoth\`ese plus g\'en\'erale que
$h$ est CR-transversale et transversale \`a l'origine, ou avec
d'autres hypoth\`eses interm\'ediaires.

\section*{\S3.~G\'eom\'etrie locale 
des paires de feuilletages minimales}

\subsection*{3.1.~Avertissement}
Nos consid\'erations seront toujours locales, centr\'ees en un
point fixe de $\C^n$ que l'on supposera \^etre l'origine, sans perte
de g\'en\'eralit\'e. Nous \'eviterons soigneusement d'employer le
langage des germes, dont l'ambigu\"{i}t\'e fondamentale qui consiste
\`a ne pas pr\'eciser dans quels petits ouverts on travaille, loin de
simplifier les \'enonc\'es et les d\'emonstrations, entretient des
impr\'ecisions qui nuisent \`a la coh\'erence de l'ensemble et
occultent le sens g\'eom\'etrique concret des concepts, localis\'es au
fur et \`a mesure des preuves\footnotemark[3].

\footnotetext[3]{
Par exemple, le fait de travailler avec des polydisques
embo\^{\i}t\'es dont les rayons d\'ecroissants sont pr\'ecis\'es
successivement en fonction de contraintes explicites est un
ingr\'edient substantiel pour d\'emontrer de mani\`ere rigoureuse et
compl\`ete que l'ensemble des automorphismes holomorphes locaux d'une
sous-vari\'et\'e locale de $\C^n$ analytique r\'eelle, g\'en\'erique
et finiment non-d\'eg\'en\'er\'ee est un {\sl groupe de Lie local de
dimension finie}, notion qui ne s'accommode gu\`ere du langage des
germes si l'on tient \`a la relier concr\`etement \`a l'objet
g\'eom\'etrique qui est stabilis\'e ({\it voir}\, le Th\'eor\`eme~4.1
dans~\cite{ gm2004}~; le th\'eor\`eme le plus g\'en\'eral dans cette
direction obtenu auparavant par les auteurs de~\cite{ za1997}
et~\cite{ ber1999a} se limite, sans n\'ecessit\'e apparente, au
sous-groupe d'isotropie d'un point fixe donn\'e \`a l'avance).
}

Dans cette section et dans celle qui suit, nous pr\'esentons les notions
de base. Elles sont d\'evelopp\'ees en partie dans d'autres
r\'ef\'erences, mais elles y sont g\'en\'eralement expos\'ees dans un
autre langage ou d'une mani\`ere parfois incompl\`ete que nous jugeons
trop peu conceptuelle ou insuffisamment \'epur\'ee. Apr\`es une
pr\'esentation progressive et g\'en\'erale des deux objets
fondamentaux que sont la paire de feuilletages invariants (Section~3)
et les jets des sous-vari\'et\'es de Segre (Section~4), les objets
analytiques avec lesquels nous travaillerons r\'eellement seront
clairement et concr\`etement pos\'es dans les deux r\'esum\'es qui
apparaissent aux sous-sections~3.31 et 4.33.

Bien que les concepts de base poss\`edent tous un sens g\'eom\'etrique
initial, en v\'erit\'e, ce qui constituera pour nous l'essence m\^eme
du Th\'eor\`eme principal~1.23, c'est l'architecture purement
alg\'ebrique des calculs formels qui appara\^{\i}tront dans sa
d\'emonstration d\'evelopp\'ee~; c'est l'encha\^{\i}nement
structur\'e, r\'egl\'e et \'epur\'e des gestes formels~; et surtout,
c'est la {\it nature duelle des calculs}, toujours absolument
sym\'etriques par conjugaison complexe. C'est pourquoi nous
exprimerons nos calculs en d\'eployant leurs deux versions
parall\`eles d'une mani\`ere simultan\'ee. Enfin, c'est gr\^ace \`a
cette compr\'ehension interne de la sym\'etrie entre les variables
holomorphes et les variables anti-holomorphes que nous serons \`a
m\^eme de transf\'erer presque directement un large pan de la
g\'eom\'etrie CR \`a l'\'etude des syst\`emes compl\`etement
int\'egrables d'\'equations aux d\'eriv\'ees partielles analytiques
({\it voir}~\cite{ me2004}).

Ainsi, dans ce m\'emoire, sera privil\'egi\'ee la concr\'etude
explicite des calculs par rapport \`a leur possible abstraction
structurale.

\subsection*{3.2.~Sous-vari\'et\'es g\'en\'eriques analytiques 
r\'eelles de $\C^n$}
Une sous-vari\'et\'e analytique r\'eelle locale $M$ de $\C^n$ passant
par l'origine est dite {\sl g\'en\'erique} si son espace tangent {\it
g\'en\`ere}\, $\C^n$, c'est-\`a-dire $T_0 M + J T_0 M = T_0 \C^n$,
o\`u $J$ d\'esigne la structure complexe standard de $\C^n$. Il en
d\'ecoule trivialement que la codimension r\'eelle $d$ de $M$
satisfait $d \leq n$. Si $t =(t_1 ,\dots ,\, t_n) \in \C^n$ sont des
coordonn\'ees centr\'ees \`a l'origine, on peut repr\'esenter
concr\`etement $M$ par $d$ \'equations cart\'esiennes ind\'ependantes
$\rho_1 (t,\, \bar t) =0, \dots, \, \rho_d (t,\, \bar t)=0$, o\`u les
s\'eries enti\`eres convergentes $\rho_j (t ,\, \bar t) \in \C\{ t,\,
\bar t \}$ s'annulent \`a l'origine et satisfont la condition de
r\'ealit\'e $\rho_j (t ,\, \bar t) \equiv \overline{ \rho }_j( \bar t,
\, t)$. L'hypoth\`ese que $M$ est une sous-vari\'et\'e (sans
singularit\'es) de $\C^n$ \'equivaut au fait que la matrice $\left(
\frac{ \partial \rho_j }{ \partial t_{k_1 }} (t,\, \bar t) \ \ \frac{
\partial \rho_j }{ \partial \bar t_{ k_1 }}( t,\, \bar t) \right)_{ 1
\leq j \leq d }^{ 1 \leq k_1,\, k_2 \leq n}$ \`a $d$ lignes et \`a
$2n$ colonnes est de rang $d$ \`a l'origine. En exprimant
analytiquement la condition g\'eom\'etrique $T_0 M + JT_0 M=T_0 \C^n$,
on voit aussi que la g\'en\'ericit\'e de $M$ \'equivaut au fait que la
matrice $\left( \frac{ \partial \rho_j }{ \partial t_k }(t,\, \bar t)
\right)_{1 \leq j\leq d}^{ 1\leq k \leq n}$ de taille $d \times n$ est
de rang $d$ \`a l'origine. On v\'erifie que cette derni\`ere
condition est ind\'ependante du choix des \'equations d\'efinissantes
et qu'elle est alors satisfaite en tout point $p \in M$ suffisamment
proche de l'origine~: la g\'en\'ericit\'e est une condition
ouverte. Sans perte de g\'en\'eralit\'e, on pourra donc supposer
qu'elle est satisfaite en tout point de $M$.

Les sous-vari\'et\'es analytiques r\'eelles dites {\sl Cauchy-Riemann}
(CR) de $\C^n$ les plus g\'en\'erales ne sont pas forc\'ement
g\'en\'eriques, n\'eanmoins, il est bien connu qu'elles sont
g\'en\'eriques dans leur complexification intrins\`eque, qui est une
sous-vari\'et\'e complexe localement biholomorphe
\`a $\C^{n'}$ pour un entier $n'\leq n$~; par cons\'equent,
ne travailler qu'avec des sous-vari\'et\'es
g\'en\'eriques n'est en 
rien restrictif.

Gr\^ace \`a la formule de la dimension pour la somme de deux
sous-espaces vectoriels $\dim_\R (E + F) = \dim_\R \, E + \dim_\R \,
F- \dim_\R (E \cap F)$, on d\'eduit que la distribution de
sous-espaces lin\'eaires r\'eels $M \ni p \mapsto T_p M \cap J T_pM$
est de rang constant $2(n-d)$. Le sous-espace $T_p M \cap J T_p M$ est
appel\'e sous-espace complexe tangent \`a $M$ et not\'e
$T_p^cM$. Puisque $J^2= -{\rm Id}$, c'est l'unique sous-espace
$J$-invariant de $T_pM$ de dimension maximale. L'entier $m:= n-d$ est
appel\'e la {\sl dimension CR}\, de $M$.

Seuls deux cas limites sont inint\'eressants du point
de vue de la g\'eom\'etrie CR locale~: lorsque la codimension $d$
s'annule, auquel cas $M$ s'identifie \`a un cube ouvert de $\C^n$, et
lorsque la dimension CR $m$ s'annule, auquel cas $M$ s'identifie \`a
un cube ouvert de $\R^n$. Par cons\'equent, {\it dans tout le
m\'emoire nous travaillerons avec des sous-vari\'et\'es g\'en\'eriques
locales de codimension et de dimension CR strictement positives}.

Apr\`es une transformation lin\'eaire inversible de $\C^n$
arbitrairement proche de l'identit\'e, on peut supposer que dans les
coordonn\'ees $t$ scind\'ees en deux groupes $t= (z,\, w) = (z_1,\,
\dots,\, z_m,\, w_1,\, \dots,\, w_d)= (x +iy,\, u +iv) \in \C^m \times
\C^d$, l'espace tangent $T_0M$ ne contient aucun vecteur r\'eel de
l'espace des $v$~; de mani\`ere \'equivalente, on a~: $T_0 M + \left(
\{0\} \times \R_v^d \right) = T_0 \C^n$. Puisque la sous-vari\'et\'e
$M$ est alors graph\'ee au-dessus de $\C_z^m \times \R_u^d$, on peut la
repr\'esenter par $d$ \'equations de la forme $v_j = \varphi_j (x,\,
y, \,u)$, $j=1, \dots,\, d$, o\`u les s\'eries analytiques r\'eelles
$\varphi_j (x,\, y,\, u) \in\R \{ x,\, y,\, u \}$ satisfont bien s\^ur
$\varphi_j (0) = 0$.

De notre point de vue, toute repr\'esentation de $M$ par des
\'equations d\'efinissantes r\'eelles poss\`ede le d\'efaut majeur de
ne pas diff\'erencier clairement les variables $t$ et $\bar t$. C'est
pourquoi nous utiliserons toujours d'autres \'equations
d\'efinissantes, dites {\sl \'equations complexes}.

Rempla\c cons $x$ par $(z +\bar z) /2$, $y$ par $(z- \bar z)/2i$, $u$
par $(w + \bar w)/2$ et $v$ par $(w- \bar w)/2i$ dans ces \'equations,
ce qui donne $w_j - \bar w_j = 2i \, \varphi_j \left( (z+ \bar z)/2,\,
(z-\bar z)/ 2i,\, (w+\bar w)/2\right)$. La matrice des d\'eriv\'ees
partielles de ces $d$ \'equations scalaires par rapport aux variables
$w_l$, calcul\'ee \`a l'origine, vaut~: $I_{ d\times d} - 2i \left(
\frac{ 1}{2} \, \frac{ \partial \varphi_j}{ \partial u_l} (0,\, 0,\,
0) \right)_{ 1\leq j\leq d}^{ 1\leq l \leq d}$~; elle est non nulle,
puisque sa partie r\'eelle ne s'annule pas. Par cons\'equent, le
th\'eor\`eme des fonctions implicites (version analytique complexe)
s'applique et il permet de r\'esoudre $w$ en fonction des variables
$z$, $\bar z$ et $\bar w$. Notons les $d$ \'equations complexes
obtenues sous la forme $w_j = \overline{ \Theta}_j (z,\, \bar z,\,
\bar w )$, $j=1 ,\dots,\, d$, o\`u les $\overline{ \Theta }_j (z ,\,
\bar z,\, \bar w) \in \C\{ z, \, \bar z ,\, \bar w\}$ sont les uniques
solutions des identit\'es analytiques\footnotemark[4]
\def\theequation{3.3}\begin{equation}
\frac{ 
\overline{\Theta}_j(z,\, \bar z,\, \bar w) -\bar w_j}{2i} \equiv
\varphi_j\left(
\frac{ z+\bar z}{2}, \, 
\frac{ z- \bar z}{2i}, \, 
\frac{ 
\overline{\Theta}(z,\, \bar z,\, \bar w) + \bar w}{2}
\right). 
\end{equation}
On v\'erifie que de telles \'equations complexes
existent avec l'hypoth\`ese (un peu plus g\'en\'erale)
$T_0M + \left( \{0\}\times \C_w^n \right) = T_0 \C^n$.
Un probl\`eme de coh\'erence surgit alors imm\'ediatement~: les
s\'eries enti\`eres $\overline{ \Theta }_j$ \'etant \`a valeurs
complexes, les parties r\'eelles et imaginaires des \'equations $w_j =
\overline{ \Theta}_j(z,\, \bar z,\, \bar w)$ fournissent en v\'erit\'e
$2d$ \'equations r\'eelles. De mani\`ere \'equivalente, il faudrait
leur ajouter les \'equations conjugu\'ees $\bar w_j = \Theta_j(\bar z,
z, w)$, $j= 1, \dots,\, d$, ce qui semble contredire le fait que $M$
est de codimension $d$. Autre ambigu\"{\i}t\'e~: on aurait pu choisir
de r\'esoudre par rapport \`a $\bar w$ plut\^ot que par rapport \`a
$w$.

\footnotetext[4]{
Nous notons ici ces solutions $\overline{ \Theta}_j$, en les
\'equipant d'embl\'ee d'une barre de conjugaison complexe, avec
l'id\'ee que dans leurs $2m+d$ arguments $(z,\, \bar z,\, \bar w)$,
ceux qui sont des conjugu\'es de variables complexes dominent. Dans la
suite, nous verrons \`a quel point le jeu alternatif entre les
s\'eries enti\`eres $\overline{ \Theta}_j( z,\bar z, \bar w)$ et leurs
conjugu\'ees $\Theta_j(\bar z,\, z,\, w)$ est crucial. }

Heureusement, on d\'emontre ({\it cf.}~\cite{ ber1999a} ou
\cite{me2003}) qu'il existe une matrice inversible $\left( \overline{
a}_{ i, j} \left( \bar t ,\, t \right) \right)_{ 1 \leq i \leq d }^{ 1
\leq j \leq d}$ de s\'eries formelles convergentes, de taille $d
\times d$, \'egale \`a $-{\rm Id }_{ d\times d}$ \`a l'origine, telle
qu'on a l'identit\'e formelle vectorielle
\def\theequation{3.4}\begin{equation}
\bar w - \Theta (\bar z,\, z, \, w) \equiv
\overline{ a} (\bar t,\, t) \, \left[
w-\overline{ \Theta} (z,\, \bar z,\, \bar w)
\right]
\end{equation}
dans $\C\{ t,\, \bar t\}^d$. R\'eciproquement, on v\'erifie que pour
toute s\'erie enti\`ere vectorielle analytique $\Theta (\bar z,\, z,\,
w) \in \C\{ \bar z,\, z,\, w \}^d$ satisfaisant une telle \'equation,
le sous-ensemble $M := \{(z,\, w)\in \C^n: \, w= \overline{ \Theta}
(z,\, \bar z ,\, \bar w)\}$ est une sous-vari\'et\'e locale de $\C^n$
analytique {\it r\'eelle}, g\'en\'erique et de codimension $d$.

Il d\'ecoule aussi visiblement de~\thetag{ 3.4} qu'on n'obtiendrait
aucune \'equation ind\'ependante nouvelle pour la repr\'esentation de
$M$ par les \'equations $w_j = \overline{ \Theta }_j (z,\, \bar z,\,
\bar w )$ en leur ajoutant les \'equations conjugu\'ees $\bar w_j =
\Theta_j ( \bar z,\, z,\, w)$.

Mais on se gardera d'en d\'eduire (comme dans~\cite{ ber1999a}) que
l'on peut choisir d\'efinitivement pour la repr\'esentation de $M$
l'une de ses deux collections d'\'equations d\'efinissantes
conjugu\'ees. En effet, nous allons constater nettement dans la
Section~7 ci-dessous qu'il est n\'ecessaire d'effectuer un jeu
alternatif permanent entre les deux syst\`emes d'\'equations, vus
comme un couple d'objets sym\'etriques articul\'es par la
relation~\thetag{ 3.4}~: {\it cette ambigu\"{\i}t\'e est
fondamentale}.

\subsection*{3.5.~Complexification extrins\`eque}
Soient maintenant $\zeta \in \C^m$ et $\xi \in \C^d$ des nouvelles
coordonn\'ees ind\'ependantes correspondant aux complexifications des
variables $\bar z$ et $\bar w$, ce que l'on peut \'ecrire
symboliquement $\zeta := (\bar z)^c$ et $\xi := (\bar w)^c$, o\`u la
lettre <<c>> est l'initiale du mot <<complexification>>. On notera
$\tau := (\zeta,\, \xi ) \in \C^n$ la complexification de $\bar t$.
Dans la suite, on utilisera les notations l\'eg\`erement abr\'eg\'ees
$\overline{\Theta} (z,\, \bar t)$ et $\Theta( \bar z,\, t)$ plus
fr\'equemment que $\overline{ \Theta} (z ,\, \bar z ,\, \bar w)$ et
$\Theta( \bar z ,\, z,\, w)$. En rempla\c cant $\bar t$ par $\tau$
dans les s\'eries enti\`eres convergentes $\overline{ \Theta}_j (z,\,
\bar t)= \sum_{ \beta \in \N^m, \, \alpha \in \N^n }\, \overline{
\Theta}_{ j,\, \beta,\, \alpha }\, z^\beta \, \bar t^\alpha$, o\`u
$\overline{ \Theta}_{j,\, \beta,\, \alpha} \in \C$, on obtient des
s\'eries enti\`eres convergentes $\overline{ \Theta}_j (z,\, \tau) :=
\sum_{ \beta \in \N^m, \, \alpha \in \N^n }\, \overline{ \Theta}_{
j,\, \beta ,\, \alpha }\, z^\beta \, \tau^\alpha$ des $m+n$ variables
ind\'ependantes $(z,\, \tau)$. La {\sl complexification extrins\`eque}
$\mathcal{ M}:= (M )^c$ de $M$ est alors la sous-vari\'et\'e
analytique complexe de $\C^n \times \C^n$ passant par l'origine et de
codimension $d$ qui est d\'efinie par l'une des deux collections de
$d$ \'equations d\'efinissantes holomorphes
\def\theequation{3.6}\begin{equation}
w_j = \overline{\Theta}_j (z,\, \tau)
\ \ \ \ \ \ \
{\rm ou}
\ \ \ \ \ \ \
\xi_j = \Theta_j(\zeta, \, t),
\end{equation}
qui sont \'evidemment \'equivalentes en vertu de la complexification
de la relation~\thetag{ 3.4}, qui s'\'ecrit~:
\def\theequation{3.7}\begin{equation}
\xi- \Theta (\zeta,\, t) \equiv
\overline{ a} (\tau,\, t) \, \left[
w-\overline{ \Theta} (z,\tau)
\right].
\end{equation}
Notons que $\dim_\C \, \mathcal{ M}=2m+d$. Notons aussi que $M$ se
plonge dans $\C^n \times \C^n$ comme l'intersection de sa
complexification $\mathcal{ M}$ avec la diagonale antiholomorphe
$\underline{ \Lambda} := \{(t,\, \tau)\in \C^n \times \C^n: \, \tau =
\bar t\}$. Cependant, dans toute la suite de ce m\'emoire, nous ne
travaillerons d\'esormais qu'avec des sous-vari\'et\'es g\'en\'eriques
complexifi\'ees. C'est pourquoi {\it notre objet g\'eom\'etrique
fondamental de d\'epart est la sous-vari\'et\'e analytique complexe
locale $\mathcal{ M}$ d\'efinie par les \'equations~\thetag{ 3.6} qui
sont articul\'ees par la relation de sym\'etrie~\thetag{ 3.7}}. Plus
analytiquement encore, {\it notre objet fondamental de d\'epart est la
collection des s\'eries enti\`eres analytiques complexes $\overline{
\Theta}_j (z,\, \tau)$ ainsi que leurs conjugu\'ees $\Theta_j
(\zeta,\, t)$}.

\subsection*{3.8.~Paire de feuilletages invariants} \
Toute biholomorphisme local de $\C^n$ la forme $t' = h(t) = (h_1(t),
\dots,\, h_n(t)) \in \C \{ t \}^n$ fixant l'origine induit par
conjugaison l'anti-biholomorphisme $\bar t'= \bar h(\bar t)$ et par
cons\'equent, il se complexifie en un biholomorphisme $(t',\, \tau') =
(h(t),\, \bar h(\tau))$ de $\C^n\times \C^n$ qui est d'une forme
particuli\`ere, puisqu'il est \`a variables
s\'epar\'ees. G\'eom\'etriquement parlant, un tel biholomorphisme du
produit $\C^n\times \C^n$ envoie les sous-ensembles $\{t= ct.\}$ et
$\{ \tau=ct.\}$ sur les sous-ensembles $\{t'=ct.\}$ et $\{ \tau' =
ct.\}$~: il stabilise la paire de feuilletages triviaux qui sont
parall\`eles aux axes de coordonn\'ees <<horizontales>> $t$ et
<<verticales>> $\tau$. Par cons\'equent, tous les concepts
analytico-g\'eom\'etriques locaux qui sont attach\'es \`a $M$ d'une
mani\`ere qui est invariante par rapport aux changements de
coordonn\'ees holomorphes de la forme $t \mapsto h(t)$ co\"{\i}ncident
avec les objets analytico-g\'eom\'etriques de $\mathcal{ M}$ qui sont
invariants par rapport au sous-groupe (infini) de transformations de
la forme $(t,\,\tau) \mapsto ( h(t) ,\, \bar h(\tau))$. En
particulier, les deux feuilletages de $\mathcal{ M}$ dont les feuilles
sont les intersections de $\mathcal{ M}$ avec les sous-ensembles $\{
t= ct. \}$ et $\{ \tau= ct. \}$ sont invariants. Analytiquement, ces
feuilles sont ce qu'on appellera les {\sl sous-vari\'et\'es de Segre
complexifi\'ees} $\mathcal{ S}_{\tau_p}$ et les {\sl sous-vari\'et\'es
de Segre complexifi\'ees conjugu\'ees}\footnotemark[5]\!, d\'efinies
par
\def\theequation{3.9}\begin{equation}
\left\{
\aligned
\mathcal{ S}_{\tau_p}
& \
:=\left\{
(t,\, \tau)\in \C^{2n}: \ 
\tau= \tau_p, \
w= \overline{\Theta}( z, \, \tau_p)
\right\}
=\mathcal{ M}\cap \{
\tau = \tau_p\} \ \ \ \ \ {\rm et}
 \\
\underline{\mathcal{ S}}_{t_p}
& \
:=\left\{ (t,\, \tau)\in \C^{2n} 
: \ 
t= t_p, \ 
\xi = 
\Theta (\zeta,\, t_p)
\right\} \ 
=\mathcal{ M}\cap \{ t =
t_p\},
\endaligned\right.
\end{equation}
o\`u $\tau_p\in \C^n$ et $t_p\in \C^n$ sont fixes.

\footnotetext[5]{
Avant complexification, on peut d\'efinir pr\'ealablement les
sous-vari\'et\'es de Segre classiques $S_{ \bar t_p } \subset \C^n$
(\cite{ se1931a}, \cite{ se1931b}, \cite{ se1932}, \cite{ we1977},
\cite{ we1977}) et conjugu\'ees $\overline{ S}_{ t_p } \subset \C^n$
(traditionnellement ignor\'ees), que l'on envisage d'embl\'ee comme un
couple sym\'etrique articul\'e par la conjugaison complexe ({\it
voir}~surtout \cite{ me1998} et \cite{ me2003}). }

\bigskip
\begin{center}
\input complexification.pstex_t
\end{center}
\bigskip

Le diagramme ci-dessus a pour objet de repr\'esenter cette paire
fondamentale de feuilletages. Cependant, nous mettons le lecteur en
garde, parce que sur cette figure bidimensionnelle, la codimension
dans $\mathcal{ M}$ de la paire de feuilletages semble nulle alors
qu'en r\'ealit\'e, elle est strictement positive~:
\def\theequation{3.10}\begin{equation}
\dim_\C \, \mathcal{ M} - \dim_\C \, 
\mathcal{ S}_{ \tau_p } - \dim_\C
\, \underline{ \mathcal{ S }}_{ t_p}
= d \geq 1.
\end{equation}
Pour se repr\'esenter intuitivement la situation g\'eom\'etrique d'une
mani\`ere plus ad\'equate, on devrait imaginer par exemple que
$\mathcal{M}$ est un cube de dimension trois \'equip\'e de deux
feuilletages par des courbes qui sont en position 
g\'en\'erale.

\subsection*{3.11.~Flots de champs CR 
complexifi\'es} On dira intuitivement que $\mathcal{ M}$ est {\sl
minimale \`a l'origine} si l'on peut recouvrir un voisinage de $0$
dans $\mathcal{ M}$ en se d\'epla\c cant alternativement le long des
sous-vari\'et\'es de Segre complexifi\'ees et le long des
sous-vari\'et\'es de Segre complexifi\'ees conjugu\'ees. Afin de
d\'efinir rigoureusement cette condition ({\it voir}~la
D\'efinition~3.19 {\it infra}), il est n\'ecessaire d'exprimer
math\'ematiquement ce que l'on entend par d\'eplacement alternatif le
long de la paire de feuilletages invariants.

Pour cela, complexifions une famille g\'en\'eratrice $L_1,\dots,\,
L_m$ de champs de vecteurs CR tangents \`a $M$ de type $(1,\, 0)$
ainsi que leurs conjugu\'es $\overline{ L}_1,\dots,\, \overline{ L}_m$
qui sont de type $(0,\, 1)$. On peut choisir explicitement les
g\'en\'erateurs $L_k := \frac{ \partial }{ \partial z_k } + \sum_{ j=1
}^d\, \frac{ \partial \overline{ \Theta }_j}{ \partial z_k} ( z,\,
\bar z,\, \bar w) \, \frac{ \partial}{ \partial w_j}$ pour $k =1,
\dots,\,m$. Les complexifications fournissent deux collections de $m$
champs de vecteurs \`a coefficients analytiques donn\'es explicitement
par
\def\theequation{3.12}\begin{equation}
\left\{
\aligned
\mathcal{L}_k:= & \ 
{\partial \over\partial z_k}+
\sum_{j=1}^d\,
{\partial\overline{\Theta}_j\over \partial z_k}
(z,\, \zeta,\, \xi)\,
\, {\partial\over\partial w_j}, \ \ \ \ 
k=1,\dots,\,m, \ \ \ \ \ \ {\rm et} \\
\underline{\mathcal{L}}_k:= & \
{\partial\over\partial \zeta_k}+\sum_{j=1}^d\,
{\partial \Theta_j\over \partial \zeta_k}(\zeta,\, z,\, w)\, \,
{\partial\over\partial \xi_j}, \ \ \ \
k=1,\dots,\,m.
\endaligned\right.
\end{equation}
On v\'erifie imm\'ediatement que $\mathcal{ L}_k \left( w_j-
\overline{ \Theta}_j (z,\, \zeta,\, \xi) \right) \equiv 0$, ce qui
montre que les champs de vecteurs $\mathcal{ L}_k$ sont tangents \`a
$\mathcal{ M}$. De mani\`ere analogue, $\underline{ \mathcal{ L}}_k
\left( \xi_j - \Theta_j (\zeta,\, z,\, w) \right) \equiv 0$, de telle
sorte que les champs de vecteurs $\underline{ \mathcal{ L}}_k$ sont
aussi tangents \`a $\mathcal{ M}$. De plus, on v\'erifie
imm\'ediatement les relations de commutation $[\mathcal{ L}_k,\,
\mathcal{ L}_{k'}] =0$ et $[\underline{ \mathcal{ L}}_k,\, \underline{
\mathcal{ L}}_{ k'}] =0$ pour tous $k,\, k'= 1,\dots,\, m$. D'apr\`es
le th\'eor\`eme de Frobenius, il d\'ecoule de ces relations de
commutation que chacune des distributions $m$-dimensionnelles
engendr\'ees par ces deux collections de champs de vecteurs est
int\'egrable~; cela n'a rien de surprenant, puisque les vari\'et\'es
int\'egrales de $\{ \mathcal{ L}_k \}_{1 \leq k \leq m}$ ne sont
autres que les vari\'et\'es de Segre complexifi\'ees, tandis que les
vari\'et\'es int\'egrales de $\{ \underline{ \mathcal{ L }}_k\}_{ 1
\leq k \leq m}$ ne sont autres que les vari\'et\'es de Segre
complexifi\'ees conjugu\'ees. Bien s\^ur, en g\'en\'eral les
$\mathcal{ L}_k$ ne commutent pas avec les $\underline{ \mathcal{ L
}}_{k'}$~: c'est justement la non-int\'egrabilit\'e de la distribution
CR complexifi\'ee engendr\'ee par les deux familles $\{ \mathcal{
L}_k, \, \underline{ \mathcal{ L}}_{k'}\}_{ 1\leq k,\, k' \leq m}$ qui
est responsable de la minimalit\'e ({\it cf.} Lemme~3.22 {\it
infra}).

Gr\^ace \`a cette paire de familles de champs de vecteurs, on peut
param\'etrer les sous-vari\'et\'es de Segre complexifi\'ees
(conjugu\'ees). En effet, introduisons les flots <<multiples>> des
deux collections $(\mathcal{ L}_k )_{1 \leq k \leq m}$ et
$(\underline{\mathcal{L }}_{ k'})_{ 1\leq k' \leq m}$. Si $p$ est un
point arbitraire de $\mathcal{ M}$ dont les coordonn\'ees $(w_p,\,
z_p,\, \zeta_p,\, \xi_p) \in \C^{ 2n}$ satisfont les
\'equations~\thetag{ 3.6} et si $z_1 := (z_{ 1,\, 1},\dots,\, z_{1,\,
m})\in\C^m$ est un param\`etre de <<multitemps>> complexe arbitraire,
d\'efinissons le <<multiflot>> de $\mathcal{ L}$ par
\def\theequation{3.13}\begin{equation}
\left\{
\aligned
\mathcal{L}_{z_1}(z_p,\, w_p,\, \zeta_p,\, \xi_p) & \ :=
\exp\left(z_1\mathcal{L})(p):=
\exp(z_{1,\, 1}\mathcal{L}_1(\cdots(\exp(z_{1,\, m}\mathcal{L}_m(
p)))\cdots)\right) \\
& \
:=\left(
z_p+z_1,\, \overline{\Theta}(z_p+z_1,\, \zeta_p,\, \xi_p),\, 
\zeta_p,\, \xi_p\right).
\endaligned\right.
\end{equation}
Bien entendu, $\mathcal{ L}_{ z_1}(p)$ appartient \`a $\mathcal{
M}$. De mani\`ere analogue, pour $p\in \mathcal{ M}$ et $\zeta_1 \in
\C^m$, le <<multiflot>> de $\underline{ \mathcal{ L}}$
se d\'efinit par
\def\theequation{3.14}\begin{equation}
\underline{\mathcal{L}}_{\zeta_1}(z_p,\,w_p,\,\zeta_p,\,\xi_p):=
(z_p,\,w_p,\,\zeta_p+\zeta_1,\,\Theta(\zeta_p+\zeta_1,\,z_p,\,w_p)), 
\end{equation}
et l'on a $\underline{ \mathcal{ L}}_{\zeta_1}(p) \in \mathcal{
M}$. Notons que les deux applications~\thetag{ 3.13} et~\thetag{ 3.14}
sont holomorphes par rapport \`a leurs variables.

\subsection*{3.15~Cha\^{\i}nes de Segre}
Ainsi, pla\c cons tout d'abord le point $p$ \`a l'origine et d\'epla\c
cons-nous le long de la vari\'et\'e de Segre complexifi\'ee
conjugu\'ee $\underline{ \mathcal{ S} }_0$ d'une hauteur de $z_1\in
\C^m$, c'est-\`a-dire consid\'erons le point $\underline{ \mathcal{ L
}}_{ z_1}(0)$, que nous noterons aussi $\underline{ \Gamma}_1(z_1)$.
Bien s\^ur, on a $\underline{ \Gamma}_1 (0)=0$. Soit $z_2 \in \C^m$. En
partant de ce point $\underline{ \Gamma}_1(z_1)$, d\'epla\c cons-nous
horizontalement le long de la vari\'et\'e de Segre complexifi\'ee
d'une longueur de $z_2 \in \C^m$, c'est \`a dire consid\'erons le
point
\def\theequation{3.16}\begin{equation}
\underline{\Gamma}_2(z_1,\,z_2):=
\mathcal{L}_{z_2}(\underline{\mathcal{L}}_{z_1}(0)).
\end{equation}
Ensuite, d\'efinissons $\underline{ \Gamma}_3 (z_1,\, z_2,\, z_3) :=
\underline{ \mathcal{ L}}_{ z_3}( \mathcal{ L}_{ z_2} (\underline{
\mathcal{ L} }_{ z_1} (0 )))$, puis 
\def\theequation{3.17}\begin{equation}
\underline{ \Gamma}_4 (z_1,\,
z_2,\, z_3,\, z_4) := \mathcal{ L}_{ z_4}( \underline{ \mathcal{ L}}_{
z_3}( \mathcal{ L}_{ z_2} (\underline{ \mathcal{ L}}_{ z_1} (0 )))),
\end{equation}
et ainsi de suite. Le diagramme suivant illustre le proc\'ed\'e~:

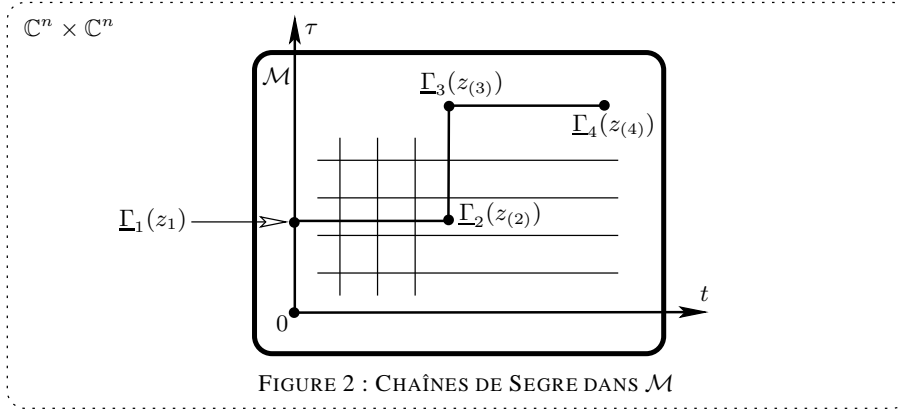
\bigskip
\begin{center}
\input orbite.pstex_t
\end{center}
\bigskip

Par r\'ecurrence, pour tout entier positif $k$, on obtient une
application holomorphe locale $\underline{ \Gamma }_k (z_1,\,
\dots,\,z_k)$ \`a valeurs dans $\mathcal{ M}$, d\'efinie pour $z_1,\,
\dots,\, z_k \in \C^m$ suffisamment petits et satisfaisant
$\underline{ \Gamma}_k(0, \,\dots, \, 0) =0$. Dans la suite, nous
utiliserons souvent l'abr\'eviation $z_{ (k)} :=(z_1,\, \dots,\, z_k)
\in \C^{ mk}$ et nous appellerons $\underline{ \Gamma}_k$ la {\sl
$k$-i\`eme cha\^{\i}ne de Segre conjugu\'ee}.

Si l'on commen\c cait cette suite d'applications holomorphes
compos\'ees par le flot de $\mathcal{ L}$ au lieu de commencer par
celui de $\underline{ \mathcal{ L }}$, on obtiendrait des applications
$\Gamma_1( z_1):= \mathcal{ L}_{ z_1}(0)$, puis $\Gamma_2 (z_{(2)}) :=
\underline{ \mathcal{ L}}_{ z_2} (\mathcal{ L}_{z_1}(0))$, {\it etc.},
et g\'en\'eralement $\Gamma_k(z_{(k)})$. Nous appellerons $\Gamma_k$
la {\sl $k$-i\`eme cha\^{\i}ne de Segre}.

Puisque $\Gamma_k (0) = \underline{ \Gamma}_k (0) = 0$, pour tout
entier strictement positif $k$, il existe un cube (polydisque)
suffisamment petit $\square_{mk}(\delta_k)$ centr\'e \`a l'origine
dans $\C^{mk}$ et de rayon $\delta_k >0$ tel que $\Gamma_k(z_{(k)})$
et $\underline{ \Gamma}_k (z_{(k)})$ appartiennent \`a $\mathcal{M}$
pour tout $z_{(k)} \in \square_{mk}( \delta_k)$.

Il existe une relation de sym\'etrie entre $\Gamma_k$ et
$\underline{\Gamma}_k$. En effet, soit $\overline{ \sigma}$
l'involution antiholomorphe de $\C^n\times \C^n$ d\'efinie par
$\overline{ \sigma} (t,\,\tau):=(\bar\tau,\,\bar t)$. Puisque l'on a
$w = \overline{ \Theta}( z,\, \zeta,\, \xi)$ si et seulement si $\xi =
\Theta( \zeta,\, z,\, w)$, cette involution envoie $\mathcal{M}$ dans
$\mathcal{M}$ et elle fixe aussi point par point la diagonale
antiholomorphe $\underline{ \Lambda}$. En appliquant $\overline{
\sigma}$ aux d\'efinitions ~\thetag{3.13} et~\thetag{3.14} des flots
de $\mathcal{ L}$ et de $\underline{ \mathcal{ L}}$, on v\'erifie
ais\'ement que $\overline{ \sigma} ( \mathcal{ L}_{ z_1} (p)) =
\underline{ \mathcal{ L }}_{ \bar z_1} (\overline{ \sigma} (p))$. Il
en d\'ecoule la relation de sym\'etrie g\'en\'erale $\overline{
\sigma} \left( \Gamma_k( z_{ (k) }) \right) = \underline{ \Gamma}_k
\left( \overline{ z_{ (k) }} \right)$. Dans la suite de ce m\'emoire,
nous travaillerons essentiellement avec les applications $\underline{
\Gamma}_k$.

\subsection*{3.18.~Minimalit\'e locale}
Observons que $\mathcal{ L}_0 (p) =p$ et que $\underline{ \mathcal{
L}}_0 (p)=p$~; autrement dit, $\mathcal{ L}_0$ et $\underline{
\mathcal{L}}_0$ co\"{\i}ncident avec l'application identit\'e. Nous
en d\'eduisons l'identit\'e~: $\Gamma_{k+1}(z_{(k)},\,0) \equiv
[\mathcal{ L} \ {\rm ou} \ \underline{ \mathcal{ L}}]_0 \left(
\Gamma_k( z_{ (k)})\right) \equiv \Gamma_k (z_{ (k)})$. Par
cons\'equent, les rangs des applications $\underline{ \Gamma}_k$
croissent avec $k$. Bien entendu, ces rangs sont born\'es par $\dim_\C
\, \mathcal{ M}= 2m+d$. On d\'eduit aussi des relations
$\Gamma_{k+1}(z_{(k)},\,0) = \Gamma_k (z_{ (k)})$ que pour tout $k\geq
2$, le rang \`a l'origine de $\underline{ \Gamma}_k$ est
invariablement \'egal \`a $2m$, mais en des points $\underline{ z}_{
(k)}^* \in \square_{ mk} (\delta_k)$ distincts de l'origine, le rang
des applications $\underline{ \Gamma}_k$ peut augmenter jusqu'\`a
atteindre $2m+d$. Nous pouvons maintenant \'enoncer la d\'efinition
pr\'ecise de la notion de minimalit\'e.

\def\thedefinition{3.19}\begin{definition}
{\rm 
La sous-vari\'et\'e complexifi\'ee $\mathcal{ M}$ d'une
sous-vari\'et\'e locale de $\C^n$ analytique r\'eelle, g\'en\'erique
$M$ est dite {\sl minimale \`a l'origine} s'il existe un entier $\mu_0
\geq 1$ et des points $\underline{ z}_{( \mu_0)}^* \in \square_{ m
\mu_0} ( \delta_{ \mu_0 })$, arbitrairement proches de l'origine
satisfaisant $\underline{ \Gamma }_{ \mu_0} \left( \underline{ z}_{ (
\mu_0 )}^* \right) =0$ tels que l'application $\underline{ \Gamma }_{
\mu_0}$ est de rang (maximal possible) \'egal \`a $2m + d$ 
en ces points $\underline{ z}_{ ( \mu_0) }^*$.
}
\end{definition}

Dans ce cas, l'image d'un petit voisinage d'un tel point $\underline{
z}_{ (\mu_0)}^*$ dans $\C^{ m\mu_0}$ contient un petit voisinage de
$0$ dans $\mathcal{ M}$. Gr\^ace \`a la relation $\overline{ \sigma}
\left( \Gamma_k( z_{ (k) }) \right) = \underline{ \Gamma}_k \left(
\overline{ z_{ (k)}} \right)$, une propri\'et\'e similaire est
satisfaite par $\Gamma_{ \mu_0}$, avec le m\^eme entier $\mu_0$. On
dira aussi que $M$ est minimale \`a l'origine si sa complexifi\'ee
$\mathcal{ M}$ l'est, au sens de cette d\'efinition.

Bien qu'il ne soit pas n\'ecessaire, {\it stricto sensu}, de commenter
cette condition de minimalit\'e pour comprendre la d\'emonstration du
Th\'eor\`eme~1.23, formulons quand m\^eme quelques \'enonc\'es
explicatifs.

\subsection*{ 3.20.~Commentaires}
Rappelons que le {\sl rang g\'en\'erique} d'une application holomorphe
entre deux cubes complexes est le maximum de son rang aux
diff\'erents points du cube \`a la source. En travaillant d'abord avec
le rang g\'en\'erique des applications $\underline{ \Gamma}_k$, qui
lui aussi cro\^{\i}t avec $k$, on d\'emontre des propri\'et\'es
\'el\'ementaires que nous r\'esumons dans l'\'enonc\'e suivant,
\'etabli dans la Section~7 de~\cite{ me1998} et dans le Chapitre~2
de~\cite{ me2003}.

\def\thetheorem{3.21}\begin{theorem}
La minimalit\'e en $0$ de la sous-vari\'et\'e g\'en\'erique
complexifi\'ee $\mathcal{ M}= (M )^c$ de $\C^n$ est une propri\'et\'e
invariante par biholomorphisme local~{\rm :} elle ne d\'epend ni du
choix d'\'equations d\'efinissantes pour $M$, ni du choix d'un
syst\`eme de coordonn\'ees holomorphes s'annulant au point de
r\'ef\'erence, ni du choix d'un syst\`eme g\'en\'erateur de champs CR
complexifi\'es {\rm (}conjugu\'es{\rm )} $(\mathcal{ L }_k)_{ 1 \leq k
\leq m}$ et $(\underline{ \mathcal{ L} }_k)_{ 1 \leq k \leq m}$. De
plus, il existe un entier invariant $\nu_0$ satisfaisant $\nu_0 \leq
d+1$, qu'on appellera le {\rm type de Segre}\, de $M$ \`a l'origine,
qui est le plus petit entier $k$ tel que les applications $\Gamma_k$
et $\underline{ \Gamma}_k$ sont de rang g\'en\'erique \'egal \`a $2m
+d$ sur le cube $\square_{ mk}( \delta_k)$, pour tout $k\geq \nu_{
0}+1$. Enfin, l'entier {\rm impair} $\mu_0 :=2 \nu_0 +1$, que l'on
appellera le {\rm type de Segre de} $\mathcal{M}$ \`a l'origine, est
le plus petit entier $k$ tel que les applications $\Gamma_k$ et
$\underline{ \Gamma }_k$ sont de rang \'egal \`a $2m+d$ en des points
$z_{ (k )}^* \in \square_{ mk} (\delta_k )$ et $\underline{ z}_{ (k
)}^* \in \square_{ mk} (\delta_k )$ satisfaisant $\Gamma_k \left(z_{
(k )}^* \right) = 0$ et $\underline{ \Gamma}_k \left( \underline{ z}_{
(k )}^*\right) = 0$ qui sont arbitrairement proches de l'origine dans
$\C^{ mk}$.
\end{theorem}

Le fait qu'il existe de tels points $z_{ (k) }^*$ et $\underline{ z}_{
(k)}^*$ arbitrairement proches de l'origine est d\^u \`a la
propri\'et\'e qu'ont les applications holomorphes locales d'atteindre
leur rang g\'en\'erique en tout point d'un ouvert de Zariski dense de
l'espace source, gr\^ace au principe du prolongement analytique ({\it
cf.}~\cite{ me1998}, \cite{ me2003}). Ce fait assez crucial sera
utilis\'e dans le \S3.31 ci-dessous.

Le proc\'ed\'e de d\'emonstration de ce th\'eor\`eme est inspir\'e de
la construction des orbites de champs de vecteurs, telle qu'elle
appara\^{\i}t dans l'article~\cite{ su1973}, dont la nouveaut\'e
principale r\'esidait dans le traitement des syst\`emes de champs de
vecteurs de classe $\mathcal{ C}^\infty$, par opposition aux champs
dont les coefficients sont analytiques. Mais dans la cat\'egorie
analytique, un th\'eor\`eme semblable \'etait connu depuis
l'article~\cite{ na1966}, o\`u l'auteur raisonne plut\^ot en
consid\'erant l'alg\`ebre de Lie engendr\'ee par un syst\`eme de
champs de vecteurs analytiques. H.J.~Sussmann d\'emontre une
proposition g\'en\'erale \'etablissant l'\'equivalence entre ces deux
proc\'ed\'es (\cite{ su1973}, Theorem~8.1 et \S9), laquelle,
sp\'ecifi\'ee \`a la g\'eom\'etrie CR analytique locale, nous offre
l'\'enonc\'e suivant~:

\def\thelemma{3.22}\begin{lemma}
{\rm (\cite{ ber1996}, \cite{ me1998})}
La sous-vari\'et\'e g\'en\'erique {\rm analytique} complexifi\'ee
$\mathcal{ M} = (M)^c$ de $\C^n$ est minimale \`a l'origine {\rm (}au
sens de la D\'efinition~3.19{\rm )} si et seulement si l'alg\`ebre de
Lie engendr\'ee par les sections locales du fibr\'e tangent complexe
$T^cM=TM\cap JTM$ engendre l'espace tangent \`a $M$ \`a
l'origine {\rm (}cette alg\`ebre de Lie est constitu\'ee de toutes les
combinaisons lin\'eaires \`a coefficients analytiques r\'eels de
crochets de Lie embo\^{\i}t\'es $[X_1 [X_2 [ X_3 [\dots [X_k,\, X_{
k+1}] \dots ]]]]$ de longueur finie $k$ arbitraire, o\`u 
les $X_l$ sont des sections de $T^cM${\rm )}.
\end{lemma}

Cette deuxi\`eme condition, plus ancienne, est classique. Certains
auteurs appellent la sous-vari\'et\'e g\'en\'erique $M$ {\sl de type
fini \`a l'origine} (au sens de T.~Bloom et I.~Graham) si les crochets
de Lie de $T^cM$ de longueur arbitraire engendrent $TM$ \`a l'origine.
Nous pr\'ef\'erons l'appellation de {\sl minimalit\'e \`a l'origine}
(au sens de J.-M.~Tr\'epreau et A.E.~Tumanov), puisque ce sont les
orbites des champs de vecteurs CR complexifi\'es qui sont les <<bons
objets>>, et non leurs crochets de Lie. En effet, ce sont les travaux
profonds de J.-M.~Tr\'epreau~\cite{ tr1986}, \cite{ tr1990} et de
A.E.~Tumanov~\cite{ tu1988}, \cite{ tu1994} qui ont fait
d\'efinitivement comprendre l'ad\'equation de la correspondance entre
les orbites CR et les wedges attach\'es pour l'extension holomorphe
des fonctions CR ({\it cf.} aussi~\cite{ me1994}).

Rappelons aussi que les deux concepts que sont les flots de champs de
vecteurs et leurs crochets de Lie sont \'etrangers l'un \`a l'autre
quant \`a la combinatoire des calculs, comme l'illustrent les
nombreux exemples que l'on trouve \`a la Section~8 de~\cite{ me1998}.
De plus, la d\'emonstration d'un th\'eor\`eme essentiellement
\'equivalent au Th\'eor\`eme~3.21 que l'on trouve dans~\cite{ ber1996}
et dans~\cite{ ber1999a} utilise fortement les crochets de champs de
vecteurs et un syst\`eme de coordonn\'ees dites <<normales>>, qui sont
adapt\'ees aux {\sl nombres de L.~H\"ormander} de la structure,
lesquels sont des invariants combinatoires d\'ecrivant les sauts de
dimension occasionn\'es par le calcul successif des crochets de Lie
({\it voir}~les chapitres~4 et~10 de~\cite{ ber1999a}). En adoptant le
point de vue <<crochets de Lie>>, les d\'emonstrations deviennent
extraordinairement techniques, et ce, sans n\'ecessit\'e interne. On
pourrait de surcro\^{\i}t s'\'etonner que dans le livre~\cite{
ber1999a} (publi\'e dans une collection prestigieuse) qui contient une
copie de la d\'emonstration du th\'eor\`eme de H.J.~Sussmann (\S3.1~:
{\it Nagano's theorem}~; \S3.2~: {\it Sussmann's theorem}), le lien
naturel entre le proc\'ed\'e de H.J.~Sussmann et ce que les auteurs
appellent <<ensembles de Segre>> (qui sont en v\'erit\'e des
projections sur $\C^n$ des images dans $\mathcal{ M}$ des applications
$\underline{ \Gamma}_k$ ou $\Gamma_k$), n'ait pas \'et\'e observ\'e.

Deux ans apr\`es l'apparition de notre pr\'epublication
\'electronique~\cite{ me1998} (travail non publi\'e), S.M.~Baouendi,
P.~Ebenfelt et L.-P.~Rothschild ont repris en partie notre point de
vue bas\'e sur l'article de H.J.~Sussmann, afin de construire plus
\'economiquement leurs <<ensembles de Segre>>~; leur pr\'epublication
fut publi\'ee dans une revue sp\'ecialis\'ee en g\'eom\'etrie
alg\'ebrique pure~: \cite{ ber2003}. N\'eanmoins, dans cette
r\'ef\'erence, le point de vue exprim\'e par les {\sc Figures~1} et~2
ci-dessus est absent~: on n'y trouve ni mention de la paire de
feuilletages invariants, ni description intrins\`eque de la
sym\'etrie par conjugaison complexe ni aucune <<vision>>
g\'eom\'etrique. De plus, ces auteurs \'evitent de mentionner que
l'extraordinaire technicit\'e des chapitres~4 et~10 de leur
livre~\cite{ ber1999a} devient caduque.

Le point de vue ensembliste qui consiste \`a consid\'erer les
projections sur $\C^n$ des images des cha\^{\i}nes de Segre, {\it
i.e.} les ensembles $\pi_t \left( \underline{ \Gamma }_k ( z_{ (k) })
\right)$ n'apporte rien du point de vue fonctionnel. Ce sont au
contraire les cha\^{\i}nes de Segre $\Gamma_k$ et
$\underline{ \Gamma}_k$ 
vues comme {\it applications holomorphes locales}\, ainsi
que leur propri\'et\'e de submersivit\'e (dans le cas minimal) qui
sont vraiment utilis\'ees dans l'\'etude du principe de r\'eflexion
analytique. Par cons\'equent, nous n'adopterons jamais la
terminologie <<ensembles de Segre>>.

En conclusion, retenons seulement la propri\'et\'e
de submersivit\'e de $\underline{ \Gamma}_k$
\'enonc\'ee dans la D\'efinition~3.19.

\subsection*{3.23.~Projections des submersions $\Gamma_k$
et $\underline{ \Gamma}_k$ sur $\mathcal{ M}$} 
Soit $\mu_0 = 2 \nu_0 + 1$ le type de Segre type de $\mathcal{ M}$ \`a
l'origine $0$, qui est toujours impair. Si $\mathcal{ M}$ est
minimale \`a l'origine, les deux applications holomorphes locales
\def\theequation{3.24}\begin{equation}
\Gamma_{\mu_0} \ {\rm et} \ \, 
\underline{ \Gamma}_{\mu_0} \ : \ \ \ \
\ \square_{m\mu_0}( \delta_{\mu_0})
\longrightarrow \mathcal{ M}
\end{equation}
satisfont $\Gamma_{ \mu_0} ( z_{ (\mu_0)}^*) =0$ et $\underline{
\Gamma}_{ \mu_0} ( \underline{ z}_{ \mu_0}^*) =0$ et elles sont
submersives en $z_{ (\mu_0)}^*$ et $\underline{ z}_{ (\mu_0)}^*$,
c'est-\`a-dire de rang (maximal) \'egal \`a $\dim_\C \, \mathcal{ M}$.

Si l'on veut quitter l'espace $\C^n \times \C^n$ o\`u vit la
complexification $\mathcal{ M}$ et revenir \`a l'espace $\C^n$ des
coordonn\'ees $t$ o\`u vit $M$ (ou \`a l'espace des coordonn\'ees
$\tau$ o\`u vit la conjugu\'ee $\overline{ M }$ de $M$), on peut aussi
utiliser la propri\'et\'e de submersivit\'e des deux
applications~\thetag{ 3.24} en les composant \`a gauche avec l'une des
deux projections $\pi_t (t,\, \tau) := t$ et $\pi_\tau ( t,\, \tau) :=
\tau$ (sur les sous-espaces de coordonn\'ees <<horizontales>> et
<<verticales>>), qui sont bien \'evidemment submersives, ce qui donne
deux couples de possibilit\'e~:
$\pi_t \left(
\Gamma_{\mu_0} (z_{(\mu_0)})
\right)$, $\pi_t \left(
\underline{ \Gamma}_{ \mu_0} (
z_{(\mu_0)})
\right)$ et $\pi_\tau \left(
\Gamma_{\mu_0} (z_{(\mu_0)})
\right)$, $\pi_\tau \left(
\underline{ \Gamma}_{ \mu_0} (
z_{(\mu_0)})
\right)$. 
Pour deux de ces quatre expressions, une l\'eg\`ere simplification
formelle intervient alors~: nous affirmons que l'on a les
deux relations
\def\theequation{3.25}\begin{equation}
\left\{
\aligned
\pi_t \left(
\underline{ \Gamma}_{2\nu_0+1}(z_{(2\nu_0+1)})
\right)
& \
\equiv
\pi_t \left(
\underline{ \Gamma}_{2\nu_0}(z_{(2\nu _0)})
\right) 
\ \ \ \ \
{\rm et} \\
\pi_\tau \left(
\Gamma_{2\nu_0+1}(z_{(2\nu_0+1)})
\right)
& \
\equiv
\pi_\tau \left(
\Gamma_{2\nu_0}(z_{(2\nu_0)})
\right).
\endaligned\right.
\end{equation}
En effet, puisque $\mu_0 = 2\nu_0 +1$ est impair, le premier terme de
flot (\`a gauche) de la cha\^{\i}ne de Segre conjugu\'ee $\underline{
\Gamma}_{ 2\nu_0 +1}(z_{(2\nu_0+1)})$ est le terme $\underline{
\mathcal{ L}}_{z_{2\nu_0+1}}$, c'est-\`a-dire que l'on peut \'ecrire
\def\theequation{3.26}\begin{equation}
\underline{ \Gamma}_{2\nu_0+1} 
\left( z_{ (2\nu_0 +1)} \right) =
\underline{ \mathcal L}_{z_{2\nu_0+1}}\left(
\underline{ \Gamma}_{2\nu_0}( z_{(2\nu_0)})\right).
\end{equation}
Si l'on note les quatre coordonn\'ees du point $\underline{ \Gamma}_{
2\nu_0} (z_{(2\nu_0)})$ par
\def\theequation{3.27}\begin{equation}
\left( z(z_{(2\nu_0)}),\, w(z_{(2\nu_0)}),\,
\zeta(z_{(2\nu_0)}),\, \xi(z_{(2\nu_0)}) \right),
\end{equation}
une application de la formule~\thetag{ 3.14} nous
donnne
\def\theequation{3.28}\begin{equation}
\left\{
\aligned
\underline{ \Gamma}_{ 2 \nu_0 + 1} 
(z_{ (2\nu_0 + 1 )}) 
& \
=
\underline{ \mathcal{ L
}}_{z_{2\nu_0+1}}\left(
\underline{ \Gamma}_{
2\nu_0}( z_{(2\nu_0)})
\right) \\
& \
=
\left(
z(z_{(2\nu_0)}),\, 
w(z_{(2\nu_0)}),\, 
z_{2\nu_0+1}+ 
\zeta(z_{(2\nu_0)}),\,
\ \ \ \ \ \ \ \ \ 
\right .\\
\ \ \ \ \ \ \ \ \
\ \ \ \ \ \ \ \ \
\ \ \ \ \ \ \ \ \ 
&
\left.
\Theta \left( z_{2\nu_0+1}+ 
\zeta(z_{(2\nu_0)}),\,
z(z_{(2\nu_0)}),\, 
w(z_{(2\nu_0)})\right)
\right).
\endaligned\right.
\end{equation}
Puisque $\pi_t(z,\, w,\, \zeta,\, \xi)= (z,\, w)$, la premi\`ere
relation~\thetag{ 3.25} est \'evidente. La seconde se v\'erifie de
mani\`ere analogue. 

Au total, on obtient deux submersions \`a valeurs dans $\C_t^n$ et
dans $\C_\tau^n$.

\def\thecorollary{3.29}\begin{corollary} 
Si $M$ est minimale \`a l'origine, il existe un entier $\nu_0 \leq d
+1$, le type de Segre de $M$ \`a l'origine, et il existe des
points $\underline{ z}_{ (2 \nu_0 )}^* \in \C^{ 2m \nu_0}$ et $z_{ (2
\nu_0 )}^* \in \C^{ 2m \nu_0}$ arbitrairement proches de l'origine
tels que les deux applications
\def\theequation{3.30}\begin{equation}
\left\{
\aligned
\C^{2m\nu_0}
\ni z_{(2\nu_0)} 
& \
\longmapsto \pi_t \left(\underline{ 
\Gamma}_{2\nu_0} ( 
z_{ (2\nu_0)})\right) \in \C^n 
\ \ \ \ \ {\rm et} \\
\C^{2m \nu_0}
\ni z_{(2\nu_0)} 
& \
\longmapsto \pi_\tau \left( 
\Gamma_{2\nu_0} (
z_{ (2\nu_0)})\right) \in \C^n
\endaligned\right.
\end{equation}
sont de rang $n$ et s'annulent aux deux points $\underline{ z}_{ (2
\nu_0)}^*$ et $z_{(2\nu_0)}^*$.
\end{corollary}

\subsection*{ 3.31.~R\'esum\'e~: utilisation concr\`ete
de la minimalit\'e} Nous sommes maintenant en mesure 
de ramener les th\'eor\`emes de convergence annonc\'es
dans l'Introduction \`a des propri\'et\'es
de convergence sur les cha\^{ \i}nes de
Segre conjugu\'ees. 

\def\thelemma{3.32}\begin{lemma}
Pour d\'emontrer le Th\'eor\`eme~1.39, il suffit d'\'etablir que pour
tout $k \in \N$, les applications formelles $z_{ (k)} \longmapsto_{
\mathcal{ F}} \, h \left( \pi_t\left( \underline{ \Gamma }_k (z_{
(k)}) \right) \right)$ sont convergentes, {\it i.e.} appartiennent \`a
$\C \{ z_{ (k)} \}^{ n'}$. De m\^eme, pour d\'emontrer le Th\'eor\`eme
principal~1.23, il suffit d'\'etablir que les applications formelles
$z_{ (k)} \longmapsto_{ \mathcal{ F }} \ \xi' - \Theta' \left(
\zeta',\, \, h \left( \pi_t \left( \underline{ \Gamma }_k( z_{ (k)})
\right) \right) \right) \in \C \{ \zeta',\, z_{ (k)} \}^{ d'}$ sont
convergentes pour tout $k \in \N$. Ces deux propri\'et\'es sont
satisfaites en rempla\c cant $\underline{ \Gamma}_k$ par $\Gamma_k$.
\end{lemma}

Ce sont effectivement ces propri\'et\'es de convergence qui
appara\^{\i}tront naturellement dans les Sections~5, 6 et~7
ci-dessous. Par souci d'\'el\'egance et de sym\'etrie, nous
travaillerons simultan\'ement avec les deux composantes $h$ et
$\overline{ h}$ de l'application complexifi\'ee $h^c = ( h,\,
\overline{ h })$~: sans composer avec les projections $\pi_t$ et
$\pi_\tau$, nous d\'emontrerons que $z_{ (k)}
\longmapsto_{ \mathcal{ F}} \ h^c \left( \underline{ \Gamma }_k
( z_{ (k) }) \right)$ converge pour tout $k$

\proof
En effet, prenons $k := 2 \nu_0$ et supposons l'application formelle
\def\theequation{3.33}\begin{equation}
z_{ (2 \nu_0) } \longmapsto_{ \mathcal{ F}} \, h \left( \pi_t \left(
\underline{ \Gamma}_{ 
2 \nu_0} (z_{ ( 2\nu_0 ) }) \right) \right)=: H
(z_{ (2\nu_0)})
\end{equation}
convergente. Puisque l'on peut choisir un point $\underline{ z }_{ (2
\nu_0) }^*$ arbitrairement proche de l'origine o\`u la premi\`ere
application~\thetag{ 3.30} est submersive, assurons-nous que ce point
$z_{ (2 \nu_0)}^*$ appartient au domaine de convergence normale de $H
( z_{ (2\nu_0)})$. Soit $s:= (s_1,\, \dots,\, s_n ) \in \C^n$ et soit
$s \mapsto \phi_{ 2\nu_0} (s) \in \C^{ 2m\nu_0}$ une application
affine satisfaisant $\phi_{ 2\nu_0} (0) = \underline{z}_{ 2\nu_0}^*$
telle que l'application holomorphe locale $s \longmapsto \pi_t \left(
\underline{ \Gamma}_{ 2\nu_0} (\phi_{ 2\nu_0}(s)) \right)$ est de rang
$n$ en $s=0$. On d\'eduit d'abord que
\def\theequation{3.34}\begin{equation}
h \left( \pi_t \left( \underline{ \Gamma}_{ 2 \nu_0} (\phi_{ ( 2\nu_0
) } (s) ) \right) \right) \equiv 
H( \phi_{
2\nu_0} (s))\in \C\{ s\}^{ n'}
\end{equation} 
est convergente. Ensuite, puisque $s \longmapsto
\pi_t \left( \underline{ \Gamma}_{ 2\nu_0} (\phi_{ 2\nu_0}(s))
\right)$ est de rang $n$ en $s=0$, on conclut que
$h (t) \in \C \{ t\}^{ n'}$ est convergente. 

La seconde assertion du Lemme~3.32 se v\'erifie
de mani\`ere similaire.
\endproof

\section*{\S4.~Jets de sous-vari\'et\'es de Segre et application
de r\'eflexion}

\subsection*{4.1.~Pr\'eliminaire}
Dans l'article fondamental~\cite{ we1978}, g\'en\'eralis\'e plus tard
en codimension sup\'erieure ({\it cf.}~\cite{ we1982}), S.M.~Webster a
introduit l'application qui, \`a un point $p$ d'une hypersurface
analytique r\'eelle $M$ de $\C^n$, associe le plan tangent complexe
\`a $M$ en $p$, c'est-\`a-dire $p\mapsto (p,\, T_p^cM)$, et il a
observ\'e que l'image de l'hypersurface $M$ est une sous-vari\'et\'e
analytique totalement r\'eelle de $\C^{2n-1}$ si et seulement si $M$
est Levi non-d\'eg\'en\'er\'ee. Puisque la sous-vari\'et\'e de Segre
$S_{\bar p}$ passant par $p\in M$ admet comme espace tangent en $p$ le
m\^eme sous-espace $T_p^cM$, cette application s'identifie avec
l'application qui, \`a un point $p$ de $M$, associe le jet d'ordre $1$
de la sous-vari\'et\'e de Segre passant par $p$. Dans~\cite{ dw1980},
K.~Diederich et S.M.~Webster ont g\'en\'eralis\'e cette id\'ee en
introduisant les jets d'ordre arbitraire des sous-vari\'et\'es de
Segre~; ils ont ainsi exhib\'e des conditions nouvelles de
non-d\'eg\'en\'erescence, plus g\'en\'erales que la Levi
non-d\'eg\'en\'erescence. Puisqu'aucun travail de fondation n'a
\'et\'e publi\'e jusqu'\`a pr\'esent pour d\'ecrire ces concepts, nous
entreprenons de r\'esumer ici les \'el\'ements d'une th\'eorie
autonome des jets de sous-vari\'et\'es de Segre. Nous \'elaborerons
ainsi un point de vue alternatif aux calculs non g\'eom\'etriques du
Chapitre~11 de~\cite{ ber1999a}.

\subsection*{4.2.~D\'efinitions fondamentales}
Soit $\mathcal{ M}$ la complexification extrins\`eque d'une
sous-vari\'et\'e analytique r\'eelle g\'en\'erique d\'efinie par les
deux jeux sym\'etriques d'\'equations~\thetag{ 3.6}. Rappelons que la
sous-vari\'et\'e de Segre complexifi\'ee $\underline{ \mathcal{ S}}_t$
peut \^etre consid\'er\'ee comme le sous-ensemble de $\C^n$ d\'efini
par $\{\tau= (\zeta,\, \xi)\in \C^n: \, \xi = \Theta (\zeta,\, t)\}$,
o\`u $t$ est consid\'er\'e comme fix\'e. Elle est param\'etr\'ee par
$\zeta \in \C^m$. Introduisons alors l'application des jets d'ordre
$k$ de $\underline{ \mathcal{ S}}_t$ en l'un de ses points $(\zeta,\,
\Theta(\zeta,\, t))$, qui est d\'efinie pr\'ecis\'ement par
\def\theequation{4.3}\begin{equation}
J_\tau^k\underline{\mathcal{S}}_t:=\left( \zeta,\
\left({1\over\beta!}\, \partial_\zeta^\beta
\Theta_j(\zeta,\, t)\right)_{1\leq j\leq d,\, \vert \beta \vert \leq
k}\right).
\end{equation}
Elle est \`a valeurs dans $\C^{m+ N_{ d,\, m,\,k}}$ o\`u $N_{
d,\,m,\,k} := d \, \frac{ ( m+k )! }{ m! \ k!}$ est le nombre de
d\'eriv\'ees partielles qui apparaissent \`a droite de la virgule du
second membre de~\thetag{ 4.3}. Remarquons que les termes $(\zeta_1
,\, \dots,\, \zeta_m)$ apparaissent comme premi\`eres composantes de
l'application~\thetag{ 4.3}. Remarquons aussi que si $k_2 \geq k_1$ et
si $\pi_{k_2,\,k_1}$ d\'esigne la projection caonique $\C^{ m+ N_{d
,\,m ,\,k_2}} \to \C^{ m+N_{d,\,m,\,k_1}}$, on a~: $\pi_{ k_2,\, k_1}(
J_\tau^{ k_2} \underline{ \mathcal{S }}_t)= J_\tau^{ k_1} \underline{
\mathcal{ S }}_t$. Nous noterons dans la suite cette application par
$(t ,\, \tau)\mapsto J_\tau^k \underline{ \mathcal{ S }}_t$, o\`u nous
sous-entendons que $(t ,\, \tau)\in \mathcal{ M}$.

De mani\`ere analogue, introduisons l'application des jets d'ordre $k$
de la sous-vari\'et\'e de Segre complexifi\'ee $\mathcal{S}_\tau$ en
l'un de ses points $\left( z,\, \overline{ \Theta} ( z, \tau)
\right)$. Elle est d\'efinie pr\'ecis\'ement par
\def\theequation{4.4}\begin{equation}
J_t^k\mathcal{S}_\tau:=\left(
z,\,\left({1\over\beta!}\, 
\partial_z^\beta 
\overline{\Theta}_j(z,\,\tau)\right)_{1\leq j\leq d,\,\, 
\vert \beta \vert \leq k}\right).
\end{equation}
Le lien de sym\'etrie entre ces deux applications
est tr\`es simple~:
\def\theequation{4.5}\begin{equation}
J_{\bar t}^k\underline{\mathcal{S}}_{\bar \tau} \equiv
\overline{
J_t^k\mathcal{S}_\tau
}
\end{equation}
En d'autres termes, le diagramme suivant est commutatif~:
$$
\diagram \mathcal{M} \rto^{\sigma} 
\dto_{J_\bullet^k\mathcal{S}_\bullet} 
& \mathcal{M} 
\dto^{J_\bullet^k\underline{\mathcal{S}}_\bullet} \\
\C^{m+N_{d,\,m,\,k}} \rto^{(\overline{\bullet})} & 
\C^{m+N_{d,\,m,\,k}}
\enddiagram,
$$
o\`u $(\overline{ \bullet})$ d\'esigne l'op\'erateur de conjugaison
complexe. Puisque ces deux applications de jets sont essentiellement
\'equivalentes, il n'\'etait pas restrictif de pr\'esenter dans le
\S1.7 les cinq conditions de non-d\'eg\'en\'erescence {\bf (nd1)},
{\bf (nd2)}, {\bf (nd3)}, {\bf (nd4)} et {\bf (nd5)} seulement \`a
partir de $J_\tau^k \underline{ \mathcal{ S }}_t$.

\subsection*{4.6.~Invariance biholomorphe de l'application de jets}
Soit $t'=h(t)$ un changement de coordonn\'ees holomorphes locales
centr\'e \`a l'origine. Soit $\mathcal{ M}':= h^c ( \mathcal{ M})$
l'image de la complexification $\mathcal{ M}$ de $M$ par ce changement
de coordonn\'ees. Puisque $h^c$ est inversible, la codimension de
$\mathcal{ M'}$ est la m\^eme que celle de $\mathcal{ M}$. Soient
$\xi_j'= \Theta_j' ( \zeta',\, t')$, $j=1,\, \dots,\,d$, $\zeta'\in
\C^m$, $\xi'\in \C^m$, $t'\in \C^n$, des \'equations complexes
graph\'ees pour $\mathcal{ M}'$ dans un syst\`eme de cordonn\'ees
analogue au syst\`eme de coordonn\'ees dans lequel $\mathcal{ M}$ est
repr\'esent\'ee. Par hypoth\`ese, $h_i(t) \in \C \{ t\}$ avec $h_i
(0)=0$, pour $i=1,\,\dots,\, n$, et il existe une matrice $b(t,\,
\tau)$ de taille $d \times d$ de s\'eries enti\`eres convergentes
telle que l'on a l'identit\'e formelle vectorielle $r' \left( h(t),\,
\overline{ h} ( \bar t) \right) \equiv b(t,\, \tau)\, r (t,\, \tau)$
dans $\C \{ t,\, \tau\}^d$, o\`u l'on a pos\'e
\def\theequation{4.7}\begin{equation}
r_j(t,\, \tau):= 
\xi_j- \Theta_j (\zeta,\, t), 
\ \ \ \ \
{\rm et} 
\ \ \ \ \
r_j' \left( t',\, \tau' \right) := 
\xi_j'- \Theta_j' \left( \zeta',\, t' \right),
\end{equation}
pour $j=1,\, \dots ,\, d$. Puisque nous scindons les coordonn\'ees
$t'= (z',\, w')\in \C^m\times \C^d$ et $\tau'= (\zeta',\, \xi')\in
\C^m \times \C^d$ en deux groupes, divisons aussi les composantes du
biholomorphisme $h$ en deux groupes $h(t) := (f(t),\, g(t))\in \C \{
t\}^m \times \C\{ t\}^d$. En rempla\c cant $\xi$ par $\Theta (\zeta,\,
t)$ dans l'identit\'e fondamentale $r' \left( h(t),\, \overline{ h }
(\bar t) \right) \equiv b(t,\, \tau)\, r (t,\, \tau)$, le membre de
droite s'annule identiquement et nous obtenons les $d$ identit\'es
formelles suivantes, valables dans $\C\{ \zeta,\, t\}$ et qui seront
notre point de d\'epart
\def\theequation{4.8}\begin{equation}
\overline{ g}_j \left( \zeta,\,\Theta(\zeta,\,t) \right)
\equiv \Theta_j'\left(\overline{ f}
(\zeta,\,\Theta(\zeta,\,t)),\,
h(t)\right), \ \ \ \ \ \ j=1,\,\dots,\,d.
\end{equation}
En diff\'erentiant ces identit\'es une infinit\'e de fois par rapport
\`a $\zeta \in \C^m$, nous allons \'etablir les relations suivantes,
dont l'apparence technique ne doit pas cacher qu'elles ont une
signification g\'eom\'etrique cruciale que nous expliquons ci-apr\`es.

\def\thelemma{4.9}\begin{lemma}
Pour tout $j=1,\, \dots,\, d$ et tout $\beta \in \N^m$, il existe une
application rationnelle $Q_{j,\, \beta}$ dont l'expression explicite,
qui ne d\'epend ni de $\mathcal{ M}$, ni de $h$, ni de $\mathcal{
M}'$, peut \^etre calcul\'ee gr\^ace \`a des formules combinatoires
universelles\footnotemark[6], telle que l'on a les identit\'es
formelles suivantes, valables dans $\C \{ \zeta,\, t\}$~{\rm :}
\def\theequation{4.10}\begin{equation}
\left\{
\aligned
{}
&
{1\over \beta!}\, 
{\partial^{\vert \beta\vert} \Theta_j'\over
\partial (\zeta')^\beta}\left(
\overline{ f} (\zeta,\,\Theta(\zeta,\,t)), \
h(t)\right)\equiv\\
&
\ \ \ \ \
\equiv
Q_{j,\,\beta}\left(
\left(\partial_\zeta^{\beta_1}\Theta_{j_1}(\zeta,\,t)
\right)_{1\leq j_1\leq d, \,
\vert \beta_1\vert \leq \vert \beta\vert}, \,
\left(\partial_\tau^{\alpha_1}
\overline{ h}_{i_1}(\zeta,\,\Theta(\zeta,\,t))
\right)_{1\leq i_1\leq n,\, 
\vert \alpha_1 \vert\leq 
\vert \beta \vert}\right) \\
&
\ \ \ \ \ 
=:
R_{j,\, \beta} \left(\zeta,\ 
\left(\partial_\zeta^{\beta_1}\Theta_{j_1}(\zeta,\,t)
\right)_{1\leq j_1\leq d, \,
\vert \beta_1\vert \leq \vert \beta\vert}
\right) \\
&
\ \ \ \ \ 
=:
q_{j,\, \beta} (\zeta,\, t),
\endaligned\right.
\end{equation}
o\`u l'avant-derni\`ere \'equation d\'efinit la s\'erie enti\`ere
$R_{j,\, \beta}$ \`a partir de $Q_{j,\, \beta}$ par simple oubli de la
d\'ependance en les jets de $\overline{ h}$, et la derni\`ere
d\'efinit la s\'erie enti\`ere $q_{j,\, \beta}$ \`a partir de $Q_{j,\,
\beta}$, par simple oubli de toute d\'ependance par rapport aux jets
de $\Theta$ ou de $\overline{ h}$. Ici, les s\'eries enti\`eres
$Q_{j,\, \beta}$ sont holomorphes au voisinage du jet constant
obtenu en posant $(\zeta,\, t) = (0,\, 0)$, {\it i.e.}
$\left(( \partial_\zeta^{ \beta_1} \Theta_{j_1} (0,\,0) )_{1 \leq j_1
\leq d,\, \, \vert \beta_1 \vert \leq \vert \beta\vert},
\,(\partial_\tau^{ \alpha_1} \overline{ h}_{ i_1}(0,\,0) )_{1 \leq
i_1\leq n,\, \vert \alpha_1 \vert \leq \vert \beta \vert}\right)$.
Des relations sym\'etriques analogues sont satisfaites en rempla\c
cant $\Theta$, $\Theta'$, $\zeta$, $t$, $\overline{ f}$, $h$ par
$\overline{ \Theta}$, $\overline{ \Theta}'$, $z$, $\tau$ $f$,
$\overline{ h}$.
\end{lemma}

\footnotetext[6]{
Mais le travail d'explicitation compl\`ete serait particuli\`erement
fastidieux, car quatre ingr\'edients se combineraient ensemble dans le
calcul~: la formule g\'en\'eralis\'ee de Fa\`a di Bruno pour la
diff\'erentiation de fonctions de plusieurs variables compos\'ees~; le
d\'eveloppement formel des d\'erivations compos\'ees $( \underline{
\mathcal{ L}}_1)^{ \beta_1} ( \underline{ \mathcal{ L}}_2)^{ \beta_2}
\cdots ( \underline{ \mathcal{ L}}_m)^{ \beta_m}$, o\`u 
$(\beta_1,\,\beta_2,\, \dots,\, 
\beta_m)\in \N^m$~; la formule universelle pour
les d\'erivations partielles d'ordre quelconque d'un quotient de deux
s\'eries enti\`eres~; les d\'eterminants de Cramer. Heureusement,
nous n'aurons pas besoin ici d'une formule explicite pour les
s\'eries enti\`eres analytiques $Q_{j,\, \beta}$. }

L'existence des s\'eries enti\`eres $R_{j,\, \beta}$ exprime que
l'application~\thetag{ 4.3} des jets d'ordre $k$ des sous-vari\'et\'es
de Segre complexifi\'ees conjugu\'ees est invariante par
biholomorphisme, c'est-\`a-dire que le diagramme suivant est
commutatif~:
$$
\diagram 
\mathcal{M} \rto^{h} 
\dto_{J_\bullet^k\underline{\mathcal{S}}_\bullet} 
& \mathcal{M}'
\dto^{J_\bullet^k\underline{\mathcal{S}}_\bullet'} \\
\C^{m+N_{d,\,m,\,k}} \rto^{ R^k(h)} & 
\C^{m+N_{d,\,m,\,k}}
\enddiagram,
$$
o\`u l'application $R^k(h)$, qui d\'epend de $h$, est d\'efinie par
ses composantes $R_{j,\,\beta}$ pour $j=1,\,\dots,\, d$ et $\vert
\beta \vert \leq k$. Comme $h$ est inversible, on v\'erifie que la
transformation associ\'ee $R^k(h)$ est elle aussi un biholomorphisme
local. Gr\^ace \`a l'inversibilit\'e de $R^k (h)$, on d\'emontre sans
difficult\'e que les cinq conditions de non-d\'eg\'en\'erescence {\bf
(nd1)}, {\bf (nd2)}, {\bf (nd3)}, {\bf (nd4)} et {\bf (nd5)}
introduites dans le \S1.7 sont invariantes par biholomorphisme local
({\it voir}~le Chapitre~3 de~\cite{ me2003} pour les d\'etails).

\proof
Nous allons diff\'erentier les \'equations~\thetag{ 4.8} par rapport
\`a $\zeta_k$, pour $k=1,\,\dots,\,m$. Rappelons ici les expressions
explicites des champs de vecteurs CR de type $(0,\, 1)$ complexifi\'es
$\underline{ \mathcal{ L }}_k$ d\'efinis pr\'ec\'edemment~:
\def\theequation{4.11}\begin{equation}
\underline{\mathcal{L}}_k=
{\partial \over \partial \zeta_k}+\sum_{j=1}^d \, 
{\partial \Theta_j\over \partial \zeta_k}(\zeta,\,t)\, 
{\partial \over \partial \xi_j},
\end{equation}
pour $k=1,\, \dots,\, m$. On voit imm\'ediatement que le fait de
diff\'erenter une s\'erie enti\`ere compos\'ee $\psi ( \zeta,\,
\xi)\vert_{\xi= \Theta( \zeta,\, t)} = \psi( \zeta,\, \Theta(
\zeta,\,t))$ par rapport \`a $\zeta_k$ \'equivaut \`a lui appliquer le
champ de vecteurs $\underline{ \mathcal{ L }}_k$, entendu comme
d\'erivation, c'est-\`a-dire que l'on a~:
\def\theequation{4.12}\begin{equation}
\frac{ \partial }{\partial \zeta_k} \, 
\psi(\zeta,\, \Theta(\zeta,\, t)) \equiv
\left[
\underline{ \mathcal{ L}}_k \psi \right] 
(\zeta,\, \Theta(\zeta,\, t)).
\end{equation}
Cette observation est triviale, mais elle a son importance~;
g\'eom\'etriquement parlant, appliquer l'op\'erateur $\underline{
\mathcal{ L }}_k$ signifie que l'on diff\'erentie le long de la
sous-vari\'et\'e de Segre complexifi\'ee conjugu\'e $\underline{
\mathcal{ S }}_t$, le tout \'etant param\'etr\'e par $t$.

Ainsi, en diff\'erentiant les relations~\thetag{ 4.8} par rapport
\`a $\zeta_k$ et en utilisant la formule de d\'erivation compos\'ee,
on obtient les relations
\def\theequation{4.13}\begin{equation}
\underline{\mathcal{L}}_k \overline{ g}_j
(\zeta,\,\Theta(\zeta,\,t))\equiv
\sum_{l=1}^m\, \underline{\mathcal{L}}_k
\overline{f}_l
(\zeta,\,\Theta(\zeta,\,t))\,
{\partial \Theta_j'\over 
\partial \zeta_l'} \left(
\overline{ f}
(\zeta,\,\Theta(\zeta,\,t)),\, h(t) \right),
\end{equation}
pour $k = 1,\, \dots,\, m$ et $j = 1,\, \dots,\, d$. 

Par ailleurs, en posant $t=0$ dans~\thetag{ 4.8}, on obtient
\def\theequation{4.14}\begin{equation}
\overline{ g}_j (\zeta,\, \Theta(\zeta,\, 0))\equiv
\Theta_j' \left(
\overline{ f}(\zeta,\, \Theta(\zeta,\, 0)), \, 0
\right). 
\end{equation}
Ces identit\'es expriment que la restriction de $h^c=(h,\, \overline{
h})$ \`a la sous-vari\'et\'e de Segre complexifi\'ee conjugu\'ee
$\underline{ \mathcal{ S}}_0 =\{ (0,\, 0,\, \zeta,\, \Theta(\zeta,\,
0)): \, \zeta \in \C^m \}$, qui est de dimension complexe $m$, est \`a
valeurs dans la sous-vari\'et\'e de Segre complexifi\'ee conjugu\'ee
$\underline{ \mathcal{ S}}_0' =\{ \left( 0,\, 0,\, \zeta',\,
\Theta '(\zeta',\, 0) \right) : \, \zeta' \in \C^m\}$, qui est de m\^eme
dimension. Puisque $h^c$ est un biholomorphisme, cette restriction
est forc\'ement inversible, c'est-\`a-dire que
l'application $\zeta \mapsto \overline{ h} ( \zeta,\, \Theta (\zeta,\,
0))$ est de rang $m$ en $\zeta = 0$. Autrement dit, la matrice
(constante) des d\'eriv\'ees partielle \`a l'origine $\left(
\underline{ \mathcal{ L}}_k \, \overline{ h}_i (0) \right)_{1 \leq i
\leq n}^{ 1 \leq k\leq m}$ est de rang $m$. Mais \`a cause des
relations~\thetag{ 4.13}, prises en $t=0$, les $d$ derni\`eres lignes
de la matrice $\left( \underline{ \mathcal{ L}}_k \, \overline{ h}_i
(0) \right)_{1 \leq i \leq n}^{ 1 \leq k\leq m}$ sont des combinaisons
lin\'eaires \`a coefficients constants de ses $m$ premi\`eres lignes.
Il en d\'ecoule finalement que le d\'eterminant suivant de taille
$m\times m$ ne s'annule pas~:
\def\theequation{4.15}\begin{equation}
{\rm det}\, 
\left(\underline{\mathcal{L}}_{k_1}
\overline{ f}_{k_2}(0)\right)_{
1\leq k_1,\,k_2\leq m}\neq 0.
\end{equation}
Par cons\'equent, nous pouvons diviser
par le d\'eterminant 
\def\theequation{4.16}\begin{equation}
\mathcal{D}(\zeta,\,t):=
{\rm det}\, 
\left(\underline{\mathcal{L}}_{k_1} 
\overline{ f}_{k_2}(\zeta,\,\Theta(\zeta,\,t))\right)_{
1\leq k_1,\,k_2\leq m},
\end{equation}
pourvu que $\zeta\in \C^m$ et $t\in \C^n$ soient suffisamment petits
pour qu'il ne s'annule pas.

Maintenant, envisageons les \'equations~\thetag{ 4.13} \`a $j$ fix\'e
comme un syst\`eme lin\'eaire non homog\`ene dont les inconnues sont
les $m$ d\'eriv\'ees partielles $\partial \Theta_j'/ \partial
\zeta_1',\, \dots,\, \partial \Theta_j'/ \partial \zeta_m'$. Gr\^ace
aux formules de Cramer, nous pouvons r\'esoudre ce syst\`eme, ce qui
donne des expressions de la forme
\def\theequation{4.17}\begin{equation}
{\partial \Theta_j'\over \partial \zeta_k'} \left(
\overline{ f}
(\zeta,\,\Theta(\zeta,\,t)),\, h(t)\right) \equiv
{T_{j,\,k}\left( \left(\underline{\mathcal{L}}_{k_1'} 
\overline{h}_{ i_1}(\zeta,\,
\Theta(\zeta,\,t)) \right)_{1\leq i_1\leq n, \, 
1\leq k_1'\leq m}\right)\over
\mathcal{D}(\zeta,\,t)}.
\end{equation}
Ici, en examinant la forme explicite des d\'eterminants de Cramer, on
v\'erifie ais\'ement que les termes $T_{j,\,k}$ sont des polyn\^omes
universels en leurs variables.

Traitons maintenant le cas $\vert \beta \vert = 2$ de~\thetag{ 4.10}.
Pour un mutiindice arbitraire $\beta = (\beta_1,\, \beta_2,\, \dots,\,
\beta_m) \in \N^m$, nous noterons $\underline{ \mathcal{ L} }^\beta$
la d\'erivation d'ordre $\vert \beta\vert$ d\'efinie par la
composition $(\underline{ \mathcal{ L }}_1 )^{ \beta_1}(\underline{
\mathcal{ L }}_2 )^{ \beta_2} \cdots ( \underline{ \mathcal{L }}_m)^{
\beta_m}$.

Diff\'erentions \`a nouveau les identit\'es~\thetag{ 4.13} par rapport
aux variables $\zeta_k$. Nous obtenons \`a nouveau des syst\`emes de
Cramer de m\^eme d\'eterminant $\mathcal{ D} (\zeta,\, t)$ que l'on
peut r\'esoudre. Aussi, pour toute paire d'entiers $(k_1,\, k_2)$ avec
$1\leq k_1,\,k_2\leq m$ et pour tout $j=1,\, \dots,\,d$, il existe des
polyn\^omes universels $T_{j,\,k_1,\,k_2}$ tels que l'on peut \'ecrire
\def\theequation{4.18}\begin{equation}
{\partial^2\Theta_j'\over 
\partial \zeta_{k_1}'\partial \zeta_{k_2}'}
\left(\overline{ f}
(\zeta,\,\Theta(\zeta,\,t)),\,h(t)\right)\equiv
{T_{j,\,k_1,\,k_2}\left( \left(
\underline{ \mathcal{ L}}^{\beta_1} 
\overline{ h}_{i_1}(\zeta,\,
\Theta(\zeta,\,t))\right)_{
1\leq i_1\leq n, \, 1\leq \vert \beta_1 \vert \leq 2}\right)\over
\mathcal{D}(\zeta,\,t)^3}.
\end{equation}
Le lecteur aura remarqu\'e l'exposant $3$ qui appara\^{\i}t au
d\'enominateur. On doit le comprendre comme
``$3$''$=$``$2$''$+$``$1$'', o\`u ``$2$'' provient de la
diff\'erentiation du quotient $T_{j,\,k }/ \mathcal{D}$
dans~\thetag{4.17} et o\`u ``$1$'' provient de la seconde application
de la r\`egle de Cramer.

Diff\'erentions successivement les relations~\thetag{ 4.8} par rapport
\`a $\zeta^\beta= \zeta_1^{ \beta_1} \zeta_2^{ \beta_2} \cdots
\zeta_m^{ \beta_m}$. En raisonnant par r\'ecurrence, 
on d\'emontre que pour tout $\beta\in\N^m$ et pour tout
$j=1,\,\dots,\,d$, il existe un polyn\^ome $T_{j,\,\beta}$
dont l'expression explicite est complexe, mais universelle, 
tel que l'identit\'e suivante est satisfaite, dans 
$\C\{ \zeta,\, t\}$~:
\def\theequation{4.19}\begin{equation}
{1\over \beta!}\, 
{\partial^{\vert \beta\vert} \Theta_j'\over
\partial (\zeta')^\beta}
\left(\overline{ 
f}(\zeta,\,\Theta(\zeta,\,t)),\,h(t)\right)
\equiv
{T_{j,\,\beta}\left( \left(\underline{\mathcal{L}}^{\beta_1}
\overline{ h}_{
i_1}(\zeta,\,\Theta(\zeta,\,t))\right)_{1\leq i_1\leq n,\,
1\leq \vert \beta_1 \vert \leq \vert \beta\vert}\right)\over
[\mathcal{D}(\zeta,\,t)]^{2\vert \beta\vert -1}}.
\end{equation} 
\`A ce stade, il est important d'observer qu'en d\'eveloppant les
d\'erivations multiples $\underline{ \mathcal{ L}}^{ \beta_1 }$, o\`u
$\beta_1 \in \N^m$ est un multiindice fix\'e, on obtient des
op\'erateurs diff\'erentiels \`a coefficients non constants. Plus
pr\'ecis\'ement, en utilisant les expressions explicites~\thetag{
4.11} des champs de vecteurs $\underline{ \mathcal{ L} }_k$, et en
raisonnant par r\'ecurrence, il est ais\'e d'\'etablir que les
coefficients des op\'erateurs $\underline{\mathcal{L}}^{\beta_1}$,
exprim\'es par rapports aux op\'erateurs $\partial_t^\alpha$, sont des
polyn\^omes universels en les d\'eriv\'ees partielles $\left(
\partial^{ \vert \beta_2 \vert} \Theta_{ j_2}( \zeta,\,t)/\partial
\zeta^{ \beta_2} \right)_{ 1 \leq j_2\leq d, \, 1 \leq \vert
\beta_2\vert \leq \vert \beta_1 \vert}$. Ainsi, le num\'erateur
de~\thetag{4.19} devient un polyn\^ome universel (dont l'expression
explicite ne nous serait pas utile) en les variables
\def\theequation{4.20}\begin{equation}
\left(
\left(\partial_\zeta^{\beta_1}\Theta_{j_1}(\zeta,\,t)
\right)_{1\leq j_1\leq d, \,
\vert \beta_1\vert \leq \vert \beta\vert}, \,
\left(\partial_t^{\alpha_1}
\overline{ h}_{i_1}(\zeta,\,\Theta(\zeta,\,t))
\right)_{1\leq i_1\leq n,\, 
\vert \alpha_1 \vert\leq 
\vert \beta \vert}\right)
\end{equation} 
De m\^eme, en d\'eveloppant le d\'eterminant qui appara\^{\i}t au
d\'enominateur, on obtient un polyn\^ome universel en les
variables~\thetag{ 4.20}, pour $\vert \beta_1 \vert \leq 1$. Concluons
que nous avons effectivement construit l'application rationnelle
souhait\'ee $Q_{j,\, \beta} := T_{j,\, \beta} / [\mathcal{ D }]^{ 2\vert
\beta \vert - 1}$.
\endproof

Il nous reste maintenant \`a \'etablir l'invariance biholomorphe
de l'application de r\'eflexion CR formelle.

\subsection*{4.21.~Inversion}
Dans cet objectif, d\'eveloppons les s\'eries enti\`eres $\Theta_j'$
par rapport aux puissances de $\zeta'$, ce qui donne $\Theta_j'
\left(\zeta',\, t' \right) =: \sum_{\beta \in \N^m}\, (\zeta')^\beta
\, \Theta_{j,\,\beta}' \left(t' \right)$, pour $j=1,\, \dots,\, d$,
o\`u les $\Theta_{j,\, \beta} ' (t')$ sont des s\'eries enti\`eres
analytiques complexes. Avec ces notations, nous pouvons d\'evelopper
le membre de gauche de~\thetag{ 4.10} en conservant pour membre de
droite la derni\`ere ligne de~\thetag{ 4.10}, ce qui donne~:
\def\theequation{4.22}\begin{equation}
\sum_{ \gamma \in \N^m}\, 
\frac{ (\beta+ \gamma) !}{ \beta ! \ \gamma !} \
\overline{ f} (\zeta,\, \Theta(\zeta,\, t))^\gamma \
\Theta_{j,\, \beta + \gamma} ' ( h(t)) \equiv 
q_{j,\, \beta} (\zeta,\, t).
\end{equation}
Puisque nous voulons \'etablir
l'invariance biholomorphe de l'application de r\'eflexion, nous 
devrions r\'esoudre les termes $\Theta_{ j,\, \beta}'(h(t))$ qui
apparaissent dans le d\'eveloppement de la premi\`ere ligne
de~\thetag{ 4.10} par rapport aux autres termes, de telle sorte que
n'apparaissent au final que des termes $\Theta_{j_1,\, \beta_1}(t)$.

Pour cela, il est utile de savoir qu'une collection 
infinie d'identit\'es formelles de la forme
\def\theequation{4.23}\begin{equation}
\sum_{\gamma\in\N^m}\, 
{(\beta+\gamma)! \over \beta! \ \gamma!}\,
{\zeta'}^\gamma \, 
\theta_{j,\,\beta+\gamma}'= q_{j,\,\beta},
\end{equation}
o\`u $j=1,\, \dots,\, d$, o\`u $\beta \in \N^m$, o\`u $\zeta' \in \C^m$
et o\`u les constantes $\theta_{j,\, \beta}' \in \C$ et $q_{j,\,
\beta}\in\C$ sont arbitraires et ind\'ependantes, peut \^etre
r\'esolue formellement par rapport aux inconnues $\theta_{j,\,
\beta}'$ au moyen d'une formule totalement similaire, \`a des
diff\'erences de signe pr\`es~:
\def\theequation{4.24}\begin{equation}
\sum_{\gamma\in\N^m}\,
(-1)^\gamma\, 
{(\beta+\gamma)! \over \beta! \ \gamma!}\, 
{\zeta'}^\gamma \,
q_{j,\,\beta+\gamma}=
\theta_{j,\,\beta}', 
\end{equation}
pour tout $j=1,\, \dots,\, d$ et tout $\beta \in \N^m$. On peut se
convaincre de cette formule d'inversion de deux mani\`eres~: ou bien
en rempla\c cant l'expression ~\thetag{ 4.24} de $\theta_{j,\,
\beta}'$ directement dans~\thetag{ 4.23} et en effectuant le calcul,
ou bien en raisonnant avec des s\'eries de Taylor convergentes en deux
points distincts et en extrapolant ensuite la relation obtenue aux
s\'eries formelles non convergentes.

L'application de cette formule
d'inversion aux identit\'es~\thetag{ 4.22} 
nous donne les relations formelles
\def\theequation{4.25}\begin{equation}
\small
\left\{
\aligned
{}
&
\Theta_{j,\,\beta}'(h(t))\equiv
\sum_{\gamma\in\N^m}\, 
(-1)^\gamma \, 
{(\beta+\gamma)!\over \beta! \ \gamma !}\ 
\overline{ f} (\zeta,\,\Theta(\zeta,\,t))^\gamma\, \cdot
q_{j,\, \beta + \gamma }(\zeta,\, t) 
\\
& \ 
\ \ \ \ \ \ \ \ \ \ \ \ \ 
\ \ \ \
\equiv
\sum_{\gamma\in\N^m}\, 
(-1)^\gamma \, 
{(\beta+\gamma)!\over \beta! \ \gamma !}\ 
\overline{ f} (\zeta,\,\Theta(\zeta,\,t))^\gamma\, \cdot
\\
&
\cdot Q_{j,\,\beta+\gamma}\left(
\left(\partial_\zeta^{\beta_1}\Theta_{j_1}(\zeta,\,t)
\right)_{1\leq j_1\leq d, \,
\vert \beta_1\vert \leq \vert \beta\vert+\vert \gamma \vert}, \,
\left(\partial_\tau^{\alpha_1}
\overline{ h}_{i_1}(\zeta,\,\Theta(\zeta,\,t))
\right)_{1\leq i_1\leq n,\, 
\vert \alpha_1 \vert\leq 
\vert \beta \vert+\vert \gamma \vert}
\right).
\endaligned\right.
\end{equation}
Ici, un argument de convergence est n\'ecessaire. Heureusement,
gr\^ace aux estim\'ees de Cauchy pour les d\'eriv\'ees partielles de
s\'eries enti\`eres analytiques, il existe des constantes $\sigma'>0$,
$\rho'>0$ et $C'>0$ telles qu'on a une estimation de la forme
\def\theequation{4.26}\begin{equation}
\left\vert
\frac{ 1}{\beta !} \, 
\frac{ \partial^{\vert \beta \vert} \Theta_j'}{
\partial (\zeta')^\beta} (\zeta',\, t') 
\right\vert <
C' \, (\rho')^{-\vert \beta \vert}, 
\end{equation}
pour tous $\zeta'$, $t'$ satisfaisant $\vert \zeta ' \vert < \sigma'$,
$\vert t' \vert < \sigma'$. Soit $\sigma >0$ une constante
suffisamment
petite pour que l'on ait $\vert h(t) \vert < \sigma'$ (d'o\`u aussi
$\left\vert \overline{ h} (\tau) \right\vert < 
\sigma'$) pour tout $t$ tel que
$\vert t \vert < \sigma$. On en d\'eduit la m\^eme estim\'ee de
Cauchy que~\thetag{ 4.26} pour le
membre de gauche de~\thetag{ 4.10}, 
avec les m\^emes constantes $\rho ' > 0$ et $C' > 0$. 
Ensuite, ces estimations se transf\`erent imm\'ediatement
aux membres de la derni\`ere ligne de~\thetag{ 4.10}~:
\def\theequation{4.27}\begin{equation}
\left\vert
q_{j,\, \beta} (\zeta,\, t) 
\right\vert < C' \, 
(\rho')^{-\vert \beta \vert}.
\end{equation}
Cette derni\`ere estimation suffit pour assurer que les
formules~\thetag{ 4.25} sont normalement convergentes lorsque
$\vert \zeta \vert < \sigma$ et $\vert t \vert < \sigma$.

\subsection*{4.28.~Formules fondamentales}
Pour terminer, posons $\zeta =0$ dans les formules~\thetag{ 4.25}, ce
qui exige quelques commentaires. Notons que les termes $\Theta_j (
0,\, t)$ s'identifient bien \'evidemment aux termes $\left. \Theta_{
j,\, \beta} (t) \right\vert_{ \beta =0}$ du d\'eveloppement $\Theta_j
( \zeta,\, t) =: \sum_{ \beta \in \N^m}\, \zeta^\beta \, \Theta_{ j,\,
\beta} (t)$ de $\Theta_j$ par rapport aux puissances de $\zeta$. Par
cons\'equent, les puissances $\overline{ f} ( 0,\, \Theta (0,\,
t))^\gamma$ et les jets $\partial_\tau^{ \alpha_1} \overline{ h }_{
i_1 }( \zeta,\, \Theta( 0,\, t))$ qui apparaissent dans~\thetag{
4.25}$\vert_{ \zeta=0}$ sont des s\'eries enti\`eres analytiques par
rapport aux $d$ variables $\Theta_{ j_1,\, 0 } (t)$, pour $j_1=1,\,
\dots,\, d$. Au total, {\it les identit\'es~\thetag{
4.25}$\vert_{\zeta = 0}$ ne font intervenir \`a droite que les
fonctions $\Theta_{j_1,\, \beta_1}(t)$}. Notons ici que l'indice
$\beta_1$ n'est pas born\'e, \`a cause de la pr\'esence de la somme
$\sum_{ \gamma \in \N^m}$.

En conclusion, si l'on note $S_{ j,\, \beta}$ ce second membre, on
obtient des formules absolument fondamentales entre les composantes
$\Theta_{ j,\, \beta}' (t')$ attach\'ees \`a l'\'equation de
$\mathcal{ M}' = h^c (\mathcal{ M})$ et les composantes $\Theta_{j,\,
\beta } (t)$ attach\'ees \`a l'\'equation de $\mathcal{ M}$, formules
que nous \'ecrirons~:
\def\theequation{4.29}\begin{equation}
\Theta_{j,\, \beta}'( h(t)) \equiv S_{j,\, \beta} \left( \left\{
\Theta_{j_1,\, \beta_1}(t) \right\}_{1\leq j_1\leq d,\, \beta_1\in
\N^m} \right).
\end{equation}
Ici, $S_{j,\, \beta}$ est une s\'erie enti\`ere qui d\'epend d'une
infinit\'e de variables. N\'eanmoins, il ne sera pas n\'ecessaire d'en
appeler \`a la th\'eorie des fonctions holomorphes d'une infinit\'e de
variables, pourvu que $S_{j,\, \beta}$ soit entendu comme d\'efini
par~\thetag{ 4.25}$\vert_{ \zeta = 0}$, avec les estim\'ees de
Cauchy~\thetag{ 4.27}$\vert_{ \zeta =0}$.

\subsection*{4.30.~Invariance biholomorphe
de la convergence de l'application de r\'eflexion CR formelle} Gr\^ace
aux relations~\thetag{ 4.29}, nous pouvons \'etablir l'assertion
qui fait imm\'ediatement suite \`a l'\'enonc\'e du Th\'eor\`eme~1.23.

En effet, comme dans les hypoth\`eses du Th\'eor\`eme~1.23, soit $h:
(M,\, 0) \longrightarrow_{ \mathcal{ F} } (M',\, 0)$ une application
CR formelle entre sous-vari\'et\'es g\'en\'eriques analytiques
r\'eelles de codimension $d$ et $d'$ dans $\C^n$ et dans $\C^{n'}$.
Soit un syst\`eme de coordonn\'ees holomorphes $(z',\, w')\in \C^{m' }
\times \C^{ d'}$ dans lequel la complexification $\mathcal{ M}'$ est
repr\'esent\'ee par les \'equations graph\'ees $\xi_{ j'}' = \Theta_{
j'} '( \zeta ',\, t')$, $j' = 1,\, \dots,\, d'$.

Supposons que l'application de r\'eflexion $\mathcal{ R}_h' (\tau',\,
t) := \xi' - \sum_{ \gamma' \in \N^{ m'}} \, (\zeta')^{ \gamma'} \,
\Theta_{ \gamma '} ' ( h(t))$ converge. De mani\`ere \'equivalente,
toutes les s\'eries formelles $\Theta_{j',\, \gamma '} '( h(t)) \in
\C\{ t\}$ conver\-gent et elles satisfont une estim\'ee de Cauchy telle
que $\left\vert \Theta_{j',\, \gamma '}' (h(t)) \right\vert < C' \,
(\rho')^{ -\vert \gamma ' \vert}$, pour $\vert t \vert < \sigma$.

Soit $t'' := \phi' (t')$ un biholomorphisme local fixant l'origine qui
transforme $M'$ en $M'' := \phi' (M')$. Apr\`es une renum\'erotation
\'eventuelle des coordonn\'ees $t''$, on peut supposer qu'elles se
scindent en $t''= (z'',\, w'') \in \C^{ m'} \times \C^{ d'}$ de telle
sorte que $T_0 M'' + \left( \{ 0 \} \times \C_{ w''}^{ d'} \right) =
T_0 \C_{ t''}^{ n'}$. Soient alors $\xi_{ j'}'' = \Theta_{ j'} ''
(\zeta'',\, t'')$, $j'= 1,\, \dots,\, d'$, des \'equations graph\'ees
pour la complexification $\mathcal{ M}''$. Le but est de d\'emontrer
que l'application de r\'eflexion transform\'ee $\mathcal{ R}_{ \phi'
\circ h}'' (\tau'',\, t):= \xi'' - \Theta '' ( \zeta'',\, \phi'(
h(t)))$ converge elle aussi.

Pour le v\'erifier, scindons aussi les composantes du changement
de coordonn\'ees comme suit~:
$\phi' (t') =: (\varphi' (t'),\, \psi' (t')) \in \C\{ t'\}^{ m'}
\times \C \{ t'\}^{ d'}$. Les relations~\thetag{ 4.29} pour la
transformation $\phi'$ s'\'ecrivent
\def\theequation{4.31}\begin{equation}
\Theta_{ j',\, \gamma'} '' \left(\phi'(t') \right) \equiv
S_{ j',\, \gamma'}'\left(
\left\{
\Theta_{ j_1',\, \gamma_1'}' (t')
\right\}_{1\leq j_1' \leq d',\, 
\gamma_1 ' \in \N^{ m'}}
\right).
\end{equation}
Rempla\c cons-y $t'$ par $h(t)$, ce qui donne~:
\def\theequation{4.32}\begin{equation}
\Theta_{ j',\, \gamma'} '' \left( \phi'(h(t)) \right) \equiv
S_{ j',\, \gamma'}'\left(
\left\{
\Theta_{ j_1',\, \gamma_1'}' (h(t))
\right\}_{1\leq j_1' \leq d',\, 
\gamma_1 ' \in \N^{ m'}}
\right).
\end{equation}
Puisque les $S_{ j',\, \gamma'}'$ sont analytiques et puisque les
$\Theta_{ j_1',\, \gamma_1'}' (h(t))$ convergent par hypoth\`ese, nous
d\'eduisons de~\thetag{ 4.32} que les composantes $\Theta_{ j',\,
\gamma'}'' \left(\phi' (h(t))\right)$ de l'application de r\'eflexion
$\mathcal{ R}_h''( \tau'',\, t) := \xi'' - \Theta ''( \zeta'',\,
h(t))$ convergent. Un examen de la construction des $S_{ j',\,
\gamma'}$ montre que toutes les estim\'ees de Cauchy n\'ecessaires
sont satisfaites. L'assertion est d\'emontr\'ee.

\subsection*{ 4.33.~R\'esum\'e}
En d\'efinitive, toutes les hypoth\`eses du Th\'eor\`eme
principal~1.23 et du Th\'eor\`eme~1.39 sont invariantes par changement
de coordonn\'ees. Sans perte de g\'en\'eralit\'e, nous pourrons donc
travailler avec un syst\`eme fix\'e de coordonn\'ees dans lesquelles
la complexification $\mathcal{ M}'$ est repr\'esent\'ee par $\xi' =
\Theta' ( \zeta',\, t')$.

\section*{\S5.~Convergence d'applications 
CR formelles finiment non-d\'eg\'en\'er\'ees} 

\subsection*{ 5.1.~Pr\'eliminaire} 
Afin de rendre plus accessible la lecture du Th\'eor\`eme~1.39 {\bf
(iii)} (Section~6) et du Th\'eor\`eme~1.23 (Section~7), nous
d\'emontrons ici le Th\'eor\`eme~1.39 {\bf (i)} (connu) en utilisant
notre propre formalisme~; cette section est exclusivement consacr\'ee
\`a cette d\'emonstration. Certaines observations combinatoires, qui
seront r\'eutilis\'ees ult\'erieurement, appara\^{\i}tront au cours
des raisonnements ({\it voir} par exemple le Lemme~5.4 ci-dessous).
Puisqu'elles sont explicit\'ees d'une mani\`ere parfois incompl\`ete
ou exag\'er\'ement r\'esum\'ee dans~\cite{ ber1998} et dans \cite{
ber1999a}, nous les d\'etaillerons soigneusement.

\subsection*{ 5.2.~R\'esolution de $h(t)$ en fonction 
des jets de $\overline{ h} (\tau)$} Soit $h^c : \, (\mathcal{ M},\, 0)
\longrightarrow_{ \mathcal{ F}} (\mathcal{ M}',\, 0)$ la
complexification de l'application CR formelle du Th\'eor\`eme~1.39
{\bf (i)}. Nous utiliserons les notations du \S1.25 et du \S1.33.

Rappelons que, d'apr\`es le \S1.33, en appliquant les d\'erivation
$\underline{ \mathcal{ L}}^\beta$ \`a la premi\`ere ligne de~\thetag{
1.27}, nous obtenons les identit\'es de r\'eflexion \'ecrites \`a la
premi\`ere ligne de~\thetag{ 1.31}, que nous r\'e\'ecrivons plus
pr\'ecis\'ement comme suit~:
\def\theequation{5.3}\begin{equation}
\underline{ \mathcal{ L}}^\beta \, 
\overline{ g}_{j'} (\zeta,\, \Theta (\zeta,\, t)) 
\equiv 
\sum_{ \gamma ' \in \N^{ m'}} \, 
\underline{ \mathcal{ L}}^\beta 
\overline{ f}^\gamma (\zeta,\, \Theta (\zeta,\, t)) \
\Theta_{j',\, \gamma'} ' (h(t)), 
\end{equation} 
pour tout $j'= 1,\, \dots,\, d'$ et 
tout $\beta \in \N^m$. Il est n\'ecessaire
d'examiner ce que donnent les termes $\underline{ \mathcal{ L}}^\beta
\overline{ g}_{ j'}$ et $\underline{ \mathcal{ L}}^\beta \overline{
f}^\gamma$.

Les lettres $\beta$ et $\gamma$ (resp. $\beta'$ et $\gamma'$)
d\'esigneront toujours des multiindices de $\N^m$ 
(resp. de~$\N^{ m'}$)~;
la lettre $\alpha$ d\'esignera toujours un multiindice de $\N^n$.
Soit $\ell \in \N$ un entier positif. Nous noterons $J_t^\ell
h(t) := \left( \partial_t^\alpha h_{ i'}(t) \right)_{ 1 \leq i' \leq
n',\, \vert \alpha \vert \leq \ell}$ le jet d'ordre $\ell$ de $h$
et aussi $J_\tau^\ell \overline{ h} (\tau) := \left(
\partial_\tau^\alpha \overline{ h}_{ i'} (\tau) \right)_{ 1\leq i'\leq
n',\, \vert \alpha \vert \leq \ell}$ le jet d'ordre $\ell$ de
$\overline{ h} (\tau)$. Le nombre de ces d\'eriv\'ees partielles est
\'egal \`a $N_{ n',\, n,\, \ell}:= n' \, \frac{ (n+ \ell) !}{ n! \
\ell !}$. Nous appellerons {\sl jet strict} de $h$ (ou de
$\overline{ h}$) les d\'eriv\'ees partielles de $h$ (ou de
$\overline{ h}$) d'ordre strictement
positif.

En observant que les coefficients des champs de vecteurs $\underline{
\mathcal{ L}}_k$ d\'efinis \`a la deuxi\`eme ligne de~\thetag{ 1.26}
sont analytiques et en raisonnant par r\'ecurrence, on d\'emontre
l'assertion suivante.

\def\thelemma{5.4}\begin{lemma}
Pour tout $i' = 1,\, \dots,\, n'$ et tout $\beta \in \N^m$, il existe
un polyn\^ome $P_{ i',\, \beta}$ en le jet $J_\tau^{ \vert \beta
\vert} \overline{ h}(\tau)$ dont les coefficients sont des s\'eries
enti\`eres convergentes par rapport aux variables $(t,\, \tau)$
qui ne d\'ependent que des s\'eries d\'efinissantes $\Theta_j
(\zeta,\, t)$ de $\mathcal{ M}$, tel que
\def\theequation{5.5}\begin{equation}
\underline{ \mathcal{ L}}^\beta
\overline{ h}_{ i'} (\tau) \equiv
P_{ i',\, \beta} \left(
t,\, \tau,\, J_\tau^{\vert \beta \vert} \, 
\overline{ h} (\tau)
\right),
\end{equation}
pour $(t,\, \tau) \in \C^n \times \C^n$.
\end{lemma}

\proof
Pour $\beta = 0$, l'existence de $P_{ i',\, \beta}$ est
\'evidente. Raisonnons par r\'ecurrence sur l'entier $\vert \beta
\vert$, en supposant plus pr\'ecis\'ement que le premier groupe
d'arguments de $P_{ i',\, \beta}$ est le jet $J_\zeta^{ \vert
\beta \vert} \Theta (\zeta,\, t)$ par rapport
\`a $\zeta$ de la s\'erie d\'efinissante
$\Theta(\zeta,\, t)$, et en supposant que
$P_{ i',\, \beta}$ est
un polyn\^ome universel
par rapport aux variables $\left(J_\zeta^{ \vert
\beta \vert} \Theta (\zeta,\, t),\,
J_\tau^{ \vert \beta \vert} \overline{ h} (\tau) 
\right)$.
Fixons $\ell \in \N$ et supposons
les relations~\thetag{ 5.5} vraies pour tout $\beta \in \N^m$ tel que
$\vert \beta \vert \leq \ell$. En appliquant la d\'erivation
$\underline{ \mathcal{ L}}_k$ \`a~\thetag{ 5.5}, et en utilisant
l'hypoth\`ese de r\'ecurrence, nous pouvons \'ecrire~:
\def\theequation{5.6}\begin{equation}
\left\{
\aligned
\underline{ \mathcal{ L}}_k 
\underline{ \mathcal{ L}}^\beta 
\overline{ h}_{ i'} (\tau) \equiv
& \
\underline{ \mathcal{ L}}_k 
\left[
P_{ i',\, \beta} \left(
J_\zeta^{ \vert
\beta \vert} \Theta (\zeta,\, t),\, 
J_\tau^{ \vert \beta \vert} \overline{ 
h} (\tau)
\right)
\right] \\
\equiv
& \
\sum_{\beta_1\in\N^m,\, \vert \beta_1 \vert \leq \vert \beta \vert} \,
\sum_{ j=1}^d\, 
\frac{ \partial P_{ i',\, \beta}}{ \partial J_j^{ \vert
\beta_1 \vert}} \cdot
\partial_{\zeta_k} \partial_\zeta^{ \beta_1} \Theta_j (\zeta,\, t) + 
\\
& \
+
\sum_{ i_1 '=1}^{ n'} \, 
\sum_{\alpha_1\in\N^n,\, 
\vert \alpha_1 \vert \leq \vert \beta \vert} \, 
\frac{ \partial P_{ i_1',\, \beta}}{ \partial J_{ i_1'}^{ 
\alpha_1}} \cdot \partial_{ \zeta_k} \partial_\tau^{ 
\alpha_1} \overline{ h}_{ i_1'} 
(\tau) + \\
& \
+
\sum_{ i_1 '= 1}^{ n'} \, 
\sum_{\alpha_1 \in \N^n,\, 
\vert \alpha_1 \vert \leq
\vert \beta \vert} \, 
\sum_{ j=1}^d \,
\frac{ \partial P_{ i',\, \beta}}{ \partial J_{ i_1'}^{ \alpha_1}}\cdot
\partial_{ \zeta_k} \Theta_j \cdot
\partial_{ \xi_j} \partial_{ \tau}^{ \alpha_1} 
\overline{ h}_{ i_1'} (\tau) \\
& \
=:
P_{ i',\, \beta + \1_k^m}
\left(
J_\zeta^{ \vert \beta \vert +1} \Theta (\zeta,\, t),\, 
J_\tau^{ \vert \beta \vert} \overline{ h} (\tau), 
\right),
\endaligned\right.
\end{equation}
o\`u le symbole $\1_k^m$ d\'esigne le multiindice $(0,\, \dots,\, 0,\,
1,\, 0,\, \dots,\, 0)\in \N^m$, avec $1$ \`a la $k$-i\`eme place et
$0$ aux autres places. Le calcul~\thetag{ 5.6} d\'emontre bien que
$P_{ i',\, \beta+ \1_k^m}$ est un polyn\^ome et qu'il se d\'eduit de
$P_{ i',\, \beta}$ par une formule alg\'ebrico-diff\'erentielle
universelle de diff\'erentiation compos\'ee. Puisque tout multiindice
de longueur $\ell +1$ peut s'\'ecrire $\beta + \1_k^m$, o\`u $\vert
\beta \vert = \ell$, ceci d\'emontre le Lemme~5.4.
\endproof

Gr\^ace \`a la formule de Leibniz pour la diff\'erentiation d'un
produit de fonctions, on observe plus g\'en\'eralement qu'un
d\'eveloppement polynomial analogue existe pour $\underline{ \mathcal{
L}}^\beta \overline{ h }^\alpha (\tau)$, lorsque $\alpha \in \N^{ n'}$
est un multiindice arbitraire. Nous n'aurons pas besoin de formules
combinatoires explicites, que du reste, nous n'avons
pas tent\'e d'\'ecrire.

Gr\^ace \`a ce pr\'eliminaire, en tenant compte du fait que $\tau =
(\zeta,\, \Theta (\zeta,\, t))$ sur $\mathcal{ M}$, nous pouvons
d\'evelopper les identit\'es~\thetag{ 5.3} comme suit~:
\def\theequation{5.7}\begin{equation}
R_{ j',\, \beta}' \left(
\zeta,\, t,\, 
J_\tau^{\vert \beta \vert} 
\overline{ h} (\zeta,\, 
\Theta (\zeta,\, t)) \, : 
h(t)
\right) \equiv 0,
\end{equation}
o\`u les $R_{ j',\, \beta}' \left( \zeta,\, t,\, J_\tau^{ \vert \beta
\vert} \overline{ h} : t' \right)$ sont des s\'eries enti\`eres
convergentes, relativement polynomiales par rapport au jet strict de
$\overline{ h}$. Le dernier argument $t'$ est plac\'e apr\`es un
signe <<:>> pour bien notifier que dans~\thetag{ 5.3}, ce sont les
composantes de $\overline{ h}$ que l'on d\'erive en appliquant les
op\'erateurs $\underline{ \mathcal{ L}}^\beta$, tandis que les
composantes de $h$ ne sont pas diff\'erenti\'ees.

En revenant \`a la d\'efinition~\thetag{ 1.34}
des applications $\Psi_{ j',\, \beta}'$ introduites
dans le \S1.33, nous avons la co\"{\i}ncidence de notations~:
\def\theequation{5.8}\begin{equation}
\Psi_{ j',\, \beta} ' \left( 0,\, 0,\, t' \right)
\equiv R_{ j',\, \beta} ' \left(0,\, 0,\, 
J_\tau^{ \vert \beta \vert } \overline{ h} (0) \, : t' \right).
\end{equation}
D'apr\`es l'hypoth\`ese principale du Th\'eor\`eme~1.39 {\bf (i)},
pour $k$ assez grand, l'application $t' \longmapsto \left( \Psi_{
j',\, \beta} ' (0,\, 0,\, t')\right)_{ 1\leq j'\leq d',\, \vert \beta
\vert \leq k}$ est de rang $n'$ en $t'= 0$. Nous noterons $\ell_0$
un tel entier $k$. De mani\`ere \'equivalente, l'application
holomorphe locale
\def\theequation{5.9}\begin{equation}
t' \longmapsto 
\left(
R_{ j',\, \beta}' \left(
0,\, 0,\, J_\tau^{ \vert \beta \vert }
\overline{ h} (0) \, : t'
\right)
\right)_{ 1\leq j' \leq d',\, 
\vert \beta \vert \leq \ell_0}
\end{equation}
est de rang $n'$ en $t'=0$.

Gr\^ace \`a cette hypoth\`ese forte, nous pouvons appliquer le
th\'eor\`eme des fonctions implicites aux identit\'es~\thetag{ 5.7},
\'ecrites seulement pour $\vert \beta \vert \leq \ell_0$, afin de
r\'esoudre $h(t)$ en fonction des autres variables, ce qui donne une
expression de la forme
\def\theequation{5.10}\begin{equation}
h(t) \equiv 
\Phi \left(
\zeta,\, t,\, J_\tau^{\ell_0} 
\overline{ h} (\zeta,\, \Theta(\zeta,\, t))
\right), 
\end{equation}
o\`u l'application $\Phi$ est \`a valeurs dans $\C^{ n'}$, et
holomorphe au voisinage de $\left( 0,\, 0,\, J_\tau^{ \ell_0}
\overline{ h } (0) \right)$ dans $\C^{ m+ n+ N_{ n',\, n,\, 
\ell_0,\, }}$.

En admettant un l\'eger \'ecart de notation, nous r\'e\'ecrirons 
l'identit\'e~\thetag{ 5.10} sous la forme
\def\theequation{5.11}\begin{equation}
h(t) = \Phi \left( t,\, \tau,\, 
J_\tau^{\ell_0} \overline{ h} (\tau) \right),
\end{equation}
\'etant entendu que $(t,\, \tau) \in \mathcal{ M}$, c'est-\`a-dire que
$\xi$ doit \^etre remplac\'e par $\Theta (\zeta,\, t)$ si l'on veut
interpr\'eter~\thetag{ 5.11} comme une v\'eritable identit\'e
formelle dans $\C \dl \zeta,\, t \dr^{ n'}$.

En partant des identit\'es de r\'eflexion \'ecrites \`a la deuxi\`eme
ligne de~\thetag{ 1.32} et en raisonnant comme ci-dessus, 
on d\'emontre que $\overline{ h}(\tau)$ satisfait
l'identit\'e suivante sur $\mathcal{ M}$~:
\def\theequation{5.12}\begin{equation}
\overline{ h} (\tau) = \overline{ \Phi} \left( \tau,\, t,\, 
J_t^{\ell_0} h ( t )\right),
\end{equation}
Pour se convaincre que c'est bien l'application conjugu\'ee
$\overline{ \Phi}$ qui appara\^{\i}t dans~\thetag{ 5.12}, on peut aussi
poser $\tau := \bar t$ dans~\thetag{ 5.11}, conjuguer de part et
d'autre du signe <<=>>, puis complexifier
la variable $(\bar t)^c =: \tau$ dans
l'identit\'e obtenue, ce qui donne~\thetag{ 5.12}.

En r\'esum\'e, nous avons obtenu un 
couple d'identit\'es formelles, 
valables pour 
$(t,\, \tau) \in \mathcal{ M}$~:
\def\theequation{5.13}\begin{equation}
\left\{
\aligned
h(t) 
& \
= \Phi \left( t,\, \tau,\, 
J_\tau^{\ell_0} \overline{ h} (\tau)\right), \\
\overline{ h} (\tau) 
& \
= \overline{ \Phi} \left( \tau,\, t,\, 
J_t^{\ell_0} h ( t )\right).
\endaligned\right.
\end{equation}
Soit $\Pi(t,\, \tau) \in \C\dl t,\, \tau \dr$ une s\'erie formelle
sans terme constant, par exemple $h_{ i'} (t) - \Phi_{ i'} \left( t,\,
\tau,\, J_\tau^{ \ell_0} \overline{ h} (\tau) \right)$, o\`u $i'$ est
fix\'e. Quand nous disons que $\Pi(t,\, \tau) = 0$ pour $(t,\, \tau)
\in \mathcal{ M}$, nous entendons que l'on a les deux identit\'es
formelles $\Pi(t,\, \zeta,\, \Theta (\zeta,\, t)) \equiv 0$ dans
$\C\dl \zeta,\, t \dr$ et $\Pi( z,\, \overline{ \Theta} (z,\, \tau),\,
\tau) \equiv 0$ dans $\C \dl z,\, \tau\dr$. De mani\`ere
\'el\'ementaire, on v\'erifie qu'en posant $r (t,\, \tau):= \xi
-\Theta (\zeta,\, t)$ et $\overline{ r} (\tau,\, t) := w - \overline{
\Theta} (z,\, \tau)$, ces identit\'es sont \'equivalentes entre elles
et qu'elles sont \'equivalentes \`a la propri\'et\'e suivante~: il
existe deux matrices de s\'eries formelles $c(t,\, \tau)$ et $d(t,\,
\tau)$ de taille $1 \times d$ telles que $\Pi(t,\, \tau) \equiv c(t,\,
\tau) \, r(t,\, \tau)$ et $\Pi(t,\, \tau) \equiv d(t,\, \tau) \,
\overline{ r}(\tau,\, t)$ dans $\C \dl t,\, \tau\dr$. Retenons donc
que l'expression $\Pi(t,\, \tau) = 0$ pour $(t,\, \tau) \in \mathcal{
M}$ signifie rigoureusement que $\Pi(t,\, \zeta,\, \Theta (\zeta,\,
t)) \equiv 0$ ou que $\Pi( z,\, \overline{ \Theta} (z,\, \tau),\,
\tau) \equiv 0$.

Le couple d'identit\'es~\thetag{ 5.13} 
est absolument fondamental~: il signifie que
$h(t)$ (resp. $\overline{ h} (\tau)$) peut \^etre r\'esolu par
rapport au jet d'ordre $\ell_0$ de $\overline{ h} (\tau)$ (resp. de
$h(t)$) et que cette r\'esolution est {\it analytique}\, (puisque
$\Phi$ l'est), bien que $h(t)$, $\overline{ h} (\tau)$ et leurs jets
soient, par hypoth\`ese, des applications {\it formelles}\, qui sont
{\it a priori}\, non convergentes. Ces identit\'es ne permettent pas
encore de d\'emontrer que $h(t)$ et $\overline{ h} (\tau)$
convergent~: en effet, les variables $t$ et $\tau$ sont encore
m\^el\'ees dans~\thetag{ 5.13}, et les deux termes $J_\tau^{ \ell_0 }
\overline{ h} (\tau)$, $J_t^{ \ell_0} h(t)$ sont {\it a priori}\, non
convergents. Mais dans le reste de cette section, 
en <<d\'em\'elangeant>> les variables $t$ et $\tau$,
nous allons
\'etablir l'\'enonc\'e suivant, qui compl\`etera la d\'emonstration
du Th\'eor\`eme~1.39 {\bf (i)}.

\def\theassertion{ 5.14}\begin{assertion}
En supposant $\mathcal{ M}$ minimale \`a l'origine, le couple
d'identit\'es~\thetag{ 5.13} implique que $h(t) \in \C \{ t\}^{ n'}$
et que $\overline{ h} (\tau) \in \C \{ \tau \}^{ n'}$ convergent.
\end{assertion}

\subsection*{ 5.15.~Passage aux jets d'ordre quelconque} 
Notre premier objectif est de g\'en\'eraliser~\thetag{ 5.13} en faisant
appara\^{\i}tre, \`a la place des membres de gauche $h (t)$ et
$\overline{ h } (\tau)$, des jets d'ordre arbitraires $J_t^\ell h(t)$
et $J_\tau^\ell \overline{ h} (\tau)$, o\`u $\ell \in \N$.
Nous verrons ci-apr\`es pourquoi cela est n\'ecessaire. 

Pour cela, introduisons les champs de vecteurs $\Upsilon_j$, $j=1,\,
\dots,\, d$, et $\underline{ \Upsilon }_j$, $j=1,\, \dots,\, d$,
d\'efinis comme suit~:
\def\theequation{5.16}\begin{equation}
\Upsilon_j := 
\frac{ \partial }{ \partial w_j} +
\sum_{ l=1}^d\, 
\frac{ \partial \Theta_l}{ \partial w_j}
(\zeta,\, t)\, 
\frac{ \partial }{ \partial \xi_l}~;
\ \ \ \ \ \ \
\ \ \ \ \ \ \
\underline{ \Upsilon}_j := 
\frac{ \partial }{ \partial \xi_j} +
\sum_{ l=1}^d \, 
\frac{ \partial \overline{ \Theta}_l}{ \partial \xi_j}
(z,\, \tau) \, 
\frac{ \partial }{ \partial w_l}.
\end{equation}
Ces champs de vecteurs sont tangents \`a $\mathcal{ M}$~: en effet, on
v\'erifie imm\'ediatement que $\Upsilon_{ j_1} \left[ \xi_{ j_2}
-\Theta_{ j_2} ( \zeta,\, t) \right] \equiv 0$ et que $\underline{
\Upsilon}_{ j_1} \left[ w_{ j_2} - \overline{ \Theta}_{ j_2} (z,\,
\tau) \right] \equiv 0$, pour $j_1,\, j_2 = 1,\, \dots,\, d$.
Observons que la collection des $2m+d$ champs de vecteurs $\mathcal{
L}_k$, $\underline{ \mathcal{ L}}_k$ et $\Upsilon_j$ engendre le
fibr\'e tangent $T\mathcal{ M}$. La m\^eme propri\'et\'e est
satisfaite par la collection $\mathcal{ L}_k$, $\underline{ \mathcal{
L}}_k$ et $\underline{ \Upsilon}_j$. Notons au passage des relations
de commutation que nous n'utiliserons pas~: $\left[ \Upsilon_j,\,
\underline{\mathcal{ L}}_k \right]= 0$ et $ \left[ \underline{
\Upsilon}_j,\, \mathcal{ L}_k \right]= 0$. Consid\'erons l'action sur
$h(t)$ des champ de vecteurs $\Upsilon_j$, interpr\'et\'es comme
d\'erivations. Soit $\delta = (\delta_1,\, \dots,\, \delta_d) \in
\N^d$ un multiindice arbitraire. Notons $\Upsilon^\delta$ la
d\'erivation compos\'ee $(\Upsilon_1 )^{ \delta_1} (\Upsilon_2 )^{
\delta_2 } \cdots (\Upsilon_d)^{ \delta_d}$, qui est d'ordre $\vert
\delta \vert$. On observe imm\'ediatement que
\def\theequation{5.17}\begin{equation}
\Upsilon^\delta h(t) \equiv
\partial_w^\delta h(t).
\end{equation}
De mani\`ere similaire, $\underline{ \Upsilon}^\delta \overline{ h}
(\tau) \equiv \partial_{ \xi}^\delta \overline{ h} (\tau)$.

En proc\'edant comme dans la d\'emonstration du Lemme~5.4, on
v\'erifie que pour tout $i'= 1,\, \dots,\, n'$ et tout multiindice
$\alpha = (\beta,\, \delta) \in \N^m \times \N^d$, il existe des
s\'eries enti\`eres $Q_{ i',\, \beta,\, \delta}$ telles que
\def\theequation{5.18}\begin{equation}
Q_{ i',\, \beta,\, \delta} 
\left( 
t,\, \tau,\,
J_t^{ \vert \beta \vert + \vert \delta \vert} 
h_{ i_1'} (t)
\right).
\end{equation}
Ces s\'eries sont polynomiales par rapport aux variables de jets et
poss\`edent des coefficients analytiques par rapport \`a $(t,\,
\tau)$, lesquels coefficients ne d\'ependent que des s\'eries
d\'efinissantes $\Theta_j (\zeta,\, t)$ de $\mathcal{ M}$. 

R\'eciproquement, puisque les $n$ champs de vecteurs $\mathcal{ L}_k$
et $\Upsilon_j$ engendrent l'espace des $(t_1,\, \dots,\, t_n)$, nous
pouvons inverser les formules $\mathcal{ L}^\beta \Upsilon^\delta h_{
i'} (t) \equiv Q_{ i',\, \beta,\, \delta}$, en consid\'erant les
d\'eriv\'ees partielles $\partial_t^{ \alpha_1} h_{ i_1'}(t)$ comme
inconnues. En examinant les termes de plus haute d\'erivation dans les
polyn\^omes $Q_{ i',\, \beta,\, \delta}$, on constate que le syst\`eme
lin\'eaire que l'on doit inverser est de type trigonal, et de
d\'eterminant $1$~: l'inversion est \'el\'ementaire et conserve la
polynomialit\'e. On en d\'eduit~:

\def\thelemma{5.19}\begin{lemma}
Pour tout $i'= 1,\, \dots,\, n'$ et tout multiindice $\alpha =
(\beta,\, \delta) \in \N^m \times \N^d$, il existe une s\'erie
enti\`ere $P_{ i',\, \alpha}$, polynomiale par rapport \`a son dernier
groupe de variables, dont les coefficients sont des
s\'eries enti\`eres convergentes par
rapport \`a $(t,\, \tau)$ qui ne d\'ependent que des s\'eries
d\'efinissantes $\Theta_j(\zeta,\, t)$ de $\mathcal{ M}$, tel que
l'identit\'e suivante est satisfaite dans $\C \dl t \dr$~{\rm :}
\def\theequation{5.20}\begin{equation}
\partial_t^\alpha h_{ i'} (t) 
\equiv P_{ i',\, \alpha} 
\left(
t,\, \tau,\, 
\left(
\mathcal{ L}^{ \beta_1} 
\Upsilon^{ \delta_1} h(t)
\right)_{ \vert \beta_1 \vert \leq
\vert \beta \vert,\, \vert \delta_1 \vert \leq \vert
\delta \vert}
\right),
\end{equation}
lorsque $(t,\, \tau) \in \C^n \times \C^n$.
\end{lemma}

Appliquons maintenant la d\'erivation $\mathcal{ L}^\beta
\Upsilon^\delta$ \`a la premi\`ere ligne de~\thetag{ 5.12} et la
d\'erivation $\underline{ \mathcal{ L}}^\beta \underline{
\Upsilon}^\delta$ \`a la seconde ligne de~\thetag{ 5.12}~: cela est
autoris\'e, puisque les champs de vecteurs $\mathcal{ L}_k$,
$\underline{ \mathcal{ L}}_k$, $\Upsilon_j$ et $\underline{
\Upsilon}_j$ sont tangents \`a $\mathcal{ M}$. Comme les coefficients
de ces d\'erivations sont analytiques, en utilisant la <<r\`egle de la
cha\^{ \i}ne>> pour diff\'erentier les membres de droite $\Phi$ et
$\overline{ \Phi}$ de~\thetag{ 5.12}, nous voyons appara\^{\i}tre des
expressions polynomiales par rapport aux variables de jets stricts,
dont les coefficients sont analytiques par rapport aux variables
$(t,\, \tau) \in \mathcal{ M}$. L'ordre de ces variables de jets monte
jusqu'\`a $\ell_0 + \vert \beta \vert + \vert \delta \vert$. En
appliquant le Lemme~5.19, nous pouvons effectuer des substitutions
lin\'eaires sur les identit\'es obtenues afin de r\'esoudre les
d\'eriv\'ees partielles $\partial_{ t}^\alpha h_{ i'} (t)$, ce qui
donne des identit\'es de la forme $\partial_{ t}^\alpha h_{ i'} (t)
\equiv \Phi_{ i',\, \alpha} \left( t,\, \zeta,\, \Theta (\zeta,\,
t),\, J_\tau^{ \ell_0 + \vert \beta \vert +
\vert \delta \vert} \overline{ h}( \zeta,\, \Theta (\zeta,\,
t)) \right)$.

Soit $\ell \in \N$ un entier arbitraire. Collectons
les identit\'es obtenues pour tous les multiindices $\alpha$
satisfaisant $\vert \alpha \vert \leq \ell$ et r\'esumons les calculs
effectu\'es.

\def\thelemma{5.21}\begin{lemma}
Pour tout entier $\ell \in \N$, il existe une s\'erie enti\`ere
$\Phi_\ell = \Phi_\ell \left( t,\, \tau,\, J^{ \ell_0 + \ell} \right)$
\`a valeurs dans $\C^{ N_{ n',\, n,\, \ell}}$, polynomiale par rapport
aux variables de jets stricts, dont les coefficients sont des s\'eries
enti\`eres convergentes par rapport \`a $(t,\, \tau)$ qui ne
d\'ependent que des s\'eries d\'efinissantes $\Theta_j(\zeta,\, t)$ de
$\mathcal{ M}$, telle que les deux identit\'es suivantes sont
satisfaites~{\rm :}
\def\theequation{5.22}\begin{equation}
\left\{
\aligned
J_t^\ell h(t) 
& \
\equiv
\Phi_\ell \left(
t,\, \zeta,\, \Theta (\zeta,\, t),\, 
J_\tau^{ \ell_0 + \ell} \overline{ h}(
\zeta,\, \Theta (\zeta,\, t))
\right) \ \ \ \ \ {\rm et} \\
J_\tau^\ell \overline{ h} (\tau) 
& \
\equiv
\overline{ \Phi}_\ell \left(
\tau,\, z,\, \overline{ \Theta} (z,\, \tau),\, 
J_t^{ \ell_0 + \ell} h( z,\, \overline{ \Theta} (z,\, \tau))
\right),
\endaligned\right.
\end{equation}
dans $\C \dl \zeta,\, t \dr^{ N_{ n',\, n,\, \ell}}$ et
dans $\C \dl z,\, \tau \dr^{ N_{ n',\, n,\, \ell}}$.
\end{lemma}

Ces formules constituent la g\'en\'eralisation d\'esir\'ee de~\thetag{
5.12}. Elles seront notre nouveau point de d\'epart.

La deuxi\`eme \'etape cruciale dans la d\'emonstration de
l'assertion~5.13 va consister \`a effectuer un grand nombre
d'<<auto-substitutions>> entre ces formules, en se <<d\'epla\c cant>>
le long des cha\^{\i}nes de Segre (conjugu\'ees). Par souci de
clart\'e et compr\'ehensibilit\'e, commen\c cons par d\'etailler les
deux premi\`eres \'etapes.

\subsection*{ 5.23.~Restrictions \`a la premi\`ere et \`a la
deuxi\`eme cha\^{ \i}ne de Segre conjugu\'ee} Rappelons que $h^c (t,\,
\tau) = \left(h(t),\, \overline{ h} (\tau) \right)$. Nous admettrons
l'abus de notation suivant~: au lieu d'\'ecrire rigoureusement $h
\left( \pi_t (t,\, \tau) \right)$ et $\overline{ h} \left( \pi_\tau
(t,\, \tau)\right)$, nous \'ecrirons $h(t,\, \tau)= h(t)$ et
$\overline{ h} (t,\, \tau) = \overline{ h} (\tau)$.

Soit $\nu \in \N$, $\nu \geq 1$, soit ${\sf x} \in \C^\nu$ et soit
${\sf Q} ( {\sf x} ) = \left( {\sf Q}_1 ( {\sf x} ),\, \dots,\, {\sf
Q}_{ 2n} ( {\sf x} ) \right) \in \C \dl {\sf x} \dr^{ 2n}$ une
application formelle quelconque satisfaisant ${\sf Q} (0)= 0$.
Puisque le flot multiple de $\mathcal{ L}$ d\'efini par~\thetag{ 3.13}
n'agit pas sur les variables $\tau$ et que le flot multiple de
$\underline{ \mathcal{ L}}$ d\'efini par~\thetag{ 3.14} n'agit pas sur
les variables $t$, nous avons deux relations simples~:
\def\theequation{5.24}\begin{equation}
\left\{ 
\aligned
J_t^\ell h
\left(
\underline{ \mathcal{ L}}_{ z_1} \left(
{\sf Q} ({\sf x} )
\right)
\right) 
& \
\equiv
J_t^\ell h 
\left(
{\sf Q} ({\sf x} )
\right) 
\ \ \ \ \ \ {\rm et} \\
J_\tau^\ell \overline{ h}
\left(
\mathcal{ L}_{ z_1} 
\left(
{\sf Q} ({\sf x} )
\right)
\right) 
& \
\equiv
J_\tau^\ell \overline{ h}
\left(
{\sf Q} ({\sf x} )
\right).
\endaligned\right.
\end{equation}
Pour g\'en\'eraliser ces observations, rappelons que lorsque $k$ est
impair, la $k$-i\`eme cha\^{\i}ne de Segre conjugu\'ee $\underline{
\Gamma }_k$ commence \`a gauche par le multiflot $\underline{
\mathcal{ L} }_{ z_k} (\cdots)$~; de m\^eme, lorsque $k$ est pair,
elle commence \`a gauche par le multiflot $\mathcal{ L}_{ z_k} (
\cdots)$. Nous d\'eduisons alors de~\thetag{ 5.24} 
deux relations triviales, mais absolument cruciales~:
\def\theequation{5.25}\begin{equation}
\left\{ 
\aligned
J_t^\ell h \left(
\underline{ \Gamma}_k \left( z_{ (k)} \right)
\right) 
& \
\equiv
J_t^\ell h \left(
\underline{ \Gamma}_{ k-1} \left( z_{ (k-1)} \right)
\right),
\ \ \ \ \ \ \
\text{\sf si} \ k \
\text{\sf est impair}~; \\
J_\tau^\ell \overline{ h} \left(
\underline{ \Gamma}_k \left( z_{ (k)} \right)
\right) 
& \
\equiv
J_\tau^\ell \overline{ h} \left(
\underline{ \Gamma}_{ k-1} \left( z_{ (k-1)} \right)
\right),
\ \ \ \ \ \ \
\text{\sf si} \ k \
\text{\sf est pair}.
\endaligned\right.
\end{equation}
Des relations similaires sont satisfaites avec la cha\^{\i}ne de Segre
$\Gamma_k$.

Restreignons maintenant les relations~\thetag{ 5.22} \`a la premi\`ere
cha\^{ \i}ne de Segre conjugu\'ee, c'est-\`a-dire rempla\c cons $(t,\,
\tau)$ par $\underline{ \Gamma}_1 \left( z_{ (1)} \right)$.
En tenant compte de~\thetag{ 5.24}, nous obtenons~:
\def\theequation{5.26}\begin{equation}
\left\{ 
\aligned
J_t^\ell h \left(
0
\right) 
& \
\equiv 
\Phi_\ell \left(
\underline{ \Gamma}_1 (
z_1
), \
J_\tau^{ \ell_0+ \ell} \overline{ h}
\left(
\underline{ \Gamma}_1 
(
z_1
)
\right)
\right) 
\ \ \ \ \ \ \ 
{\rm et} \\
J_\tau^\ell \overline{ h} \left(
\underline{ \Gamma}_1 (
z_1 )
\right) 
& \
\equiv
\overline{ \Phi}_\ell
\left(
\underline{ \Gamma}_1
(
z_1
), \
J_t^{\ell_0 + \ell} h (0)
\right).
\endaligned\right.
\end{equation}
La premi\`ere relation ne sera pas utile~; interpr\'etons la seconde~:
{\it puique le jet $J_t^\ell h(0)$ est constant et puisque $\overline{
\Phi}_\ell$ est analytique, le membre de droite converge par rapport
\`a la variable $z_1 \in \C^m$}~; par cons\'equent, le membre de
gauche est convergent par rapport \`a $z_1$, ce qui signifie que {\it
la restriction \`a la premi\`ere cha\^{\i}ne de Segre conjugu\'ee du
jet d'ordre quelconque $J_\tau^\ell \overline{ h}$ converge}.

Passons maintenant \`a la deuxi\`eme cha\^{ \i}ne de Segre
conjugu\'ee. Rempla\c cons $(t,\, \tau)$ par $\underline{ \Gamma}_2 (
z_{ (2)}) = \underline{ \Gamma}_2 (z_1,\, z_2)$ dans~\thetag{ 5.22}.
En tenant compte de~\thetag{ 5.24}, nous obtenons~:
\def\theequation{5.27}\begin{equation}
\left\{ 
\aligned
J_t^\ell h \left(
\underline{ \Gamma}_2 ( z_{ (2)}) 
\right) 
& \
\equiv 
\Phi_\ell \left(
\underline{ \Gamma}_2 (
z_{ (2)}
), \
J_\tau^{ \ell_0+ \ell} \overline{ h}
\left(
\underline{ \Gamma}_1 
(
z_1
)
\right)
\right) 
\ \ \ \ \ \ \ 
{\rm et} \\
J_\tau^\ell \overline{ h} \left(
\underline{ \Gamma}_1 (
z_1 )
\right) 
& \
\equiv
\overline{ \Phi}_\ell
\left(
\underline{ \Gamma}_2
(
z_{ (2)}
), \
J_t^\ell h \left(
\underline{ \Gamma}_2 ( z_{ (2)}) 
\right)
\right).
\endaligned\right.
\end{equation}
La seconde relation ne sera pas utile. Quant \`a la premi\`ere, on
voit tout de suite qu'il faut y remplacer le terme $J_\tau^{ \ell_0 +
\ell} \overline{ h} \left( \underline{ \Gamma}_1 ( z_1) \right)$ par
l'expression que nous venons d'obtenir dans~\thetag{ 5.26}, sans oublier
de remplacer $\ell$ par $\ell_0 + \ell$, ce qui nous donne~:
\def\theequation{5.28}\begin{equation}
\left\{
\aligned
J_t^\ell h \left(
\underline{ \Gamma}_2 ( z_{ (2)}) 
\right) 
& \
\equiv 
\Phi_\ell \left(
\underline{ \Gamma}_2 (
z_{ (2)}
), \
J_\tau^{ \ell_0+ \ell} \overline{ h}
\left(
\underline{ \Gamma}_1 
(
z_1
)
\right)
\right) \\
& \
\equiv
\Phi_\ell 
\left(
\underline{ \Gamma}_2 (
z_{ (2)}
), \
\overline{ \Phi}_{ \ell_0 +\ell} 
\left(
\underline{ \Gamma}_1 
(z_1),\, J_t^{2\ell_0 + \ell} h( 0)
\right)
\right) \\
& \
=:
\Pi_{\ell,\, 2} 
\left(
z_{ (2)},\, 
J_t^{2 \ell_0 + \ell} h(0)
\right).
\endaligned\right.
\end{equation}
Pour pouvoir effectuer une telle substitution, il \'etait important
de formuler \`a l'avance les identit\'es~\thetag{ 5.22} qui
g\'en\'eralisent aux jets d'ordre quelconque les identit\'es~\thetag{
5.12}. Interpr\'etons maintenant le r\'esultat~: l'expression de la
seconde ligne de~\thetag{ 5.28} est convergente, puisque $\Phi_\ell$
et $\overline{ \Phi}_{ \ell_0+ \ell}$ sont analytiques et puisque le
jet $J_t^{ 2\ell_0 +\ell} h(0)$ de l'application formelle $h$ est
constant. On en d\'eduit que {\it la restriction \`a la seconde cha\^{
\i}ne de Segre conjugu\'ee du jet d'ordre quelconque $J_t^\ell h$
converge}.

Le lecteur aura devin\'e qu'en poursuivant ces substitutions le long
des cha\^{\i}nes de Segre conjugu\'ees de longueur arbitraire, nous
pourrons conclure la d\'emonstration de l'assertion~5.13, gr\^ace au
Lemme~3.32. \'Ecrivons donc la d\'emonstration finale.

\subsection*{ 5.29.~Substitutions sur les cha\^{\i}nes de
Segre conjugu\'ees de longueur arbitraire} 
L'assertion suivante g\'en\'eralise les substitutions
pr\'ec\'edentes.

\def\thelemma{5.30}\begin{lemma}
Pour tout entier $k\in \N$ et tout entier $\ell \in \N$, il existe une
s\'erie formelle $\Pi_{ \ell,\, k} = \Pi_{ \ell,\, k} \left( z_{
(k)},\, J^{k \ell_0 + \ell} \right)$ \`a valeurs dans $\C^{ N_{ n',\,
n,\, \ell}}$, polynomiale par rapport aux variables de jets stricts,
dont les coefficients sont des s\'eries enti\`eres convergentes par
rapport \`a $z_{ (k)}$ qui ne d\'ependent que des s\'eries
d\'efinissantes $\Theta_j (\zeta,\, t)$ de $\mathcal{ M}$, telle que
les deux identit\'es formelles suivantes~{\rm :}
\def\theequation{5.31}\begin{equation}
\left\{
\aligned
J_\tau^\ell
\overline{ h} \left(
\underline{ \Gamma}_k (z_{ (k)})
\right) 
& \
\equiv
\overline{ \Pi}_{ \ell,\, k} 
\left(
z_{ (k)},\, 
J_t^{ k\ell_0+ \ell} 
h(0)
\right), 
\ \ \ \ \ \ \
\text{ \sf si} \ k
\text{ \sf est impair}~; \\
J_t h \left(
\underline{ \Gamma}_k (z_{ (k)}) 
\right)
& \
\equiv
\Pi_{ \ell,\, k}
\left(
z_{ (k)},\, 
J_t^{ k\ell_0 +\ell} h (0)
\right) 
\ \ \ \ \ \ \
\text{ \sf si} \ k
\text{ \sf est pair},
\endaligned\right.
\end{equation}
sont satisfaites dans $\C \dl z_{ (k)} \dr^{ N_{ n',\, n,\, \ell}}$.
\end{lemma}

Ce lemme nous permet de conclure la d\'emonstration de
l'assertion~5.13, et par l\`a-m\^eme celle du Th\'eor\`eme~1.39 {\bf
(iii)}. En effet, en posant $\ell = 0$ dans la seconde
identit\'e~\thetag{ 5.31}, gr\^ace \`a la convergence de $\Pi_{ 0,\,
k}$, l'hypoth\`ese du Lemme~3.32 est satisfaite.

\proof
Pour $k= 1,\, 2$, le travail est achev\'e. Raisonnons par r\'ecurrence
en supposant le lemme vrai pour un entier $k$. On supposera $k$
impair\,--\,le cas o\`u il est pair se traitera
de mani\`ere similaire.

Rempla\c cons $(t,\, \tau)$ par $\underline{ \Gamma}_{ k+1} (z_{
(k+1)} )$ dans la premi\`ere ligne de~\thetag{ 5.22}~; utilisons la
seconde relation~\thetag{ 5.25} pour simplifier~; appliquons
l'hypoth\`ese de r\'ecurrence, c'est-\`a-dire la premi\`ere ligne
de~\thetag{ 5.31}~: ces trois actions reviennent \`a effectuer le
calcul suivant~:
\def\theequation{5.32}\begin{equation}
\left\{
\aligned
J_t^\ell h \left(
\underline{ \Gamma}_{ k+1} (z_{ (k+1)}) 
\right)
& \
\equiv
\Phi_\ell 
\left(
\underline{ \Gamma}_{ k+1} (z_{ (k+1)}),\, 
J_\tau^{ \ell_0 + \ell} \overline{ h} 
\left(
\underline{ \Gamma}_{ k+1} (z_{ (k+1)})
\right)
\right) \\
& \
\equiv
\Phi_\ell 
\left(
\underline{ \Gamma}_{ k+1} (z_{ (k+1)}),\, 
J_\tau^{ \ell_0 + \ell} \overline{ h} 
\left(
\underline{ \Gamma}_k (z_{ (k)})
\right)
\right) \\
& \
\equiv
\Phi_\ell 
\left(
\underline{ \Gamma}_{ k+1} (z_{ (k+1)}),\
\overline{ \Pi}_{ \ell_0+\ell,\, k} 
\left(
z_{ (k)},\, 
J_t^{ k\ell_0 + \ell_0 + \ell} h(0)
\right)
\right) \\
& \
=:
\Pi_{ \ell,\, k+1} \left(
z_{ (k+1)},\, 
J_t^{ (k+1) \ell_0+ \ell} h (0)
\right).
\endaligned\right.
\end{equation}
Le r\'esultat obtenu n'est autre que la formule~\thetag{ 5.31} au
niveau $k+1$.
\endproof

\'Etudions \`a pr\'esent la convergence de $h$ avec des hypoth\`eses
moins simples que la non-d\'eg\'en\'erescence finie {\bf (h2)}.

\section*{\S6.~Convergence d'applications CR formelles 
Segre non-d\'eg\'en\'er\'ees}

\subsection*{ 6.1.~Pr\'eliminaire} 
Cette section est enti\`erement consacr\'ee \`a d\'emontrer le
Th\'eor\`eme~1.39 {\bf (iii)}, valable pour des applications CR
formelles satisfaisant la condition de non-d\'eg\'en\'e\-res\-cence {\bf
(h4)}. D'apr\`es le \S5.2, nous pouvons d\'evelopper les identit\'es
de r\'eflexion \thetag{ 5.3} sous la forme~\thetag{ 5.6}. En revenant
\`a la d\'efinition~\thetag{ 1.34} des s\'eries $\Psi_{ j',\,
\beta}'(z,\, w,\, \zeta,\, \xi,\, t')$ introduites dans le \S1.33,
nous avons la co\"{\i}ncidence de notations~:
\def\theequation{6.2}\begin{equation}
\left\{ 
\aligned
0
& \
\equiv
\Psi_{ j',\, \beta}'
\left(
z,\, \overline{ 
\Theta} (z,\, 0),\, 0,\, 0,\, 
h\left(
z,\, \overline{ \Theta} (z,\, 0)
\right)
\right) \\
& \
\equiv
R_{ j',\, \beta} ' \left(
z,\, \overline{ \Theta} (z,\, 0),\, 
0,\, 0,\, 
J_{ \tau}^{ \vert \beta \vert} 
\overline{ h} (0) \, :
h \left(
z,\, \overline{ \Theta} (z,\, 0)
\right)
\right).
\endaligned\right.
\end{equation}
Une fois pour toutes, choisissons des entiers $j'(1),\, \dots,\,
j'(n')$ et des multiindices distincts $\beta( 1),\, \dots,\, \beta
(n')$ tels que le d\'eterminant~\thetag{ 1.36} ne s'annule pas. Pour
all\'eger la pr\'esentation, nous noterons les s\'eries $R_{ j'(i_1'),
\, \beta( i_1')}'$ par $R_{ i_1'} '$, pour $i_1'=1,\, \dots,\, n'$.
D\'efinissons $\ell_0 := \max_{ 1\leq i_1 ' \leq n'} \, \vert
\beta(i_1') \vert$.

\subsection*{ 6.3.~Proposition principale}
Nous ramenons ainsi la d\'emonstration du Th\'eor\`eme~1.39 {\bf
(iii)} \`a la proposition suivante, l\'eg\`erement plus g\'en\'erale,
puisqu'on ne suppose pas que $h$ est une application CR formelle \`a
valeurs dans une sous-vari\'et\'e g\'en\'erique analytique r\'eelle
complexifi\'ee
$\mathcal{ M}'$, mais seulement que c'est une application formelle
satisfaisant certaines \'equations analytiques.

\def\theproposition{6.4}\begin{proposition}
Soient $R_{ i'}' \left( t,\, \tau,\, J^{\ell_0} \, : t'\right)$, $i'
=1,\, \dots,\, n'$, des s\'eries enti\`eres holomorphes dans un
voisinage de $\left( 0,\, 0,\, J_\tau^{ \ell_0} \overline{ h} (0) \, :
0 \right)$ dans $\C^n \times \C^n \times \C^{ N_{ n',\, n,\,
\ell_0}}$, polynomiales par rapport aux variables de jets stricts.
Soit $\mathcal{ M}$ la complexification d'une sous-vari\'et\'e
analytique r\'eelle g\'en\'erique repr\'esent\'ee par les
\'equations~\thetag{ 3.6} et soit $h(t) \in \C \dl t \dr^{ n'}$ une
application formelle, avec $h(0) = 0$. Supposons que pour $(t,\,
\tau) \in \mathcal{ M}$, les deux identit\'es formelles conjugu\'ees
qui suivent sont satisfaites~{\rm :}
\def\theequation{6.5}\begin{equation}
\left\{
\aligned
0 
& \
= R_{ i'} ' \left(
t,\, \tau,\, J_\tau^{ \ell_0} \overline{ h} (\tau) \, : 
h(t)
\right), \\ 
0 
& \
= 
\overline{ R}_{ i'} ' 
\left(
\tau,\, t,\, J_t^{ \ell_0} h (t)\, : \overline{ h} (\tau)
\right).
\endaligned\right.
\end{equation}
pour $i'= 1,\, \dots,\, n'$. Si $\mathcal{ M}$ est minimale en $0$ et
si le d\'eterminant suivant~{\rm :}
\def\theequation{6.6}\begin{equation}
{\rm det} \left(
\frac{ \partial R_{ i'
}'}{ \partial t_{ i_1'}}
\left(
z,\, \overline{ \Theta} 
(z,\, 0),\, 0,\, 0,\, 
J_\tau^{ \ell_0} 
\overline{ h} (0)\, : h\left(
z,\, \overline{ \Theta} (z,\, 0)
\right)
\right)
\right)_{ 1\leq i',\, i_1' \leq n'} \not \equiv 0,
\end{equation}
ne s'annule pas identiquement
dans $\C\{ z\}$, l'application $h(t) \in \C\{ t\}^{ n'}$
est convergente.
\end{proposition}

\proof
En travaillant avec les cha\^{\i}nes de Segre $\Gamma_k$ (alors que
nous avons travaill\'e avec les cha\^{\i}nes conjugu\'ees $\underline{
\Gamma}_k$ dans la Section~5), \'etablissons que $h \left( \Gamma_ k
(z_{ (k)}) \right) \in \C\{ z_{ (k) }\}^{ n'}$ pour tout $k \in
\N$. D'apr\`es le Lemme~3.32, cela suffira pour obtenir la conclusion
de la Proposition~6.4. Commen\c cons par \'etablir un crit\`ere
simplifi\'e pour la convergence des jets de $h$ et de $\overline{ h}$
sur les cha\^{\i}nes de Segre.

\subsection*{6.7.~Jets transversaux de $h$}
Rappelons que les champs de vecteurs $\Upsilon_j$ et $\underline{
\Upsilon}_j$ tangents \`a $\mathcal{ M}$ sont d\'efinis par~\thetag{
5.16}. Notons les relations
\def\theequation{6.8}\begin{equation}
\left\{ 
\aligned
J_\tau^\ell \overline{ h} \left(
\Gamma_k \left( z_{ (k)} \right)
\right) 
& \
\equiv
J_\tau^\ell \overline{ h} \left(
\Gamma_{ k-1} \left( z_{ (k-1)} \right)
\right),
\ \ \ \ \ \ \
\text{\sf si} \ k \
\text{\sf est impair}, \\
J_t^\ell h \left(
\Gamma_k \left( z_{ (k)} \right)
\right) 
& \
\equiv
J_t^\ell h \left(
\Gamma_{ k-1} \left( z_{ (k-1)} \right)
\right),
\ \ \ \ \ \ \
\text{\sf si} \ k \
\text{\sf est pair},
\endaligned\right.
\end{equation}
qui sont analogues et \'equivalentes par conjugaison 
aux relations~\thetag{ 5.25}.

\def\thelemma{6.9}\begin{lemma}
Soit $\ell \in \N$,
soit $k \in \N$ et soit $\Gamma_k ( z_{ (k)})$ la $k$-i\`eme cha\^{
\i}ne de Segre. Les propri\'et\'es suivantes sont \'equivalentes~{\rm
:}
\begin{itemize}
\item[{\bf (i)}]
Le jet d'ordre
$\ell$ de $h$ ou de $\overline{ h}$ converge sur
la $k$-i\`eme cha\^{\i}ne de Segre, {\it i.e.}
en tenant compte des simplifications~\thetag{ 6.8}~{\rm :}
\def\theequation{6.10}\begin{equation}
\left\{
\aligned
J_t^\ell h \left(
\Gamma_k ( z_{ (k)}) 
\right) 
& \
\in 
\C\{ z_{ (k)} \}^{ N_{ n',\, n,\, \ell}} \ \ \ \ \
\text{\sf si} \ k \
\text{\sf est impair},
\\
J_\tau^\ell \overline{ h}\left(
\Gamma_k ( z_{ (k)}) 
\right) 
& \
\in 
\C\{ z_{ (k)} \}^{ N_{ n',\, n,\, \ell}} \ \ \ \ \
\text{\sf si} \ k \
\text{\sf est pair}~;
\endaligned\right.
\end{equation}
\item[{\bf (ii)}]
pour tout $\beta \in \N^m$ et tout $\delta \in \N^d$ satisfaisant
$\vert \beta \vert + \vert \delta \vert \leq \ell$, les d\'eriv\'ees
suivantes de $h$ ou de $\overline{ h}$ convergent~{\rm :}
\def\theequation{6.11}\begin{equation}
\left\{
\aligned
\left[
\mathcal{ L}^\beta
\Upsilon^\delta h
\right] \left(
\Gamma_k ( z_{ (k)})
\right) 
& \
\in \C \{ z_{ (k)} \}^{ n'} \ \ \ \ \
\text{\sf si} \ k \
\text{\sf est impair}, \\
\left[
\underline{\mathcal{ L}}^\beta
\underline{ \Upsilon}^\delta \overline{ h}
\right] \left(
\Gamma_k ( z_{ (k)})
\right) 
& \
\in \C \{ z_{ (k)} \}^{ n'} \ \ \ \ \
\text{\sf si} \ k \
\text{\sf est pair}~;
\endaligned\right.
\end{equation}
\item[{\bf (iii)}]
pour tout $\delta \in \N^d$ satisfaisant $\vert \delta \vert \leq
\ell$, les d\'eriv\'ees
suivantes de $h$ ou de $\overline{ h}$ convergent~{\rm :}
\def\theequation{6.12}\begin{equation}
\left\{
\aligned
\left[
\Upsilon^\delta h
\right] \left(
\Gamma_k ( z_{ (k)})
\right) 
& \
\in \C \{ z_{ (k)} \}^{ n'} \ \ \ \ \
\text{\sf si} \ k \
\text{\sf est impair},
\\
\left[
\underline{ \Upsilon}^\delta \overline{ h}
\right] \left(
\Gamma_k ( z_{ (k)})
\right) 
& \
\in \C \{ z_{ (k)} \}^{ n'} \ \ \ \ \
\text{\sf si} \ k \
\text{\sf est pair}.
\endaligned\right.
\end{equation}
\end{itemize}
\end{lemma}

Interpr\'etons la derni\`ere propri\'et\'e {\bf (iii)}~: on doit se
figurer que les diff\'erentiations $\mathcal{ L }^\beta$ et
$\underline{ \mathcal{ L} }^\beta$ sont <<horizontales>>, puisqu'elles
agissent le long des sous-vari\'et\'es de Segre complexifi\'ees
(conjugu\'ees)~; au contraire, les diff\'erentiations
$\Upsilon^\delta$ et $\underline{ \Upsilon }^\delta$ doivent \^etre
consid\'er\'ees comme <<transversales>>, puisque les 
$2m+ d$ champs de vecteurs $\mathcal{ L}_k$,
$\underline{ \mathcal{ L }}_k$ et $\Upsilon_j$ (ou $\underline{
\Upsilon}_j$) engendrent $T \mathcal{ M}$ ({\it cf.} le~\S5.15). Ainsi,
le Lemme~6.9 ram\`ene la convergence {\bf (i)} des jets de $h$ et de
$\overline{ h}$ \`a celle {\bf (iii)} de leurs seuls <<jets
transversaux>>.

\proof
Rappelons qu'en d\'eveloppant $\mathcal{ L}^\beta
\Upsilon^\delta h_{ i'}(t)$ (ou $\underline{ \mathcal{ L}}^\beta
\underline{ \Upsilon}^\delta \overline{ h}_{ i'} (\tau)$) par rapport
aux jets $\partial_t^\alpha h_{ i'} (t)$, on obtient les
formules~\thetag{ 5.18} et leurs inverses~\thetag{ 5.20} (ainsi que
leurs conjugu\'ees). Nous en d\'eduisons imm\'ediatement
l'\'equivalence entre {\bf (i)} et {\bf (ii)}.

{\bf (ii)} implique {\bf (iii)} trivialement. R\'eciproquement,
supposons $k$ impair pour fixer les id\'ees. Observons que la
$k$-i\`eme cha\^{\i}ne de Segre $\Gamma_k (z_{ (k) })$
s'\'ecrit $\mathcal{ L}_{ z_k} \left( \Gamma_{ k-1} (z_{ (k-1)})
\right)$. Soit $\phi(t) \in \C \dl t \dr^{ n'}$ une application
formelle arbitraire. Soit $z \in \C^m$. D'apr\`es la
d\'efinition~\thetag{ 3.13} du flot de $\mathcal{ L}$, on $\frac{
\partial }{\partial z_l} \left[ \phi (\mathcal{ L}_z (0)) \right]
\equiv \left[ \mathcal{ L}_l \phi \right] (\mathcal{ L}_z (0))$, pour
$l= 1,\, \dots,\, m$. Plus g\'en\'eralement, en notant $z_{ k;\, l}$,
$l=1,\, \dots,\, m$, les composantes de $z_k \in \C^m$, on a~:
\def\theequation{6.13}\begin{equation}
\frac{ \partial}{\partial z_{ k;\, l}}
\left[
\phi \left( \mathcal{
L}_{ z_k} \left( \Gamma_{k-1} 
(z_{ (k-1)}) \right) \right)\right] 
\equiv
\left[ \mathcal{ L}_l \phi \right]
\left( \Gamma_k (z_{ (k)}) \right).
\end{equation}
Pour un multiindice $\beta \in \N^m$ arbitraire, 
nous en d\'eduisons l'identit\'e~:
\def\theequation{6.14}\begin{equation}
\partial_{ z_k}^\beta \left[ \phi \left(
\Gamma_k (z_{ (k)})
\right) \right] \equiv
\left[
\mathcal{ L}^\beta \phi \right]
\left(
\Gamma_k (z_{ (k)})
\right).
\end{equation}

Supposons maintenant que $\left[ \Upsilon^\delta h \right] \left(
\Gamma_k (z_{ (k)}) \right)\in \C\{ 
z_{ (k)}\}^{ n'}$ est convergent. Puisque la propri\'et\'e
de convergence d'une s\'erie formelle est stable par d\'erivation
quelconque par rapport \`a ses variables, nous en d\'eduisons, en
appliquant la diff\'erentiation $\partial_{ z_k}^\beta$ \`a la s\'erie
$\left[ \Upsilon^\delta h \right] \left( \Gamma_k (z_{ (k)}) \right)$
et en utilisant au passage la formule~\thetag{ 6.14}, que $\left[
\mathcal{ L}^\beta \Upsilon^\delta \right] \left( \Gamma_k (z_{ (k)})
\right)$ est aussi convergent, pour tout $\beta \in \N^m$. Ceci
d\'emontre que {\bf (iii)} implique {\bf (ii)} dans le cas o\`u $k$
est impair. Le cas o\`u $k$ est pair se traite de mani\`ere analogue.
\endproof

\subsection*{ 6.15~Convergence des jets de
$h$ et de $\overline{ h}$ sur la premi\`ere cha\^{\i}ne de Segre} Nous
pouvons maintenant amorcer la d\'emonstration de la Proposition~6.4.
Rappelons que pour $k=1$, on a $z_{ (1)}\equiv z_1 \in \C^m$.
Rempla\c cons $(t,\, \tau)$ par $\Gamma_1 (z_1)$ dans la premi\`ere
ligne de~\thetag{ 6.5}, ce qui nous donne les $n'$ identit\'es
formelles~:
\def\theequation{6.16}\begin{equation}
0 \equiv
R_{ i'} ' 
\left(
\Gamma_1 (z_1),\, 
J_\tau^{\ell_0} \overline{ h}
(0) \, : h( \Gamma_1 (z_1))
\right), 
\end{equation}
pour $i' = 1,\, \dots,\, n'$. Au passage, nous avons utilis\'e la
relation connue $J_\tau^{ \ell_0} \overline{ h} \left(\Gamma_1 (z_1)
\right) \equiv J_\tau^{ \ell_0} \overline{ h} (0)$. G
r\^ace \`a l'analyticit\'e des $R_{ i'}'$,
gr\^ace au fait que $J_\tau^{\ell_0} \overline{ h} (0)$ est constant,
gr\^ace au Corollaire~2.9 et enfin, gr\^ace \`a l'hypoth\`ese analytique
principale~\thetag{ 6.6} de la Proposition~6.4, nous d\'eduisons
imm\'ediatement de~\thetag{ 6.16} que la restriction de $h$ \`a la
premi\`ere cha\^{\i}ne de Segre, {\it i.e.} $h( \Gamma_1 (z_1)) \equiv
h \left( z_1,\, \overline{ \Theta} (z_1,\, 0) \right) \in \C \{ z_1
\}^{ n'}$, est convergente.

La prochaine \'etape consiste \`a g\'en\'eraliser aux jets de $h$
d'ordre quelconque cette propri\'et\'e de convergence, apr\`es
composition par $\Gamma_1 (z_1)$. Nous verrons dans le \S6.29
ci-apr\`es pourquoi cela est n\'ecessaire. 

\def\thelemma{6.17}\begin{lemma}
Pour tout $\ell \in \N$, on 
a $J_t^\ell h \left(
\Gamma_1 (z_1)
\right) \in \C\{ z_1\}^{ N_{ n',\, n,\, \ell}}$.
\end{lemma}

\proof
Gr\^ace au Lemme~6.9, il
suffit d'\'etablir que
\def\theequation{6.18}\begin{equation}
\left[
\Upsilon^\delta h\right] \left(
\Gamma_1 (z_1)
\right) \in \C\{ 
z_1 \}^{ n'},
\end{equation}
pour tout $\delta \in \N^d$. Pour cela, nous aurons besoin de
formules analogues \`a~\thetag{ 6.13} et~\thetag{ 6.14}, en rempla\c
cant $\mathcal{ L}$ par $\Upsilon$ et $\underline{ \mathcal{ L}}$ par
$\underline{ \Upsilon}$. Soient $p = \left( z_p,\, w_p,\, \zeta_p,\,
\xi_p \right) \in \mathcal{ M}$, $w \in \C^d$ et $\xi \in \C^d$. Les
multiflots des champs de vecteurs $\Upsilon_j$ et $\underline{
\Upsilon}_j$ introduits dans~\thetag{ 5.16} ont pour expression
explicite
\def\theequation{6.19}\begin{equation}
\left\{
\aligned
\Upsilon_w (p) 
& \
:= 
\left(
z_p,\, w_p+ w,\, 
\zeta_p,\, 
\Theta (\zeta_p,\, z_p,\, w_p+ w)
\right) \ \ \ \ \
{\rm et} \\
\underline{ \Upsilon}_\xi ( p) 
& \
:= 
\left(
z_p,\, \overline{ \Theta} ( z_p,\, \zeta_p,\, 
\xi_p+ \xi),\, \zeta_p,\, \xi_p+\xi\right).
\endaligned\right.
\end{equation} 
Il est clair que $\Upsilon_w (p) \in \mathcal{ M}$ et que $\underline{
\Upsilon}_\xi (p) \in \mathcal{ M}$. Soit $\phi(t) \in \C \dl t \dr^{
n'}$. L'argument qui nous a permis d'\'etablir~\thetag{ 6.14}
permet de v\'erifier que pour tout $\delta \in \N^d$, la
relation formelle suivante est 
satisfaite~:
\def\theequation{6.20}\begin{equation}
\partial_{ w}^\delta \left[
\phi \left(
\Upsilon_w \left( \Gamma_1
(z_1)
\right)
\right)
\right] \equiv
\left[ \Upsilon^\delta \phi \right]
\left(
\Upsilon_w \left( \Gamma_1
(z_1)
\right)
\right).
\end{equation}

Rempla\c cons maintenant $(t,\, \tau)$ par $\Upsilon_w \left( \Gamma_1
(z_1)\right) \in \mathcal{ M}$ dans la premi\`ere ligne de~\thetag{
6.5}, ce qui donne les identit\'es formelles~:
\def\theequation{6.21}\begin{equation}
0 \equiv
R_{ i'} '
\left(
\Upsilon_w \left(
\Gamma_1 (z_1) 
\right),\, 
J_\tau^{ \ell_0} \overline{ h}
\left(\Upsilon_w \left(
\Gamma_1 (z_1) 
\right)
\right) \, : 
h \left(\Upsilon_w \left(
\Gamma_1 (z_1) 
\right)
\right)
\right),
\end{equation}
dans $\C \dl z_1,\, w \dr$, pour $i'=1,\, \dots,\, n'$. Dans ces
identit\'es, on ne peut effectuer aucune simplification en utilisant
les relations~\thetag{ 6.8}. N\'eanmoins, en diff\'erentiant ces
identit\'es par rapport \`a $w_j$, $j=1,\, \dots,\, d$, et en posant
$w=0$, on pourra appliquer la simplification cruciale $J_\tau^{
\ell_0} \overline{ h} (\Gamma_1 (z_1)) \equiv J_\tau^{ \ell_0}
\overline{ h}(0)$. Pour exprimer pr\'ecis\'ement ce que fournit
l'application de l'op\'erateur $\partial_{ w_j} (\cdot) \vert_{ w=0}$
\`a~\thetag{ 6.21}, groupons les variables dont d\'ependent les
s\'eries convergentes $R_{ i'}'$ en trois sous-groupes~: $(t,\,
\tau)$~; $J^{ \ell_0}$~; et $t'$. En appliquant l'op\'erateur
$\partial_{ w_j} (\cdot) \vert_{ w=0}$ \`a~\thetag{ 6.21}, gr\^ace \`a
la <<r\`egle de la cha\^{\i}ne>> pour la diff\'erentiation de s\'eries
compos\'ees et gr\^ace \`a~\thetag{ 6.20}, nous obtenons une somme
massive constitu\'ee de trois sous-groupes de polyn\^omes
diff\'erentiels, que nous \'ecrivons en d\'etail comme suit, juste
avant de formuler des commentaires
explicatifs~:
\def\theequation{6.22}\begin{equation}
\left\{
\aligned
{}
& \
0 \equiv
\frac{ \partial R_{ i'}'}{\partial w_j}
\left(
\Gamma_1 (z_1),\, J_\tau^{ \ell_0} \overline{ h}
(0)\, : h \left(
\Gamma_1 (z_1)
\right)
\right)+ \\
& \
+ 
\sum_{ l=1}^n\, 
\frac{\partial \Theta_l}{\partial w_j} 
(0,\, z_1,\, \overline{ \Theta}(z_1,\, 0)) \, 
\frac{ \partial R_{ i'}'}{\partial \xi_l}
\left(
\Gamma_1 (z_1),\, J_\tau^{ \ell_0} \overline{ h}
(0)\, : h \left(
\Gamma_1 (z_1)
\right)
\right) + \\
& \
+ \sum_{ i_1'=1}^{ n'}\, 
\sum_{ \alpha_1 \in \N^n,\, \vert
\alpha_1 \vert \leq \ell_0}\,
\frac{ \partial R_{ i'} '}{ \partial \left(
\partial_\tau^{ \alpha_1} \overline{ h}_{ i_1'}
\right)}\left(
\Gamma_1 (z_1),\, J_\tau^{ \ell_0} \overline{ h}
(0)\, : h \left(
\Gamma_1 (z_1)
\right)
\right) \cdot \Upsilon_j \partial_{\tau}^{ \alpha_1} 
\overline{ h}_{ i'} (0) + \\
& \
+
\sum_{ i_1 ' = 1}^{ n'} \, 
\frac{ \partial R_{ i'} '}{\partial t_{ i_1'}'}
\left(
\Gamma_1 (z_1),\, J_\tau^{ \ell_0} \overline{ h}
(0)\, : h \left(
\Gamma_1 (z_1)
\right)
\right) \cdot
\left[ \Upsilon_j h_{ i_1'}\right] (\Gamma_1 (z_1)),
\endaligned\right.
\end{equation}
pour $i'= 1,\, \dots,\, n'$ et $j=1,\, \dots,\, d$. Les deux
premi\`eres lignes explicitent les d\'eriv\'ees partielles par rapport
au premier groupe de variables $(t,\, \tau)$~; la troisi\`eme ligne
collecte les d\'eriv\'ees partielles par rapport au second groupe
$J_\tau^{\ell_0} \overline{ h}$~; la quatri\`eme ligne collecte les
d\'eriv\'ees partielles par rapport au troisi\`eme groupe
$t'$. D'apr\`es le d\'ebut du \S6.15, nous savons d\'ej\`a que $h(
\Gamma_1 (z_1)) \in \C\{ z_1 \}^{ n'}$ converge. En examinant tous
les termes de la premi\`ere, de la deuxi\`eme et de la troisi\`eme
ligne de~\thetag{ 6.22}, nous voyons qu'ils sont convergents. De
plus, les termes de la quatri\`eme ligne situ\'es avant le signe de
multiplication <<$\cdot$>> sont eux aussi convergents. Par
cons\'equent, nous pouvons r\'e\'ecrire~\thetag{ 6.22} sous la forme
contract\'ee~:
\def\theequation{6.23}\begin{equation}
0 \equiv
b_{ i',\, j}' (z_1) + 
\sum_{ i_1' =1}^{ n'} \, 
r_{i',\, i_1',\, j}' (z_1) \cdot
\left[
\Upsilon_j h_{ i_1'} \right] \left(
\Gamma_1 (z_1)
\right),
\end{equation}
pour $i'=1,\, \dots,\, n'$ et $j=1,\, \dots,\, d$. Ici, comme
nous venons de le remarquer, 
les s\'eries enti\`eres $b_{ i',\, j}'(z_1)$ et $r_{
i',\, i_1',\, j}'( z_1)$ sont convergentes. En fixant $j$ pour
appliquer le lemme suivant, nous d\'eduisons de~\thetag{ 6.23} que
$\left[ \Upsilon_j h_{ i_1'} \right] \left( \Gamma_1 (z_1) \right) \in
\C\{ z_1 \}$ converge, pour tout $i_1'=1,\, \dots,\, n'$ et tout
$j=1,\, \dots,\, d$.

\def\thelemma{6.24}\begin{lemma}
Soit $\nu \in \N$, $\nu \geq 1$, soit ${\sf x} \in \C^\nu$, soit $n'
\in \N$, $n' \geq 1$, soit $r_{ i',\, i_1'}' ( {\sf x} ) \in \C\{ {\sf
x} \}$ pour $i',\, i_1'= 1,\, \dots,\, n'$ et 
soit $b_{ i'}' ({\sf x}) \in
\C\{ {\sf x} \}$ pour $i' = 1,\, \dots,\, n'$. Si $n'$ s\'eries {\rm a
priori} formelles $a_{ i_1'}' ( {\sf x} ) \in \C \dl {\sf x} \dr$,
$i_1'=1,\, \dots,\, n'$, satisfont formellement le syst\`eme lin\'eaire
\def\theequation{6.25}\begin{equation}
0 \equiv 
b_{ i'}' ( {\sf x}) +
\sum_{ i_1'=1}^{ n'} \, 
r_{ i',\, i_1'} '( {\sf x} ) \, a_{ i_1'}' ( {\sf x}),
\end{equation}
et si
\def\theequation{6.26}\begin{equation}
{\rm det} \left(
r_{ i',\, i_1'} ' ({\sf x})
\right)_{ 1\leq i',\, i_1' \leq n'}
\not\equiv 0
\end{equation}
dans $\C \{ {\sf x} \}$, les s\'eries $a_{ i_1'} ' ( {\sf x}) \in \C\{
{\sf x}\}$ sont en fait convergentes.
\end{lemma}

On peut d\'emontrer ce lemme de mani\`ere \'el\'ementaire en
appliquant les formules de Cramer~: pour chaque $a_{ i_1'}' ({\sf
x})$, on obtient alors un quotient de s\'eries convergentes, dont le
d\'enominateur est constitu\'e du d\'eterminant~\thetag{ 6.26}~; par
unicit\'e des s\'eries formelles $a_{ i_1'}' ({\sf x})$, le
d\'enominateur est n\'ecessairement absorb\'e par le num\'erateur~; en
appliquant le th\'eor\`eme de division de K.~Weierstrass (version
formelle et version analytique), on en d\'eduit que les s\'eries $a_{
i_1'}' ( {\sf x}) \in \C\{ {\sf x}\}$ convergent. Observons aussi que
ce lemme est un cas particulier du Corollaire~2.9, si l'on consid\`ere
les \'equations analytiques $r_{ i'}' ({\sf x},\, {\sf y}) := b_{ i'}'
({\sf x}) + \sum_{ i_1'= 1}^{ n'} \, r_{ i',\, i_1'}' ({\sf x}) \,
{\sf y}_{ i_1'} = 0$.

\smallskip

En r\'esum\'e, nous avons \'etabli que $\left[ \Upsilon_j h \right]
\left( \Gamma_1 (z_1) \right) \in \C\{ z_1 \}^{ n'}$ pour $j=1,\,
\dots,\, d$, ce qui nous donne le Lemme~6.17 pour $\ell =1$. 

\smallskip

Pour traiter le cas g\'en\'eral, nous allons d\'emontrer~\thetag{
6.18} par r\'ecurrence. Soit $\ell \in \N$ tel que $\ell \geq 1$~;
supposons la propri\'et\'e~\thetag{ 6.18} satisfaite pour tout
multiindice $\delta_2 \in \N^d$ tel que $\vert \delta_2 \vert \leq
\ell$. Soit $j$ tel que $1\leq j \leq d$. Soit $\delta_1 \in \N^d$ un
multiindice arbitraire tel que $\vert \delta_1 \vert = \ell$.
Rappelons qu'un multiindice $\delta \in \N^d$ arbitraire tel que
$\vert \delta \vert = \ell+1$ s'\'ecrit $\delta = \delta_1 + \1_j^d$,
o\`u $\1_j^d$ d\'esigne le multiindice $(0,\,\dots,\,0,\, 1,\, 0,\,
\dots,\, 0)\in \N^d$, o\`u l'entier $1$ se situe \`a la $j$-i\`eme
place et o\`u $0$ se trouve aux autres places. En appliquant la
d\'erivation $\partial_{ w_j} \partial_w^{\delta_1}$ \`a~\thetag{
6.21}, en utilisant la <<r\`egle de la cha\^{\i}ne>> pour la
diff\'erentiation de s\'eries enti\`eres compos\'ees et en utilisant
la formule de Leibniz pour la diff\'erentiation des divers produits
qui apparaissent, nous pourrions \'ecrire sous forme d\'evelopp\'ee
une formule (compliqu\'ee) qui g\'en\'eraliserait~\thetag{ 6.22}, mais
un tel travail formel n'est pas r\'eellement n\'ecessaire. En effet,
en tenant compte de l'hypoth\`ese de r\'ecurrence, nous affirmons
qu'il existe des s\'eries enti\`eres convergentes $b_{ i',\,
\delta_1,\, j} '( z_1) \in \C\{ z_1 \}$ telles que l'application des
d\'erivations $\partial_{ w_j} \partial_w^{\delta_1}$\`a~\thetag{
6.21} fournit les identit\'es formelles suivantes~:
\def\theequation{6.27}\begin{equation}
0 \equiv
b_{ i',\, \delta_1,\, j}' (z_1) +
\sum_{ i_1'=1}^{ n'} \, 
\frac{ \partial R_{ i'}'}{\partial t_{ i_1'}'}
\left(
\Gamma_1 (z_1),\, J_\tau^{ \ell_0} \overline{ h}
(0)\, : h \left(
\Gamma_1 (z_1)
\right)
\right) \cdot
\left[ \Upsilon_j \Upsilon^{\delta_1} 
h_{ i_1'}\right] (\Gamma_1 (z_1)),
\end{equation}
pour $i'=1,\, \dots,\, n'$, et pour $j=1,\,
\dots,\, d$. En effet, le terme $b_{ i' ,\, \delta_1,\, j} '(z_1)$
collecte toutes les d\'eriv\'ees partielles de $R_{ i'}'$ (par rapport
aux trois groupes de variables $(t,\, \tau)$~; $J^{\ell_0}
\overline{ h}$~; et $t'$) autres que celles \'ecrites dans~\thetag{
6.27}~; ces d\'eriv\'ees partielles peuvent \^etre multipli\'ees par
des coefficients entiers et par des termes de la forme $\Upsilon^{
\delta_2} h_{ i_2'} \left( \Gamma_1 (z_1) \right)$, pour $i_2'= 1,\,
\dots,\, n'$ et pour $\vert \delta_2 \vert \leq \ell$. Tous ces termes
convergent, gr\^ace \`a l'hypoth\`ese de r\'ecurrence~: c'est pourquoi
$b_{ i',\, \delta_1,\, j}'( z_1) \in \C\{ z_1\}$.

Pour terminer, r\'e\'ecrivons~\thetag{ 6.27} sous la forme suivante
\def\theequation{6.28}\begin{equation}
0 \equiv
b_{ i',\, \delta_1,\, j}' (z_1) +
\sum_{ i_1'=1}^{ n'}\, 
r_{ i',\, i_1',\, j}(z_1) \cdot 
\left[ \Upsilon_j \Upsilon^{\delta_1} 
h_{ i_1'}\right] (\Gamma_1 (z_1)),
\end{equation}
qui est analogue \`a~\thetag{ 6.23}, avec les m\^emes coefficients
convergents $r_{ i',\, i_1',\, j} (z_1) \in \C\{ z_1\}$. Une
application directe du Lemme~6.24 permet de conclure que $\left[
\Upsilon_j \Upsilon^{ \delta_1} h_{ i_1'} \right] (\Gamma_1 (z_1))\in
\C\{ z_1 \}$ pour tout $j=1,\, \dots,\, d$, tout $\vert \delta_1 \vert
= \ell$ et tout $i_1'= 1,\, \dots,\, n'$. Ceci compl\`ete la
d\'emonstration de~\thetag{ 6.18} par r\'ecurrence sur $\vert \delta
\vert$. Le Lemme~6.17 est d\'emontr\'e.
\endproof

\subsection*{6.29.~Convergence des jets
de $h$ et de $\overline{ h}$ sur la seconde cha\^{ \i}ne de Segre}
Dans le pr\'ec\'edent paragraphe, nous n'avons pas insist\'e sur le
fait que les jets de $\overline{ h}$ convergent sur la premi\`ere
cha\^{\i}ne de Segre, parce que cette propri\'et\'e est triviale~: en
effet, gr\^ace \`a~\thetag{ 6.8}, le jet $J_\tau^\ell \overline{ h}
\left( \Gamma_1 (z_1) \right) \equiv J_\tau^\ell \overline{ h}(0)$ est
constant, donc convergent. Un ph\'enom\`ene analogue pr\'evaut sur la
deuxi\`eme cha\^{\i}ne de Segre~: puisque l'on a $J_t^\ell h \left(
\Gamma_2 (z_{ (2)}) \right) \equiv J_t^\ell h (\Gamma_1 (z_1))$
gr\^ace \`a~\thetag{ 6.8}, il d\'ecoule du Lemme~6.17 d\'emontr\'e
pr\'ec\'edemment que les jets de $h$ convergent (gratuitement) sur la
deuxi\`eme cha\^{\i}ne de Segre.

Pour d\'emontrer que la s\'erie enti\`ere $\overline{ h}\left(
\Gamma_2 (z_{ (2)}) \right)$ converge par rapport \`a $z_{ (2)} \in
\C^{ 2m}$, \'ecrivons la seconde ligne de~\thetag{ 6.5}
en y rempla\c cant $(t,\, \tau)$ par
$\Gamma_2 (z_{ (2)})$, ce qui nous donne les
$n'$ identit\'es formelles~:
\def\theequation{6.30}\begin{equation}
0 \equiv 
\overline{ R}_{ i'} ' \left(
\sigma \left(
\Gamma_2 (z_{ (2)})
\right),\,
J_t^{ \ell_0} h \left(
\Gamma_1 (z_1)
\right) \, : 
\overline{ h} \left(
\Gamma_2 (z_{ (2)})
\right)
\right),
\end{equation}
pour $i'=1,\, \dots,\, n'$. Ici, puisque le premier groupe $(\tau,\,
t)$ d'arguments de $\overline{ R}_{ i'}'$ dans~\thetag{ 6.5} est
\'ecrit dans l'ordre inverse de l'ordre habituel $(t,\, \tau)$, pour
\^etre rigoureux sur le plan notationnel, nous avons compos\'e la
cha\^{\i}ne $\Gamma_2 (z_{ (2)})$ avec l'involution holomorphe $\sigma
(t,\, \tau):= (\tau,\, t)$ introduite juste avant le~\S3.18. Mais
dans la Section~7 ci-dessous, nous admettrons le l\'eger \'ecart de
notation qui consiste \`a sous-entendre $\sigma$. Bien entendu, nous
avons utilis\'e~\thetag{ 6.8} pour simplifier $J_t^{ \ell_0} h \left(
\Gamma_2 (z_{(2)}) \right) \equiv J_t^{ \ell_0} h \left( \Gamma_1
(z_1) \right)$. Gr\^ace au Lemme~6.18, le terme $J_t^{ \ell_0} h
\left( \Gamma_1 (z_1) \right)$ est convergent~: c'est bien parce qu'il
appara\^{\i}t dans~\thetag{ 6.30} que nous avions besoin de
d\'emontrer \`a l'avance que les jets de $h$ convergent sur $\Gamma_1
(z_1)$, au moins jusqu'\`a l'ordre $\ell_0$. Par cons\'equent, les
\'equations~\thetag{ 6.30} sont convergentes par rapport \`a leurs
deux premiers groupes de variables. Pour appliquer le Corollaire~2.9,
il suffit \`a pr\'esent de v\'erifier que le d\'eterminant suivant~:
\def\theequation{6.31}\begin{equation}
{\rm det} \left(
\frac{ \partial \overline{ R}_{ i'}'}{\partial \bar
t_{ i_1'}'}
\left(
\sigma \left(
\Gamma_2 (z_{ (2)})
\right),\,
J_t^{ \ell_0} h \left(
\Gamma_1 (z_1)
\right) \, : 
\overline{ h} \left(
\Gamma_2 (z_{ (2)})
\right)
\right)
\right)_{ 1\leq i',\, i_1' \leq n'}
\end{equation}
ne s'annule pas identiquement dans $\C \dl z_{ (2)} \dr = \C\dl z_1,\,
z_2 \dr$. Or, en posant $z_1 =0$ et en tenant compte de la relation
$\Gamma_2 (0,\, z_2) = \underline{ \mathcal{ L}}_{ z_2} ( \mathcal{
L}_{ z_1} (0))\vert_{ z_1 =0} = \underline{ \mathcal{ L}}_{ z_2}
(\mathcal{ L }_0 (0)) = \underline{ \Gamma }_1 (z_2)$, ce
d\'eterminant se simplifie comme suit~:
\def\theequation{6.32}\begin{equation}
\left\{
\aligned
{\rm det} \left(
\frac{ \partial \overline{
R}_{ i'}'}{\partial \bar t_{ i_1'}'}
\left(
\sigma \left(
\underline{ \Gamma}_1 (z_2 )
\right),\,
J_t^{ \ell_0} h (0) \, : 
\overline{ h} \left(
\underline{ \Gamma}_1 (z_2)
\right)
\right)
\right)_{ 1\leq i',\, i_1' \leq n'} \\
\equiv
{\rm det} \left(
\frac{ \partial \overline{
R}_{ i'}'}{\partial \bar t_{ i_1'}'}
\left(
z_2,\, \Theta (z_2,\, 0),\, 0,\, 0,\, 
J_t^{ \ell_0} h (0) \, : 
\overline{ h} 
\left(
z_2,\, \Theta (z_2,\, 0)
\right)
\right)
\right)_{ 1\leq i',\, i_1' \leq n'}.
\endaligned\right.
\end{equation} 
Gr\^ace \`a l'hypoth\`ese principale~\thetag{ 6.6} (retrouv\'ee ici
\`a une conjugaison pr\`es), ce dernier d\'eterminant ne s'annule pas
identiquement dans $\C \dl z_2 \dr$, donc le d\'eterminant~\thetag{
6.31} ne s'annule pas identiquement dans $\C \dl z_{ (2)} \dr$ et le
corollaire~2.9 nous permet de conclure que $\overline{ h} \left(
\Gamma_2 (z_{ (2)}) \right) \in \C\{ z_{ (2)} \}^{ n'}$ converge.

Nous pourrions poursuivre le raisonnement et \'etablir une version du
Lemme~6.17 sur la deuxi\`eme cha\^{\i}ne de Segre, {\it i.e.}
d\'emontrer que $J_\tau^{ \ell} \overline{ h} \left( \Gamma_2 ( z_{
(2)}) \right) \in \C \{ z_{ (2)} \}^{ N_{ n',\, n\, \ell}}$ pour tout
$\ell \in \N$. Mais puisque nous avons maintenant pass\'e en revue et
comment\'e en d\'etail tous les points techniques d\'elicats, nous
pouvons d'ores et d\'ej\`a entamer la d\'emonstration finale de la
Proposition~6.4 en raisonnant par double r\'ecurrence sur la longueur
$k$ des cha\^{\i}nes de Segre et sur l'ordre $\ell$ des jets de $h$.

\subsection*{ 6.33.~D\'emonstration finale}
Pour conclure la Proposition~6.4, nous savons qu'il 
suffit d'\'etablir le lemme suivant, qui 
g\'en\'eralise le Lemme~6.17.

\def\thelemma{6.34}\begin{lemma}
Pour tout $k\in \N$ et tout $\ell \in \N$, on a
\def\theequation{6.35}\begin{equation}
\left\{
\aligned
J_t^\ell h \left(
\Gamma_k (z_{ (k)})
\right) 
& \
\in \C \{ z_{ (k)} \}^{ N_{ n',\, n,\, \ell}},
\ \ \ \ \
\text{\sf si} \ k \
\text{\sf est impair}~; \\
J_\tau^\ell \overline{ h} \left(
\Gamma_k (z_{ (k)})
\right) 
& \
\in \C \{ z_{ (k)} \}^{ N_{ n',\, n,\, \ell}}, 
\ \ \ \ \
\text{\sf si} \ k \
\text{\sf est pair}.
\endaligned\right.
\end{equation}
\end{lemma}

\proof
Pour $k=1$, ce lemme est d\'ej\`a d\'emontr\'e. Supposons
que la propri\'et\'e \thetag{ 6.35} est satisfaite 
au niveau $k$ et d\'emontrons-la au
niveau $k+1$. Pour fixer les id\'ees, nous supposerons $k$ 
impair\,--\,le cas o\`u $k$ est pair se traitant de mani\`ere similaire.

L'objectif 
est de d\'emontrer que $J_\tau^\ell \overline{ h}
\left( \Gamma_{ k +1} (z_{
(k+1)}) \right) \in \C \{ z_{ (k+1)} \}^{ N_{ n',\, n,\, \ell}}$
converge. Gr\^ace au Lemme~6.9, il suffit de d\'emontrer que
\def\theequation{6.36}\begin{equation}
\left[
\underline{ \Upsilon}^\delta 
\overline{ h}
\right] \left(
\Gamma_{k+1} (z_{ (k+1)})
\right) \in 
\C \dl z_{ (k)} \dr^{ n'},
\end{equation}
pour tout $\delta \in \N^d$. Nous allons raisonner par r\'ecurrence
sur l'entier $\vert \delta \vert$. 

Traitons d'abord le cas $\delta =0$, d\'ej\`a vu ci-dessus pour $k
+1=2$. Pour cela, posons $(t,\, \tau) := \Gamma_{ k+1 } (z_{ (k
+1)})$ dans la seconde ligne de~\thetag{ 6.5}, ce qui nous donne~:
\def\theequation{6.37}\begin{equation}
0 \equiv
\overline{ R}_{ i'} '
\left(
\sigma \left(
\Gamma_{ k+1} (z_{ (k+1)})
\right),\, 
J_t^{\ell_0} h \left(
\Gamma_k (z_{ (k)})
\right) \, : 
\overline{ h} \left(
\Gamma_{ k+1} (z_{ (k+1)})
\right)
\right),
\end{equation}
pour $i'= 1,\, \dots,\, n'$. Gr\^ace \`a l'hypoth\`ese de
r\'ecurrence, le deuxi\`eme argument de $\overline{ R}_{ i'}'$
est convergent. Pour appliquer le Corollaire~2.9, il suffit donc de
v\'erifier que le d\'eterminant suivant~:
\def\theequation{6.38}\begin{equation}
{\rm det} \left(
\frac{ \partial \overline{ R}_{ i'}'}{\partial \bar
t_{ i_1'}'}
\left(
\sigma \left(
\Gamma_{k+1} (z_{ (k+1)})
\right),\,
J_t^{ \ell_0} h \left(
\Gamma_k (z_k)
\right) \, : 
\overline{ h} \left(
\Gamma_{ k+1} (z_{ (k+1)})
\right)
\right)
\right)_{ 1\leq i',\, i_1' \leq n'}
\end{equation}
ne s'annule pas identiquement dans $\C \dl z_{ (k+1 )} \dr = \C \dl
z_{ (k)},\, z_{ k+ 1} \dr$. Or, en posant $z_{ (k)}= 0$ et en tenant
compte de la relation $\Gamma_{ k+1} (0,\, z_{ k+1 }) = \underline{
\Gamma}_1 (z_{ k+1})$, ce d\'eterminant se simplifie et co\"{\i}ncide
avec le d\'eterminant~\thetag{ 6.32} dans lequel $z_2$ est rempla\c
c\'e par $z_{ k+1}$. Gr\^ace \`a l'hypoth\`ese principale~\thetag{
6.6}, il ne s'annule pas identiquement dans $\C \dl z_{ k+1} \dr$,
donc le d\'eterminant~\thetag{ 6.38} ne s'annule pas identiquement
dans $\C \dl z_{ (k+1)} \dr$ et le corollaire~\thetag{ 2.9} nous
permet de conclure que $\overline{ h} \left( \Gamma_{ k+1} (z_{
(k+1)}) \right) \in \C \{ z_{ (k+1)}\}^{ n'}$ converge.

Soit $\ell \in\N$ avec $\ell \geq 0$~; supposons la
propri\'et\'e~\thetag{ 6.36} vraie pour tout multiindice $\delta_2 \in
\N^d$ tel que $\vert \delta_2 \vert \leq \ell$. Soit $\delta_1 \in
\N^d$ un multiindice arbitraire tel que $\vert \delta_1 \vert =
\ell$. Consid\'erons le multiindice $\delta_1 + \1_j^d$.

En rempla\c cant $(t,\, \tau)$ par $\underline{ \Upsilon}_\xi \left(
\Gamma_{ k+1} (z_{ (k+1)}) \right)$ dans la deuxi\`eme ligne
de~\thetag{ 6.5}, o\`u $\xi \in \C^d$, nous obtenons les
identit\'es
formelles~:
\def\theequation{6.39}\begin{equation}
\small
0 \equiv
\overline{ R}_{ i'} ' \left(
\sigma \left(
\underline{ \Upsilon}_\xi 
\left(
\Gamma_{ k+1} (z_{ (k+1)})
\right)
\right),\, 
J_t^{ \ell_0} h\left(
\underline{ \Upsilon}_\xi 
\left(
\Gamma_{ k+1} (z_{ (k+1)})
\right)
\right) \, : 
\overline{ h}
\left(
\underline{ \Upsilon}_\xi 
\left(
\Gamma_{ k+1} (z_{ (k+1)})
\right)
\right)
\right),
\end{equation}
dans $\C \dl z_{ (k+1)},\, \xi\dr$ pour $i'= 1,\, \dots,\, n'$.
Appliquons la d\'erivation $\left.
\partial_{ \xi_j}\partial_\xi^{\delta_1}
(\cdot) \right\vert_{ \xi=0}$ \`a~\thetag{ 6.39}, 
ce qui revient \`a appliquer
$\underline{ \Upsilon}_j \underline{
\Upsilon}^{ \delta_1}(\cdot) \vert_{ \xi=0}$. 
Nous
obtenons un polyn\^ome (compliqu\'e) dont nous extrayons seulement le
terme que nous consid\'erons comme principal~; ainsi nous
r\'esumons le r\'esultat du calcul sous la forme~:
\def\theequation{6.40}\begin{equation}
\left\{
\aligned
0 
& \
\equiv
b_{ i',\, \delta_1,\, j}' (z_{ (k+1)})+ \\
& \
\ \ \ \ \
+
\sum_{ i_1'= 1}^{ n'} \, 
\frac{ \partial \overline{ R}_{ i'}}{\partial 
\bar t_{ i_1'}'} 
\left(
\Gamma_{ k+1}( z_{ (k+1)}),\, 
J_t^{\ell_0} h \left(
\Gamma_k (z_{ (k)})
\right) \, : 
\overline{ h} \left(
\Gamma_{ k+1} (z_{ (k+1)})
\right)
\right)\cdot \\
& \ 
\ \ \ \ \ \ \ \ \ \ \ \ \ \ \ \ \ \
\ \ \ \ \ \ \ \ \ \ \ \ \ \ \ \ \ \
\ \ \ \ \ \ \ \ \ \ \ \ \ \ \ \ \ \
\ \ \ \ \ \ \ \ \ \ \ \ \ \ \ \ \ \
\cdot
\left[
\Upsilon_j \Upsilon^{ \delta_1} \overline{ h} \right]
\left(
\Gamma_{ k+1} (z_{ (k+1)})
\right), 
\endaligned\right.
\end{equation}
pour $i'= 1,\, \dots,\, n'$, o\`u le terme $b_{ i',\, \delta_1,\, j}'
(z_{ (k+1)})$ est le <<reste>> qui incorpore tous les termes qui
apparaissent, except\'e celui qui est \'ecrit \`a droite du signe
<<$+$>>. Nous affirmons que $b_{ i',\, \delta_1,\, j}'
(z_{ (k+1)})$ est convergent.

En effet, le terme $b_{ i',\, \delta_1,\, j}' (z_{ (k+1)})$ est un
polyn\^ome \`a coefficients entiers qui contient comme
mon\^omes deux cat\'egories de termes. La
premi\`ere cat\'egorie est constitu\'ee des
d\'eriv\'ees <<transversales>> suivantes~:
\def\theequation{6.41}\begin{equation}
\left[
\underline{ \Upsilon}^{\delta_2} h_{ i_2'} \right]
\left(
\Gamma_{ k+1} (z_{ (k+1)})
\right),
\end{equation}
o\`u $i_2' = 1,\, \dots,\, n'$ et o\`u, de mani\`ere cruciale, la
longueur de $\delta_2$ satisfait l'in\'egalit\'e $\vert \delta_2 \vert
\leq \ell$~: en effet, le seul terme qui incorpore une d\'eriv\'ee
<<transversale>> d'ordre $\ell+ 1$ est celui qui est \'ecrit \`a
droite du signe <<$+$>> dans~\thetag{ 6.40}. Gr\^ace \`a
l'hypoth\`ese de r\'ecurrence, tous les termes~\thetag{ 6.41}
convergent. La deuxi\`eme cat\'egorie de termes est consitut\'ee des
d\'eriv\'ees partielles suivantes de $\overline{ R}_{ i'}$~:
\def\theequation{6.42}\begin{equation}
\left\{
\aligned
{}
& \
\left[
\frac{ \partial^{\eta_2}}{
\partial t_{i''(1)}' \cdots
\partial t_{i''(\eta_2)}'} \
\frac{\partial^{ \varepsilon_2}}{
\partial \left(
\partial_t^{ \alpha(1)} h_{ i'(1)}
\right)\cdots\partial \left(
\partial_t^{ \alpha(\varepsilon_2)} h_{ i'(\varepsilon_2)}
\right)}
\underline{ \Upsilon}^{ \delta_2}
\left(
\overline{ R}_{ i'}
\right)\right] \\
& \ 
\ \ \ \ \ \ \ \ 
\ \ \ \ \ \ \ \ 
\ \ \ \ \ \ \ \
\left(
\Gamma_{ k+1}( z_{ (k+1)}),\, 
J_t^{\ell_0} h \left(
\Gamma_k (z_{ (k)})
\right) \, : 
\overline{ h} \left(
\Gamma_{ k+1} (z_{ (k+1)})
\right)
\right)
\endaligned\right.
\end{equation}
(dans cette expression, la deuxi\`eme ligne d\'esigne l'argument de la
premi\`ere ligne), o\`u $1\leq i''(1),\, \dots,\, i''(\eta_2) \leq
n'$, o\`u $\alpha(1),\, \dots,\, \alpha (\varepsilon_2) \in \N^n$,
o\`u $1\leq i'(1),\, \dots,\, i'(\varepsilon_2) \leq n'$ et o\`u
$\eta_2+\varepsilon_2 + \delta_2 + \leq \ell+1$. Gr\^ace \`a
l'hypoth\`ese de r\'ecurrence (il s'agit en fait de l'hypoth\`ese de
r\'ecurrence qui porte sur $k$ et non sur $\ell$), les termes~\thetag{
6.42} convergent aussi.

Puisque $b_{ i',\, \delta_1,\, j}' (z_{ (k+1)})$ converge et puisque
le d\'eterminant~\thetag{ 6.38} est convergent et ne s'annule pas
identiquement, le Lemme~\thetag{ 6.24} s'applique. Nous en d\'eduisons
que $\underline{ \Upsilon}_j \underline{ \Upsilon}^{ \delta_1}
\overline{ h}_{ i_1'} \left( \Gamma_{ k+1} (z_{ (k+1)}) \right) \in
\C\{ z_{ (k+1)}\}$ pour $j=1,\, \dots,\, d$, $\vert \delta_1 \vert
\leq \ell$ et $i'=1,\, \dots,\, n'$~: c'est la propri\'et\'e~\thetag{
6.36} au niveau $\vert \delta \vert = \ell+1$.

Ceci compl\`ete la d\'emonstration du Lemme~6.34.
\endproof

La Proposition~6.4 est d\'emontr\'ee.
\endproof

\section*{\S7.~Convergence de l'application de r\'eflexion 
CR formelle}

\subsection*{ 7.1.~Pr\'eliminaire}
Cette derni\`ere section est consacr\'ee \`a d\'emontrer le
Th\'eor\`eme~1.23. Contrairement aux deux pr\'ec\'edentes
d\'emonstrations, nous ne pourrons pas ramener les raisonnements \`a
la seule consid\'eration des composantes de $h$~: il nous faudra
\'etudier simultan\'ement le nombre {\it infini} des composantes
$\Theta_{ j',\, \gamma'}' (h(t))$ de l'application de
r\'eflexion~\thetag{ 1.28}. Avertissons d\`es maintenant le lecteur~:
la d\'emonstration du Th\'eor\`eme~1.23 {\it n'est ni courte ni
simple}~; ce th\'eor\`eme constitue en effet une avanc\'ee
substantielle par rapport aux travaux de S.M.~Baouendi, P.~Ebenfelt et
L.-P.~Rothschild, dont le principe \'el\'ementaire a \'et\'e expos\'e
dans le \S5~; comme nous l'avons expliqu\'e dans l'Introduction, cette
avanc\'ee ne se limite pas \`a \`a la simple greffe de techniques
connues. Pour toutes ces raisons, nous devrons faire preuve dans cette
section d'un r\'eel souci d'exposition, afin de pr\'esenter nos
calculs majeurs de la mani\`ere la plus accessible qui soit.

\smallskip

Soit $\chi (t)$ une application formelle et soit $\ell \in \N$.
Rappelons que nous notons $J_t^\ell \chi(t)$ son jet d'ordre
$\ell$. Par exemple, le jet d'ordre $\ell$ de l'application formelle
$\Theta_{ \gamma '} '( h(t))$ de $\C^n$ dans $\C^{ d'}$ est
constitu\'e de toutes les d\'eriv\'ees partielles par rapport \`a $t$
de l'application formelle compos\'ee $t \longmapsto_{ \mathcal{ F}}
\Theta_{ \gamma'} '( h(t))$. De m\^eme, le jet $J_t^\ell \mathcal{
R}_h' (\tau',\, t)$ s'identifie \`a~: $J_t^\ell \xi'-\sum_{ \gamma'
\in \N^{ m'}} \, (\zeta')^{ \gamma'} \, J_t^\ell \Theta_{ \gamma'} '
(h(t))$~; ici, pour $\alpha =0$, la d\'eriv\'ee partielle
$\partial_t^\alpha \xi'$ s'identifie \`a $\xi'$, mais dispara\^{\i}t
pour $\vert \alpha \vert \geq 1$.

Pour commencer, \'enon\c cons ce que fournit la diff\'erentiation de
l'application formelle compos\'ee $t \longmapsto_{ 
\mathcal{ F}} \Theta_{ \gamma'} '( h(t))$.

\def\thelemma{7.2}\begin{lemma}
Soit $\alpha \in \N^n$. Il existe une application polynomiale
universelle $Q_\alpha$ \`a valeurs dans $\C^{ d'}$, lin\'eaire par
rapport \`a son second groupe de variables et ne
d\'ependant que du jet strict d'ordre $\vert \alpha \vert$
de $h$, qui satisfait~{\rm :}
\def\theequation{7.3}\begin{equation}
\left\{
\aligned
\partial_t^\alpha 
\left[
\Theta_{ \gamma'} ' (h(t))
\right] 
& \
\equiv
Q_\alpha
\left(
J_t^{ \vert \alpha \vert} h(t), \
\left[
J_{ t'}^{ \vert \alpha \vert}
\Theta_{ \gamma'}' \right] (h(t))
\right), \\
0 
& \
\equiv
Q_\alpha \left(J_t^{ \vert \alpha \vert} h(t), \
0
\right).
\endaligned\right.
\end{equation}
\end{lemma}
La propri\'et\'e d'annulation \'ecrite \`a la seconde ligne d\'ecoule
de la lin\'earit\'e par rapport au second groupe de variables~; elle
sera utilis\'ee dans la d\'emonstration du Lemme~7.7 ci-dessous. Pour
\'etablir cet \'enonc\'e, on proc\`ede en raisonnant par
r\'ecurrence~; la d\'emonstration est \'el\'ementaire, puisqu'on ne
demande pas l'expression explicite des applications $Q_\alpha$,
inutile pour la suite. Mentionnons toutefois qu'une application de la
formule de Fa\`a di Bruno \`a plusieurs variables fournirait cette
expression.

\subsection*{ 7.4.~Estim\'ees de Cauchy}
Dans la conclusion du Th\'eor\`eme~1.23, la convergence de
l'application de r\'eflexion~\thetag{ 1.28} \'equivaut \`a deux
propri\'et\'es bien distinctes~:
\begin{itemize}
\item[{\bf (i)}]
pour tout $\gamma' \in \N^{ m'}$, la s\'erie formelle vectorielle
$\Theta_{ \gamma'} ' (h(t)) \in \C\{ t \}^{ d'}$ converge~;
\item[{\bf (ii)}]
il existe des constantes $\sigma >0$, $C '>$ et $\rho' >0$ telles que
\def\theequation{7.5}\begin{equation}
\vert t \vert < \sigma \Longrightarrow
\left\vert
\Theta_{ \gamma'} ' (h(t))
\right\vert < C' \, (\rho')^{ - \vert \gamma ' \vert},
\end{equation}
pour tout $\gamma ' \in \N^{m'}$.
\end{itemize}

La seconde propri\'et\'e est une {\sl estim\'ee de Cauchy}. Nous
souhaitons \'enoncer d\`es maintenant qu'elle d\'ecoule de la
premi\`ere (Lemme~7.7 {\it infra}). Gr\^ace \`a cette remarque, nous
pourrons nous dispenser dans la suite de toute consid\'eration au
sujet des estim\'ees de croissance et nous contenter d'\'etablir
seulement des propri\'et\'es de convergence ({\it voir}~la Proposition
principale~7.18 ci-dessous).

Comme dans la Section~6, nous travaillerons avec la
cha\^{\i}ne de Segre $\Gamma_k$, sans utiliser la cha\^{\i}ne
conjugu\'ee $\underline{ \Gamma}_k$. Soit une
s\'erie formelle $\Psi' (t') \in \C \dl t'
\dr$ sans terme constant. Nous noterons
\def\theequation{7.6}\begin{equation}
\left[
J_t^\ell \Psi'(h)
\right]\left(
\Gamma_k (z_{ (k)})
\right) := 
\left(
\left.
\partial_t^\alpha 
\left[
\Psi' (h(t))
\right] \right\vert_{ t = \Gamma_k (z_{ (k)})}
\right)_{ \vert \alpha \vert \leq \ell}
\end{equation}
le jet de l'application formelle compos\'ee $t\longmapsto \Psi'(h(t))$,
dans lequel on remplace l'argument $t$ par $\Gamma_k (z_{ (k)})$.
Dans la suite, nous utiliserons r\'eguli\`erement des crochets et des
grandes parenth\`eses, en vue d'une clart\'e maximale~; en effet, loin
d'alourdir les calculs, le respect d'un parenth\'esage rigoureux
permettra d'\'eliminer toute ambigu\"{\i}t\'e notationnelle. Voici
donc l'\'enonc\'e qui nous dispensera de mentionner les estim\'ees de
Cauchy. Nous travaillerons avec une cha\^{\i}ne de Segre de longueur
$k$ arbitraire et avec des jets d'ordre $\ell$ arbitraire.

\def\thelemma{7.7}\begin{lemma}
Soit $k\in\N$ et $\ell \in \N$. Si $k$ est {\sf impair}, 
les propri\'et\'es suivantes sont \'equivalentes~{\rm :}
\begin{itemize}
\item[{\bf (1)}]
$\left[
J_t^\ell \mathcal{ R}_h'
\right] \left(
\tau',\, \Gamma_k (z_{ (k)})
\right) \in \C \{ 
\tau',\, z_{ (k)} \}^{ N_{ d',\, n,\, \ell}}$~{\rm ;}
\item[{\bf (2)}]
pour tout $\gamma' \in \N^{ m'}$, on a~{\rm :} $\left[ J_t^\ell \Theta_{
\gamma'} ' (h) \right] \left( \Gamma_k (z_{ (k)}) \right) \in \C \{
z_{ (k)}\}^{ N_{ d',\, n,\, \ell}}$, et il existe des constantes
$\sigma >0$, $C' >0$ et $\rho'
>0$ telles que
\def\theequation{7.8}\begin{equation}
\vert z_{ (k)} \vert < \sigma
\Longrightarrow 
\left\vert
J_t^\ell \Theta_{ \gamma'} ' \left(h\left(
\Gamma_k (z_{ (k)}) \right) \right)
\right\vert < 
C ' \, (\rho')^{ - \vert \gamma '\vert}~;
\end{equation}
\item[{\bf (3)}]
pour tout $\gamma' \in \N^{ m'}$, on a~{\rm :} $\left[ J_t^\ell \Theta_{
\gamma'} ' (h) \right] \left( \Gamma_k (z_{ (k)}) \right) \in \C \{
z_{ (k)}\}^{ N_{ d',\, n,\, \ell}}$.
\end{itemize}
Si $k$ est {\sf pair}, on a trois propri\'et\'es \'equivalentes analogues
conjugu\'ees, obtenues en rempa\c cant $t$ par $\tau$, $h$ par
$\overline{ h}$, $\Theta_{ \gamma'} '$ par $\overline{ \Theta }_{
\gamma'}'$ et $\mathcal{ R}_h' \left( \tau',\, t\right)$ par
$\overline{ \mathcal{ R}}_h' 
\left( t',\, \tau \right) := w' - \overline{
\Theta} ' \left( z',\, \overline{ h} (\tau) \right)$.

\end{lemma}

La troisi\`eme propri\'et\'e exprime seulement que le jet d'ordre
$\ell$ de chaque composante $\Theta_{ \gamma'} ' (h(t))$ de
l'application de r\'eflexion converge sur la $k$-i\`eme cha\^{\i}ne de
Segre, sans estimation de Cauchy.

\proof
Puisque 
\def\theequation{7.9}\begin{equation}
J_t^\ell \mathcal{ R}_h' \left(
\tau',\, \Gamma_k (z_{ (k)})
\right)
=
J_t^\ell \xi'
- \sum_{ \gamma' \in \N^{m'}} \, 
(\zeta')^{ \gamma'} \,
\left[ 
J_t^\ell
\Theta_{ \gamma'} ' (h) \right] \left( 
\Gamma_k (z_{ (k)}) \right),
\end{equation} 
l'\'equivalence entre les propri\'et\'es {\bf (1)} et {\bf (2)}
provient de la th\'eorie \'el\'ementaire des fonctions analytiques de
plusieurs variables complexes.

L'implication {\bf (2)} $\Rightarrow$ {\bf (3)} est triviale.
R\'eciproquement, supposons que la propri\'et\'e de convergence
$\left[ J_t^\ell \Theta_{ \gamma'}' (h) \right] \left( \Gamma_k (z_{
(k)}) \right) \in \C\{ z_{ (k)}\}^{ N_{ d',\, n,\, \ell}}$ est
satisfaite pour tout $\gamma' \in \N^{ m'}$. Notons $\theta_{
\gamma',\, \ell}' ( z_{ (k)}) \in \C\{ z_{ (k)}\}^{ N_{ d',\, n,\,
\ell}}$ ces s\'eries vectorielles convergentes. L'objectif est
d'\'etablir qu'elles satisfont une estim\'ee de Cauchy.

Pour cela, partons des identit\'es formelles vectorielles~:
\def\theequation{7.10}\begin{equation}
0 \equiv
\left[ J_t^\ell \Theta_{ \gamma'}' 
(h) \right] \left( \Gamma_k (z_{
(k)}) \right) - \theta_{
\gamma',\, \ell}' ( z_{ (k)}),
\end{equation}
valables dans $\C\dl z_{ (k)} \dr^{
N_{ d',\, n,\, \ell}}$.
Gr\^ace au Lemme~7.2, nous pouvons d\'evelopper le premier terme~:
\def\theequation{7.11}\begin{equation}
0 \equiv 
Q_\ell 
\left(
J_t^\ell h 
\left( \Gamma_k (z_{ (k)}) \right), \
\left[
J_{ t'}^{ \ell}
\Theta_{ \gamma'}'
\right] \left( h \left(
\Gamma_k (z_{ (k)})\right) \right)
\right) -
\theta_{
\gamma',\, \ell}' ( z_{ (k)}).
\end{equation}
Ces identit\'es sont satisfaites par l'application formelle $J_t^\ell
h \left( \Gamma_k (z_{ (k)}) \right) \in \C\dl z_{ (k)} \dr^{ N_{
n',\, n,\, \ell}}$. Notons que les identit\'es~\thetag{ 7.11} sont en
nombre infini, puisque $\gamma'$ varie dans $\N^{ m'}$~; heureusement,
le Th\'eor\`eme d'approximation~2.4 s'applique dans une telle
situation, gr\^ace \`a la remarque qui suit son \'enonc\'e. Dans le
Th\'eor\`eme~2.4, on suppose aussi que toutes les s\'eries formelles
consid\'er\'ees s'annulent \`a l'origine. Ici, $h(0)=0$, mais le {\sl
jet strict} de $h$ d'ordre $\ell$ \`a lorigine ({\it i.e.} toutes les
d\'eriv\'ees partielles $\partial_t^{ \alpha_1 } h_{ i_1'} (0)$, pour
$1 \leq \vert \alpha_1 \vert \leq \ell$ et pour $1\leq i_1 \leq n'$)
ne s'annule pas forc\'ement. Heureusement, en soustrayant de
l'\'equation vectorielle~\thetag{ 7.11} l'\'equation
vectorielle~\thetag{ 7.11} prise en $z_{ (k)}= 0$, on se ram\`ene \`a
des s\'eries formelles qui sont toutes sans terme constant.

En appliquant donc le Th\'eor\`eme~2.4 \`a~\thetag{ 7.11} avec l'ordre
d'approximation $N = 1$, nous d\'eduisons qu'il existe une application
convergente ${\sf H}_\ell \left( z_{ (k )} \right) \in \C\{ z_{ (k
)}\}^{ N_{ n',\, n,\, \ell}}$, que nous \'ecrirons sous forme
d\'evelopp\'ee comme suit~: 
${\sf H}_\ell \left( z_{ (k)} \right) = \left( {\sf
H}_\alpha \left( z_{ (k)} \right) \right)_{ \vert \alpha \vert \leq
\ell}$, telle que ${\sf H}_{ \ell} (0) = J_t^\ell h (0)$ et telle que
les \'equations formelles suivantes sont satisfaites~:
\def\theequation{7.12}\begin{equation}
0 \equiv 
Q_\ell 
\left(
{\sf H}_\ell \left( z_{ (k)} \right), \
\left[
J_{ t'}^\ell 
\Theta_{ \gamma'}'
\right] \left( {\sf H}_0 \left(
z_{ (k)} \right) \right)
\right) -
\theta_{
\gamma',\, \ell}' ( z_{ (k)}).
\end{equation}
Utilisons maintenant les estim\'ees de Cauchy satisfaites par le jet
d'ordre $\ell$, par rapport \`a $t'$, de la s\'erie d\'efinissante
$\Theta' (\zeta',\, t') = \sum_{ \gamma' \in \N^{ m'}} \, (\zeta')^{
\gamma'} \, \Theta_{ \gamma'} ' (t')$~: il existe des constantes
$\sigma' > 0$, $\widetilde {C} ' >0$ et $\rho' 
>0$ telles que
\def\theequation{7.13}\begin{equation}
\vert t' \vert < \sigma' 
\Longrightarrow 
\left\vert
J_{ t'}^\ell \Theta_{ \gamma'} ' (t') 
\right\vert < \widetilde{ C} ' \, 
(\rho')^{ - \vert \gamma ' \vert}.
\end{equation}
Soit $\sigma >0$ suffisamment petit pour que ${\sf H}_\ell\left(
\Gamma_k (z_{ (k)}) \right)$ converge dans le cube $\square_{
\sigma}^{ km} = \{z_{ (k)} \in \C^{ km} : \, \left\vert z_{ (k)}
\right\vert < \sigma$, et pour que $\vert z_{ (k)} \vert < \sigma$
implique $\left\vert {\sf H}_\ell \left( z_{ (k)}\right) -J_t^\ell
h(0) \right\vert < \sigma'$, et en particulier $\left\vert {\sf H}_0
\left( z_{ (k)}\right) \right\vert < \sigma '$. On en d\'eduit~:
\def\theequation{7.14}\begin{equation}
\vert z_{ (k)} \vert < \sigma 
\Longrightarrow 
\left\vert
\left[
J_{ t'}^\ell \Theta_{ \gamma'} '\right] \left( {\sf H}_0 
\left( z_{ (k)} 
\right) \right)
\right\vert < \widetilde{ C} ' \, 
(\rho')^{ - \vert \gamma ' \vert}.
\end{equation}
Enfin, en utilisant le fait que l'application polynomiale $Q_\ell$ est
lin\'eaire par rapport \`a son second groupe de variables et en
majorant les coefficients de cette expression lin\'eaire, qui sont des
polyn\^omes par rapport aux
variables ${\sf H}_{ i_1',\, \alpha_1} 
\left( z_{ (k)} \right)$, pour
$1 \leq \vert \alpha_1 \vert \leq
\vert \alpha \vert$ et pour $i_1' = 1,\, \dots,\, n'$, on
v\'erifie qu'il existe une constante $C' \geq \widetilde{ C}'$ telle que
\def\theequation{7.15}\begin{equation}
\vert z_{ (k)} \vert < \sigma 
\Longrightarrow 
\left\vert
Q_\ell 
\left(
{\sf H}_\ell \left( z_{ (k)} \right), \
\left[
J_{ t'}^\ell 
\Theta_{ \gamma'}'
\right] \left( {\sf H}_0 \left(
z_{ (k)} \right) \right)
\right)
\right\vert <
C ' \, 
(\rho')^{ - \vert \gamma ' \vert}.
\end{equation}
Gr\^ace \`a l'identit\'e~\thetag{ 7.12}, on en d\'eduit l'estimation
de Cauchy d\'esir\'ee~:
\def\theequation{7.16}\begin{equation}
\vert z_{ (k)} \vert < \sigma 
\Longrightarrow 
\left\vert
\theta_{ \gamma',\, \ell} '( z_{ (k)}) 
\right\vert <
C ' \, 
(\rho')^{ - \vert \gamma ' \vert}.
\end{equation}
Ceci compl\`ete la d\'emonstration du 
Lemme~7.7.
\endproof

\subsection*{ 7.17.~\'Enonc\'e de la proposition principale}
Gr\^ace au Lemme~7.7, nous ramenons la d\'emonstration du
Th\'eor\`eme~1.23 \`a des \'enonc\'es de convergence pour les jets des
composantes de l'application de r\'eflexion.

\def\theproposition{7.18}\begin{proposition}
Pour tout $k\in \N$, tout $\ell \in \N$, les propri\'et\'es de
convergence suivantes sont satisfaites~{\rm :}
\begin{itemize}
\item[{\bf (1)}]
{\sf si $k$ est impair}, pour tout $\gamma' \in \N^{ m'}$,
on a~{\rm :}
\def\theequation{7.19}\begin{equation}
\left[J_t^\ell \Theta_{ \gamma'}' (h) \right]
\left(
\Gamma_k \left(
z_{ (k)} \right) 
\right) \in \C\{ z_{ (k)} \}^{ N_{ d',\, n,\, \ell}}~;
\end{equation}
\item[{\bf (2)}]
{\sf si $k$ pair}, pour tout $\gamma' \in \N^{ m'}$,
on a~{\rm :} 
\def\theequation{7.20}\begin{equation}
\left[J_\tau^\ell \overline{\Theta}_{ 
\gamma'}' ( \overline{ h}) \right]
\left(
\Gamma_k \left(
z_{ (k)} \right) 
\right) \in \C\{ z_{ (k)} \}^{ N_{ d',\, n,\, \ell}}.
\end{equation}
\end{itemize}
\end{proposition}

Cette proposition implique le Th\'eor\`eme~1.23, gr\^ace au Lemme~7.7
et gr\^ace au Lem\-me~3.32. Notons que pour $k$ et $\ell$ fix\'es, il y
a une infinit\'e de propri\'et\'es de convergence dans~\thetag{ 7.19}
et dans~\thetag{ 7.20} ({\it cf.}~le~\S1.42). De plus, rappelons qu'en
vertu de~\thetag{ 6.8}, on a~:
\def\theequation{7.21}\begin{equation}
\left\{
\aligned
\left[
J_\tau^\ell 
\overline{ \Theta}_{ \gamma'}' \left(
\overline{ h}
\right)\right] 
\left(
\Gamma_k \left(
z_{ (k)}
\right)
\right) 
& \
\equiv 
\left[
J_\tau^\ell \overline{ \Theta}_{ \gamma'} '
\left(
\overline{ h}
\right)
\right] 
\left(
\Gamma_{ k-1} \left(
z_{ (k-1)}
\right)
\right),
\ \ \ \ \
\text{\sf si} \ k \
\text{\sf est impair}~; \\
\left[
J_t^\ell \Theta_{ \gamma'}' \left( h
\right)\right] 
\left(
\Gamma_k \left(
z_{ (k)}
\right)
\right) 
& \
\equiv 
\left[
J_t^\ell \Theta_{ \gamma'} '
\left( h
\right)
\right]
\left(
\Gamma_{ k-1} \left(
z_{ (k-1)}
\right)
\right),
\ \ \ \ \
\text{\sf si} \ k \
\text{\sf est pair}.
\endaligned\right.
\end{equation}
Par cons\'equent, il n'est pas restrictif de supposer $k$ impair
dans~\thetag{ 7.19} et $k$ pair dans~\thetag{ 7.20}.

Nous allons d'abord \'etablir la Proposition~7.18 dans le cas $k=1$,
$\ell =0$ (\S7.65), puis dans le cas $k=1$, 
$\ell$ quelconque (\S7.88) et puis dans le cas $k=2$, $\ell=0$
(\S7.122). Ensuite, dans le \S7.135 {\it infra}, nous \'ecrirons la
d\'emonstration finale, pour $k$ et $\ell$ quelconques. Pour
commencer, en guise de pr\'eliminaire, formulons deux \'enonc\'es
cruciaux~: un principe d'unicit\'e reli\'e \`a la CR-transversalit\'e
de $h$ (\S7.22 juste ci-dessous) et un principe de sym\'etrie entre
les identit\'es de r\'eflexion conjugu\'ees ({\it cf.} le~\S1.25 {\it
supra} et le~ Lemme~7.57 {\it infra}).

\subsection*{ 7.22.~CR-transversalit\'e}
Rappelons que l'application CR formelle $h$ est suppos\'ee
CR-transversale \`a l'origine. D'apr\`es~\thetag{ 1.22}, cela signifie
que si une s\'erie formelle $\overline{ F}' (\zeta') \in \C \dl \zeta'
\dr$ satisfait l'identit\'e
\def\theequation{7.23}\begin{equation}
\overline{ F} ' \left(
\overline{ f}_1 \left(
\zeta,\, \Theta (\zeta,\, 0)
\right),\, \dots,\, 
\overline{ f}_{ m'} \left(
\zeta,\, \Theta (\zeta,\, 0)
\right)
\right)\equiv 0, 
\end{equation}
dans $\C \dl \zeta \dr$, elle doit s'annuler identiquement. 
Analysons cette hypoth\`ese.

Pour tout $\gamma ' \in \N^{ m'}$, d\'eveloppons la puissance
$\overline{ f}^{ \gamma'}$, restreinte \`a la premi\`ere cha\^{\i}ne
de Segre conjugu\'ee $\underline{ \Gamma}_1 (\zeta)= \left( 0,\, 0,\,
\zeta,\, \Theta (\zeta,\, 0) \right)$, par rapport aux puissances de
$\zeta$, ce qui nous donne~:
\def\theequation{7.24}\begin{equation}
\overline{ f}^{ \gamma'}
(\zeta,\, \Theta (\zeta,\, 0)) =: 
\sum_{ \beta \in \N^m}\, 
\overline{ f}_{ \gamma',\, \beta} \, 
\frac{ \zeta^\beta}{\beta !}~, 
\end{equation} 
avec certains coefficients $\overline{ f}_{ \gamma',\, \beta} \in
\C$. Puisque $\overline{ f} (0)= 0$, la puissance $\overline{ f}^{
\gamma'}$ s'annule \`a l'ordre $\vert \gamma ' \vert$ \`a l'origine~;
nous en d\'eduisons~:
\def\theequation{7.25}\begin{equation}
\overline{ f}_{ \gamma',\, \beta} = 0, 
\ \ \ \ \
\text{ \rm pour tout} \
\gamma' \
\text{ \rm tel que} \
\vert \gamma ' \vert > \vert \beta \vert.
\end{equation}

D\'eveloppons aussi la s\'erie $\overline{ F}' (\zeta') =: \sum_{
\gamma' \in \N^{ m'}} \, \overline{ F}_{ \gamma'} '\cdot (\zeta')^{
\gamma'}$, o\`u $\overline{ F}_{ \gamma'}' \in \C$. Rempla\c cons
$(\zeta')^{ \gamma'}$ par~\thetag{ 7.24} dans le d\'eveloppement de
$\overline{ F}' (\zeta')$, ce qui donne la double somme~:
\def\theequation{7.26}\begin{equation}
\sum_{ \beta \in \N^m} \, 
\frac{ \zeta^\beta}{ \beta!} \left(
\sum_{ \gamma' \in \N^{ m'}} \, 
\overline{ F}_{ \gamma'} ' \, 
\overline{ f}_{ \gamma',\, \beta}
\right).
\end{equation}
Gr\^ace \`a la propri\'et\'e~\thetag{ 7.25}, pour tout $\beta \in
\N^m$, la somme entre parenth\`eses dans~\thetag{ 7.26} est finie~:
ces coefficients existent et la double somme a un sens. Nous en
d\'eduisons que l'hypoth\`ese de CR-transversalit\'e se traduit par
l'implication concr\`ete~:
\def\theequation{7.27}\begin{equation}
\left\{
\aligned
{}
& \
\left(
\forall \, \beta \in \N^m \ :
\ \ \ \ 
\sum_{ \gamma' \in \N^{ m'}}\, 
\overline{ F}_{ \gamma'} ' \, 
\overline{ f}_{ \gamma',\, \beta} ' = 0
\right) \\
& \ 
\ \ \ \ \ \ \ \ \ \ \ \ \ 
\ \ \ \ \ \ \ \ \ \ \ \ \ 
\Downarrow \\
& \
\ \ \ \ \ \ \ \ \ \ 
\left( 
\overline{ F}_{ \gamma'}' = 0, 
\ \ \ \ \
\forall \ \gamma ' \in \N^{ m'}
\right).
\endaligned\right.
\end{equation}
Cette implication peut s'interpr\'eter en disant que le syst\`eme
lin\'eaire homog\`ene infini \'ecrit \`a la premi\`ere ligne
de~\thetag{ 7.27} est de noyau nul.

\subsection*{ 7.28.~Notation pour les d\'erivations CR}
Soit $\underline{ \mathcal{ L}}^\beta$ la d\'erivation CR d\'efinie
par~\thetag{ 1.30}. Rappelons que diff\'erentier~\thetag{ 7.24} par
rapport \`a $\zeta_k$ revient \`a appliquer la d\'erivation
$\underline{ \mathcal{ L}}_k$, prise en $(\zeta,\, t)= (\zeta,\,
0)$. En posant $\zeta =0$, nous en d\'eduisons~:
\def\theequation{7.29}\begin{equation}
\underline{ \mathcal{ L}}^\beta \overline{ f}^{ \gamma'} (0)
= \overline{ f}_{ \gamma',\, \beta}.
\end{equation}
Avant de poursuivre, effectuons une remarque importante.
L'application formelle $\overline{ f}(\tau)$ ne d\'epend que de
$\tau$. Mais lorsqu'on lui applique la d\'erivation $\underline{
\mathcal{ L}}^\beta$, dont les coefficients d\'ependent de $t$ et de
$\tau$ ({\it voir}~\thetag{ 3.12}), on obtient une expression qui
d\'epend de $t$ et de $\tau$, \`a moins que $\beta =0$. Aussi, pour
\^etre pr\'ecis et seulement lorsque cela sera n\'ecessaire, nous
n'\'ecrirons pas $\underline{ \mathcal{ L}}^\beta \overline{
f}^{\gamma'} (\tau)$, mais $\underline{ \mathcal{ L}}^\beta \overline{
f}^{ \gamma'} (t,\, \tau)$, ou, de mani\`ere \'equivalente
$\underline{ \mathcal{ L}}^\beta \overline{ f}^{ \gamma'} (z,\, w,\,
\zeta,\, \xi)$. Voici l'information utile concernant le
d\'eveloppement de ces expressions. 

\def\thelemma{7.30}\begin{lemma}
Pour tout $\beta \in \N^m$, il existe un polyn\^ome universel
$\overline{ Q}_\beta$, lin\'eaire par rapport \`a son second groupe de
variables, tel que l'identit\'e formelle suivante est satisfaite~{\rm
:}
\def\theequation{7.31}\begin{equation}
\underline{ \mathcal{ L}}^\beta
\overline{ f}^{ \gamma'} (t,\, \tau)
\equiv
\overline{ Q}_\beta
\left(
J_\zeta^{\vert \beta \vert} \Theta
(\zeta,\, t), \
J_\tau^{\vert \beta \vert}
\overline{ f}^{\gamma'}
(\tau)
\right),
\end{equation}
dans $\C\dl t,\, \tau \dr$, pour tout $\gamma ' \in \N^{ m'}$.
\end{lemma}

On \'etablit ce lemme gr\^ace \`a une r\'ecurrence \'el\'ementaire sur
le multiindice $\beta \in\N^m$.

\subsection*{ 7.32.~Principe d'unicit\'e}
Restreignons~\thetag{ 7.31} \`a la premi\`ere cha\^{\i}ne de Segre
$\Gamma_1 (z_1) = \left( z_1,\, \overline{ \Theta} (z_1,\, 0),\, 0,\,
0 \right)$, o\`u $z_1 \in \C^m$, ce qui donne~:
\def\theequation{7.33}\begin{equation}
\underline{ \mathcal{ L}}^\beta
\overline{ f}^{ \gamma'}
\left(
z_1,\, \overline{ \Theta} (z_1,\, 0),\, 0,\, 0
\right)
\equiv
\overline{ Q}_\beta 
\left(
J_\zeta^{\vert \beta \vert} \Theta
\left(0,\, z_1,\, \overline{ \Theta} (z_1,\, 0)\right), \
J_\tau^{\vert \beta \vert} \overline{ f}^{ \gamma'}
(0)
\right). 
\end{equation}
Puisque le polyn\^ome $\overline{ Q }_\beta$ est lin\'eaire par
rapport \`a son second groupe de variables, on a~{\rm :} $\overline{
Q}_\beta \left( J_\zeta^{ \vert \beta \vert} \Theta, \ 0 \right) =
0$. En outre, puisque $\overline{ f}(0) =0$, on a~: $J_\tau^{\vert
\beta \vert} \overline{ f}^{ \gamma'} (0)= 0$ pour tout $\gamma'$ tel
que $\vert \gamma' \vert > \vert \beta \vert$. Il d\'ecoule de ces
deux observations et de~\thetag{ 7.33} que~:
\def\theequation{7.34}\begin{equation}
\underline{ \mathcal{ L}}^\beta
\overline{ f}^{ \gamma'}
\left(
z_1,\, \overline{ \Theta} (z_1,\, 0),\, 0,\, 0
\right)
\equiv 0,
\ \ \ \ \
\text{ \rm pour tout} \
\gamma' \
\text{ \rm tel que} \
\vert \gamma ' \vert > \vert \beta \vert.
\end{equation} 
Pour tout $\gamma' \in \N^{ m'}$, supposons donn\'ee une s\'erie
formelle arbitraire $\overline{ F}_{ \gamma '} (z_1) \in \C\dl z_1
\dr$. Gr\^ace \`a~\thetag{ 7.34}, pour tout $\beta \in \N^m$, la somme
infinie de produits de s\'eries formelles d\'efinie par~:
\def\theequation{7.35}\begin{equation}
\sum_{ \gamma' \in \N^{ m'}}\, 
\underline{ \mathcal{ L}}^\beta
\overline{ f}^{ \gamma'}
\left(
z_1,\, \overline{ \Theta} (z_1,\, 0),\, 0,\, 0
\right) \cdot
\overline{ F}_{ \gamma'} (z_1),
\end{equation}
est en fait une somme finie~; par cons\'equent, elle poss\`ede un sens
en tant que s\'erie formelle. Nous pouvons maintenant formuler le
principe d'unicit\'e qui nous sera pr\'ecieux pour la d\'emonstration
de la Proposition~7.18. Rappelons que $h$ est CR-transversale \`a
l'origine, ce qui se traduit par l'implication~\thetag{ 7.27}.

\def\thelemma{7.36}\begin{lemma}
L'implication plus g\'en\'erale suivante est satisfaite~{\rm :}
\def\theequation{7.37}\begin{equation}
\left\{
\aligned
{}
& \
\left(
\forall \, \beta \in \N^m \ :
\ \ \ \ 
0 \equiv
\sum_{ \gamma' \in \N^{ m'}}\,
\underline{ \mathcal{ L}}^\beta
\overline{ f}^{ \gamma'}
\left(
z_1,\, \overline{ \Theta} (z_1,\, 0),\, 0,\, 0
\right) \cdot
\overline{ F}_{ \gamma'} (z_1)
\right) \\
& \ 
\ \ \ \ \ \ \ \ \ \ \ \ \ 
\ \ \ \ \ \ \ \ \ \ \ \ \ 
\Downarrow \\
& \
\ \ \ \ \ \ \ \ \ \ 
\left( 
\overline{ F}_{ \gamma'}'(z_1) \equiv 0, 
\ \ \ \ \
\forall \ \gamma ' \in \N^{ m'}
\right).
\endaligned\right.
\end{equation}
\end{lemma}

\proof
Les arguments sont \'el\'ementaires. En revenant \`a~\thetag{ 7.29},
observons que $\underline{ \mathcal{ L}}^\beta \overline{ f}^{
\gamma'} \left( 0,\, 0,\, 0,\, 0 \right) = \overline{ f}_{ \gamma',\,
\beta}$. Par cons\'equent, on peut \'ecrire
\def\theequation{7.38}\begin{equation}
\underline{ \mathcal{ L}}^\beta
\overline{ f}^{ \gamma'}
\left(
z_1,\, \overline{ \Theta} (z_1,\, 0),\, 0,\, 0
\right) = 
\overline{ f}_{ \gamma',\, \beta} + 
\sum_{ \gamma \in \N^m, \, \vert \gamma \vert \geq 1} \, 
\overline{ f}_{ \gamma',\, \beta,\, \gamma} \, z_1^\gamma,
\end{equation}
avec certains coefficients $\overline{ f}_{ \gamma',\, \beta,\,
\gamma } \in \C$. De m\^eme, d\'eveloppons
chaque s\'erie formelle 
\def\theequation{7.39}\begin{equation}
\overline{ F}_{ \gamma'} (z_1) =: 
\sum_{ \gamma \in \N^m} \, 
\overline{ F}_{ \gamma',\, \gamma} \, z_1^\gamma,
\end{equation}
o\`u $\overline{ F}_{ \gamma',\, \gamma} \in \C$, et r\'e\'ecrivons la
premi\`ere ligne de~\thetag{ 7.37}~:
\def\theequation{7.40}\begin{equation}
\forall \, \beta \in \N^m \ : 
\ \ \ \ \
0 \equiv
\sum_{ \gamma' \in \N^{ m'}}
\left(
\overline{ f}_{ \gamma',\, \beta} + 
\sum_{ \gamma \in \N^m, \, \vert \gamma \vert \geq 1} \, 
\overline{ f}_{ \gamma',\, \beta,\, \gamma} \, z_1^\gamma
\right) \cdot 
\left(
\sum_{ \gamma \in \N^m} \, 
\overline{ F}_{ \gamma',\, \gamma} \, z_1^\gamma
\right).
\end{equation}
En posant tout d'abord $z_1=0$, on obtient les \'equations
lin\'eaires~:
\def\theequation{7.41}\begin{equation}
\forall \, \beta \in \N^m \ : 
\ \ \ \ \
0 =
\sum_{ \gamma' \in \N^{ m'}}
\overline{ f}_{ \gamma',\, \beta} \cdot
\overline{ F}_{ \gamma',\, 0},
\end{equation}
qui impliquent imm\'ediatement que $\overline{ F}_{ \gamma',\, 0}= 0$,
pour tout $\gamma ' \in \N^{ m'}$, gr\^ace \`a l'implication~\thetag{
7.27}. Soit $\ell \in \N$. En raisonnant par r\'ecurrence sur la
longueur des multiindices $\gamma$, supposons que l'on a~: $\overline{
F}_{ \gamma',\, \gamma_2 } =0$, pour tout multiindice $\gamma_2\in
\N^m$ tel que $\vert \gamma_2 \vert \leq \ell$ et pour tout multindice
$\gamma ' \in \N^{ m'}$. Dans~\thetag{ 7.40}, la somme du second terme
du produit commence alors pour $\vert \gamma \vert \geq \ell$. Soit
$\gamma_1 \in \N^m$ un multiindice arbitraire de longueur $\ell
+1$. En annulant le coefficient de $z_1^{ \gamma_1}$ dans~\thetag{
7.40}, on trouve les \'equations lin\'eaires
\def\theequation{7.42}\begin{equation}
\forall \, \beta \in \N^m \ : 
\ \ \ \ \
0 =
\sum_{ \gamma' \in \N^{ m'}}
\overline{ f}_{ \gamma',\, \beta} \cdot
\overline{ F}_{ \gamma',\, \gamma_1},
\end{equation}
qui impliquent imm\'ediatement que $\overline{ F}_{ \gamma',\,
\gamma_1}= 0$, pour tout $\gamma ' \in \N^{ m'}$, en utilisant de
nouveau~\thetag{ 7.27}. Ceci compl\`ete la d\'emonstration.
\endproof

\subsection*{ 7.43.~Combinatoire d'identit\'es conjugu\'ees}
Avant d'amorcer la d\'emonstration de la Proposition~7.18, notre
objectif sera d'analyser le lien qui existe entre les quatre
identit\'es de r\'eflexion~\thetag{ 1.31} et~\thetag{ 1.32}. \`A
partir de maintenant, nous devrons en effet utiliser de mani\`ere
exhaustive les relations qui existent entre les objets
analytico-g\'eom\'etriques introduits jusqu'\`a pr\'esent, et leurs
conjugu\'es.

Comme dans le \S1.15, posons $r(t,\, \tau) := \xi - \Theta (\zeta,\,
t)$ et $\overline{ r} (\tau,\, t) := w - \overline{ \Theta} ( z,\,
\tau)$. De m\^eme, posons $r' \left(t',\, \tau' \right) := \xi' -
\Theta ' \left( \zeta',\, t' \right)$ et $\overline{ r} ' \left(
\tau',\, t' \right) := w' - \overline{ \Theta} ' \left( z',\, \tau'
\right)$. Rappelons ({\it cf.}~\thetag{ 3.7}) qu'il existe deux
matrices inversibles $a(t,\, \tau)$ et $a' \left( t',\, \tau' \right)$
de s\'eries formelles convergentes, de taille $d\times d$ et $d'
\times d'$, telles qu'on a les quatre identit\'es formelles
vectorielles suivantes, conjugu\'ees par paires verticales~:
\def\theequation{7.44}\begin{equation}
\left\{
\aligned
\overline{ r} (\tau,\, t) 
& \
\equiv
a (t,\, \tau) \, r( t,\, \tau), 
\ \ \ \ \ \ \ \ \ \ \ \ \ \ \ \ 
\overline{ r}' \left( \tau',\, t' \right) 
&
\equiv
a' \left( t',\, \tau' \right) \, r' \left(t',\, \tau' \right), \\
r(t,\, \tau) 
& \
\equiv
\overline{ a} (\tau,\, t) \, \overline{ r} (\tau,\, t), 
\ \ \ \ \ \ \ \ \ \ \ \ \ \ \ \ 
r'(t',\, \tau') 
&
\equiv
\overline{ a}'\left(\tau',\, t' \right) \, 
\overline{ r}' \left(\tau',\, t'\right),
\endaligned\right.
\end{equation}
dans $\C\dl t,\, \tau \dr^{ d}$ et dans $\C\dl t',\, \tau ' \dr^{
d'}$. Par hypoth\`ese ({\it cf.} le~\S1.15), 
il existe une matrice $b(t,\, \tau)$ de taille
$d'\times d$ de s\'eries formelles telle que $r' \left( h(t),\,
\overline{ h} (\tau) \right) \equiv b(t,\, \tau) \, r(t,\, \tau)$.
Respectons la combinatoire sous-jacente~: ce n'est pas une, mais
quatre identit\'es fondamentales, toutes \'equivalentes entre elles, 
qu'il nous faut \'ecrire~:
\def\theequation{7.45}\begin{equation}
\left\{
\aligned
r' \left( h(t),\,
\overline{ h} (\tau) \right) 
& \
\equiv b(t,\, \tau) \, r(t,\, \tau), 
\ \ \ \ \ \ \ \ \ \ \ 
\overline{ r}' \left( \overline{ h}(\tau),\,
h (t) \right) 
&
\equiv c(t,\, \tau) \, r (t,\, \tau), 
\\
\overline{ r}' \left(
\overline{ h} (\tau),\, h(t) \right)
& \
\equiv
\overline{ b} (\tau,\, t) \, 
\overline{ r} (\tau,\, t),
\ \ \ \ \ \ \ \ \ \ \ 
r' \left( h(t),\,
\overline{ h} (\tau) \right) 
&
\equiv \overline{ c}(\tau,\, t) \, \overline{ r} (\tau,\, t).
\endaligned\right.
\end{equation}
Ici, nous avons pos\'e $c(t,\, \tau) := \overline{ b}( \tau,\, t) \, a
(t,\, \tau)$. Pour $(t,\, \tau) \in \mathcal{ M}$, on a $r(t,\, \tau)
=0$ et $\overline{ r}(\tau,\, t) =0$. En partant des quatre
identit\'es~\thetag{ 7.45}, on pourrait penser qu'il ne reste alors
que deux \'egalit\'es
\def\theequation{7.46}\begin{equation}
r' \left(
h(t),\, \overline{ h}(\tau)
\right) = 0 
\ \ \ \ \ \ \
{\rm et}
\ \ \ \ \ \ \
\overline{ r}' \left(
\overline{ h}(\tau),\, h (t),
\right) =0.
\end{equation}
lorsque $(t,\, \tau) \in \mathcal{ M}$, mais on se m\'eprendrait alors
sur la v\'eritable nature des choses~: en effet, dire que $(t,\, \tau)
\in \mathcal{ M}$ revient \`a dire que $w= \overline{ \Theta} (z,\,
\tau)$ ou que $\xi = \Theta (\zeta,\, t)$, et par cons\'equent, en
rempla\c cant $w$ et $\xi$ par ces deux valeurs possibles dans les
deux \'egalit\'es~\thetag{ 7.46}, on obtient quatre identit\'es
formelles vectorielles~:
\def\theequation{7.47}\begin{equation}
\left\{
\aligned
r' \left( h(t),\,
\overline{ h} (\zeta,\, \Theta (\zeta,\, t)) \right) 
& \
\equiv 0, 
\ \ \ \ \ \ \ \ \ \ \
\overline{ r}' \left( 
\overline{ h}(\zeta,\, \Theta (\zeta,\, t)),\,
h (t) \right) 
&
\equiv 0,
\\
\overline{ r}' \left(
\overline{ h} (\tau),\, h\left(
z,\, \overline{ \Theta} (z,\, \tau) \right) \right)
& \
\equiv
0,
\ \ \ \ \ \ \ \ \ \ \ 
r' \left( h\left(z,\, 
\overline{ \Theta}(z,\, \tau)\right),\,
\overline{ h} (\tau) \right) 
&
\equiv 0, 
\endaligned\right.
\end{equation}
qui ne sont autres, dans le m\^eme ordre, que les quatre
relations~\thetag{ 7.45} restreintes \`a $\mathcal{ M}$. \`A nouveau,
elles sont conjugu\'ees par paires verticales.

Venons-en maintenant aux identit\'es de r\'eflexion. Rappelons que
les champs de vecteurs $\mathcal{ L}_k$ et $\underline{ \mathcal{
L}}_k$ d\'efinis par~\thetag{ 3.12} satisfont $\mathcal{ L}_k \,
\overline{ r} (\tau,\, t) \equiv 0$ et $\underline{ \mathcal{ L} }_k
\, r (t,\, \tau) \equiv 0$. En appliquant ces d\'erivations \`a la
premi\`ere colonne de~\thetag{ 7.44}, nous en d\'eduisons que
$\underline{ \mathcal{ L} }_k \, \overline{ r} (\tau,\, t) \equiv
\underline{ \mathcal{ L} }_k \, a(t,\, \tau) \cdot r(t,\, \tau)$ et
que $\mathcal{ L}_k \, r (t,\, \tau) \equiv \mathcal{ L}_k \,
\overline{ a} (\tau,\, t) \cdot \overline{ r} (\tau,\, t)$. Par
cons\'equent, nous avons $\underline{ \mathcal{ L}}_k \, \overline{ r}
(\tau,\, t) =0$ et $\mathcal{ L}_k \, r (t,\, \tau) =0$, pour
$(t,\, \tau) \in \mathcal{ M}$, mais pas d'identit\'e formelle dans
$\C \dl t,\, \tau\dr^d$. Plus g\'en\'eralement, pour tout
multindice $\beta \in
\N^m$ et pour $(t,\, \tau) \in\mathcal{ M}$, 
nous avons les quatre relations~:
\def\theequation{7.48}\begin{equation}
\left\{
\aligned
\mathcal{ L}^\beta \, r (t,\, \tau) 
& \
= 0, 
\ \ \ \ \ 
\ \ \ \ \ 
\ \ \ \ \ 
\mathcal{ L}^\beta \,
\overline{ r} (\tau,\, t) 
&
= \ 0, \\
\underline{ \mathcal{ L}}^\beta \,
\overline{ r } (\tau,\, t) 
& \
= 0, 
\ \ \ \ \ 
\ \ \ \ \ 
\ \ \ \ \ 
\underline{ \mathcal{ L}}^\beta \, 
r (t,\, \tau) 
&
= \ 0,
\endaligned\right.
\end{equation}
satisfaites pour $(t,\, \tau ) \in \mathcal{ M}$. Lorsque le
multiindice $\beta$ est de longueur strictement positive, les
relations de la deuxi\`eme colonne sont des identit\'es formelles
vectorielles dans $\C \dl t,\, \tau \dr^{ d}$. Cette remarque ne sera
pas utilis\'ee dans la suite.

En appliquant les d\'erivations $\mathcal{ L}^\beta$ et $\underline{
\mathcal{ L}}^\beta$ aux deux relations~\thetag{ 7.46}, nous obtenons
les quatre identit\'es de r\'eflexion~\thetag{ 1.31} et~\thetag{ 1.32}
de l'Introduction~:
\def\theequation{7.49}\begin{equation}
\left\{
\aligned
0 
& \ 
=
\underline{ \mathcal{ L}}^\beta \,
r' \left(
h(t),\, \overline{ h}(\tau)
\right), 
\ \ \ \ \ \ \ \ \ \ \ \ \
0 
&
= 
\underline{ \mathcal{ L}}^\beta \, 
\overline{ r}' \left(
\overline{ h} (\tau),\, h(t)
\right), \\
0 
& \
=
\mathcal{ L}^\beta \, 
\overline{ r}' \left(
\overline{ h}(\tau),\, 
h(t)
\right), 
\ \ \ \ \ \ \ \ \ \ \ \ \
0 
&
= 
\mathcal{ L}^\beta \,
r' \left(
h(t),\, \overline{ h}(\tau)
\right),
\endaligned\right.
\end{equation}
o\`u $\beta$ parcourt $\N^m$, avec $(t,\, \tau) \in \mathcal{ M}$,
comme \`a l'accoutum\'ee. La premi\`ere ligne de~\thetag{ 7.49}
co\"{\i}ncide avec~\thetag{ 1.31}~; la deuxi\`eme ligne, avec~\thetag{
1.32}.

Pour obtenir ces quatre identit\'es de r\'eflexion, on peut proc\'eder
d'une mani\`ere tr\`es l\'eg\`erement diff\'erente en partant des
quatre identit\'es formelles~\thetag{ 7.47}. En effet, rappelons que
le fait de diff\'erentier une s\'erie $\overline{ \psi} \left( z,\,
\overline{ \Theta} (z,\, \tau) \right)$ par rapport \`a $z^\beta$
revient \`a appliquer la d\'erivation $\mathcal{ L}^\beta$ \`a
$\overline{ \psi}$ ({\it cf.} aussi~\thetag{ 4.12}). Ainsi, en
diff\'erentiant les quatre identit\'es formelles~\thetag{ 7.47} par
rapport \`a $z^\beta$ ou par rapport \`a $\zeta^\beta$, nous obtenons
quatre familles infinies d'identit\'es formelles~:
\def\theequation{7.50}\begin{equation}
\left\{
\aligned
0 
& \ 
\equiv
\underline{ \mathcal{ L}}^\beta \,
r' \left(
h(t),\, \overline{ h}(\zeta,\, \Theta (\zeta,\, t))
\right), 
\ \ \ \ \ \ \ \ \ \ \ \ \ \,
0 
\equiv
\underline{ \mathcal{ L}}^\beta \, 
\overline{ r}' \left(
\overline{ h} (\zeta,\, \Theta (\zeta,\, t)),\, h(t)
\right), \\
0 
& \
\equiv
\mathcal{ L}^\beta \, 
\overline{ r}' \left(
\overline{ h}(\tau),\, 
h\left( z,\, \overline{ \Theta} (z,\, 0) \right)
\right), 
\ \ \ \ \ \ \ \ \ \ \
0 
\equiv
\mathcal{ L}^\beta \,
r' \left(
h\left( z,\, \overline{ \Theta} (z,\, 0) 
\right),\, \overline{ h}(\tau)
\right),
\endaligned\right.
\end{equation}
satisfaites pour tout $\beta \in \N^m$.
Elles co\"{\i}ncident avec~\thetag{ 7.49}.

\subsection*{ 7.51.~\'Equivalence entre les identit\'es de r\'eflexion}
Maintenant, v\'erifions que ces quatre familles infinies d'identit\'es
sont \'equivalentes entre elles. Les deux membres de chaque colonne
de~\thetag{ 7.49} (ou de~\thetag{ 7.50}) sont conjugu\'es entre eux,
donc \'equivalents. Pour d\'emontrer que les deux familles infinies
\'ecrites \`a la premi\`ere ligne de~\thetag{ 7.49} (ou de~\thetag{
7.50}) sont \'equivalentes, partons des deux
identit\'es formelles~:
\def\theequation{7.52}\begin{equation}
\left\{
\aligned
r'\left(
h(t),\, \overline{ h}(\tau)
\right) 
& \
\equiv
\overline{ a} ' \left(
\overline{ h} (\tau),\, h (t)
\right) \, 
\overline{ r}' \left(
\overline{ h} (\tau),\, h(t)
\right), \\
\overline{ r}' \left(
\overline{ h} (\tau),\, 
h (t)
\right) 
& \
\equiv
a' \left(
h (t),\, \overline{ h} (\tau)
\right) \, 
r' \left(
h(t),\, \overline{ h} (\tau)
\right),
\endaligned\right.
\end{equation}
obtenues en rempla\c cant les variables $t'$ et $\tau'$ par $h (t)$ et
par $\overline{ h} (\tau)$ dans la deuxi\`eme colonne de~\thetag{
7.44}. Si $\beta$ et $\gamma$ sont deux multiindices de $\N^m$, on
notera $\gamma \leq \beta$ si $\gamma_k \leq \beta_k$ pour
$k=1,\,\dots,\, m$. Appliquons aux deux identit\'es~\thetag{ 7.52} la
d\'erivation $\underline{ \mathcal{ L}}^\beta$ en tenant compte de la
formule de Leibniz, ce qui donne, en n'\'ecrivant pas les arguments
$(t,\, \tau)$ pour all\'eger les formules~:
\def\theequation{7.53}\begin{equation}
\left\{
\aligned
\underline{ \mathcal{ L}}^\beta \, 
r'\left( h,\, \overline{ h}
\right) 
& \
\equiv
\
\sum_{ \gamma \leq \beta}\, 
\frac{ \beta !}{\gamma ! \
(\beta - \gamma)!} \cdot
\underline{ \mathcal{ L}}^{ \beta - \gamma}
\left(
\overline{ a} ' \left(
\overline{ h},\, h
\right) \right)\cdot
\underline{ \mathcal{ L}}^\gamma
\left(
\overline{ r}' \left(
\overline{ h},\, h
\right) \right), \\
\underline{ \mathcal{ L}}^\beta \, 
\overline{ r}' \left(
\overline{ h},\, 
h
\right) 
& \
\equiv\
\sum_{ \gamma \leq \beta}\, 
\frac{ \beta !}{\gamma ! \
(\beta - \gamma)!} \cdot
\underline{ \mathcal{ L}}^{ \beta - \gamma}
\left(
a' \left(
h,\, \overline{ h}
\right)\right)\cdot \underline{ \mathcal{ L}}^\gamma
\left(
r' \left(
h,\, \overline{ h}
\right)\right).
\endaligned\right.
\end{equation}
Gr\^ace \`a ces relations, nous d\'eduisons imm\'ediatement
l'\'equivalence d\'esir\'ee~:
\def\theequation{7.54}\begin{equation}
\left\{
\aligned
{}
& \
\left(
\forall \, \beta \in \N^m \ :
\underline{ \mathcal{ L}}^\beta \,
r' \left( h (t),\, 
\overline{ h} (\tau)
\right)
=0 
\right)
\\
& \ 
\ \ \ \ \ \ \ \ \ \ \ \ \ 
\ \ \ \ \ \ \ \ \ \ \ \ \ 
\Updownarrow \\
& \
\left( 
\forall \, \beta \in \N^m \ :
\underline{ \mathcal{ L}}^\beta \, 
\overline{ r}' \left(
\overline{ h} (\tau),\, 
h (t)
\right) =0 
\right).
\endaligned\right.
\end{equation}
Ainsi, les quatre familles infinies d'identit\'es de
r\'eflexions~\thetag{ 7.49} sont \'equivalentes \`a une seule
d'entre elles.

\subsection*{ 7.55.~Sp\'eculation cruciale}
Une telle observation semble maintenant contredire notre affirmation
au sujet de la n\'ecessit\'e d'utiliser r\'eellement les deux familles
infinies d'identit\'es de r\'eflexion~\thetag{ 1.31} (ou leurs
conjugu\'ees~\thetag{ 1.32}), comme nous l'avons avanc\'e dans le
\S1.25. En effet, gr\^ace \`a l'\'equivalence~\thetag{ 7.54}, on peut
fort bien soutenir qu'il est redondant de mentionner plus d'une
famille infinie d'identit\'es de r\'eflexion, les trois autres
n'ajoutant aucune relation nouvelle. C'est le point de vue qui est
adopt\'e dans toutes les r\'ef\'erences cit\'ees apr\`es
l'\'equation~\thetag{ 1.30}, bien que dans aucun de ces articles, on
ne puisse d\'eceler le moindre souci d'\'etablir de telles
\'equivalences. Nous admettrons bien s\^ur que le point de vue
classique est justifi\'e lorsque des hypoth\`eses de finitude telles
que {\bf (h1)}, {\bf (h2)}, {\bf (h3)} ou {\bf (h4)}, ou encore,
telles que l'alg\'ebricit\'e de $M'$, sont ajout\'ees au probl\`eme
pour en simplifier l'\'etude. Mais en revanche, lorsqu'aucune
hypoth\`ese particuli\`ere de finitude n'est suppos\'ee sur la
sous-vari\'et\'e g\'en\'erique image $M'$, le fait que l'on
diff\'erentie $\overline{ g}$ et les puissances $\overline{ f}^{
\gamma'}$ dans la premi\`ere ligne de~\thetag{ 1.31}, alors qu'on
diff\'erentie les composantes $\overline{ \Theta}_{ \gamma'} ' \left(
\overline{ h} \right)$ de l'application de r\'eflexion CR conjugu\'ee
dans la deuxi\`eme ligne de~\thetag{ 1.31}, constitue une diff\'erence
formelle que nous consid\'erons comme majeure, bien qu'elle soit
occult\'ee par l'\'equivalence~\thetag{ 7.54} et qu'on puisse, pour
cette raison, la juger comme infime et n\'egligeable. Mentionnons
qu'\`a ce jour, nous ne sommes pas parvenu \`a construire une seconde
d\'emonstration du Th\'eor\`eme~1.23 qui n'utilise pas un jeu
alternatif permanent entre chaque couple d'identit\'es de
r\'eflexions~\thetag{ 1.31} et~\thetag{ 1.32}\footnotemark[7]. Le
reste de cette section est consacr\'e \`a expliquer comment doit
fonctionner ce jeu alternatif, afin d'\'etablir la Proposition~7.18.

\footnotetext[7]{
Dans l'article~\cite{ bmr2002}, cette strat\'egie est emprunt\'ee
librement, sans r\'ef\'erence \`a nos travaux ant\'erieurs~\cite{
me2000}, \cite{ me2001b} et~\cite{ me2001c}. Cependant, dans ledit
article, l'approche originale que nous exposons ici est submerg\'ee
par une technicit\'e quasi-cryptique~; il en d\'ecoule un manque
certain de clart\'e conceptuelle qui rend cet article herm\'etique,
m\^eme du point de vue d'un sp\'ecialiste rompu au calcul formel.
}

\subsection*{ 7.56.~Principe d'\'equivalence polaris\'ee}
Commen\c cons par polariser l'\'equivalence~\thetag{ 7.54}. Soit $\nu
\in \N$, $\nu \geq 1$, soit ${\sf x} \in \C^\nu$ et soient ${\sf H}
({\sf x}) \in \C \dl {\sf x} \dr^n$ et ${\sf Q} ({\sf x})\in \C\dl
{\sf x} \dr^{ 2n}$ deux applications formelles sans terme constant.

\def\thelemma{7.57}\begin{lemma}
L'\'equivalence suivante est satisfaite~:
\def\theequation{7.58}\begin{equation}
\left\{
\aligned
{}
& \
\left( 
\forall \, \beta \in \N^m \ :
\left.
\underline{ \mathcal{ L}}^\beta 
\left[
\overline{ r}' \left(
\overline{ h} (\tau),\, 
{\sf H} ({\sf x})
\right)\right]
\right\vert_{ (t,\, \tau) = {\sf Q}({\sf x})} 
\equiv 0,
\ \ \
\text{\rm dans} \
\C\dl {\sf x} \dr^{ d'} 
\right)
\\
& \ 
\ \ \ \ \ \ \ \ \ \ \ \ \ 
\ \ \ \ \ \ \ \ \ \ \ \ \ 
\Updownarrow \\
& \
\left(
\forall \, \beta \in \N^m \ :
\left.
\underline{ \mathcal{ L}}^\beta
\left[
r' \left(
{\sf H} ({\sf x}),\, 
\overline{ h} (\tau)
\right) \right]
\right\vert_{ (t,\, \tau) = {\sf Q}({\sf x})}
\equiv 0,
\ \ \
\text{\rm dans} \
\C\dl {\sf x} \dr^{ d'} 
\right).
\endaligned\right.
\end{equation}
\end{lemma}

Puisque les coefficients des d\'erivations $\underline{ \mathcal{
L}}^\beta$ d\'ependent de $t$ et de $\tau$, et que l'on remplace
$(t,\, \tau)$ par ${\sf Q} ({\sf x})$, ces identit\'es formelles sont
\`a interpr\'eter dans $\C \dl {\sf x} \dr^{ d'}$. Une \'equivalence
analogue est v\'erifi\'ee, en rempla\c cant $\underline{ \mathcal{
L}}^\beta$ par $\mathcal{ L}^\beta$ et $\overline{ h}(\tau)$ par
$h(t)$.

\proof
En adaptant l'argument du \S7.51 ci-dessus, partons
des deux identit\'es formelles conjugu\'ees~:
\def\theequation{7.59}\begin{equation}
\left\{
\aligned
\overline{ r}' \left(
\overline{ h} (\tau),\, 
{\sf H} ( {\sf x})
\right) 
& \
\equiv
a' \left({\sf H} ( {\sf x}),\, 
\overline{ h} (\tau)
\right) \, 
r' \left({\sf H} ( {\sf x}),\, 
\overline{ h} (\tau)
\right), \\
r'\left(
{\sf H} ( {\sf x}),\, \overline{ h}(\tau)
\right) 
& \
\equiv
\overline{ a} ' \left(
\overline{ h} (\tau),\, {\sf H} ( {\sf x})
\right) \, 
\overline{ r}' \left(
\overline{ h} (\tau),\, {\sf H} ( {\sf x})
\right), 
\endaligned\right.
\end{equation}
valables dans $\C \dl \tau,\, {\sf x} \dr^{ d'}$, que
nous obtenons
en rempla\c cant les variables $t'$ et $\tau'$ par ${\sf H} (
{\sf x})$ et par $\overline{ h} (\tau)$ dans la deuxi\`eme colonne
de~\thetag{ 7.44}. Appliquons les d\'erivations $\underline{ \mathcal{
L}}^\beta$ et la formule de Leibniz, ce qui donne les identit\'es
suivantes, dans $\C \dl t,\, \tau,\, {\sf x}\dr^{ d'}$~:
\def\theequation{7.60}\begin{equation}
\left\{
\aligned 
\underline{ \mathcal{ L}}^\beta \, 
\overline{ r}' \left(
\overline{ h},\, 
{\sf H}
\right) 
& \
\equiv\
\sum_{ {\beta_1} \leq \beta}\, 
\frac{ \beta !}{{\beta_1} ! \
(\beta - {\beta_1})!} \cdot
\underline{ \mathcal{ L}}^{ \beta - {\beta_1}}
\left(
a' \left(
{\sf H},\, \overline{ h}
\right)\right)\cdot \underline{ \mathcal{ L}}^{\beta_1}
\left(
r' \left(
{\sf H},\, \overline{ h}
\right)\right), \\
\underline{ \mathcal{ L}}^\beta \, 
r'\left( {\sf H},\, \overline{ h}
\right) 
& \
\equiv
\
\sum_{ {\beta_1} \leq \beta}\, 
\frac{ \beta !}{{\beta_1} ! \
(\beta - {\beta_1})!} \cdot
\underline{ \mathcal{ L}}^{ \beta - {\beta_1}}
\left(
\overline{ a} ' \left(
\overline{ h},\, {\sf H}
\right) \right)\cdot
\underline{ \mathcal{ L}}^{\beta_1}
\left(
\overline{ r}' \left(
\overline{ h},\, {\sf H}
\right) \right).
\endaligned\right.
\end{equation}
Gr\^ace \`a ces relations, nous d\'eduisons 
imm\'ediatement l'\'equivalence
d\'esir\'ee.
\endproof 

\subsection*{ 7.61.~\'Ecriture d\'efinitive des quatre
familles d'identit\'es de 
r\'eflexion} Dor\'enavant, pour \'ecrire les identit\'es de
r\'eflexion~\thetag{ 1.13} et~\thetag{ 1.32}, nous adopterons le
principe de notation qui est expliqu\'e au \S7.28. Plus
pr\'ecis\'ement, nous les \'ecrirons en respectant un parenth\'esage
rigoureux, en utilisant le signe <<$\cdot$>> pour la multiplication,
et en faisant appara\^{\i}tre les deux arguments $(t,\, \tau)$ des
termes auxquels s'appliquent les diff\'erentiations $\underline{
\mathcal{ L}}^\beta$ et $\mathcal{ L}^\beta$, ce qui donne~:
\def\theequation{7.62}\begin{equation}
\left\{
\aligned
\left[
\underline{ \mathcal{ L}}^\beta \, 
\overline{ g}_{ j'}\right] (\tau,\, t)
& \
=
\sum_{ \gamma' \in \N^{m'}}\,
\left[
\underline{ \mathcal{ L}}^\beta 
\overline{ f}^{\gamma '}\right] 
(\tau,\, t) \cdot
\Theta_{j',\, \gamma' } ' (h(t)), \\
0
& \
=
\sum_{\gamma '\in \N^{ m'}} \, 
f(t)^{\gamma '} \cdot
\left[
\underline{ \mathcal{ L} }^\beta
\overline{ \Theta}_{ j',\, 
\gamma'} ' \left(\overline{ h} 
\right) \right] (\tau,\, t), \\
\left[
\mathcal{ L}^\beta 
g_{ j'} \right] (t,\, \tau) 
& \
=
\sum_{\gamma' \in \N^{ m'}} \,
\left[
\mathcal{ L}^\beta 
f^{\gamma '}\right] (t,\, \tau) \cdot 
\overline{ \Theta}_{ j',\, 
\gamma'} ' \left( \overline{ h} (\tau) \right), \\
0
& \
=
\sum_{ \gamma' \in \N^{m'}}\, 
\overline{ f} (\tau)^{\gamma '} \cdot
\left[
\mathcal{ L }^\beta
\Theta_{j',\, \gamma' } ' (h)\right] (t,\, \tau).
\endaligned\right.
\end{equation}
Ici~: $(t,\, \tau) \in \mathcal{ M}$~; $j'$ satisfait $1\leq j' \leq
d'$~; $\beta$ varie dans $\N^m$. Nous sous-entendons que pour $\beta
=0$, la deuxi\`eme ligne de~\thetag{ 7.62} co\"{\i}ncide avec la
relation non diff\'erenti\'ee $g_{ j'} (t) = \sum_{ \gamma' \in \N^{
m'}} \, f(t)^{ \gamma'} \, \overline{ \Theta}_{j',\, \gamma'} '
(h(t))$~; de m\^eme, nous sous-entendrons que la quatri\`eme ligne
de~\thetag{ 7.62} co\"{\i}ncide avec $\overline{ g}_{ j'} (\tau) =
\sum_{ \gamma ' \in \N^{ m'}} \, \overline{ f}^{ \gamma'} (\tau) \,
\Theta_{ j',\, \gamma'} '( h(t))$ pour $\beta =0$.

Le Lemme~7.30 fournit l'expression d\'evelopp\'ee de $\left[
\underline{ \mathcal{ L}}^\beta \overline{ f}^{ \gamma'} \right] (t,\,
\tau)$. En raisonnant par r\'ecurrence, on \'etablit aussi
l'\'enonc\'e suivant, utilse dans le \S7.65 ci-dessous.

\def\thelemma{7.63}\begin{lemma}
Pour tout $\beta \in \N^m$, il existe un polyn\^ome universel
$\overline{ P}_\beta'$, lin\'eaire par rapport \`a son second groupe
de variables, tel que l'identit\'e formelle suivante est
satisfaite~{\rm :}
\def\theequation{7.64}\begin{equation}
\left[
\underline{ \mathcal{ L}}^\beta
\overline{ \Theta}_{ j',\, \gamma'} '
\left(
\overline{ h}
\right)
\right] (t,\, \tau) 
\equiv 
\overline{ P}_{\beta} ' 
\left(
J_\zeta^{\vert \beta \vert} \Theta
(\zeta,\, t), \
J_\tau^{ \vert \beta \vert}
\left[
\overline{ \Theta}_{j',\, \gamma' } ' 
\left(\overline{ h} (\tau) \right) \right]
\right),
\end{equation}
dans $\C \dl t, \, \tau \dr^{ d'}$, pour
tout $j'=1,\, \dots,\, d'$ et tout
$\gamma ' \in \N^{ m'}$.
\end{lemma}

Pour pr\'evenir une confusion possible, notons que le second
argument de $\overline{ P}_\beta'$ dans \thetag{ 7.64} n'est pas le jet
$\left[ J_{ \tau'}^{ \vert \beta \vert} \overline{ \Theta}_{ j',\,
\gamma'}' \right] \left( \overline{ h}(\tau) \right)$, mais le jet
d'ordre $\vert \beta \vert$ de l'application formelle compos\'ee $\tau
\longmapsto_{ \mathcal{ F}} \, \overline{ \Theta}_{ j',\, \gamma'}
\left( \overline{ h} (\tau) \right)$.

\subsection*{ 7.65.~Convergence des composantes
l'application de r\'eflexion sur la premi\`ere cha\^{\i}ne de Segre}
Apr\`es ces pr\'eliminaires \'etendus, nous pouvons enfin d\'ebuter la
d\'emonstration de la Proposition~7.18. Nous allons d'abord
\'etablir~\thetag{ 7.19} dans le cas $k=1$, $\ell=0$.
Nous affirmons en effet que
pour tout $j'=1,\, \dots,\, d'$ et tout $\gamma ' \in \N^{ m'}$, on
a~:
\def\theequation{7.66}\begin{equation}
\Theta_{ j',\, \gamma' } ' \left( h \left( z_1,\, \overline{
\Theta} (z_1,\, 0) \right) \right) \in \C\{ z_1 \}.
\end{equation}

\proof
Pour cela, posons $(t,\, \tau) = \left( z_1,\, \overline{ \Theta}_1,\,
0,\, 0\right)$ dans la seconde ligne de~\thetag{ 7.62}, o\`u $\vert
\beta \vert \geq 1$, ainsi que dans l'identit\'e qui est sous-entendue
pour $\beta =0$, ce qui nous donne~:
\def\theequation{7.67}\begin{equation}
\left\{ 
\aligned
g_{ j'} \left(
z_1,\, \overline{ \Theta} (z_1,\, 0)
\right) 
&
\equiv
\sum_{ \gamma ' \in \N^{ m'}}\, 
\left[
f^{\gamma'}
\right]
\left(
z_1,\, \overline{ \Theta} (z_1,\, 0)
\right) \cdot
\overline{ \Theta}_{ j',\, \gamma '} '
\left(\overline{ h} (0) \right), \\
0
& 
\equiv
\sum_{ \gamma ' \in \N^{ m'}}\, 
\left[
f^{\gamma'}
\right]
\left(
z_1,\, \overline{ \Theta} (z_1,\, 0)
\right) \cdot \\
& \
\ \ \ \ \ \ \ 
\cdot
\overline{ P}_{\beta} ' 
\left(
J_\zeta^{\vert \beta \vert} 
\Theta \left(
0,\, z_1,\, 
\overline{ \Theta} (z_1,\, 0)
\right), \
J_\tau^{\vert \beta \vert} \, 
\overline{ \Theta}_{j',\, \gamma'}' \left(
\overline{ h}
\right)(0)
\right).
\endaligned\right.
\end{equation}
Ces \'equations sont analytiques par rapport \`a $z_1$, puisque le
second groupe d'arguments de $\overline{ P}_\beta'$ est constant.
Elles sont satisfaites par l'application formelle $z_1 \longmapsto_{
\mathcal{ F}} h \left( z_1,\, \overline{ \Theta} (z_1,\,
0)\right)$. En appliquant le Th\'eor\`eme~2.4 avec l'ordre
d'approximation $N=1$, nous d\'eduisons qu'il existe une application
convergente ${\sf H} (z_1) =: ({\sf F} (z_1),\, {\sf G} (z_1)) \in
\C\{ z_1 \}^{ m'} \times \C\{ z_1\}^{ d'}$, avec ${\sf H} (0)=0$,
qui satisfait les m\^emes \'equations analytiques, {\it i.e.}~:
\def\theequation{7.68}\begin{equation}
\left\{ 
\aligned
{\sf G}_{ j'} (z_1) 
&
\equiv
\sum_{ \gamma ' \in \N^{ m'}}\, 
{\sf F}^{\gamma'}(z_1) \cdot
\overline{ \Theta}_{ j',\, \gamma '} '
\left(\overline{ h} (0) \right), \\
0
& 
\equiv
\sum_{ \gamma ' \in \N^{ m'}}\, 
{\sf F}^{\gamma'}(z_1) \cdot
\overline{ P}_{ \beta} ' 
\left(
J_\zeta^{\vert \beta \vert} 
\Theta \left(
0,\, z_1,\, 
\overline{ \Theta} (z_1,\, 0)
\right), \
J_\tau^{\vert \beta \vert} \, 
\overline{ \Theta}_{j',\, \gamma'}' \left(
\overline{ h}
\right)(0)
\right).
\endaligned\right.
\end{equation} 
Notons que ces derni\`eres identit\'es~\thetag{ 7.68} co\"{\i}ncident
avec les identit\'es de la premi\`ere ligne de~\thetag{ 7.58}, si on
les \'ecrit avec ${\sf x} := z_1$, avec ${\sf Q}( {\sf x}) := \left(
z_1,\, \overline{ \Theta} (z_1,\, 0),\, 0,\, 0 \right)$ et avec ${\sf
H} ( {\sf x}) := {\sf H} (z_1)$. 

{\it Gr\^ace au principe d'\'equivalence~\thetag{ 7.58}, nous pouvons
donc faire basculer ces identit\'es vers la seconde ligne
de~\thetag{ 7.58}}, ce qui nous donne, pour
tout $j'= 1,\, \dots,\, d'$ et tout $\beta \in \N^{m'}$~:
\def\theequation{7.69}\begin{equation}
\left\{
\aligned
\left[
\underline{ \mathcal{ L}}^\beta \, 
\overline{ g}_{ j'} \right]
\left(
z_1,\, \overline{ \Theta} (z_1,\, 0),\, 0,\, 0
\right) 
& 
\equiv
\sum_{ \gamma ' \in \N^{ m'}} \, 
\left[
\underline{ \mathcal{ L}}^\beta 
\overline{ f}^{ \gamma'} 
\right]
\left(
z_1,\, \overline{ \Theta} (z_1,\, 0),\, 0,\, 0
\right) \cdot \\
& \
\ \ \ \ \ \ \ \ \ \ \
\ \ \ \ \ \ \ \ \ \ \
\ \ \ \ \ \ \ \ \ \ \ 
\ \ \ \
\cdot
\Theta_{ j',\, \gamma '} '( {\sf H} (z_1)).
\endaligned\right.
\end{equation}
Presque miraculeusement, nous voyons appara\^{\i}tre dans~\thetag{
7.69} les applications convergentes $\Theta_{ j',\, \gamma'} '(
{\sf H} (z_1)) \in \C\{ z_1\}$. Puisque ces termes constituent la seule
diff\'erence par rapport aux identit\'es de r\'eflexion \'ecrites \`a
la premi\`ere ligne de~\thetag{ 7.62} et prises en $(t,\, \tau) =
\left( z_1,\, \overline{ \Theta} (z_1,\, 0),\, 0,\, 0 \right)$, que
nous explicitons comme suit~:
\def\theequation{7.70}\begin{equation}
\left\{
\aligned
\left[
\underline{ \mathcal{ L}}^\beta \, 
\overline{ g}_{ j'} \right]
\left(
z_1,\, \overline{ \Theta} (z_1,\, 0),\, 0,\, 0
\right) 
& 
\equiv
\sum_{ \gamma ' \in \N^{ m'}} \, 
\left[
\underline{ \mathcal{ L}}^\beta 
\overline{ f}^{ \gamma'} 
\right]
\left(
z_1,\, \overline{ \Theta} (z_1,\, 0),\, 0,\, 0
\right) \cdot \\
& \
\ \ \ \ \ \ \ \ \ \ \
\ \ \ \ \ \ \ \ \ \ \
\ \ \ \ \ 
\cdot
\Theta_{ j',\, \gamma '} '\left( h\left(z_1,\, 
\overline{ \Theta} (z_1,\, 0) \right)\right),
\endaligned\right.
\end{equation}
soustrayons imm\'ediatement~\thetag{ 7.69} de~\thetag{ 7.70}~:
\def\theequation{7.71}\begin{equation}
\left\{
\aligned
0
& 
\equiv
\sum_{ \gamma ' \in \N^{ m'}} \, 
\left[
\underline{ \mathcal{ L}}^\beta 
\overline{ f}^{ \gamma'} 
\right]
\left(
z_1,\, \overline{ \Theta} (z_1,\, 0),\, 0,\, 0
\right) \cdot \\
& \
\ \ \ \ \ \ \ \ \ \ 
\ \ \ \ \ \ \ \ \ \ 
\cdot
\left\{
\Theta_{ j',\, \gamma '} '\left( h\left(z_1,\, 
\overline{ \Theta} (z_1,\, 0) \right)\right)
- \Theta_{ j',\, \gamma'} '( {\sf H} (z_1))
\right\}.
\endaligned\right.
\end{equation}
Pour terminer la d\'emonstration, il ne reste plus qu'\`a appliquer le
principe d'unicit\'e~\thetag{ 7.37} (pour chaque $j'$) \`a 
ces derni\`eres identit\'es, ce qui fournit les relations~:
\def\theequation{7.72}\begin{equation}
\Theta_{ j',\, \gamma'} '
\left( h
\left(
z_1,\, \overline{ \Theta} (z_1,\, 0)
\right)
\right) \equiv 
\Theta_{ j',\, \gamma'} '( {\sf H} (z_1)).
\end{equation}
Elles \'etablissent la propri\'et\'e de convergence
annonc\'ee~\thetag{ 7.66}, pour tout $j'= 1,\, \dots,\, d'$ et tout
$\gamma' \in \N^{ m'}$.
\endproof

Ainsi, en passant de~\thetag{ 7.68} \`a~\thetag{ 7.69}, nous avons
initi\'e un jeu alternatif entre deux familles d'identit\'es de
r\'eflexion, et ce jeu crucial se poursuivra ult\'erieurement.
Puisque les composantes $\Theta_{ j',\, \gamma'} ' ({\sf H} (z_1))$
n'\'etaient pas visibles dans~\thetag{ 7.68}, alors qu'elles sont
visibles dans~\thetag{ 7.69}, nous voyons bien qu'il \'etait
n\'ecessaire d'utiliser les deux familles infinies d'identit\'es de
r\'eflexion~\thetag{ 7.67} et~\thetag{ 7.70}, bien qu'elles soient
<<\'equivalentes>> d'apr\`es le \S7.51. Dans la suite de la
d\'emonstration de la Proposition~7.18, nous aurons besoin de
g\'en\'eraliser encore l'\'equivalence~\thetag{ 7.58} ({\it voir} le
Lemme~7.110 ci-dessous).

\subsection*{ 7.73.~Convergence des jets des composantes
de l'application de r\'eflexion sur la premi\`ere cha\^{\i}ne de Segre}
Cette propri\'et\'e de convergence constitue le point technique le
plus complexe de toute la d\'emonstration. Pour l'expliquer
soigneusement, nous ferons alterner l'\'etude sp\'ecifique des
calculs formels effectu\'es sur la premi\`ere cha\^{\i}ne de Segre avec
l'\'enonciation de propri\'et\'es formelles g\'en\'erales,
valables sur une cha\^{ \i}ne de Segre de longueur arbitraire. Ces
derni\`eres seront utilis\'ees dans la d\'emonstration finale de la
Proposition~7.18, o\`u nous raisonnerons par r\'ecurrence sur la
longueur des cha\^{\i}nes ({\it voir} le~\S7.135), ce qui
se fera sans mal, gr\^ace \`a tous nos pr\'eparatifs.
 
Commen\c cons par \'enoncer un crit\`ere pour la convergence des jets
de toutes les composantes de l'application de r\'eflexion $\mathcal{
R}_h'$ ou de sa conjugu\'ee $\overline{ \mathcal{ R}}_h'$ sur les
cha\^{\i}nes de Segre. La d\'emonstration, analogue \`a celle du
Lemme~6.9, ne sera pas d\'etaill\'ee.

\def\thelemma{7.74}\begin{lemma}
Soit $\ell \in\N$, soit $k \in \N$ et soit $\Gamma_k ( z_{ (k)})$ la
$k$-i\`eme cha\^{ \i}ne de Segre. Les propri\'et\'es suivantes sont
\'equivalentes~{\rm :}
\begin{itemize}
\item[{\bf (i)}]
Le jet d'ordre $\ell$ de chaque composante de l'application de
r\'eflexion $\mathcal{ R}_h'$ ou de sa conjugu\'ee $\overline{
\mathcal{ R}}_h'$ converge sur la $k$-i\`eme cha\^{\i}ne de Segre,
{\it i.e.}, en tenant compte des simplifications~\thetag{ 
7.21}, on a pour tout $\gamma ' \in \N^{ m'}$~{\rm :}
\def\theequation{7.75}\begin{equation}
\left\{
\aligned
\left[
J_t^\ell \Theta_{ \gamma'}' (h) \right] \left(
\Gamma_k ( z_{ (k)}) 
\right) 
& \
\in 
\C\{ z_{ (k)} \}^{ N_{ d',\, n,\, \ell}} \ \ \ \ \
\text{\sf si} \ k \
\text{\sf est impair},
\\
\left[
J_\tau^\ell 
\overline{ \Theta}_{ \gamma'}'\left(
\overline{ h} \right)
\right]
\left(
\Gamma_k ( z_{ (k)}) 
\right) 
& \
\in 
\C\{ z_{ (k)} \}^{ N_{ d',\, n,\, \ell}} \ \ \ \ \
\text{\sf si} \ k \
\text{\sf est pair}~;
\endaligned\right.
\end{equation}
\item[{\bf (ii)}]
pour tout $\beta \in \N^m$ et tout $\delta \in \N^d$ satisfaisant
$\vert \beta \vert + \vert \delta \vert \leq \ell$, les d\'eriv\'ees
suivantes de $\Theta_{ \gamma'} ' (h)$ ou de $\overline{ \Theta}_{
\gamma '} ' \left( \overline{ h} \right)$ convergent pour tout $\gamma
' \in \N^{m'}$~{\rm :}
\def\theequation{7.76}\begin{equation}
\left\{
\aligned
\left[
\mathcal{ L}^\beta
\Upsilon^\delta 
\Theta_{ \gamma '} ' (h)
\right] \left(
\Gamma_k ( z_{ (k)})
\right) 
& \
\in \C \{ z_{ (k)} \}^{ d'} \ \ \ \ \
\text{\sf si} \ k \
\text{\sf est impair}, \\
\left[
\underline{\mathcal{ L}}^\beta
\underline{ \Upsilon}^\delta 
\overline{ \Theta}_{ \gamma'} ' 
\left(
\overline{ h}
\right)
\right] \left(
\Gamma_k ( z_{ (k)})
\right) 
& \
\in \C \{ z_{ (k)} \}^{ d'} \ \ \ \ \
\text{\sf si} \ k \
\text{\sf est pair}~;
\endaligned\right.
\end{equation}
\item[{\bf (iii)}]
pour tout $\delta \in \N^d$ satisfaisant $\vert \delta \vert \leq
\ell$, les d\'eriv\'ees suivantes de $\Theta_{ \gamma'} ' (h)$ ou de
$\overline{ \Theta}_{ \gamma '} ' \left( \overline{ h} \right)$
convergent pour tout $\gamma ' \in \N^{m'}$~{\rm :}
\def\theequation{7.77}\begin{equation}
\left\{
\aligned
\left[
\Upsilon^\delta 
\Theta_{ \gamma '} ' (h)
\right] \left(
\Gamma_k ( z_{ (k)})
\right) 
& \
\in \C \{ z_{ (k)} \}^{ d'} \ \ \ \ \
\text{\sf si} \ k \
\text{\sf est impair}, \\
\left[
\underline{ \Upsilon}^\delta 
\overline{ \Theta}_{ \gamma'} ' 
\left(
\overline{ h}
\right)
\right] \left(
\Gamma_k ( z_{ (k)})
\right) 
& \
\in \C \{ z_{ (k)} \}^{ d'} \ \ \ \ \
\text{\sf si} \ k \
\text{\sf est pair}.
\endaligned\right.
\end{equation}
\end{itemize}
\end{lemma}

Ainsi, notre objectif est de g\'en\'eraliser la propri\'et\'e de
convergence~\thetag{ 7.66} aux jets des composantes de l'application
de r\'eflexion. Soit $w_1 \in \C^d$. En nous inspirant de la preuve
du Lemme~6.17, rempla\c cons $(t,\, \tau)$ dans la seconde ligne
de~\thetag{ 7.62} par la composition de flots $(z_1,\, w_1)
\longmapsto \Upsilon_{ w_1} \left( \Gamma_1 (z_1) \right)$, ce qui
donne~:
\def\theequation{7.78}\begin{equation}
\left\{
\aligned
g_{ j'} \left(
\Upsilon_{ w_1} (\Gamma_1 (z_1))
\right)
&
\equiv
\sum_{ \gamma ' \in \N^{ m'}} \, 
\left[
f^{ \gamma'} 
\right]
\left(
\Upsilon_{ w_1} (\Gamma_1 (z_1))
\right) \cdot
\overline{ \Theta}_{ j',\, \gamma'} 
\left( \overline{ h} 
\left(
\Upsilon_{ w_1} (\Gamma_1 (z_1))
\right)\right), \\
0 
& \
\equiv
\sum_{ \gamma' \in \N^{ m'}} \, 
\left[
f^{ \gamma'} 
\right]
\left(
\Upsilon_{ w_1} (\Gamma_1 (z_1))
\right)\cdot
\left[
\underline{ \mathcal{ L}}^\beta
\overline{ \Theta}_{ j',\, \gamma'} '
\left( \overline{ h} \right)
\right] \left(
\Upsilon_{ w_1} (\Gamma_1 (z_1))
\right).
\endaligned\right.
\end{equation}
La seconde ligne est valable pour tout $\beta \in \N^m$ satisfaisant
$\vert \beta \vert \geq 1$. Soit $\delta \in \N^d$. Appliquons la
diff\'erentiation $\left. \partial_{ w_1}^\delta \right\vert_{ w_1
=0}$ \`a ces identit\'es. Gr\^ace \`a la propri\'et\'e $\Upsilon_0 =
{\rm Id}$, gr\^ace aux relations~\thetag{ 6.20} et gr\^ace \`a la
formule de Leibniz, nous obtenons~:
\def\theequation{7.79}\begin{equation}
\left\{
\aligned
\left[
\Upsilon^\delta g_{ j'}\right] \left(
\Gamma_1 (z_1) \right)
& 
\equiv
\sum_{ \gamma' \in\N^{ m'}} \, 
\sum_{ \delta_1 \leq \delta} \, 
\frac{ \delta !}{ \delta_1 ! \ 
(\delta -\delta_1)!} \ 
\left[
\Upsilon^{ \delta - \delta_1} f^{ \gamma'}
\right](\Gamma_1 (z_1)) \cdot \\
& \
\ \ \ \ \ \ \ \ \ \ \ \ \
\ \ \ \ \ \ \ \ \ \ \ \ \ 
\ \ \ \ \ \ \ \
\cdot
\left[
\Upsilon^{ \delta_1}
\overline{ 
\Theta}_{ j',\, \gamma'}
\left(
\overline{ h}
\right)\right] \left(
\Gamma_1 (z_1)
\right), \\
0
& 
\equiv
\sum_{ \gamma' \in\N^{ m'}} \, 
\sum_{ \delta_1 \leq \delta} \, 
\frac{ \delta !}{ \delta_1 ! \ 
(\delta -\delta_1)!} \ 
\left[
\Upsilon^{ \delta - \delta_1} f^{ \gamma'}
\right](\Gamma_1 (z_1)) \cdot \\
& \
\ \ \ \ \ \ \ \ \ \ \ \ \
\ \ \ \ \ \ \ \ \ \ \ \ \ 
\ \ \ \ \ \ \ \
\cdot
\left[
\Upsilon^{ \delta_1}
\underline{ \mathcal{ L}}^\beta \overline{ 
\Theta}_{ j',\, \gamma'} 
\left(
\overline{ h}
\right)\right] \left(
\Gamma_1 (z_1)
\right).
\endaligned\right.
\end{equation}
Dans la deuxi\`eme identit\'e, on suppose $\vert \beta \vert \geq 1$.
Pour d\'evelopper les termes qui apparaissent \`a la deuxi\`eme et \`a
la quatri\`eme ligne de~\thetag{ 7.79}, voici une g\'en\'eralisation
du Lemme~7.63. Nous formulons aussi quatre autres d\'eveloppements, qui
seront utiles par la suite.

\def\thelemma{7.80}\begin{lemma}
Pour tout tout $\delta \in \N^d$, tout $\beta \in \N^m$ et tout
$\gamma' \in \N^{m'}$, il existe cinq polyn\^omes scalaires universels
$N_{ \gamma',\, \delta}$, $\overline{ Q}_{ \gamma',\, \beta,\,
\delta}$, $\overline{ P}_{ \beta,\, \delta}'$, $M_\delta'$ et
$\overline{ K}_{ \beta,\, \delta }$ tels que les quatre identit\'es
formelles suivantes sont satisfaites~{\rm :}
\def\theequation{7.81}\begin{equation}
\left\{
\aligned
\left[
\Upsilon^{\delta}
f^{ \gamma'}\right] (t)
&
\equiv
N_{\gamma',\, \delta}
\left(
\left(
\Upsilon^{\delta_1} f_{ k'} (t)
\right)_{ 1\leq k' \leq m',\, \delta_1 \leq \delta}
\right), 
\\
\left[
\Upsilon^\delta
\underline{ \mathcal{ L}}^\beta
\overline{ f}^{ \gamma'}\right] (t,\, \tau)
&
\equiv
\overline{ Q}_{ \gamma',\, \beta,\, \delta}
\left(
J_{\zeta,\, t}^{ \vert \beta \vert+ 
\vert \delta \vert} \Theta( \zeta, \, t), \
J_\tau^{\vert \beta \vert + 
\vert \delta \vert} \overline{ f}
(\tau)
\right), 
\\ 
\left[
\Upsilon^\delta
\underline{ \mathcal{ L}}^\beta
\overline{ \Theta}_{ j',\, \gamma'} '
\left(
\overline{ h}
\right)
\right] (t,\, \tau) 
&
\equiv 
\overline{ P}_{\beta,\, \delta} ' 
\left(
J_{\zeta,\, t}^{\vert \beta \vert+ \vert \delta \vert} \Theta
(\zeta,\, t), \
J_\tau^{ \vert \beta \vert + \vert \delta \vert}
\left[ 
\overline{ \Theta}_{j',\, \gamma' } ' 
\left(\overline{ h}(\tau) \right)\right]
\right), 
\\
\left[
\Upsilon^\delta
\Theta_{ j',\, \gamma'} '
\left(
h
\right)
\right] (t) 
&
\equiv 
M_\delta'
\left(
\left(
\Upsilon^{ \delta_1} h_{ i_1'} (t)
\right)_{ 1\leq i_1' \leq n',\, \delta_1 \leq \delta},
\
\left[
J_{t'}^{ \vert \delta \vert}
\Theta_{j',\, \gamma' } ' \right](h (t))
\right), 
\\
\left[
\Upsilon^\delta
\underline{ \mathcal{ L}}^\beta
\overline{ g}_{ j'}\right] (t,\, \tau)
&
\equiv
\overline{ K}_{ \beta,\, \delta}
\left(
J_{\zeta,\, t}^{ \vert \beta \vert+ 
\vert \delta \vert} \Theta( \zeta, \, t), \
J_\tau^{\vert \beta \vert + 
\vert \delta \vert} \overline{ g}_{ j'}
(\tau)
\right), 
\endaligned\right.
\end{equation}
dans $\C \dl t, \, \tau \dr$, pour tout $j'=1,\, \dots,\, d'$ et
tout $\gamma ' \in \N^{ m'}$.
\end{lemma}

Les d\'emonstrations, analogues \`a celle du Lemme~5.4, seront
omises. Cependant, formulons quelques commentaires explicatifs. Bien
que l'on ait $\Upsilon^\delta \phi (t) \equiv \partial_w^\delta \phi
(t)$, pour toute s\'erie formelle $\phi( t)$ de la seule variable $t$
({\it cf.}~\thetag{ 5.17}), nous utiliserons exclusivement la notation
$\Upsilon^\delta$ par souci d'uniformit\'e. Les coefficients des
champs de vecteurs $\Upsilon_j$ font intervenir des d\'eriv\'ees
partielles des s\'eries d\'efinissantes $\Theta_{ j_2} (\zeta,\, t)$
par rapport aux variables $w_{ j_1}$. Pour cette raison, nous
\'ecrivons dans la deuxi\`eme, dans la troisi\`eme et dans la
cinqui\`eme ligne de~\thetag{ 7.81} que le premier groupe d'arguments
de $\overline{ Q}_{\gamma', \beta,\, \delta}$, de $\overline{ P}_{
\beta,\, \delta}'$ et de $\overline{ K}_{ \beta,\, \delta}$ fait
intervenir le jet d'ordre $\vert \beta \vert + \vert \delta \vert$ de
l'application $\Theta (\zeta,\, t)$ par rapport \`a toutes les
variables $(\zeta,\, t)$. Remarquons que les polyn\^omes $\overline{
P}_{ \beta,\, \delta} '$ et $M_\delta'$ sont ind\'ependants de $j'$ et
de $\gamma'$. Enfin, notons que le deuxi\`eme groupe d'arguments de
$M_\delta'$ est le jet d'ordre $\vert \delta \vert$ par rapport \`a
$t'$ de l'application convergente $t' \longmapsto \Theta_{ j',\,
\gamma'} ' (t')$, dans lequel on remplace $t'$ par $h(t)$.

\smallskip

Revenons maintenant aux termes de la deuxi\`eme et de la quatri\`eme
ligne de~\thetag{ 7.79}. Plus g\'en\'eralement, soit $k \in \N$ et
consid\'erons la $k$-i\`eme cha\^{\i}ne de Segre $\Gamma_k (z_{
(k)})$. Gr\^ace au d\'eveloppement donn\'e par la troisi\`eme ligne
de~\thetag{ 7.81} et surtout, gr\^ace aux relations cruciales~\thetag{
7.21}, nous obtenons un abaissement d'un cran sur les cha\^{\i}nes de
Segre dans le deuxi\`eme groupe de variables de $\overline{
P}_{\beta,\, \delta_1}'$~: 
\def\theequation{7.82}\begin{equation}
\left\{
\aligned
{}
& 
\text{ \sf si} \ k \
\text{ \sf est impair}~: \ \ \ \ \ \
\left[
\Upsilon^{ \delta_1}
\underline{ \mathcal{ L}}^\beta
\overline{ \Theta}_{ j',\, \gamma'} '
\left(
\overline{ h}
\right)
\right] \left(
\Gamma_k (z_{ (k)})
\right)
\equiv 
\\
& \
\ \ \ \ \ \ \ \ \ 
\equiv
\overline{ P}_{ \beta,\, {\delta_1}} ' 
\left(
J_{\zeta,\, t}^{\vert \beta \vert+ \vert {\delta_1} \vert} \Theta
\left(\Gamma_k (z_{ (k)}) \right), \
J_\tau^{ \vert \beta \vert + \vert {\delta_1} \vert}
\left[ 
\overline{ \Theta}_{j',\, \gamma' } ' 
\left( 
\overline{ h} \right) \right]
\left(\Gamma_k (z_{ (k)})\right)
\right), 
\\
& \
\ \ \ \ \ \ \ \ \ 
\equiv
\overline{ P}_{ \beta,\, {\delta_1}} ' 
\left(
J_{\zeta,\, t}^{\vert \beta \vert+ \vert {\delta_1} \vert} \Theta
\left(\Gamma_k (z_{ (k)}) \right), \
J_\tau^{ \vert \beta \vert + \vert {\delta_1} \vert}
\left[
\overline{ \Theta}_{j',\, \gamma' } ' 
\left(\overline{ h}\right)\right]
\left(\Gamma_{k-1} (z_{ (k-1)})\right)
\right); \\
{}
& 
\text{ \sf si} \ k \
\text{ \sf est pair}~: \ \ \ \ \ \
\left[
\underline{ \Upsilon}^{ \delta_1}
\mathcal{ L}^\beta
\Theta_{ j',\, \gamma'} '
\left(
h
\right)
\right] \left(
\Gamma_k (z_{ (k)})
\right)
\equiv \\
& \
\ \ \ \ \ \ \ \ \ 
\equiv
P_{ \beta,\, {\delta_1}} ' 
\left(
J_{z,\, \tau}^{\vert \beta \vert+ \vert {\delta_1} \vert} 
\overline{ \Theta}
\left(\Gamma_k (z_{ (k)}) \right), \
J_t^{ \vert \beta \vert + \vert {\delta_1} \vert}
\left[ 
\Theta_{j',\, \gamma' } ' 
\left( 
h \right) \right]
\left(\Gamma_k (z_{ (k)})\right)
\right), \\
& \
\ \ \ \ \ \ \ \ \ 
\equiv
P_{ \beta,\, {\delta_1}} ' 
\left(
J_{z,\, \tau}^{\vert \beta \vert+ \vert {\delta_1} \vert} 
\overline{ \Theta}
\left(\Gamma_k (z_{ (k)}) \right), \
J_t^{ \vert \beta \vert + \vert {\delta_1} \vert}
\left[
\Theta_{j',\, \gamma' } ' 
\left(h\right)\right]
\left(\Gamma_{k-1} (z_{ (k-1)})\right)
\right).
\endaligned\right.
\end{equation}
En particulier, pour $k=1$~:
\def\theequation{7.83}\begin{equation}
\left\{
\aligned
{}
& 
\left[
\Upsilon^{ \delta_1}
\underline{ \mathcal{ L}}^\beta
\overline{ \Theta}_{ j',\, \gamma'} '
\left(
\overline{ h}
\right)
\right] \left(
\Gamma_1 (z_1)
\right)
\equiv \\
& \
\ \ \ \ \ \ \ \ \ 
\equiv
\overline{ P}_{ \beta,\, {\delta_1}} ' 
\left(
J_{\zeta,\, t}^{\vert \beta \vert+ \vert {\delta_1} \vert} \Theta
\left(\Gamma_1 (z_1) \right), \
J_\tau^{ \vert \beta \vert + \vert {\delta_1} \vert}
\left[
\overline{ \Theta}_{j',\, \gamma' } ' 
\left(\overline{ h}\right)\right]
\left(0)\right)
\right), 
\endaligned\right.
\end{equation}
et puisque le deuxi\`eme groupe d'arguments est {\it constant}, nous
en d\'eduisons que ces termes sont tous des s\'eries {\it
convergentes}\, par rapport \`a $z_1$.
Utilisons alors la premi\`ere ligne
de~\thetag{ 7.81} et interpr\'etons le r\'esultat
ci-apr\`es~:
\def\theequation{7.84}\begin{equation}
\left\{
\aligned
\left[
\Upsilon^\delta g_{ j'}\right] \left(
\Gamma_1 (z_1) \right)
\equiv
& \
\sum_{ \gamma' \in\N^{ m'}} \, 
\sum_{ \delta_1 \leq \delta} \, 
\frac{ \delta !}{ \delta_1 ! \ 
(\delta -\delta_1)!} \, 
\left[
\Upsilon^{ \delta_1} \overline{ \Theta}_{ j',\, \gamma'}'
(\overline{ h}) \right]
\left(
\Gamma_1 (z_1)
\right)
\cdot \\
& \
\ \ \ \ \ \ \ \ \ \ \ \ \
\cdot
N_{\gamma',\, \delta-\delta_1} 
\left(
\left(
\Upsilon^{ \delta_2} f_{ k'} (\Gamma_1 (z_1))
\right)_{ 1\leq k' \leq m',\, \delta_2 \leq \delta -\delta_1}
\right), 
\\
0 
\equiv
& \
\sum_{ \gamma' \in\N^{ m'}} \, 
\sum_{ \delta_1 \leq \delta} \, 
\frac{ \delta !}{ \delta_1 ! \ 
(\delta -\delta_1)!} \, 
\left[
\Upsilon^{ \delta_1}\underline{ \mathcal{ L}}^\beta 
\overline{ \Theta}_{ j',\, \gamma'}'
(\overline{ h}) \right]
\left(
\Gamma_1 (z_1)
\right)
\cdot \\
& \
\ \ \ \ \ \ \ \ \ \ \ \ \
\cdot
N_{\gamma',\, \delta-\delta_1} 
\left(
\left(
\Upsilon^{ \delta_2} f_{ k'} (\Gamma_1 (z_1))
\right)_{ 1\leq k' \leq m',\, \delta_2 \leq \delta -\delta_1}
\right).
\endaligned\right.
\end{equation}
Nous obtenons des \'equations analytiques satisfaites par
l'application formelle 
\def\theequation{7.85}\begin{equation}
z_1 \longmapsto_{ \mathcal{ F}} \, 
\left(
\Upsilon^{\delta_1} \, h_{ i_1'}
\left(
\Gamma_1 (z_1)
\right)
\right)_{ 1\leq i_1' \leq n',\, 
\delta_1 \leq \delta}.
\end{equation}
Cette application formelle ne s'annule pas forc\'ement en
$z_1=0$. Heureusement, puisque $h(0)=0$ et puisque ces \'equations
sont relativement polynomiales par rapport \`a toutes les variables
$\Upsilon^{ \delta_1} h_{ i_1'}$, telles que $\vert \delta_1 \vert
\geq 1$, nous pouvons appliquer le Th\'eor\`eme
d'approximation~\thetag{ 2.4} avec l'ordre d'approximation $N=1$ ({\it
cf.} la d\'emonstration du Lemme~7.7). Nous en d\'eduisons qu'il
existe une application convergente
\def\theequation{7.86}\begin{equation}
\left\{
\aligned
z_1 \longmapsto 
& \
\left(
{\sf H}_{ i_1',\, \delta_1}(z_1) 
\right)_{ 1\leq i_1 ' \leq n',\, \delta_1 \leq \delta}
\\
& \
=:
\left(
\left(
{\sf F}_{ k',\, \delta_1} (z_1)
\right)_{1\leq k'\leq m',\, \delta_1 \leq \delta}, \
\left(
{\sf G}_{ j',\, \delta_1} (z_1)
\right)_{1\leq j'\leq d',\, \delta_1 \leq \delta}
\right),
\endaligned\right.
\end{equation}
avec ${\sf H}_{ i_1',\, \delta_1} (0) =
\Upsilon^{ \delta_1} h_{ i_1'} (0)$, qui satisfait les
\'equations analytiques~\thetag{ 7.84}, c'est-\`a-dire~:
\def\theequation{7.87}\begin{equation}
\left\{
\aligned
{\sf G}_{ j',\, \delta}(z_1)
& 
\equiv
\sum_{ \gamma' \in\N^{ m'}} \, 
\sum_{ \delta_1 \leq \delta} \, 
\frac{ \delta !}{ \delta_1 ! \ 
(\delta -\delta_1)!} \ 
\left[
\Upsilon^{ \delta_1} \overline{ \Theta}_{ j',\, \gamma'}'
(\overline{ h}) \right]
\left(
\Gamma_1 (z_1)
\right)
\cdot \\
& \
\ \ \ \ \ \ \ \ \ \ \ \ \
\ \ \ \ \ \ \ \ \ \ \ \ \ 
\ \ \ \ \ \ \ \
\cdot 
N_{\gamma',\, \delta-\delta_1}
\left(
\left(
{\sf F}_{ k',\, \delta_2}(z_1)
\right)_{ 1\leq k' \leq m',\, \delta_2 \leq \delta -\delta_1}
\right), 
\\
0
& 
\equiv
\sum_{ \gamma' \in\N^{ m'}} \, 
\sum_{ \delta_1 \leq \delta} \, 
\frac{ \delta !}{ \delta_1 ! \ 
(\delta -\delta_1)!} \
\left[
\Upsilon^{ \delta_1}\underline{ \mathcal{ L}}^\beta 
\overline{ \Theta}_{ j',\, \gamma'}'
(\overline{ h}) \right]
\left(
\Gamma_1 (z_1)
\right)
\cdot \\
& \
\ \ \ \ \ \ \ \ \ \ \ \ \
\ \ \ \ \ \ \ \ \ \ \ \ \ 
\ \ \ \ \ \ \ \
\cdot 
N_{\gamma',\, \delta-\delta_1}
\left(
\left(
{\sf F}_{ k',\, \delta_2}(z_1)
\right)_{ 1\leq k' \leq m',\, \delta_2 \leq \delta -\delta_1}
\right),
\endaligned\right.
\end{equation}
pour tout $j'=1,\, \dots,\, d'$ et tout $\beta \in \N^m$ satisfaisant
$\vert \beta \vert \geq 1$. Reportons-nous maintenant \`a la
d\'emonstration de la propri\'et\'e~\thetag{ 7.66}. Gr\^ace \`a
l'\'equivalence~\thetag{ 7.58}, nous avons pu faire basculer les
identit\'es~\thetag{ 7.68} satisfaites par l'application convergente
${\sf H} (z_1)$, dans lesquelles on ne voit pas formellement
appara\^{\i}tre les composantes de l'application de r\'eflexion, vers
les <<meilleures>> identit\'es~\thetag{ 7.69}, dans lesquelles on voit
appara\^{\i}tre les modifications convergentes $\Theta_{ j',\,
\gamma'} ' ({\sf H} (z_1))$ des composantes de l'application de
r\'eflexion. Nous affirmons qu'il est aussi possible de {\it
transformer les identit\'es~\thetag{ 7.87}, satisfaites par
l'application convergente~\thetag{ 7.86}, en des identit\'es qui font
appara\^{\i}tre les d\'eriv\'ees} $\left[ \Upsilon^{ \delta_1}
\Theta_{ j',\, \gamma '} ' \right] \left( \left( {\sf H}_{ i_1',\,
\delta_2 }( z_1) \right)_{ 1 \leq i_1',\, \delta_2 \leq \delta_1}
\right)$~: expliquons soigneusement ce point-cl\'e.

\subsection*{ 7.88.~Comparaison avec la seconde famille
d'identit\'es de r\'eflexion}
Au lieu de consid\'erer les identit\'es de r\'eflexion 
\'ecrites \`a la seconde ligne de~\thetag{ 7.62}, consid\'erons
celles qui sont \'ecrites \`a la premi\`ere ligne et
rempla\c cons-y les variables $(t,\, \tau)\in \mathcal{ M}$ par
$\Upsilon_{ w_1} (\Gamma_1 (z_1))$, ce qui donne~:
\def\theequation{7.89}\begin{equation}
\left\{
\aligned
{}
& \
\left[
\underline{ \mathcal{ L}}^\beta
\overline{ g}_{ j'} \right] 
\left(
\Upsilon_{ w_1} (\Gamma_1 (z_1))
\right) 
\equiv \\
& \
\ \ \ \ \ \ \ \ \ \ \ 
\equiv
\sum_{ \gamma ' \in \N^{ m'}}\, 
\left[
\underline{ \mathcal{ L}}^\beta \overline{ f}^{ \gamma'}
\right] \left(
\Upsilon_{ w_1} (\Gamma_1 (z_1))
\right) \cdot
\Theta_{ j',\, \gamma'} ' 
\left(
h \left(
\Upsilon_{ w_1} (\Gamma_1 (z_1))
\right)
\right).
\endaligned\right.
\end{equation}
Appliquons l'op\'erateur $\left. \partial_{ w_1 }^\delta \right\vert{
w_1 =0}$, en tenant compte de la formule de Leibniz~:
\def\theequation{7.90}\begin{equation}
\left\{
\aligned
\left[
\Upsilon^\delta
\underline{ \mathcal{ L}}^\beta
\overline{ g}_{ j'} \right](\Gamma_1 (z_1)) \equiv
& \
\sum_{ \gamma ' \in \N^{ m'}} \sum_{ \delta_1 \leq \delta} \, 
\frac{ \delta! }{\delta_1 ! \ (\delta- \delta_1)!} \
\left[
\Upsilon^{\delta- \delta_1} 
\underline{ \mathcal{ L}}^\beta
\overline{ f}^{ \gamma'}
\right](\Gamma_1 (z_1)) \cdot \\
& \ 
\ \ \ \ \ \ \ \ \ \ \ \ \ \ \ \ \
\cdot
\left[
\Upsilon^{ \delta_1}
\Theta_{ j',\, \gamma'}' (h)
\right](\Gamma_1 (z_1)).
\endaligned\right.
\end{equation}
Gr\^ace aux d\'eveloppements de la deuxi\`eme et de la cinqui\`eme
ligne de~\thetag{ 7.81}, en utilisant au passage la relation
cruciale~\thetag{ 6.8}, nous voyons que les termes $ \left[
\Upsilon^\delta \underline{ \mathcal{ L}}^\beta \overline{ g}_{ j'}
\right](\Gamma_1 (z_1))$ et $\left[ \Upsilon^{\delta- \delta_1}
\underline{ \mathcal{ L}}^\beta \overline{ f}^{ \gamma'}
\right](\Gamma_1 (z_1))$ sont des s\'eries convergentes par rapport
\`a $z_1$. Enfin, utilisons la quatri\`eme ligne de~\thetag{ 7.91},
pour d\'evelopper la deuxi\`eme ligne de~\thetag{ 7.90}, ce qui
donne~:
\def\theequation{7.91}\begin{equation}
\left\{
\aligned
{}
& 
\left[
\Upsilon^\delta
\underline{ \mathcal{ L}}^\beta
\overline{ g}_{ j'} \right](\Gamma_1 (z_1)) 
\equiv
\sum_{ \gamma ' \in \N^{ m'}} \sum_{ \delta_1 \leq \delta} \, 
\frac{ \delta! }{\delta_1 ! \ (\delta- \delta_1)!} \
\left[
\Upsilon^{\delta- \delta_1} 
\underline{ \mathcal{ L}}^\beta
\overline{ f}^{ \gamma'}
\right](\Gamma_1 (z_1)) 
\cdot \\
& \ 
\ \ \ \ \ \ \
\cdot
M_{ \delta_1} '
\left(
\left(
\Upsilon^{ \delta_2} h_{ i_1'}(\Gamma_1(z_1))
\right)_{ 1\leq i_1' \leq n',\, \delta_2 \leq \delta_1}, \
\left[
J_{ t'}^{ \vert \delta_1 \vert}
\Theta_{ j',\, \gamma'}'
\right] \left(h(\Gamma_1 (z_1))\right)
\right),
\endaligned\right.
\end{equation}
pour tout $j' = 1,\, \dots,\, d'$ et tout $\beta \in \N^m$. La m\^eme
application formelle~\thetag{ 7.85} est donc solution des \'equations
analytiques~\thetag{ 7.91}, qui sont visiblement distinctes des
\'equations~\thetag{ 7.84}. Or nous avons introduit l'application {\it
convergente}~\thetag{ 7.86} qui satisfait les \'equations~\thetag{
7.87}. {\it Il nous faut maintenant g\'en\'eraliser
l'\'equivalence~\thetag{ 7.58} pour v\'erifier que cette m\^eme
application convergente~\thetag{ 7.86} est aussi solution des
\'equations~\thetag{ 7.91}}, c'est-\`a-dire qu'il nous faut
\'etablir que~: 
\def\theequation{7.92}\begin{equation}
\left\{
\aligned
{}
& 
\left[
\Upsilon^\delta
\underline{ \mathcal{ L}}^\beta
\overline{ g}_{ j'} \right](\Gamma_1 (z_1)) 
\equiv
\sum_{ \gamma ' \in \N^{ m'}} \sum_{ \delta_1 \leq \delta} \, 
\frac{ \delta! }{\delta_1 ! \ (\delta- \delta_1)!} \
\left[
\Upsilon^{\delta- \delta_1} 
\underline{ \mathcal{ L}}^\beta
\overline{ f}^{ \gamma'}
\right](\Gamma_1 (z_1)) 
\cdot \\
& \ 
\ \ \ \ \ \ \
\cdot
M_{ \delta_1} '
\left(
\left(
{\sf H}_{ i_1',\, \delta_2} (z_1)
\right)_{ 1\leq i_1' \leq n',\, \delta_2 \leq \delta_1}, \
\left[
J_{ t'}^{ \vert \delta \vert}
\Theta_{ j',\, \gamma'}'
\right] \left({\sf H}_0 (z_1)\right)
\right),
\endaligned\right.
\end{equation}
pour tout $j' = 1,\, \dots,\, d'$ et tout $\beta \in \N^m$, o\`u nous
avons pos\'e ${\sf H}_0 := \left( {\sf H}_{ i_1',\, 0} \right)_{ 1\leq
i_1' \leq n'}$. Dans le \S7.98 ci-dessous, nous \'etablirons en effet
le Lemme~7.111, qui implique l'assertion suivante.

\def\thelemma{7.93}\begin{lemma}
Fixons $\ell \in \N$ et supposons qu'une application convergente
\def\theequation{7.94}\begin{equation}
z_1 \longmapsto 
\left(
{\sf H}_{ i_1',\, \delta}(z_1)
\right)_{ 1\leq i_1 ' \leq n',\, 
\vert \delta\vert \leq \ell} 
\end{equation}
telle que ${\sf H}_{ i_1',\, \delta} (0) = \Upsilon^\delta h_{ i_1'}
(0)$, $i_1' =1,\, \dots,\, n'$, $\vert \delta \vert \leq \ell$,
satisfait les \'equations~\thetag{ 7.87}, pour tout $j'= 1,\, \dots,\,
d'$ et tout $\delta$ tel que $\vert \delta \vert \leq \ell$. Alors
elle satisfait les \'equations~\thetag{ 7.92}, pour tout $j'= 1,\,
\dots,\, d'$ et tout $\delta$ tel que $\vert \delta \vert \leq \ell$.
\end{lemma}

Gr\^ace \`a ce lemme, nous pouvons enfin \'etablir que
pour tout $j' =1,\, \dots,\, d'$, tout $\delta \in \N^d$,
et tout $\gamma ' \in \N^{ m'}$, on a la propri\'et\'e
de convergence~:
\def\theequation{7.95}\begin{equation}
\left[
\Upsilon^\delta 
\Theta_{ j',\, \gamma'} ' 
(h) \right]
\left(
z_1, \, \overline{ \Theta} (z_1,\, 0)
\right)\in \C\{ z_1\},
\end{equation}
qui g\'en\'eralise~\thetag{ 7.66}.

\proof
Soit $\ell \in \N$ arbitraire. D\'emontrons que la
propri\'et\'e~\thetag{ 7.95} est vraie pour tout $\delta$ tel que
$\vert \delta \vert \leq \ell$. En appliquant le Th\'eor\`eme
d'approximation~2.4 aux \'equations \thetag{ 7.84}, \'ecrites pour
tout $j'= 1,\, \dots,\, d'$ et tout $\delta$ tel que $\vert \delta
\vert \leq \ell$, et en appliquant le Lemme~7.93, on trouve une
application convergente~\thetag{ 7.94} qui satisfait les
\'equations~\thetag{ 7.92} pout tout $j'= 1,\, \dots,\, d'$ et tout
$\delta$ tel que $\vert \delta \vert \leq \ell$. Nous allons
d\'emontrer par r\'ecurrence que pour tout $\delta$ tel que $\vert
\delta \vert \leq \ell$, on a~:
\def\theequation{7.96}\begin{equation}
\left\{
\aligned
{}
& \
M_{ \delta}' 
\left(
\left(
\Upsilon^{ \delta_1} h_{ i_1'}(\Gamma_1(z_1))
\right)_{ 1\leq i_1' \leq n',\, \delta_1 \leq \delta}, \
\left[
J_{ t'}^{ \vert \delta \vert}
\Theta_{ j',\, \gamma'}
\right] \left(h(\Gamma_1 (z_1))\right)
\right) \equiv \\
& \
\ \ \ \ \ 
\equiv
M_{ \delta}' 
\left(
\left(
{\sf H}_{ i_1',\, \delta_1} (z_1)
\right)_{ 1\leq i_1' \leq n',\, \delta_1 \leq \delta}, \
\left[
J_{ t'}^{ \vert \delta \vert}
\Theta_{ j',\, \gamma'}
\right] \left({\sf H}_0 (z_1)\right)
\right),
\endaligned\right.
\end{equation}
pour tout $j'= 1,\, \dots,\, d'$ et tout $\gamma ' \in \N^{ m'}$, ce
qui \'etablira la convergence d\'esir\'ee, puisque
l'application~\thetag{ 7.94} converge.

Pour $\delta =0$, les relations~\thetag{ 7.96} ont d\'ej\`a \'et\'e
d\'emontr\'ees en~\thetag{ 7.72}. Soit $\ell_1 \in \N$ tel que $0
\leq \ell_1 \leq \ell- 1$ et supposons les identit\'es~\thetag{ 7.96}
vraies pour tout $\delta$ tel que $\vert \delta \vert \leq \ell_1$, et
bien s\^ur, pour tout $j'=1,\, \dots,\, d'$ et tout $\gamma ' \in
\N^{m'}$. Soit $\delta \in \N^d$ un multiindice arbitraire tel que
$\vert \delta \vert = \ell_1 +1$. Soustrayons alors les
identit\'es~\thetag{ 7.92} des identit\'es~\thetag{ 7.91}, \'ecrites
avec ce multiindice $\delta$. Dans la somme $\sum_{ \delta_1 \leq
\delta}$, on a ou bien $\delta_1 = \delta$ ou bien $\vert \delta_1
\vert \leq \ell_1$. Gr\^ace \`a l'hypoth\`ese de r\'ecurrence, il ne
reste donc que le terme $\delta_1= \delta$ dans cette somme apr\`es
soustraction, ce qui donne~:
\def\theequation{7.97}\begin{equation}
\left\{
\aligned
{}
0 \equiv 
&
\sum_{ \gamma' \in\N^{m'}} \, 
\left[
\underline{ \mathcal{ L}}^\beta \overline{ f}^{ \gamma'} \right]
(z_1,\, \overline{ \Theta} (z_1,\, 0),\, 0,\, 0) \cdot \\
&
\cdot \left(
\aligned
{}
&
M_{ \delta} '
\left(
\left(
\Upsilon^{ \delta_1} h_{ i_1'}(\Gamma_1(z_1))
\right)_{ 1\leq i_1' \leq n',\, \delta_1 \leq \delta_1}, \
\left[
J_{ t'}^{ \vert \delta_1 \vert}
\Theta_{ j',\, \gamma'}'
\right] \left(h(\Gamma_1 (z_1))\right)
\right) - \\
&
-
M_{ \delta} '
\left(
\left(
{\sf H}_{ i_1',\, \delta_1} (z_1)
\right)_{ 1\leq i_1' \leq n',\, \delta_1 \leq \delta_1}, \
\left[
J_{ t'}^{ \vert \delta \vert}
\Theta_{ j',\, \gamma'}'
\right] \left({\sf H}_0 (z_1)\right)
\right)
\endaligned
\right).
\endaligned\right.
\end{equation}
Gr\^ace au principe d'unicit\'e~\thetag{ 7.37}, nous d\'eduisons que
les identit\'es~\thetag{ 7.96} sont satisfaites pour ce multiindice
$\delta$ de longueur $\ell_1+1$, ce qui ach\`eve la d\'emonstration du
Lemme~7.93.
\endproof

\subsection*{ 7.98.~Transformation des identit\'es de r\'eflexion}
Pour formuler la g\'en\'eralisation de l'\'equivalence~\thetag{ 7.58},
quittons provisoirement la premi\`ere cha\^{\i}ne de Segre et
travaillons avec des variables $(t,\, \tau) \in \mathcal{ M}$
quelconques, ce qui sera utile pour la d\'emonstration finale de la
Proposition~7.18 ({\it voir} le~\S7.135 ci-dessous). Soit $w_1 \in
\C^d$. Puisque le multiflot de $\Upsilon$ stabilise $\mathcal{ M}$, on
a $\Upsilon_{ w_1} (t,\, \tau) \in \mathcal{M}$ et~: 
\def\theequation{7.99}\begin{equation}
\left\{
\aligned
0 = 
& \
\underline{ \mathcal{ L}}^\beta
\overline{ r}_{ j'} ' \left( \overline{ h} \left( \Upsilon_{ w_1}
(t,\, \tau) \right), \, h \left( \Upsilon_{ w_1} ( t,\, \tau) \right)
\right), \\
0 =
& \
\underline{ \mathcal{ L}}^\beta 
r_{ j'} ' \left( h \left( \Upsilon_{ w_1}
(t,\, \tau) \right), \, \overline{ h} 
\left( \Upsilon_{ w_1} ( t,\, \tau) \right)
\right),
\endaligned\right.
\end{equation} 
pour tout $j'= 1,\, \dots,\, d'$ et tout $\beta \in \N^m$. 
Polarisons ces identit\'es en y rempla\c cant $h$ par $t'$ et
d\'efinissons ainsi deux familles infinies de s\'eries 
comme suit~:
\def\theequation{7.100}\begin{equation}
\left\{
\aligned
\overline{ R}_{ j',\, 0}' \left( w_1,\, t,\, \tau : \, 
t' \right) := 
& \
\overline{ r}_{ j'} ' \left(
\overline{ h} \left(
\Upsilon_{ w_1} ( t,\, \tau)
\right), \, t'
\right) \\
=
& \
w_{j'} '
- \sum_{ \gamma' \in \N^{ m'}} \,
{z'}^{\gamma '}
\cdot
\left[
\overline{ \Theta}_{ j',\, \gamma'}' \left(
\overline{ h}
\right)
\right] 
\left(
\Upsilon_{ w_1} (t,\,\tau)
\right); \\
\overline{ R}_{ j',\, \beta}' \left( w_1,\, t,\, \tau : \, 
t' \right) := 
& \
\underline{ \mathcal{ L}}^\beta 
\overline{ r}_{ j'} \left(
\overline{ h} \left(
\Upsilon_{ w_1} (t,\, \tau)
\right),\, t'
\right)
\\
= 
& \
- \sum_{ \gamma' \in \N^{ m'}} \,
{z'}^{\gamma '}
\cdot
\left[
\underline{ \mathcal{ L}}^\beta 
\overline{ \Theta}_{ j',\, \gamma'}' \left(
\overline{ h}
\right)
\right] 
\left(
\Upsilon_{ w_1} (t,\,\tau)
\right); \\
S_{ j',\, \beta}'
\left( w_1,\, t,\, \tau : \, t'\right) := 
& \
\underline{ \mathcal{L}}^\beta 
r_{ j'} ' \left(
t',\, \overline{ h} \left(
\Upsilon_{ w_1} (t,\, \tau)
\right)
\right) \\
=
& \
\left[
\underline{ \mathcal{ L}}^\beta
\overline{ g}_{ j'}
\right] 
\left(
\Upsilon_{ w_1} (t,\,\tau)
\right) - \\
& \
\ \ \ \
- \sum_{ \gamma' \in \N^{ m'}} \, 
\left[
\underline{ \mathcal{ L}}^\beta
\overline{ f}^{\gamma '}
\right] 
\left(\Upsilon_{ w_1} (t,\,\tau)\right) \cdot
\Theta_{ j',\, \gamma'}' (t').
\endaligned\right.
\end{equation}
\`A la troisi\`eme ligne, on suppose $\vert \beta \vert \geq 1$. Avec
ces notations, les \'egalit\'es~\thetag{ 7.99} s'\'ecrivent 
maintenant~:
\def\theequation{7.101}\begin{equation}
\left\{
\aligned
0 =
& \ 
\overline{ R}_{ j',\, \beta}' \left(
w_1,\, t,\, \tau : \, 
h\left(
\Upsilon_{ w_1}( t,\, \tau)
\right)
\right), \\
0 =
& \ 
S_{ j',\, \beta}' \left(
w_1,\, t,\, \tau : \, 
h\left(
\Upsilon_{ w_1}( t,\, \tau)
\right)
\right),
\endaligned\right.
\end{equation}
toujours avec $(t,\, \tau) \in \mathcal{ M}$ et $w_1 \in \C^d$. Dans
la suite, on remplacera $(t,\, \tau)$ par $\Gamma_k (z_{ (k)})$, et
alors ces \'egalit\'es deviendront des identit\'es formelles dans $\C
\dl w_1,\, z_{ (k)}\dr$. 

Notons que $\overline{ R}_{ j',\, 0}$ s'identifie avec $\overline{
r}_{ j'}' \left( \overline{ h } \left( \Upsilon_{ w_1} (t,\, \tau)
\right),\, t' \right)$ et que $S_{ j',\, 0}$ s'identifie avec $r_{ j'}
\left(t',\, \overline{ h} \left( \Upsilon_{ w_1} (t,\, \tau) \right)
\right)$. En posant $\tau' = \overline{ h} \left( \Upsilon_{ w_1}
(t,\, \tau) \right)$ dans la deuxi\`eme colonne de~\thetag{ 7.44},
nous obtenons les identit\'es suivantes~:
\def\theequation{7.102}\begin{equation}
\left\{
\aligned
\overline{ R}_{ j',\, 0} ' 
\left(
w_1,\, 
t,\, \tau : \, t'
\right)
\equiv
& \
\sum_{ j_1 '= 1}^{ d'} \, 
A_{ j'}^{ j_1'} 
\left(
w_1,\, 
t,\, \tau : \, t'
\right) \cdot
S_{ j_1',\, 0}' 
\left(
w_1,\, 
t,\, \tau : \, t'
\right), \\
S_{ j',\, 0} ' 
\left(
w_1,\, 
t,\, \tau : \, t'
\right) \equiv
& \
\sum_{ j_1 '= 1}^{ d'} \, 
B_{ j'}^{ j_1'} 
\left(
w_1,\, 
t,\, \tau : \, t'
\right) \cdot
\overline{ R}_{ j_1',\, 0}' 
\left(
w_1,\, 
t,\, \tau : \, t'
\right),
\endaligned\right.
\end{equation}
o\`u nous avons pos\'e $A_{ j'}^{ j_1'} := {a'}_{ j'}^{ j_1'} \left(
t',\, h \left( \Upsilon_{ w_1} (t,\, \tau) \right) \right)$ et $B_{
j'}^{ j_1'} := {\overline{ a}'}_{ j'}^{ j_1'} \left( h \left(
\Upsilon_{ w_1} (t,\, \tau) \right),\, t'\right)$.
En appliquant les d\'erivations $\underline{ \mathcal{ L}}^\beta$
\`a ce couple d'identit\'es, nous obtenons les 
combinaisons lin\'eaires suivantes, que nous \'ecrivons
sans les arguments $(w_1,\, t,\, \tau : \, t')$~:
\def\theequation{7.103}\begin{equation}
\left\{
\aligned
\underline{ \mathcal{ L}}^\beta
\overline{ R}_{ j',\, 0} ' \equiv
\overline{ R}_{ j',\, \beta} ' \equiv
& \
\sum_{ j_1 '= 1}^{ d'} \,
\sum_{ \beta_1 \leq \beta} \, 
A_{ j',\, \beta}^{ j_1',\, \beta_1} \cdot
S_{ j_1',\, \beta_1}', \\
\underline{ \mathcal{ L}}^\beta
S_{ j',\, 0} ' \equiv
S_{ j',\, \beta} ' \equiv 
& \
\sum_{ j_1 '= 1}^{ d'} \,
\sum_{ \beta_1 \leq \beta} \,
B_{ j',\, \beta}^{ j_1',\, \beta_1} \cdot
\overline{ R}_{ j_1',\, \beta_1}',
\endaligned\right.
\end{equation}
o\`u nous avons pos\'e~:
\def\theequation{7.104}\begin{equation}
\left\{
\aligned
A_{ j',\, \beta}^{ j_1',\, \beta_1} :=
& \ 
\frac{ \beta !}{\beta_1! \ 
(\beta- \beta_1)!} \, 
\underline{ \mathcal{ L}}^{ 
\beta - \beta_1} 
\left(
A_{ j'}^{ j_1'}
\right), \\
B_{ j',\, \beta}^{ j_1',\, \beta_1} :=
& \ 
\frac{ \beta !}{\beta_1! \ 
(\beta- \beta_1)!} \, 
\underline{ \mathcal{ L}}^{ 
\beta - \beta_1} 
\left(
B_{ j'}^{ j_1'}
\right).
\endaligned\right.
\end{equation} 
Rappelons que les identit\'es~\thetag{ 7.103} nous ont 
servi pour \'etablir l'\'equivalence~\thetag{ 7.58}. 
Pour g\'en\'eraliser cette \'equivalence (et obtenir
la preuve du Lemme~7.93), nous allons les diff\'erentier
par rapport \`a $w_1$.

Pour cela, \'ecrivons $w_1 = (w_{1;1},\, \dots,\, w_{ 1; d}) \in \C^d$,
choisissons un entier $j$ compris entre $1$ et $d$, et appliquons la
d\'erivation $\frac{ \partial }{\partial_{ w_{1;j }}}$ aux
identit\'es~\thetag{ 7.101}, sans \'ecrire les arguments, ce qui nous
donne~:
\def\theequation{7.105}\begin{equation}
\left\{
\aligned
0 =
& \ 
\frac{ \partial \overline{ R}_{ j',\, \beta} '}{\partial w_{ 1;j}} +
\sum_{ i_1' = 1}^{ n'} \, 
\frac{ \partial \overline{ R}_{ j',\, \beta} '}{
\partial t_{ i_1'}'} \cdot
\left[\Upsilon_j h_{ i_1'}\right], \\
0 =
& \ 
\frac{ \partial S_{ j',\, \beta} '}{\partial w_{ 1;j}} +
\sum_{ i_1' = 1}^{ n'} \, 
\frac{ \partial S_{ j',\, \beta} '}{
\partial t_{ i_1'}'} \cdot
\left[\Upsilon_j h_{ i_1'}\right].
\endaligned\right.
\end{equation}
Introduisons de nouvelles variables ind\'ependantes $T_{ i_1',\, j}'$,
qui correspondent aux d\'eriv\'ees $\left[ \Upsilon_j h_{ i_1'}
\right]$ apparaissant dans~\thetag{ 7.105}, et d\'efinissons~:
\def\theequation{7.106}\begin{equation}
\left\{
\aligned
\overline{ R}_{ j',\, \beta,\, j}' 
\left(
w_1,\, 
t,\, \tau : t', \left( 
T_{ i_1', \, j}' \right)_{ 1\leq i_1' \leq n'} \right)
:=
& \ 
\frac{ \partial \overline{ R}_{ j',\, \beta} '}{\partial w_{ 1;j}} +
\sum_{ i_1' = 1}^{ n'} \, 
\frac{ \partial \overline{ R}_{ j',\, \beta} '}{
\partial t_{ i_1'}'} \cdot
T_{ i_1',\, j} '
, \\
S_{ j',\, \beta,\, j}' 
\left(
w_1,\, 
t,\, \tau : t', \left( 
T_{ i_1', \, j}' \right)_{ 1\leq i_1' \leq n'} \right)
:=
& \ 
\frac{ \partial S_{ j',\, \beta} '}{\partial w_{ 1;j}} +
\sum_{ i_1' = 1}^{ n'} \, 
\frac{ \partial S_{ j',\, \beta} '}{
\partial t_{ i_1'}'} \cdot
T_{ i_1',\, j}'.
\endaligned\right.
\end{equation}
Plus g\'en\'eralement, soit $\delta \in \N^d$ un multiindice
arbitraire. En appliquant l'op\'erateur $\partial_w^\delta$ aux
identit\'es~\thetag{ 7.101}, nous obtenons des \'egalit\'es qui
impliquent deux familles infinies de s\'eries $\overline{ R}_{ j',\,
\beta,\, \delta}'$ et $S_{ j',\, \beta,\, \delta}'$, qui sont
polynomiales par rapport aux variables diff\'erenti\'ees strictement
$\left[ \Upsilon^{ \delta_1} h_{ i_1'}\right]$, o\`u $1\leq i_1 ' \leq
n'$ et $0 \neq \delta_1 \leq \delta$, que nous \'ecrirons~:
\def\theequation{7.107}\begin{equation}
\left\{
\aligned
0 = 
& \
\overline{ R}_{ j',\, \beta,\, \delta}' 
\left(
w_1,\, 
t,\, \tau : \left( 
\left[ \Upsilon^{ \delta_1} h_{ i_1'} \right]
\left(
\Upsilon_{ w_1} (t,\, \tau)
\right)
\right)_{ 1\leq i_1' \leq n',\, \delta_1 \leq \delta} \right), 
\\
0 =
& \
S_{ j',\, \beta,\, \delta}' 
\left(
w_1,\, 
t,\, \tau : \left( 
\left[ \Upsilon^{ \delta_1} h_{ i_1'} \right]
\left(
\Upsilon_{ w_1} (t,\, \tau)
\right)
\right)_{ 1\leq i_1' \leq n',\, \delta_1 \leq \delta} \right).
\endaligned\right.
\end{equation}
Introduisons des variables ind\'ependantes $\left( T_{ i_1',\,
\delta_1}' \right)_{ 1\leq i_1 ' \leq n',\, \delta_1 \leq \delta}$
correspondant \`a ces d\'eriv\'ees $\Upsilon^{ \delta_1 } h_{ i_1' }$,
en convenant que $T_{ i_1',\, 0} '$ s'identifie \`a $t_{ i_1'}'$, et
d\'efinissons les nouvelles s\'eries~:
\def\theequation{7.108}\begin{equation}
\left\{
\aligned
{}
& \
\overline{ R}_{ j',\, \beta,\, \delta}' 
\left(
w_1,\, 
t,\, \tau : \left( 
T_{ i_1', \, \delta_1}' 
\right)_{ 1\leq i_1' \leq n',\, \delta_1\leq \delta} 
\right), \\ 
& \
S_{ j',\, \beta,\, \delta}' 
\left(
w_1,\, 
t,\, \tau : \left( 
T_{ i_1', \, \delta_1}' 
\right)_{ 1\leq i_1' \leq n',\, \delta_1 \leq \delta} 
\right).
\endaligned\right.
\end{equation}
Alors les identit\'es~\thetag{ 7.79} et~\thetag{ 7.90} (o\`u 
$\delta \in \N^d$ est fix\'e et o\`u 
$j'=1,\,\dots,\, d'$ et $\beta \in\N^m$ varient)
co\"{\i}ncident avec les deux familles d'identit\'es~:
\def\theequation{7.109}\begin{equation}
\left\{
\aligned
0 \equiv
& \
\overline{ R}_{ j',\, \beta,\, \delta} '
\left(
0,\, \Gamma_1 (z_1) : 
\left(
\Upsilon^{ \delta_1} h_{ i_1'} 
\left(
\Gamma_1 (z_1)
\right)
\right)_{ 1\leq i_1 ' \leq n',\, \delta_1 \leq \delta}
\right), \\
0 \equiv
& \
S_{ j',\, \beta,\, \delta} '
\left(
0,\, \Gamma_1 (z_1) : 
\left(
\Upsilon^{ \delta_1} h_{ i_1'} 
\left(
\Gamma_1 (z_1)
\right)
\right)_{ 1\leq i_1 ' \leq n',\, \delta_1 \leq \delta}
\right).
\endaligned\right.
\end{equation}

Apr\`es ces pr\'eliminaires, nous pouvons enfin \'enoncer la
g\'en\'eralisation attendue du Lemme~7.57. Rappelons que nous avons
pos\'e $w_1 =0$ dans~\thetag{ 7.79} et dans~\thetag{ 7.90}~; c'est
pourquoi nous poserons aussi $w_1 =0$ dans le lemme suivant.

\def\thelemma{7.110}\begin{lemma} 
Soit $\nu \in \N$, soit ${\sf x} \in \C^\nu$, soit ${\sf Q} ({\sf x})
\in \C\dl {\sf x} \dr^{ 2n}$ avec ${\sf Q} (0)= 0$, soit $\delta \in
\N^d$, et soit ${\sf x} \longmapsto_{ \mathcal{ F} } \, \left( {\sf
T}_{ i_1',\, \delta_1}' ({\sf x}) \right)_{ 1\leq i_1' \leq n',\,
\delta_1 \leq \delta}$ une application formelle telle que ${\sf T}_{
i_1',\, \delta_1} ' (0) = \left[ \Upsilon^{ \delta_1 } h_{ i_1'}
\right] (0)$. L'\'equivalence suivante est satisfaite~:
\def\theequation{7.111}\begin{equation}
\left\{
\aligned
{}
& \
\left(
\aligned
{}
& \
\forall \, j'=1,\, \dots,\, d', \ 
\forall \, \beta \in \N^m, \ 
\forall \, \delta_1 \leq \delta, 
{\text{ \rm on a dans}} \
\C\dl {\sf x} \dr \ : \\
& \
\left.
\overline{ R}_{ j',\, \beta,\, \delta_1} '
\left(
0,\, t,\, \tau : \, 
\left(
{\sf T}_{ i_1',\, \delta_2}' ({\sf x})
\right)_{ 1\leq i_1' \leq n',\, \delta_2 \leq \delta_1}
\right)
\right\vert_{ (t,\, \tau) = {\sf Q}({\sf x})}
\equiv 0
\endaligned
\right)
\\
& \ 
\ \ \ \ \ \ \ \ \ \ \ \ \ 
\ \ \ \ \ \ \ \ \ \ \ \ \ 
\ \ \ \ \ \ \ \ \ \ \ \ \ 
\ \ \ \ \ \ \ \ \ \ \ \ \ 
\Updownarrow \\
& \
\left(
\aligned
{}
& \
\forall \, j'=1,\, \dots,\, d', \ 
\forall \, \beta \in \N^m, \ 
\forall \, \delta_1 \leq \delta, 
{\text{ \rm on a dans}} \
\C\dl {\sf x} \dr \ : \\
& \
\left.
S_{ j',\, \beta,\, \delta_1} '
\left(
0,\, t,\, \tau : \, 
\left(
{\sf T}_{ i_1',\, \delta_2}' ({\sf x})
\right)_{ 1\leq i_1' \leq n',\, \delta_2 \leq \delta_1}
\right)
\right\vert_{ (t,\, \tau) = {\sf Q}({\sf x})}
\equiv 0
\endaligned
\right).
\endaligned\right.
\end{equation}
\end{lemma} 

Gr\^ace \`a l'\'equivalence ~\thetag{ 7.111} ci-dessus, le Lemme~7.93
est d\'emontr\'e~: il suffit de poser ${\sf x} = z_1$, ${\sf Q} ({\sf
x}) := \Gamma_1 (z_1)$ et ${\sf T}_{ i_1',\, \delta_1} ({\sf x}) :=
{\sf H}_{ i_1',\, \delta_1} ({\sf x})$.

\proof
Nous affirmons que les s\'eries~\thetag{ 7.108} peuvent \^etre
d\'efinies par r\'ecurrence de la mani\`ere suivante. En effet, si
l'on note $\1_j^d$ le multiindice $(0,\, \dots,\, 0,\, 1,\, 0,\,
\dots,\, 0) \in \N^d$, avec $1$ \`a la $j$-i\`eme place et $0$
aux autres places, on a les relations~:
\def\theequation{7.112}\begin{equation}
\left\{
\aligned
\overline{ R}_{ j',\, \beta,\, \delta+\1_j^d}' 
:=
& \ 
\frac{ \partial \overline{ R}_{ j',\, \beta,\, \delta}'}{
\partial w_{ 1; j}} +
\sum_{ i_1 '= 1}^{ n'} \, \sum_{ \delta_1 \leq \delta}\, 
\frac{ \partial \overline{ R}_{ j',\, \beta,\, \delta}'}{
\partial T_{ i_1',\, \delta_1}'} \cdot
T_{ i_1 ',\, \delta+ \1_j^d}', \\
S_{ j',\, \beta,\, \delta+ \1_j^d} '
:=
& \ 
\frac{ \partial S_{ j',\, \beta,\, \delta}'}{
\partial w_{ 1; j}} +
\sum_{ i_1 '=1}^{ n'} \, \sum_{ \delta_1 \leq \delta}\, 
\frac{ \partial S_{ j',\, \beta,\, \delta}'}{
\partial T_{ i_1',\, \delta_1}'} \cdot
T_{ i_1 ',\, \delta + \1_j^d}',
\endaligned\right.
\end{equation}
de telle sorte que 
\def\theequation{7.113}\begin{equation}
\left\{
\aligned
{}
& \
\partial_{ w_1}^\delta
\left[
\overline{ R}_{ j',\, \beta}' 
\left(
w_1,\, t,\, \tau: \, h 
\left(
\Upsilon_{ w_1} (t,\, \tau)
\right)
\right)
\right]
\equiv \\
& \ 
\ \ \ \ \ \ \
\equiv
\overline{ R}_{ j',\, \beta,\, \delta}' 
\left(
w_1,\, t,\, \tau : \, 
\left(
\Upsilon^{ \delta_1} h_{ i_1'} 
\left(
\Upsilon_{w_1}( t,\, \tau)
\right)
\right)_{ 
1\leq i_1' \leq n',\, \delta_1 \leq \delta}
\right), \\
& \
\partial_{ w_1}^\delta
\left[
S_{ j',\, \beta}' 
\left(
w_1,\, t,\, \tau: \, h 
\left(
\Upsilon_{ w_1} (t,\, \tau)
\right)
\right)
\right]
\equiv \\
& \ 
\ \ \ \ \ \ \
\equiv
S_{ j',\, \beta,\, \delta}' 
\left(
w_1,\, t,\, \tau : \, 
\left(
\Upsilon^{ \delta_1} h_{ i_1'} 
\left(
\Upsilon_{w_1}( t,\, \tau)
\right)
\right)_{ 
1\leq i_1' \leq n',\, \delta_1 \leq \delta}
\right). 
\endaligned\right.
\end{equation}
Nous affirmons que pour tout $\delta \in \N^d$, 
il existe deux collections
de s\'eries formelles~:
\def\theequation{7.114}\begin{equation}
\left\{
\aligned
{}
& \
A_{ j',\, \beta,\, \delta}^{ j_1',\, \beta_1,\, \delta_1}
\left(
w_1,\, t,\, \tau: \, 
\left(
T_{ i_1',\, \delta_2} '
\right)_{1\leq i_1' \leq n',\, 
\delta_2 \leq \delta- \delta_1}
\right), 
\\
& \
B_{ j',\, \beta,\, \delta}^{ j_1',\, \beta_1,\, \delta_1}
\left(
w_1,\, t,\, \tau: \, 
\left(
T_{ i_1',\, \delta_2} '
\right)_{1\leq i_1' \leq n',\, 
\delta_2 \leq \delta - \delta_1}
\right),
\endaligned\right.
\end{equation}
telles que les identit\'es formelles suivantes, 
que nous \'ecrivons d'abord
sans les arguments,
sont satisfaites~: 
\def\theequation{7.115}\begin{equation}
\left\{
\aligned
{}
&
\overline{ R}_{ j',\, \beta,\, \delta}' 
\equiv
\sum_{ j_1'= 1}^{ d'} \, 
\sum_{ \beta_1 \leq \beta} \, 
\sum_{ \delta_1 \leq \delta} \, 
A_{ j',\, \beta,\, \delta}^{ j_1',\, \beta_1,\, \delta_1}
\cdot
S_{ j_1',\, \beta_1,\, \delta_1} ', \\
&
S_{ j',\, \beta,\, \delta}' 
\equiv
\sum_{ j_1'= 1}^{ d'} \, 
\sum_{ \beta_1 \leq \beta} \, 
\sum_{ \delta_1 \leq \delta} \, 
B_{ j',\, \beta,\, \delta}^{ j_1',\, \beta_1,\, \delta_1}
\cdot
\overline{ R}_{ j_1',\, \beta_1,\, \delta_1} ';
\endaligned\right.
\end{equation}
par souci de compl\'etude, \'ecrivons-les ensuite avec leurs arguments~:
\def\theequation{7.116}\begin{equation}
\left\{
\aligned
{}
&
\overline{ R}_{ j',\, \beta,\, \delta}' 
\left(
w_1,\, t,\, \tau: \, 
\left(
T_{ i_1',\, \delta_1} '
\right)_{1\leq i_1' \leq n',\, 
\delta_1 \leq \delta}
\right) \equiv
\sum_{ j_1'= 1}^{ d'} \, 
\sum_{ \beta_1 \leq \beta} \, 
\sum_{ \delta_1 \leq \delta} \, \\
& \
\ \ \ \ \ \ \ \ \ \ \ \ \ \ \ \
\ \ \ \ \ \ \ \ \ \ \ \ \ \ \ \
\ \ \ 
A_{ j',\, \beta,\, \delta}^{ 
j_1',\, \beta_1,\, \delta_1}
\left(
w_1,\, t,\, \tau: \, 
\left(
T_{ i_1',\, \delta_2} '
\right)_{1\leq i_1' \leq n',\, 
\delta_2 \leq \delta - \delta_1}
\right) \cdot \\ 
& \
\ \ \ \ \ \ \ \ \ \ \ \ \ \ \ \
\ \ \ \ \ \ \ \ \ \ \ \ \ \ \ \
\ \ \ 
\ \ \ 
\cdot
S_{ j_1',\, \beta_1,\, \delta_1} '
\left(
w_1,\, t,\, \tau: \, 
\left(
T_{ i_1',\, \delta_2} '
\right)_{1\leq i_1' \leq n',\, 
\delta_2 \leq \delta_1}
\right), \\
&
S_{ j',\, \beta,\, \delta}' 
\left(
w_1,\, t,\, \tau: \, 
\left(
T_{ i_1',\, \delta_1} '
\right)_{1\leq i_1' \leq n',\, 
\delta_1 \leq \delta}
\right) \equiv
\sum_{ j_1'= 1}^{ d'} \, 
\sum_{ \beta_1 \leq \beta} \, 
\sum_{ \delta_1 \leq \delta} \, \\
& \
\ \ \ \ \ \ \ \ \ \ \ \ \ \ \ \
\ \ \ \ \ \ \ \ \ \ \ \ \ \ \ \
\ \ \ 
B_{ j',\, \beta,\, \delta}^{ j_1',\, \beta_1,\, \delta_1}
\left(
w_1,\, t,\, \tau: \, 
\left(
T_{ i_1',\, \delta_2} '
\right)_{1\leq i_1' \leq n',\, 
\delta_2 \leq \delta - \delta_1}
\right) \cdot \\
& \
\ \ \ \ \ \ \ \ \ \ \ \ \ \ \ \
\ \ \ \ \ \ \ \ \ \ \ \ \ \ \ \
\ \ \ 
\ \ \ 
\cdot
\overline{ R}_{ j_1',\, \beta_1,\, \delta_1} '
\left(
w_1,\, t,\, \tau: \, 
\left(
T_{ i_1',\, \delta_2} '
\right)_{1\leq i_1' \leq n',\, 
\delta_2 \leq \delta_1}
\right).
\endaligned\right.
\end{equation}
Avant d'\'etablir cette affirmation, notons que les
identit\'es~\thetag{ 7.116} impliquent imm\'ediatement l'\'equivalence
d\'esir\'ee~\thetag{ 7.111}~: il suffit de poser $w_1=0$, de remplacer
$(t,\, \tau)$ par ${\sf Q}( {\sf x})$ et de remplacer $T_{ i_1',\,
\delta_1}'$ par ${\sf T}_{ i_1',\, \delta_1} '({\sf x})$ dans~\thetag{
7.116}.

Pour $\delta =0$, les identit\'es~\thetag{ 7.115} ont d\'ej\`a \'et\'e
vues dans~\thetag{ 7.103}. \'Etablissons les identit\'es~\thetag{
7.115} pour $\delta = \1_j^d$. Pour cela, partons de~\thetag{ 7.103}, 
que nous r\'e\'ecrivons~:
\def\theequation{7.117}\begin{equation}
\left\{
\aligned
\overline{ R}_{ j',\, \beta} ' \equiv
& \
\sum_{ j_1 '= 1}^{ d'} \,
\sum_{ \beta_1 \leq \beta} \, 
A_{ j',\, \beta}^{ j_1',\, \beta_1} \cdot
S_{ j_1',\, \beta_1}', \\
S_{ j',\, \beta} ' \equiv 
& \
\sum_{ j_1 '= 1}^{ d'} \,
\sum_{ \beta_1 \leq \beta} 
B_{ j',\, \beta}^{ j_1',\, \beta_1} \cdot
\overline{ R}_{ j_1',\, \beta_1}'.
\endaligned\right.
\end{equation}
Appliquons tout d'abord l'op\'erateur $\partial_{ w_{ 1;j}}$ de part
et d'autre de~\thetag{ 7.117}, ce qui donne~:
\def\theequation{7.118}\begin{equation}
\left\{
\aligned
\frac{ \partial
\overline{ R}_{ j',\, \beta} '}{\partial w_{ 1;j}} 
\equiv
& \
\sum_{ j_1 '= 1}^{ d'} \,
\sum_{ \beta_1 \leq \beta} \left( 
\frac{ \partial
A_{ j',\, \beta}^{ j_1',\, \beta_1}}{\partial
w_{ 1;j}} \cdot
S_{ j_1',\, \beta_1}'+ 
A_{ j',\, \beta}^{ j_1',\, \beta_1} \cdot
\frac{ \partial S_{ j_1',\, \beta_1}'}{\partial w_{ 1;j}}\right),
\\
\frac{ \partial S_{ j',\, \beta} '}{
\partial w_{ 1;j}} \equiv 
& \
\sum_{ j_1 '= 1}^{ d'} \,
\sum_{ \beta_1 \leq \beta} \left( 
\frac{ \partial
B_{ j',\, \beta}^{ j_1',\, \beta_1}}{\partial
w_{ 1;j}} \cdot
\overline{ R}_{ j_1',\, \beta_1}'+ 
B_{ j',\, \beta}^{ j_1',\, \beta_1} \cdot
\frac{ \partial \overline{ R}_{ j_1',\, \beta_1}'}{\partial w_{ 1;j}}
\right).
\endaligned\right.
\end{equation}
Ensuite, appliquons l'op\'erateur $\sum_{ i_1' =1}^{ n'} \, \partial_{
t_{i_1'}'} (\cdot) \, T_{ i_1',\, j}'$ de part
et d'autre de~\thetag{ 7.117}, ce qui donne~:
\def\theequation{7.119}\begin{equation}
\left\{
\aligned
\sum_{ i_1'= 1}^{ n'} \,
\frac{ \partial 
\overline{ R}_{ j',\, \beta} '}{
\partial t_{ i_1'}'} \cdot T_{ i_1',\, j}' \equiv
& \
\sum_{ i_1'= 1}^{ n'} \,
\sum_{ j_1 '= 1}^{ d'} \,
\sum_{ \beta_1 \leq \beta} 
\left(
\frac{ \partial
A_{ j',\, \beta}^{ j_1',\, \beta_1}}{
\partial t_{ i_1'}'} \cdot T_{ i_1',\, j}'
\cdot
S_{ j_1',\, \beta_1}'+ \right. \\ 
& \
\ \ \ \ \ \ \ \ \ \ \ \ \ \ \ \
\ \ \ \ \ \ \ \ \ \ \ \ \ \ \ \
\ \ \ 
\left.
+ 
A_{ j',\, \beta}^{ j_1',\, \beta_1} \cdot
\frac{ \partial 
S_{ j_1',\, \beta_1}'}{
\partial t_{ i_1'}'} \cdot T_{ i_1',\, j}' 
\right) 
\\
\sum_{ i_1 ' =1}^{ n'}\,
\frac{ \partial 
S_{ j',\, \beta} '}{
\partial t_{ i_1'}'} \cdot T_{ i_1',\, j}' \equiv 
& \
\sum_{ i_1 ' =1}^{ n'}\,
\sum_{ j_1 '= 1}^{ d'} \,
\sum_{ \beta_1 \leq \beta} \left(
\frac{ \partial 
B_{ j',\, \beta}^{ j_1',\, \beta_1}}{
\partial t_{ i_1'}'} \cdot T_{ i_1',\, j}' \cdot
\overline{ R}_{ j_1',\, \beta_1}'+ \right. \\ 
& \
\ \ \ \ \ \ \ \ \ \ \ \ \ \ \ \
\ \ \ \ \ \ \ \ \ \ \ \ \ \ \ \
\ \ \
\left.
+ 
B_{ j',\, \beta}^{ j_1',\, \beta_1} \cdot
\frac{ \partial 
\overline{ R}_{ j_1',\, \beta_1}'}{
\partial t_{ i_1'}'} \cdot T_{ i_1',\, j}' 
\right).
\endaligned\right.
\end{equation}
Additionnons maintenant~\thetag{ 7.118} avec~\thetag{ 7.119}~; en
tenant compte de la d\'efinition~\thetag{ 7.106}, nous
obtenons l'expression d\'esir\'ee~\thetag{ 7.115}, pour 
$\delta = \1_j^d$~:
\def\theequation{7.120}\begin{equation}
\left\{
\aligned
\overline{ R}_{ j',\, \beta,\, \1_j^d} ' \equiv 
\sum_{ j_1 '= 1}^{ d'} \, 
\sum_{ \beta_1 \leq \beta} 
\left( 
A_{ j',\, \beta,\, \1_j^d}^{ j_1',\, \beta_1,\, \1_j^d} \cdot
S_{ j_1',\, \beta_1,\, \1_j^d} ' +
A_{ j',\, \beta,\, \1_j^d}^{ j_1',\, \beta_1,\, 0} \cdot
S_{ j_1',\, \beta_1,\, 0}'
\right),
\\
S_{ j',\, \beta,\, \1_j^d} ' \equiv 
\sum_{ j_1 '= 1}^{ d'} \, 
\sum_{ \beta_1 \leq \beta} 
\left(
B_{ j',\, \beta}^{ j_1',\, \beta_1} \cdot
\overline{ R}_{ j_1',\, \beta_1,\, \1_j^d} ' + 
B_{ j',\, \beta,\, \1_j^d}^{ j_1',\, \beta_1,\, 0} \cdot
\overline{ R}_{ j_1',\, \beta_1,\, 0}'
\right),
\endaligned\right.
\end{equation}
o\`u nous avons pos\'e~:
\def\theequation{7.121}\begin{equation}
\left\{
\aligned
A_{ j',\, \beta,\, \1_j^d}^{ j_1',\, \beta_1,\, \1_j^d} :=
& \ 
A_{ j',\, \beta}^{ j_1',\, \beta_1}, \\
A_{ j',\, \beta,\, \1_j^d}^{ j_1',\, \beta_1,\, 0} :=
& \
\frac{ \partial A_{ j',\, \beta}^{ j_1',\, \beta_1}}{
\partial w_{ 1;j}} +
\sum_{ i_1 '= 1}^{ n'} 
\frac{ \partial A_{ j',\, \beta}^{ j_1,\, \beta_1}}{
\partial t_{ i_1'}'} \cdot T_{ i_1',\, j}',
\\
B_{ j',\, \beta,\, \1_j^d}^{ j_1',\, \beta_1,\, \1_j^d} :=
& \ 
B_{ j',\, \beta}^{ j_1',\, \beta_1}, \\
B_{ j',\, \beta,\, \1_j^d}^{ j_1',\, \beta_1,\, 0} :=
& \
\frac{ \partial B_{ j',\, \beta}^{ j_1',\, \beta_1}}{
\partial w_{ 1;j}} +
\sum_{ i_1 '= 1}^{ n'} 
\frac{ \partial B_{ j',\, \beta}^{ j_1,\, \beta_1}}{
\partial t_{ i_1'}'} \cdot T_{ i_1',\, j}'.
\endaligned\right.
\end{equation}
En appliquant des d\'erivations d'ordre arbitraire par rapport \`a
$w_1$, en utilisant les formules~\thetag{ 7.112} et en raisonnant par
r\'ecurrence, on \'etablit les identit\'es~\thetag{ 7.115}. Puisque
nous avons expliqu\'e le principe de la d\'emonstration pour $\delta =
\1_j^d$, nous ne croyons pas utile de d\'evelopper ces
calculs formels suppl\'ementaires.

Le Lemme~7.110 est d\'emontr\'e.
\endproof

\subsection*{7.122.~Convergence des 
composantes de l'application de r\'eflexion conjugu\'ee sur la
deuxi\`eme cha\^{\i}ne de Segre} Nous savons maintenant que pour tout
$\ell \in\N$, les propri\'et\'es de convergence~:
\def\theequation{7.123}\begin{equation}
\left[
J_t^\ell \Theta_{ j',\, \gamma'}'
(h) \right] \left(
\Gamma_1 (z_1)
\right)\in \C\{ z_1\}^{ N_{ d',\, n,\, \ell}},
\end{equation} 
sont satisfaites, pour tout $j' = 1,\, \dots,\, d'$ et 
tout $\gamma' \in
\N^{ m'}$. Gr\^ace \`a cette propri\'et\'e, nous pouvons \'etablir le
Lemme~7.18 pour $k=2$ et $\ell =0$, c'est-\`a-dire~:
\def\theequation{7.124}\begin{equation}
\left[
\overline{ \Theta}_{ j',\, \gamma'}'
(\overline{ h} ) \right] \left(
\Gamma_2 (z_{(2)})
\right)\in \C\{ z_{ (2)}\}^{ d'},
\end{equation}
pour tout $j' = 1,\, \dots,\, d'$ et tout $\gamma' \in \N^{ m'}$.

\proof
Posons $(t,\, \tau) = \Gamma_2 (z_{ (2)})$ dans la quatri\`eme
famille d'identit\'es de r\'eflexion~\thetag{ 7.62}~:
\def\theequation{7.125}\begin{equation}
\left\{
\aligned
\overline{ g}_{ j'} \left(\Gamma_2 (z_{ (2)})\right) \equiv
& \
\sum_{ \gamma ' \in \N^{ m'}} \, 
\left[
\overline{ f}^{\gamma'} \right] \left(
\Gamma_2 (z_{ (2)})
\right)\cdot
\left[
\Theta_{ j',\, \gamma'} '
(h)
\right] \left(
\Gamma_2 (z_{ (2)}) \right), \\
0 \equiv
& \
\sum_{ \gamma ' \in \N^{ m'}} \, 
\left[
\overline{ f}^{\gamma'} \right] \left(
\Gamma_2 (z_{ (2)})
\right)\cdot
\left[
\mathcal{ L}^\beta \Theta_{ j',\, \gamma'} '
(h)
\right] \left(
\Gamma_2 (z_{ (2)}) \right).
\endaligned\right.
\end{equation}
Dans la deuxi\`eme ligne, on suppose $\vert \beta \vert \geq
1$. D'apr\`es l'identit\'e conjugu\'ee de la troisi\`eme ligne
de~\thetag{ 7.81}, \'ecrite pour $\delta =0$ et $\beta \in \N^m$
arbitraire, nous voyons que tous les termes
\def\theequation{7.126}\begin{equation}
\left\{
\aligned
\left[
\mathcal{ L}^\beta \Theta_{ j',\, \gamma'} '
(h)
\right] \left(
\Gamma_2 (z_{ (2)}) \right) \equiv
& \
P_{ \beta,\, 0}' 
\left(
J_{ z,\, \tau}^{ \vert \beta \vert} \overline{ \Theta}
\left(\Gamma_2 (z_{(2)})\right), \
\left[
J_t^{ \vert \beta \vert} 
\Theta_{ j',\, \gamma'} ' 
\left(
h
\right)
\right]\left(
\Gamma_2 (z_{ (2)})
\right)
\right) \\
\equiv 
& \
P_{ \beta,\, 0}' 
\left(
J_{ z,\, \tau}^{ \vert \beta \vert} \overline{ \Theta}
\left(\Gamma_2 (z_{(2)})\right), \ 
\left[
J_t^{ \vert \beta \vert}
\Theta_{ j',\, \gamma'} ' 
\left(
h
\right)
\right]\left(
\Gamma_1 (z_1)
\right)
\right).
\endaligned\right.
\end{equation}
sont convergents par rapport \`a $(z_1, \, z_2 ) \in \C^{ 2m}$,
gr\^ace \`a~\thetag{ 7.21} et gr\^ace ce que nous venons d'\'etablir
en \thetag{ 7.123}. Par cons\'equent, les \'equations~\thetag{ 7.125}
sont analytiques par rapport \`a $(z_{ (2)},\, \tau')$, et
l'application formelle $z_{ (2)} \longmapsto_{ \mathcal{ F}} \,
\overline{ h} \left( \Gamma_2 (z_{ (2)}) \right)$ les satisfait
formellement. Appliquons le Th\'eor\`eme d'approximation~\thetag{ 2.4}
avec l'ordre d'approximation $N=1$. Nous en d\'eduisons l'existence
d'une application convergente $z_{ (2)} \longmapsto \overline{ {\sf
H}} (z_{(2)}) = \left( \overline{ {\sf F} } (z_{(2)}),\, \overline{
{\sf G} } (z_{(2)}) \right) \in \C\{ z_{(2)}\}^{ m'} \times \C \{
z_{(2)} \}^{ d'}$ telle que
\def\theequation{7.127}\begin{equation}
\left\{
\aligned
\overline{ {\sf G} }_{ j'} (z_{ (2)}) \equiv
& \
\sum_{ \gamma ' \in \N^{ m'}} \, 
\left[
\overline{ {\sf F} }^{\gamma'} \right] (z_{ (2)})
\cdot
\Theta_{ j',\, \gamma'} '
(h (
\Gamma_1 (z_1) ), \\
0 \equiv
& \
\sum_{ \gamma ' \in \N^{ m'}} \, 
\left[
\overline{ {\sf F} }^{\gamma'} \right] (z_{ (2)})
\cdot
\left[
\mathcal{ L}^\beta \Theta_{ j',\, \gamma'} '
(h)
\right] \left(
\Gamma_2 (z_{ (2)}) \right).
\endaligned\right.
\end{equation}
Les \'equations~\thetag{ 7.127} co\"{\i}ncident avec $\mathcal{
L}^\beta r_{ j'}' \left( h \left( \Gamma_2 (z_{ (2)})\right),\,
\overline{ {\sf H}} (z_{ (2)}) \right) \equiv 0$, pour tout $j'= 1,\,
\dots,\, d'$ et tout $\beta \in \N^m$. En appliquant la conjugaison de
l'\'equivalence~\thetag{ 7.58}, nous pouvons faire basculer ces
identit\'es vers $\mathcal{ L}^\beta \overline{ r}_{ j'}' \left(
\overline{ {\sf H}} (z_{ (2)}),\, h \left( \Gamma_2 (z_{ (2)})\right)
\right) \equiv 0$, pour tout $j'= 1,\, \dots,\, d'$ et tout $\beta \in
\N^m$, ou, de mani\`ere d\'evelopp\'ee~:
\def\theequation{7.128}\begin{equation}
\left[
\mathcal{ L}^\beta g_{ j'} \right] \left(
\Gamma_2 ( z_{ (2)})
\right) \equiv
\sum_{ \gamma ' \in \N^{ m'}} \, 
\left[
\mathcal{ L}^\beta f^{ \gamma'}
\right]
\left(
\Gamma_2 (z_{ (2)})
\right)
\cdot
\overline{ \Theta}_{ j',\, \gamma'} ' \left(
\overline{ {\sf H}} (z_{ (2)})
\right),
\end{equation}
pour tout $j'= 1,\, \dots,\, d'$ et tout $\beta \in \N^m$.
Soustrayons ces identit\'es des identit\'es de r\'eflexion \'ecrites
\`a la troisi\`eme ligne de~\thetag{ 7.62} avec $(t,\, \tau) =\Gamma_2
( z_{ (2)})$~; nous obtenons~:
\def\theequation{7.129}\begin{equation}
0 \equiv
\sum_{ \gamma ' \in \N^{ m'}} \, 
\left[
\mathcal{ L}^\beta f^{ \gamma'}\right]
\left(
\Gamma_2 (z_{ (2)})
\right) \cdot
\left(
\overline{ \Theta}_{ j',\, \gamma'} ' \left(
\overline{ h} \left(
\Gamma_2 ( z_{ (2)})
\right)
\right) -
\overline{ \Theta}_{ j',\, \gamma'} ' \left(
\overline{ {\sf H} } (z_{ (2)})
\right)
\right).
\end{equation}
G\'en\'eralisons maintenant le principe d'unicit\'e~\thetag{ 7.37}.
Soit $k\in\N$. Pour tout $\gamma ' \in \N^{ m'}$, supposons donn\'ee
une s\'erie formelle vectorielle $\overline{ F}_{ \gamma '} (z_{ (k)})
\in \C \dl z_{ (k)} \dr^{ d'}$ et supposons que l'application formelle
$h$ est CR-transversale, ce qui se traduit par 
l'implication~\thetag{ 7.27}.

\def\thelemma{7.130}\begin{lemma}
Supposons que pour tout $\beta \in \N^m$, 
l'identit\'e formelle suivante est satisfaite~{\rm :}
\def\theequation{7.131}\begin{equation}
\left\{
\aligned
0 \equiv 
& \
\sum_{ \gamma ' \in \N^{ m'}} \, 
\left[
\underline{ \mathcal{ L}}^\beta 
\overline{ f}^{ \gamma'}
\right]\left(
\Gamma_k (z_{ (k)})
\right) \cdot
\overline{ F}_{\gamma '} ( z_{ (k)})
\ \ \ \ \ \ \
\text{ \sf si} \ k \ 
\text{ \sf est impair}, \\
0 \equiv 
& \
\sum_{ \gamma ' \in \N^{ m'}} \, 
\left[
\mathcal{ L}^\beta 
f^{ \gamma'}
\right]\left(
\Gamma_k (z_{ (k)})
\right) \cdot
\overline{ F}_{\gamma '} ( z_{ (k)})
\ \ \ \ \ \ \
\text{ \sf si} \ k \ 
\text{ \sf est pair},
\endaligned\right.
\end{equation}
dans $\C \dl z_{ (k)} \dr$. Alors~{\rm :} 
\def\theequation{7.132}\begin{equation}
\overline{ F}_{ \gamma '} ( z_{ (k)}) \equiv 0, 
\ \ \ \ \ \ \
\forall \, \gamma' \in \N^{ m'}.
\end{equation}
\end{lemma}

\proof
Traitons seulement le cas o\`u $k$ est pair. Puisque $\left[ \mathcal{
L}^\beta f^{ \gamma'}\right] (0) = f_{ \gamma',\, \beta}$ d'apr\`es
\thetag{ 7.29}, et puisque $\Gamma_k (0)= 0$, il existe des
coefficients $f_{ \gamma',\, \beta,\, \gamma_k} \in \C$ tels que l'on
peut \'ecrire~:
\def\theequation{7.133}\begin{equation}
\left[
\mathcal{ L}^\beta f^{ \gamma'} \right]
\left(
\Gamma_k ( z_{ (k)})
\right) \equiv
f_{ \gamma ',\, \beta} + 
\sum_{ \gamma_k \in \N^{ km}, \, 
\vert \gamma_k \vert \geq 1} \, 
f_{ \gamma ',\, \beta,\, \gamma_k} \cdot
z_{ (k)}^{ \gamma_k}.
\end{equation}
D\'eveloppons aussi $\overline{ F}_{ \gamma '} (z_{ (k)}) = \sum_{
\gamma_k \in \N^{ km}} \, \overline{ F}_{ \gamma',\, \gamma_k} \cdot
z_{ (k)}^{ \gamma_k}$. En raisonnant exactement comme dans la
d\'emonstration du Lemme~7.36, on obtient l'annulation de tous les
coefficients $\overline{ F}_{ \gamma ',\, \gamma_k}$, ce qui ach\`eve
la preuve.
\endproof

Une application directe de ce lemme aux identit\'es~\thetag{ 7.129}
fournit les identit\'es
\def\theequation{7.134}\begin{equation}
\overline{ \Theta}_{ j',\, \gamma'} ' \left(
\overline{ h} \left(
\Gamma_2 ( z_{ (2)})
\right)
\right) \equiv
\overline{ \Theta}_{ j',\, \gamma'} ' \left(
\overline{ {\sf H} } (z_{ (2)})
\right),
\end{equation}
pour tout $j'= 1,\, \dots,\, d'$ et tout $\gamma ' \in \N^{ m'}$, o\`u
le membre de droite est convergent, ce qui ach\`eve la d\'emonstration
de~\thetag{ 7.124}.
\endproof

\subsection*{ 7.135.~D\'emonstration finale}
Gr\^ace \`a tous ces pr\'eparatifs, nous pouvons enfin pr\'esenter la
d\'emonstration de la Proposition~7.18, en raisonnant par r\'ecurrence
sur la longueur $k$ des cha\^{\i}nes de Segre. Nous traiterons
seulement le cas o\`u $k$ est impair\,--\,le cas o\`u $k$ est pair
\'etant similaire.

Ainsi, nous supposons que pout tout $\ell \in \N$, les propri\'et\'es
de convergence~:
\def\theequation{7.136}\begin{equation}
\left[
J_t^\ell
\Theta_{ \gamma'} ' (h) \right]
\left(
\Gamma_k (z_{ (k)})
\right) \in \C\{ z_{ (k)}\}^{ N_{ d',\, n,\, \ell}}, 
\end{equation}
sont satisfaites, pour tout $\gamma ' \in \N^{ m'}$. L'objectif est de
d\'emontrer que pour tout $\ell \in \N$, les propri\'et\'es de
convergence au rang $k+1$~:
\def\theequation{7.137}\begin{equation}
\left[
J_\tau^\ell
\overline{ \Theta}_{ \gamma'} \left(
\overline{ h}
\right)
\right] \left(
\Gamma_{ k+1} (z_{ (k+1)})
\right) \in 
\C\{ z_{ (k+1)}\}^{ N_{ d',\, n,\, \ell}},
\end{equation}
sont satisfaites, pour tout $\gamma ' \in \N^{ m'}$.

Pour cela, commen\c cons par poser $(t,\, \tau) = \underline{
\Upsilon}_\xi \left( \Gamma_{ k+1} (z_{ (k+1)}) \right)$ dans la
quatri\`eme ligne de~\thetag{ 7.62}, o\`u $\xi \in\C^d$, sans oublier
d'\'ecrire l'identit\'e qui est sous-entendue pour $\beta =0$~:
\def\theequation{7.138}\begin{equation}
\left\{
\aligned
\overline{ g}_{ j'} 
\left(
\underline{ \Upsilon}_\xi \left(
\Gamma_{ k+1} ( z_{ (k+1)})
\right)
\right) \equiv
& \
\sum_{ \gamma ' \in\N^{ m'}} \,
\left[
\overline{ f}^{ \gamma'} 
\right]
\left(
\underline{ \Upsilon}_\xi \left(
\Gamma_{ k+1} ( z_{ (k+1)})
\right)
\right) \cdot \\
& \
\ \ \ \ \ \ \ \ 
\cdot
\Theta_{ j',\, \gamma'} '
\left( h
\left(
\underline{ \Upsilon}_\xi \left(
\Gamma_{ k+1} ( z_{ (k+1)})
\right)
\right)
\right), \\
0 \equiv 
& \
\sum_{ \gamma ' \in\N^{ m'}} \,
\left[
\overline{ f}^{ \gamma'} 
\right]
\left(
\underline{ \Upsilon}_\xi \left(
\Gamma_{ k+1} ( z_{ (k+1)})
\right)
\right) \cdot \\
& \
\ \ \ \ \ \ \ \ 
\cdot
\left[
\mathcal{ L}^\beta \Theta_{ j',\, \gamma'}' 
(h)
\right]\left(
\underline{ \Upsilon}_\xi \left(
\Gamma_{ k+1} ( z_{ (k+1)})
\right)
\right).
\endaligned\right.
\end{equation}
Dans la deuxi\`eme identit\'e, on suppose
$\vert \beta \vert \geq 1$. 

Fixons maintenant $\ell \in \N$. Pour d\'emontrer~\thetag{ 7.137},
appliquons les d\'erivations $\left. \partial_{ \xi}^\delta
\right\vert_{ \xi=0}$ aux identit\'es~\thetag{ 7.138} pour tout
multiindice $\delta \in \N^d$ tel que $\vert \delta \vert \leq \ell$,
ce qui donne~:
\def\theequation{7.139}\begin{equation}
\left\{
\aligned
\left[
\underline{ \Upsilon}^\delta \overline{ g}_{ j'}
\right] 
\left(
\Gamma_{ k+1} ( z_{ (k+1)})
\right) \equiv
& \
\sum_{ \gamma ' \in\N^{ m'}} \, 
\sum_{ \delta_1 \leq \delta}\,
\frac{ \delta! }{\delta_1 ! \ (\delta -\delta_1)!} \,
\left[
\underline{ \Upsilon}^{ \delta - \delta_1} \overline{ f}^{ \gamma'}
\right] 
\left(
\Gamma_{ k+1} ( z_{ (k+1)})
\right) \cdot \\
& \ 
\ \ \ \ \ \ \ \ \ \ 
\cdot
\left[
\underline{ \Upsilon}^{ \delta_1} \Theta_{ j',\, \gamma'} '
(h)
\right]
\left(
\Gamma_{ k+1} ( z_{ (k+1)})
\right), \\
0 \equiv
& \
\sum_{ \gamma ' \in\N^{ m'}} \, 
\sum_{ \delta_1 \leq \delta}\,
\frac{ \delta! }{\delta_1 ! \ (\delta -\delta_1)!} \,
\left[
\underline{ \Upsilon}^{ \delta - \delta_1} \overline{ f}^{ \gamma'}
\right] 
\left(
\Gamma_{ k+1} ( z_{ (k+1)})
\right) \cdot \\
& \ 
\ \ \ \ \ \ \ \ \ \ 
\cdot
\left[
\underline{ \Upsilon}^{ \delta_1} \mathcal{L}^\beta
\Theta_{ j',\, \gamma'} ' ( h)
\right]
\left(
\Gamma_{ k+1} ( z_{ (k+1)})
\right), 
\endaligned\right.
\end{equation}
pour tout $j'=1,\, \dots,\, d'$, tout
$\beta \in \N^m$ tel que $\vert \beta \vert \geq 1$ et
tout $\delta \in\N^d$ tel que $\vert \delta \vert \leq \ell$.
V\'erifions que tous les termes \'ecrits \`a la
deuxi\`eme et \`a la quatri\`eme ligne de~\thetag{ 7.139} 
sont des s\'eries convergentes par rapport \`a 
$z_{ (k+1)}$. En effet, en rempla\c cant
$(t,\, \tau)$ par $\Gamma_{ k+1} (z_{ (k+1)})$ dans
l'identit\'e conjugu\'ee de~\thetag{ 7.82}, nous
obtenons~:
\def\theequation{7.140}\begin{equation}
\left\{
\aligned
\left[
\underline{ \Upsilon}^{ \delta_1} 
\mathcal{L}^\beta
\Theta_{ j',\, \gamma'} ' ( h)
\right]
\left(
\Gamma_{ k+1} (z_{ (k+1)})
\right) \equiv
& \
P_{ \beta,\, \delta_1} '
\left(
J_{ z,\, \tau}^{ \vert \beta \vert+ 
\vert \delta_1 \vert} \overline{ \Theta}
\left(
\Gamma_{ k+1} (z_{ (k+1)})
\right), \right. \\
& \
\ \ \ \ \ \ \ \ \ \
\left.
\left[
J_t^{ \vert \beta \vert + \vert \delta_1 \vert}
\Theta_{ j',\, \gamma'} '( h)
\right] \left(
\Gamma_k (z_ {(k)})
\right)
\right), 
\endaligned\right.
\end{equation}
et tous ces termes sont convergents, gr\^ace \`a l'hypoth\`ese de
r\'ecurrence~\thetag{ 7.136}. Appliquons maintenant le
Th\'eor\`eme~2.4 avec l'ordre d'approximation $N=1$~: nous en
d\'eduisons qu'il existe une application convergente
\def\theequation{7.141}\begin{equation}
\left\{
\aligned
z_{(k+1)} \longmapsto 
& \
\left(
\overline{{\sf H}}_{ i_1',\, \delta}(z_{ (k+1)}) 
\right)_{ 1\leq i_1 ' \leq n',\, \vert \delta \vert \leq \ell}
\\
& \
=:
\left(
\left(
\overline{ {\sf F}}_{ k',\, \delta} (z_{ (k+1)})
\right)_{1\leq k'\leq m',\, \vert \delta \vert \leq \ell}, \
\left(
\overline{ {\sf G}}_{ j',\, \delta} (z_{ (k+1)})
\right)_{1\leq j'\leq d',\, \vert \delta \vert \leq \ell}
\right),
\endaligned\right.
\end{equation}
avec ${\sf H}_{ i_1',\, \delta_1} (0) = \Upsilon^{ \delta_1} h_{ i_1'}
(0)$, qui satisfait les \'equations analytiques~\thetag{ 7.139},
c'est-\`a-dire telle que~:
\def\theequation{7.142}\begin{equation}
\left\{
\aligned
\overline{ {\sf G}}_{ j',\, \delta}(z_{ (k+1)})
& 
\equiv
\sum_{ \gamma' \in\N^{ m'}} \, 
\sum_{ \delta_1 \leq \delta} \, 
\frac{ \delta !}{ \delta_1 ! \ 
(\delta -\delta_1)!} \ 
\left[
\underline{ \Upsilon}^{ 
\delta_1} \Theta_{ j',\, \gamma'}'
(h) \right]
\left(
\Gamma_{k+1} (z_{ (k+1)})
\right)
\cdot \\
& \
\ \ \ \ \ \ \ \ \ \ \ \ \
\ \ \ \ \ \ \ \ \ \ \ \ \ 
\ \ \ \ \ \ \ \
\cdot 
\overline{ N}_{\gamma',\, \delta-\delta_1}
\left(
\left(
\overline{ {\sf F}}_{ k',\, \delta_2}(z_{ (k+1)})
\right)_{ 1\leq k' \leq m',\, \delta_2 \leq \delta -\delta_1}
\right), 
\\
0
& 
\equiv
\sum_{ \gamma' \in\N^{ m'}} \, 
\sum_{ \delta_1 \leq \delta} \, 
\frac{ \delta !}{ \delta_1 ! \ 
(\delta -\delta_1)!} \
\left[
\underline{ \Upsilon}^{ \delta_1}\mathcal{ L}^\beta 
\Theta_{ j',\, \gamma'}'
(h) \right]
\left(
\Gamma_{k+1} (z_{ (k+1)})
\right)
\cdot \\
& \
\ \ \ \ \ \ \ \ \ \ \ \ \
\ \ \ \ \ \ \ \ \ \ \ \ \ 
\ \ \ \ \ \ \ \
\cdot 
\overline{ N}_{\gamma',\, \delta-\delta_1}
\left(
\left(
\overline{ {\sf F}}_{ k',\, \delta_2}(z_{ (k+1)})
\right)_{ 1\leq k' \leq m',\, \delta_2 \leq \delta -\delta_1}
\right),
\endaligned\right.
\end{equation}
pour tout $j'=1,\, \dots,\, d'$, tout $\beta \in \N^m$ satisfaisant
$\vert \beta \vert \geq 1$ et tout $\delta \in\N^d$ tel que $\vert
\delta \vert \leq \ell$. Utilisons maintenant
l'\'equivalence~\thetag{ 7.111} avec ${\sf x} := z_{ (k+1)}$, avec
${\sf Q}( {\sf x}):= \Gamma_{ k+1} (z_{ (k+1)})$, avec ${\sf T}_{ i_1
',\, \delta_1}' ({\sf x}) := \overline{ {\sf H}}_{ 
i_1',\, \delta_1 } ({\sf x})$ et
pour tout $\delta \in \N^d$ tel que $\vert \delta \vert \leq \ell$~:
nous en d\'eduisons que la m\^eme application formelle~\thetag{ 7.141}
satisfait les \'equations analytiques d\'evelopp\'ees suivantes~:
\def\theequation{7.143}\begin{equation}
\left\{
\aligned
\left[
\underline{ \Upsilon}^\delta
\mathcal{ L}^\beta g_{ j'}
\right] 
\left(
\Gamma_{ k+1} (z_{ (k+1)})
\right) \equiv
& \
\sum_{ \gamma ' \in\N^{ m'}} \, 
\sum_{ \delta_1 \leq \delta} \, 
\frac{ \delta ! }{\delta_1 ! \ 
(\delta- \delta_1)!} \cdot 
\\
& \
\ \ \ \ \
\cdot
\left[
\underline{ \Upsilon}^{ \delta- \delta_1}
\mathcal{ L}^\beta
f^{ \gamma'} 
\right]
\left(
\Gamma_{ k+1} (z_{ (k+1)})
\right)\cdot 
\\
& \
\ \ \ \ \
\cdot
\overline{ M}_{ \delta_1}' 
\left(
\left(
\overline{ {\sf H}}_{ i_1',\, \delta_2} ( z_{ (k+1)})
\right)_{ 1\leq i_1'\leq n',\, \delta_2 \leq \delta_1}, 
\right.
\\
& \
\ \ \ \ \ \ \ \ \ \ \ \ \ \ \ \ \
\ \ \ \ \ \ \ \ \ 
\left.
\left[
J_{\tau'}^{ \vert \delta_1 \vert}
\overline{ \Theta}_{ j',\, \gamma'} '
\right] 
\left(
\overline{ {\sf H}}_0 (z_{ (k+1)})
\right)
\right), 
\endaligned\right.
\end{equation}
pour tout $j'= 1,\, \dots,\, d'$, tout $\beta \in\N ^m$ et tout
$\delta \in \N^d$ tel que $\vert \delta \vert \leq \ell$, o\`u nous
avons pos\'e $\overline{ {\sf H}}_0 := 
\left( \overline{ {\sf H}}_{ i_1',\, 0} \right)_{ 1\leq
i_1' \leq n'}$. Comparons ces identit\'es aux identit\'es de la
troisi\`eme ligne de~\thetag{ 7.62} dans lesquelles on pose $(t,\,
\tau) = \underline{ \Upsilon}_\xi \left( \Gamma_{ k+1} (z_{ (k+1)})
\right)$ et auxquelles on applique la d\'erivation $\left.
\partial_\xi^\delta ( \cdot ) \right\vert_{ \xi= 0}$, ce qui donne~:
\def\theequation{7.144}\begin{equation}
\left\{
\aligned
\left[
\underline{ \Upsilon}^\delta
\mathcal{ L}^\beta g_{ j'}
\right] 
\left(
\Gamma_{ k+1} (z_{ (k+1)})
\right) \equiv
& \
\sum_{ \gamma ' \in\N^{ m'}} \, 
\sum_{ \delta_1 \leq \delta} \, 
\frac{ \delta ! }{\delta_1 ! \ 
(\delta- \delta_1)!} \cdot 
\\
& \
\ \ \ \ \
\cdot
\left[
\underline{ \Upsilon}^{ \delta- \delta_1}
\mathcal{ L}^\beta
f^{ \gamma'} 
\right]
\left(
\Gamma_{ k+1} (z_{ (k+1)})
\right)\cdot 
\\
& \
\ \ \ \ \
\cdot
\overline{ M}_{ \delta_1}' 
\left(
\left(
\underline{ \Upsilon}^{ \delta_2} \overline{ h}_{ i_1'}
\left(
\Gamma_{ k+1} (z_{ (k+1)})
\right)
\right)_{ 1\leq i_1'\leq n',\, \delta_2 \leq \delta_1}, 
\right.
\\
& \
\ \ \ \ \ \ \ \ \ \ \ \ \ \ \ \ \
\ \ \ \ \ \ \ \ \ 
\left.
\left[
J_{\tau'}^{ \vert \delta_1 \vert}
\overline{ \Theta}_{ j',\, \gamma'} '
\left(
\overline{ h} \right)
\right] \left( \Gamma_{ k+1} (z_{ (k+1)})\right)
\right), 
\endaligned\right.
\end{equation}
Soustrayons~\thetag{ 7.143} de~\thetag{ 7.144}, ce qui donne~:
\def\theequation{7.145}\begin{equation}
\left\{
\aligned
0 \equiv
& \
\sum_{ \gamma ' \in \N^{ m'}} \, 
\sum_{ \delta_1 \leq \delta} \, 
\frac{ \delta !}{\delta_1 ! \ 
(\delta- \delta_1)!} \, 
\cdot
\left[
\underline{ \Upsilon}^{ \delta- \delta_1}
\mathcal{ L}^\beta
f^{ \gamma'} 
\right]
\left(
\Gamma_{ k+1} (z_{ (k+1)})
\right)\cdot 
\\
& \
\ \ \ \ \
\left(
\aligned
& \
\overline{ M}_{ \delta_1}' 
\left(
\left(
\underline{ \Upsilon}^{ \delta_2} \overline{ h}_{ i_1'}
\left(
\Gamma_{ k+1} (z_{ (k+1)})
\right)
\right)_{ 1\leq i_1'\leq n',\, \delta_2 \leq \delta_1}, 
\right.
\\
& \
\ \ \ \ \ \ \ \ \ \ \ \ \ \ \ \ \
\ \ \ \ \ \ \ \ \ 
\left.
\left[
J_{\tau'}^{ \vert \delta_1 \vert}
\overline{ \Theta}_{ j',\, \gamma'} ' 
\left(
\overline{ h} \right)
\right] \left( \Gamma_{ k+1} (z_{ (k+1)})\right)
\right) -
\\
& \
-
\overline{ M}_{ \delta_1}' 
\left(
\left(
\overline{ {\sf H}}_{ i_1',\, \delta_2} ( z_{ (k+1)})
\right)_{ 1\leq i_1'\leq n',\, \delta_2 \leq \delta_1}, 
\right.
\\
& \
\ \ \ \ \ \ \ \ \ \ \ \ \ \ \ \ \
\ \ \ \ \ \ \ \ \ 
\left.
\left[
J_{\tau'}^{ \vert \delta_1 \vert}
\overline{ \Theta}_{ j',\, \gamma'} ' 
\right] 
\left(
\overline{ {\sf H}}_0 (z_{ (k+1)})
\right)
\right) 
\endaligned
\right), 
\endaligned\right.
\end{equation}
pour tout $j'=1,\, \dots,\, d'$, tout $\beta \in \N^m$ et tout $\delta
\in \N^d$ tel que $\vert \delta \vert \leq \ell$. D'apr\`es la
d\'efinition m\^eme de $\overline{ M}_{ \delta_1}'$ donn\'ee \`a la
quatri\`eme ligne de~\thetag{ 7.81}), le terme d\'evelopp\'e \`a la
deuxi\`eme et \`a la troisi\`eme ligne de~\thetag{ 7.145} s'identifie
\`a $\left[ \underline{ \Upsilon}^{ \delta_1} \overline{ \Theta}_{
j',\, \gamma'}' \left( \overline{ h} \right) \right] \left( \Gamma_{
k+1} (z_{ (k+1)}) \right)$. Notons que le terme d\'evelopp\'e \`a la
quatri\`eme et \`a la cinqui\`eme ligne de~\thetag{ 7.145} est
convergent. Ainsi, pour conclure que les propri\'et\'es de
convergence~\thetag{ 7.137} sont satisfaites (en utilisant au passage
le Lemme~7.74), il ne nous reste plus qu'\`a \'etablir que les termes
contenus dans les grandes parent\`eses s'annulent tous identiquement.

\'Ecrivons~\thetag{ 7.145} pour $\delta =0$~: la somme $\sum_{
\delta_1 \leq \delta}$ contient alors un seul terme et le Lemme~7.130
s'applique. En raisonnant par r\'ecurrence sur la longueur de
$\delta$, jusqu'\`a la longueur $\ell$, en \'ecrivant les
identit\'es~\thetag{ 7.145} pour $\vert \delta \vert$ croissant, et en
appliquant le Lemme~7.130 \`a chaque \'etape, nous v\'erifions que les
termes contenus dans les grandes parent\`eses s'annulent tous
identiquement.

Les d\'emonstrations de la Proposition~7.18 et du Th\'eor\`eme~1.23
sont achev\'ees. 
\qed

\subsection*{ 7.146.~Existence d'applications holomorphes}
Pour terminer cet article, d\'emontrons le corollaire suivant du
Th\'eor\`eme~1.23, qui implique le Corollaire~1.45.

\def\thecorollary{7.147}\begin{corollary}
Sous les hypoth\`eses du Th\'eor\`eme~1.23, pour tout entier $N \geq
1$, il existe une application convergente ${\sf H}^N (t) \in \C\{
t\}^{ n'}$ avec ${\sf H}^N (t) \equiv h(t) \
{\rm mod} \, (\mathfrak{ m} (t))^N$
{\rm (}d'o\`u $H(0)=0${\rm )}, qui induit une application
holomorphe locale de $M$ \`a valeurs dans $M'$.
\end{corollary}

\proof
D'apr\`es le Th\'eor\`eme~1.23, les composantes $\Theta_{ j',\,
\gamma'} '( h(t))$ de l'application de r\'eflexion s'identifient \`a
des s\'eries convergentes $\theta_{ j',\, \gamma'} ' (t) \in \C\{
t\}$, pour tout $j' = 1,\, \dots,\, d'$ et tout $\gamma ' \in \N^{
m'}$. Appliquons le Th\'eor\`eme d'approximation~2.4 aux \'equations
analytiques
\def\theequation{7.148}\begin{equation}
0 \equiv 
\Theta_{ j',\, \gamma'} '( h( t)) - \theta_{ j',\, \gamma'} '
(t). 
\end{equation} 
Nous en d\'eduisons que pour tout entier $N \geq 1$, il existe une
application convergente ${\sf H}^N (t) \in \C\{ t\}^{ n'}$ 
avec ${\sf H}^N (t) \equiv h(t) \ {\rm mod} \,
(\mathfrak{ m} (t))^N$ qui satisfait ces \'equations analytiques, {\it
i.e.}~:
\def\theequation{7.149}\begin{equation}
0 \equiv 
\Theta_{ j',\, \gamma'} ' \left( {\sf H}^N (t) \right) 
- \theta_{ j',\, \gamma'} ' (t). 
\end{equation} 
Nous affirmons que {\it toute telle ${\sf H}^N$ induit
n\'ecessairement une application holomorphe locale de $M$ \`a valeurs
dans $M'$}. En effet, il d\'ecoule d'abord de~\thetag{ 7.148} et
de~\thetag{ 7.149} que $\Theta_{ j',\, \gamma'} ' ( h(t)) \equiv
\Theta_{ j',\, \gamma'} ' ({\sf H}^N (t))$ et ensuite~:
\def\theequation{7.150}\begin{equation}
\left\{
\aligned
\overline{ g}_{ j'} (\zeta,\, \Theta (\zeta,\, t)) 
\equiv 
& \
\sum_{ \gamma ' \in \N^{ m'}} \, 
\left[
\overline{ f}^{ \gamma'}
\right](\zeta,\, \Theta (\zeta,\, t)) \cdot
\Theta_{ j',\, \gamma'} '( h(t)) \\ 
\equiv 
& \
\sum_{ \gamma ' \in \N^{ m'}} \, 
\left[
\overline{ f}^{ \gamma'}
\right](\zeta,\, \Theta (\zeta,\, t)) \cdot
\Theta_{ j',\, \gamma'} '( {\sf H}^N(t)),
\endaligned\right.
\end{equation}
pour tout $j'= 1,\, \dots,\, d'$. Autrement dit, nous avons~: $r_{
j'}' \left({\sf H}^N (t), \, \overline{ h} (\zeta,\, \Theta (\zeta,\,
t)) \right)\equiv 0$, pour tout $j'= 1,\, \dots,\, d'$. En rempla\c
cant $\tau'$ par $\overline{ h} (\zeta,\, \Theta (\zeta,\, t))$ et
$t'$ par ${\sf H}^N (t)$ dans la premi\`ere ligne de la deuxi\`eme
colonne de~\thetag{ 7.44}, nous en d\'eduisons~:
\def\theequation{7.151}\begin{equation}
\left\{
\aligned
0 \equiv 
& \
\overline{ r}_{ j'}' \left(
\overline{ h} (\zeta,\, \Theta (\zeta,\, t)),\, 
{\sf H}^N (t)
\right) \\
\equiv
& \
{\sf G}_{ j'}^N (t) -
\sum_{ \gamma ' \in\N^{ m'}} \,
{\sf F}^N (t)^{ \gamma'} \cdot
\overline{ \Theta}_{ j',\, \gamma'} ' \left(
\overline{ h} \left(
\zeta,\, \Theta (\zeta,\, t)
\right)\right),
\endaligned\right.
\end{equation}
pour tout $j'= 1,\, \dots,\, d'$. Rempla\c cons alors $w$ par
$\overline{ \Theta} (z,\, \tau)$ dans ces identit\'es, en tenant
compte de la relation $\xi \equiv \Theta \left( \zeta,\, z,\,
\overline{ \Theta} (z,\, \tau) \right)$, qui d\'ecoule imm\'ediatement
de~\thetag{ 3.7}, ce qui donne~:
\def\theequation{7.152}\begin{equation}
{\sf G}_{ j'}^N
\left(
z,\, \overline{ \Theta} (z,\, \tau)
\right) \equiv
\sum_{ \gamma' \in \N^{ m'}} \,
{\sf F}^N 
\left(z,\, \overline{ \Theta} (z,\, \tau)\right)^{ \gamma'}
\cdot
\overline{ \Theta}_{ j',\, \gamma'} ' \left(
\overline{ h} (\tau)
\right),
\end{equation}
pour tout $j'= 1,\, \dots,\, d'$. 

Enfin, pour terminer, rempla\c cons
$\overline{ \Theta}_{ j',\, \gamma'} ' \left( \overline{ h} (\tau)
\right)$ par $\overline{ \Theta}_{ j',\, \gamma'} ' \left( \overline{
H}^N (\tau) \right)$ dans ces identit\'es (cela est possible, gr\^ace
aux conjugu\'ees des relations~\thetag{ 7.148} et~\thetag{ 7.149}), ce
qui donne~:
\def\theequation{7.153}\begin{equation}
{\sf G}_{ j'} 
\left(
z,\, \overline{ \Theta} (z,\, \tau)
\right) \equiv
\sum_{ \gamma' \in \N^{ m'}} \,
{\sf F}^N 
\left(z,\, \overline{ \Theta} (z,\, \tau)\right)^{ \gamma'}
\cdot
\overline{ \Theta}_{ j',\, \gamma'} ' \left(
\overline{ {\sf H}}^N (\tau)
\right),
\end{equation}
pour tout $j'= 1,\, \dots,\, d'$. Ces derni\`eres identit\'es
d\'emontrent clairement que ${\sf H}^N$ induit une application de $M$
\`a valeurs dans $M'$. Le corollaire~\thetag{ 7.147} est d\'emontr\'e.
\endproof

\newpage

\bigskip

\begin{center}
\begin{minipage}[t]{11cm}
\baselineskip =0.35cm
{\footnotesize

\noindent
{\sc Abstract}.
Searching normal forms\footnotemark[1] for real analytic submanifolds
of $\C^n$ involves convergence problems. In 1983, J.K.~Moser and
S.M.~Webster provided examples of real analytic surfaces in $\C^2$
having an isolated hyperbolic (in the sense of E.~Bishop) complex
tangency, which are formally but not holomorphically normalizable
(because of the presence of small divisors), even if the normal form
is itself real analytic or algebraic. On the contrary, it appears that
such a nonconvergence phenomenon does not appear for submanifolds of
$\C^n$ whose CR dimension is locally constant, in view of recent
results by S.M.~Baouendi, P.~Ebenfelt and L.P.~Rothschild. These
results hold true with hypotheses which are relatively simple, but
satisfied at a Zariski-generic point\footnotemark[2]. Notably, these
authors establish that every invertible formal CR mapping between two
submanifolds of $\C^n$ which are real analytic, generic, finitely
nondegenerate and minimal (in the sense of J.-M.~Tr\'epreau and
A.E.~Tumanov) is convergent. In this paper, we establish a more
general convergence theorem, which is valid without any nondegeneracy
condition, and which confirms the rigidity of the CR category ({\it
see}~Theorem~1.23). This result may be interpreted as a formal Schwarz
reflection principle for CR mappings. We deduce that every formal CR
equivalence between two submanifolds of $\C^n$ which are real
analytic, generic and minimal is convergent if and only if both
submanifolds are holomorphically nondegenerate (in the sense of
N.~Stanton). Finally, we establish that two submanifolds of $\C^n$
which are real analytic, generic and minimal are formally CR
equivalent if and only if they are biholomorphically equivalent.

}

\end{minipage}
\end{center}

\vfill
\end{document}

%% file: complexification.pstex_t
\begin{picture}(0,0)%
\includegraphics{complexification.pstex}%
\end{picture}%
\setlength{\unitlength}{3947sp}%
\begingroup\makeatletter\ifx\SetFigFont\undefined%
\gdef\SetFigFont#1#2#3#4#5{%
  \reset@font\fontsize{#1}{#2pt}%
  \fontfamily{#3}\fontseries{#4}\fontshape{#5}%
  \selectfont}%
\fi\endgroup%
\begin{picture}(5724,2349)(56,-1625)
\put(866,-1099){\makebox(0,0)[lb]{\smash{{\SetFigFont{8}{9.6}{\familydefault}{\mddefault}{\updefault}{\color[rgb]{0,0,0}$0$}%
}}}}
\put(1480,417){\makebox(0,0)[lb]{\smash{{\SetFigFont{8}{9.6}{\familydefault}{\mddefault}{\updefault}{\color[rgb]{0,0,0}$\{t=t_p\}$}%
}}}}
\put(1617,224){\makebox(0,0)[lb]{\smash{{\SetFigFont{8}{9.6}{\familydefault}{\mddefault}{\updefault}{\color[rgb]{0,0,0}$\underline{\mathcal{S}}_{t_p}$}%
}}}}
\put(2375,503){\makebox(0,0)[lb]{\smash{{\SetFigFont{8}{9.6}{\familydefault}{\mddefault}{\updefault}{\color[rgb]{0,0,0}$\underline{\Lambda}$}%
}}}}
\put(2183,-381){\makebox(0,0)[lb]{\smash{{\SetFigFont{8}{9.6}{\familydefault}{\mddefault}{\updefault}{\color[rgb]{0,0,0}$\{\tau=\tau_p\}$}%
}}}}
\put(1625,-1081){\makebox(0,0)[lb]{\smash{{\SetFigFont{8}{9.6}{\familydefault}{\mddefault}{\updefault}{\color[rgb]{0,0,0}$t_p$}%
}}}}
\put(2315,-1060){\makebox(0,0)[lb]{\smash{{\SetFigFont{8}{9.6}{\familydefault}{\mddefault}{\updefault}{\color[rgb]{0,0,0}$t$}%
}}}}
\put(1859,-839){\makebox(0,0)[lb]{\smash{{\SetFigFont{8}{9.6}{\familydefault}{\mddefault}{\updefault}{\color[rgb]{0,0,0}$\mathcal{L}$}%
}}}}
\put(631,-304){\makebox(0,0)[lb]{\smash{{\SetFigFont{8}{9.6}{\familydefault}{\mddefault}{\updefault}{\color[rgb]{0,0,0}$\tau_p$}%
}}}}
\put(528,220){\makebox(0,0)[lb]{\smash{{\SetFigFont{8}{9.6}{\familydefault}{\mddefault}{\updefault}{\color[rgb]{0,0,0}$\mathcal{M}$}%
}}}}
\put(869,123){\makebox(0,0)[lb]{\smash{{\SetFigFont{8}{9.6}{\familydefault}{\mddefault}{\updefault}{\color[rgb]{0,0,0}$\underline{\mathcal{L}}$}%
}}}}
\put(1738,-267){\makebox(0,0)[lb]{\smash{{\SetFigFont{8}{9.6}{\familydefault}{\mddefault}{\updefault}{\color[rgb]{0,0,0}$\mathcal{S}_{\tau_p}$}%
}}}}
\put(1074,-1496){\makebox(0,0)[lb]{\smash{{\SetFigFont{9}{10.8}{\familydefault}{\mddefault}{\updefault}{\color[rgb]{0,0,0}{\sc Figure~1: G\'eom\'etrie de la  complexification $\mathcal{M}$}}%
}}}}
\put(3016,219){\makebox(0,0)[lb]{\smash{{\SetFigFont{8}{9.6}{\familydefault}{\mddefault}{\updefault}{\color[rgb]{0,0,0}La complexification d'une sous-vari\'et\'e}%
}}}}
\put(3016,-68){\makebox(0,0)[lb]{\smash{{\SetFigFont{8}{9.6}{\familydefault}{\mddefault}{\updefault}{\color[rgb]{0,0,0}feuilletages invariants qui sont les}%
}}}}
\put(3016,-212){\makebox(0,0)[lb]{\smash{{\SetFigFont{8}{9.6}{\familydefault}{\mddefault}{\updefault}{\color[rgb]{0,0,0}sous-vari\'et\'es int\'egrales des }%
}}}}
\put(3016,-355){\makebox(0,0)[lb]{\smash{{\SetFigFont{8}{9.6}{\familydefault}{\mddefault}{\updefault}{\color[rgb]{0,0,0}complexifi\'es des champs de}%
}}}}
\put(3016,-498){\makebox(0,0)[lb]{\smash{{\SetFigFont{8}{9.6}{\familydefault}{\mddefault}{\updefault}{\color[rgb]{0,0,0}vecteurs de types $(1,\, 0)$ et $(0,\, 1)$}%
}}}}
\put(3016,-641){\makebox(0,0)[lb]{\smash{{\SetFigFont{8}{9.6}{\familydefault}{\mddefault}{\updefault}{\color[rgb]{0,0,0}et qui s'identifient aussi aux}%
}}}}
\put(3016,-784){\makebox(0,0)[lb]{\smash{{\SetFigFont{8}{9.6}{\familydefault}{\mddefault}{\updefault}{\color[rgb]{0,0,0}sous-vari\'et\'es de Segre}%
}}}}
\put(3016,-927){\makebox(0,0)[lb]{\smash{{\SetFigFont{8}{9.6}{\familydefault}{\mddefault}{\updefault}{\color[rgb]{0,0,0}complexifi\'ees.}%
}}}}
\put(3016, 75){\makebox(0,0)[lb]{\smash{{\SetFigFont{8}{9.6}{\familydefault}{\mddefault}{\updefault}{\color[rgb]{0,0,0}analytique r\'eelle porte une paire de}%
}}}}
\end{picture}%

%% file: orbite.pstex_t
\begin{picture}(0,0)%
\includegraphics{orbite.pstex}%
\end{picture}%
\setlength{\unitlength}{4144sp}%
\begingroup\makeatletter\ifx\SetFigFont\undefined%
\gdef\SetFigFont#1#2#3#4#5{%
  \reset@font\fontsize{#1}{#2pt}%
  \fontfamily{#3}\fontseries{#4}\fontshape{#5}%
  \selectfont}%
\fi\endgroup%
\begin{picture}(5424,2454)(249,-1963)
\put(1760,-1847){\makebox(0,0)[lb]{\smash{{\SetFigFont{9}{10.8}{\familydefault}{\mddefault}{\updefault}{\color[rgb]{0,0,0}{\sc Figure~2~: Cha\^{\i}nes de Segre dans $\mathcal{ M}$}}%
}}}}
\put(1870,-1491){\makebox(0,0)[lb]{\smash{{\SetFigFont{9}{10.8}{\familydefault}{\mddefault}{\updefault}{\color[rgb]{0,0,0}$0$}%
}}}}
\put(2041,291){\makebox(0,0)[lb]{\smash{{\SetFigFont{9}{10.8}{\familydefault}{\mddefault}{\updefault}{\color[rgb]{0,0,0}$\tau$}%
}}}}
\put(4405,-1312){\makebox(0,0)[lb]{\smash{{\SetFigFont{9}{10.8}{\familydefault}{\mddefault}{\updefault}{\color[rgb]{0,0,0}$t$}%
}}}}
\put(2733,-40){\makebox(0,0)[lb]{\smash{{\SetFigFont{9}{10.8}{\familydefault}{\mddefault}{\updefault}{\color[rgb]{0,0,0}$\underline{\Gamma}_3 (z_{(3)})$}%
}}}}
\put(2964,-826){\makebox(0,0)[lb]{\smash{{\SetFigFont{9}{10.8}{\familydefault}{\mddefault}{\updefault}{\color[rgb]{0,0,0}$\underline{\Gamma}_2 (z_{(2)})$}%
}}}}
\put(931,-861){\makebox(0,0)[lb]{\smash{{\SetFigFont{9}{10.8}{\familydefault}{\mddefault}{\updefault}{\color[rgb]{0,0,0}$\underline{\Gamma}_1 (z_1)$}%
}}}}
\put(1786, -9){\makebox(0,0)[lb]{\smash{{\SetFigFont{9}{10.8}{\familydefault}{\mddefault}{\updefault}{\color[rgb]{0,0,0}$\mathcal{ M}$}%
}}}}
\put(360,279){\makebox(0,0)[lb]{\smash{{\SetFigFont{9}{10.8}{\familydefault}{\mddefault}{\updefault}{\color[rgb]{0,0,0}$\C^n\times\C^n$}%
}}}}
\put(3642,-289){\makebox(0,0)[lb]{\smash{{\SetFigFont{9}{10.8}{\familydefault}{\mddefault}{\updefault}{\color[rgb]{0,0,0}$\underline{\Gamma}_4 (z_{(4)})$}%
}}}}
\end{picture}%